\input amstex




\input colordvi
\long\def\Coloration#1#2{\textColor {#1}{#2}\textColor{.0 .0 .0 1.0}}

\def\red#1{\Coloration{.0 .99 .87 .14}#1}


\comment

olive-green  7,0,95,35       green (51,0,96,39)       red (0,100,87,14)
    darkblue (94,100,0,32)   blue (98,30,0,28)      darkbrown (0,32,61,82)
    brown (0,43,83,58)      purplered (5,98,0,29)         purple (44,97,0,44)
     verydarkblue(96,97,0,84)     pink (0,28,6,4)         orange (0,53,99,31)
\endcomment




\font\rm=cmr10 \rm

\font\bf=cmb10
\font\Rm=cmr9 at 11pt
\rm
\font\it=cmsl9 at 10pt
 at 7pt

\font\Rrm=cmr17 at 16pt
   \font\Rm=cmr12 at 11.5pt

\long\def\Pf{\par\noindent {\it Proof.} }
\def\({\left(}
\def\){\right)}
\def\st{such that }
\def\qed{\hfill$\bullet$\vskip 4pt}
\def\quotes#1{{\lq\lq #1\rq\rq}}
\def\brcs#1{\left\{ #1\right\}}
\def\sgn {\text{\rm sgn\,}}

\def\comp{\circ}
\def\Log{\text{Log\,}}
\def\iso{\cong}
\def\wrt{with respect to }
\def\:{\,:}

\def\supp{\text{supp}\,}

\def\ker{\text{ker\,}}

\def\Aut{\text{Aut}\,}

\def\rr{{\Cal V}}

\def\A{{\Cal A}}

\def\G{{\Cal G}}

\def\C{\text{\bf C}}

\def\M{{\text M}}

\def\Im{{\text{Im}\,}}

\def\I{\text{I\,}}

\def\tr{\text{\,tr}}

\def\R{\text{\bf R}}
\def\N{\text{\bf N}}
\def\Z{\text{\bf Z}}
\def\Q{\text{\bf Q}}

\def\Arrow #1;#2.{#1\:#2 \to }

\def\Set#1#2{\brcs{#1 \left|\vphantom{#1 #2} \right.#2}}

\def\Oh#1{{\pmb O}\(#1\)}

\def\oh#1{{\pmb o}\(#1\)}

 \def\supp#1{\text{supp\,}#1}
\def\LL#1#2{{\Cal M}_{#1,#2}}
\def\Rrr#1,#2{{\Cal J}_{#1,#2}}
\def\slfrac#1#2{{\raise -.07 ex\hbox{$^{#1}$}}\!/\raise .35 ex \hbox{${}_{#2}$}}
\def\ssf #1/#2{\slfrac {#1}{#2}}

\def\pd #1,#2.{\frac {\partial #1}{\partial #2}}

   \long\def\Lem
#1.#2\par{\vskip4pt{\baselineskip=13pt\font\it=cmsl12 at
11.5pt\Rm
   \noindent {\rm \uppercase{#1}} #2\vskip3pt

   }} 

\long\def\Proclaim #1.#2 \endproclaim{\vskip4pt{\baselineskip=13pt\font\it=cmsl12 at
11.5pt\Rm
   \noindent {\rm \uppercase{#1}} #2\vskip3pt

   }} 

\long\def\remark #1\endremark{\vskip 2pt \noindent {\it Remark\/} #1\par}

\long\def\Sectionhead #1.#2:\par #3{\vskip 4pt \noindent {\bf #1 #2}vskip 2pt\noindent\nospace #3}

\long\def\Title #1\par {\noindent{\Rrm #1}\vskip 9pt}

 \long\def\SubT #1.{\noindent {\it #1\/} } 
 
 \long\def\SecT
#1\par{\vskip 3pt \noindent {\bf #1}\vglue1pt
   \noindent}

\long\def\subtitle #1.{\vskip 2pt \noindent {\it #1}}

\long\def\Rmk#1\par{\vskip 1pt \noindent {\it
Remark.} #1\vskip2pt}

\long\def\Abstract #1\par{{\leftskip= 3 true cm \rightskip = 3 true cm \font\it=cmsl10 \font\rm=cmr10 \baselineskip = 10pt
\parindent=.35 true cm\rm\noindent 
{\it Abstract} #1\vskip 8pt

}}

\long\def\Author #1 \par{\noindent{\it #1}}


\def\thmxx{ fill me in!! }

\def\ceil#1{\lceil #1 \rceil}

\def\flo #1{\lfloor #1 \rfloor}
\def\PP{{\Cal P}}
\def\diag{\text{{\rm diag}}\,}

\def\Q{\text{\bf Q}}
\def\Aff{\text{Aff\,}}
\def\Inf{\text{Inf\,}}

\def\paren #1{\/{\rm(}#1\/{\rm)}}

\let\sign=\sgn %
\def\gl{\text{GL}}
\def\Gd #1{\text{Gd}\(#1\)}

\def\op{^{\text{op}}}

 %

\def\Leg#1,#2.{\(\frac{#1}{{#2}}\)}
\def\Ip #1,#2.{\langle\!\langle #1,#2\rangle\! \rangle}
\def\Ipd #1,#2.{\Ip #1,#2._d}

\def\paren #1{{\rm(}#1\/{\rm)}}
 
\def\mset #1{\left[\hglue-2.95 pt\left| \vphantom{$A_A^A$}#1 \right|\hglue-2.95 pt\right]}

\let\mset=\Lt 
\def\Spec{\text{Spec\,}}
\let\spec=\Spec
\def\1{{\pmb 1}}
\def\fillmein{{\red{{fill me in }}}}
\def\LLL{{\Cal L}}
\def\LL{\hat l}
\def\pplus{^{\prime+}}


\let\hat=\widehat


\def\Aff{\text{Aff\,}}
\NoBlackBoxes

\def\End{\text{End\,}}
\def\G{{\Cal G}}

\def\Bf#1{\text{\bf #1}}
\def\I{\text{I\,}}
\def\tr{\text{tr\,}}

\def\flo #1{\lfloor #1 \rfloor}

\let\iso=\cong
\def\SS{{\Cal S}}
\def\Supp{\text{supp}\,}

\def\First{1}
\def\firstone{\First.1}
\def\firsttwo{\First.2}
\def\firstthr{\First.3}
\def\firstfou{\First.4}

\def\firsteig{\First.8}

\def\WC{2}
\def\wcone{\WC.1}
\def\wctwo{\WC.2}
\def\wcthr{\WC.3}
\def\wcfou{\WC.4}
\def\wcfiv{\WC.5}
\def\wcsix{\WC.6}
\def\wcsev{\WC.7}
\def\wceig{\WC.8}
\def\wcnin{\WC.9}
\def\wcten{\WC.10}
\def\wcele{\WC.11}
\def\wctwe{\WC.12}

\def\Sty{11}

\def\One{12}
\def\Ends{13}
\def\oneone{\One.1}
\def\onetwo{\One.2}
\def\onethr{\One.3}
\def\onefou{\One.4}
\def\onefiv{\One.5}
\def\onesix{\One.6}
\def\onesev{\One.7}
\def\oneeig{\One.8}
\def\onenin{\One.9}
\def\oneten{\One.11}
\def\oneele{\One.12}
\def\onetwe{\One.13}
\def\onethi{\One.14}
\def\oneftn{\One.15}
\def\oneffn{\Ends.1}

\def\onetty{\Sty.1}
\def\onetto{\Ends.2}
\def\onettt{\Ends.3}
\def\onetth{\Ends.4}
\def\onettf{\Ends.5}

\def\onettw{\Ends.9}

\def\onetho{\Ends.6}

\def\onetht{\Ends.7}
\def\onethe{\Ends.8}

\def\endsone{\Ends.1}
\def\endstwo{\Ends.2}
\def\endsthr{\Ends.3}

\def\Two{9}
\def\twoone{\Two.1}

\def\twothr{\Two.3}
\def\twofou{\Two.4}
\def\twofiv{\Two.5}
\def\twosix{\Two.2}
\def\twosev{\Two.6}
\def\twoeig{\Two.7}
\def\twonin{\Two.9}
\def\twoten{\Two.10}
\def\twoele{\Two.11}
\def\twotwe{\Two.12}
\def\twothi{\Two.13}
\def\twoftn{\Two.14}
\def\twoffn{\Two.15}
\def\twosxn{\Two.16}
\def\twosvn{\Two.17}

\def\twontn{\Two.8}

\def\X{3}
\def\Xone{\X.1}
\def\Xtwo{\X.2}
\def\Xthr{\X.3}
\def\Xfou{\X.4}
\def\Xfiv{\X.5}
\def\Xsix{\X.6}
\def\Xsev{\X.7}
\def\Xeig{\X.8}
\def\Xnin{\X.10}
\def\Xten{\X.11}
\def\Xele{\X.12}
\def\Xtwe{\X.13}
\def\Xthi{\X.14}
\def\Xftn{\X.15}

\def\Xthy{\X.9}

\def\M{4}
\def\Mone{\M.1}
\def\Mtwo{\M.2}

\def\B{5}
\def\Bone{\B.1}
\def\Btwo{\B.2}
\def\Bthr{\B.3}
\def\Bfou{\B.4}
\def\Bfiv{\B.5}
\def\Bsix{\B.6}
\def\Bsev{\B.7}
\def\Beig{\B.8}
\def\Bnin{\B.9}
\def\Bten{\B.10}

\def\IP{6}
\def\IPone{\IP.1}
\def\IPtwo{\IP.2}
\def\IPthr{\IP.3}

\def\TR{7}
\def\TRone{\TR.1}
\def\TRtwo{\TR.2}
\def\TRthr{\TR.3}
\def\TRfou{\TR.4}
\def\TRfiv{\TR.5}

\def\MI{8}
\def\MIone{\MI.1}
\def\MItwo{\MI.2}
\def\MIthr{\MI.3}
\def\MIfou{\MI.4}
\def\MIfiv{\MI.5}
\def\MIsix{\MI.6}

\def\MIeig{\MI.7}

\def\IMI{10}
\def\IMIone{\IMI.1}
\def\IMItwo{\IMI.2}
\def\IMIthr{\IMI.3}
\def\IMIfou{\IMI.4}
\def\IMIfiv{\IMI.5}
\def\IMIsix{\IMI.6}
\def\IMIsev{\IMI.7}

\def\Y{A1}
\def\Yone{\Y.1}
\def\Ytwo{\Y.2}
\def\Ythr{\Y.3}
\def\Yfou{\Y.4}
\def\Yfiv{\Y.5}
\def\Ysix{\Y.6}
\def\Ysev{\Y.7}

\def\FTN{14}
\def\FTNone{\FTN.1}
\def\FTNtwo{\FTN.2}
\def\FTNthr{\FTN.3}
\def\FTNfou{\FTN.4}
\def\FTNfiv{\FTN.6}
\def\FTNsix{\FTN.7}
\def\FTNsev{\FTN.8}
\def\FTNeig{\FTN.9}
\def\FTNnin{\FTN.5}
\def\FTNten{\FTN.10}


\Title Boundary of the boundary for random walks on groups 

\Abstract We study fine structure related to finitely supported random walks on infinite finitely generated discrete groups, largely motivated by dimension group techniques. The unfaithful extreme harmonic functions (defined only on proper space-time cones), aka unfaithful pure traces,  can be represented on systems of finite support, avoiding dead ends. This motivates properties of the random walk (WC) and of the group (SWC) which become of interest in their own right. While all abelian groups satisfy WC, the do not satisfy SWC; however some abelian by finite groups do satisfy the latter, and we characterize when this occurs. 
{\par}In general, we determine the maximal order ideals, aka, maximal proper space-time subcones of that generated by the group element $1$ at time zero), and show that the corresponding quotients are stationary simple dimension groups, and that all such can occur for the free group on two generators. 
{\par } We conclude with a case study of the discrete Heisenberg group, determining among other things, the pure traces (these are the unfaithful ones, not arising from characters). 


\noindent {\it David Handelman}%
\plainfootnote{*}{\rm Mathematics Department, University of Ottawa, Ottawa ON, Canada; rochelle2\@sympatico.ca}

{}\plainfootnote{}{AMS(MOS) classification:  46A55,  60G50  secondary 06F20, 22E25, 46A22, 46N99, 52B11, 54F05, 60B10; key words and phrases:  random walks on groups,  space-time cone, harmonic function, pure trace, perfidious trace, discrete Heisenberg group, dimension group, order ideal, dead end, nilpotent group, group ring, weight function} 

 \SecT Introduction 

 We discuss fine and very fine structure of objects associated to random walks on groups. Let $G$ be a finitely generated discrete group, and let $f$ be an element of the group ring (either over the reals, $\R G$, or over the integers, $\Z G$) \st all of its coefficients are positive, $1$ belongs to its support, and its support generates $G$ as a semigroup; such an element will be called {\it admissible,} as will its support $S$. We do not require that $S$ be symmetric, nor that the coefficients add to $1$.

 Form the smallest space-time cone $\Cal C \subseteq G \times \Z^+$ that contains $1 \in G$ at time zero. Now consider the pure (extremal) space-time harmonic functions (traces) on $\Cal C$. As is well-known, the faithful ones extend to space-time harmonic functions on all of $G \times \Z^+$ and are given by eigenvectors for the multiplication operator obtained from $f$. The unfaithful ({\it perfidious\/}) extremal ones, however, do not extend, and there are always plenty of them; but they are not well-understood, except in the case that $G$ is abelian, or for a few special cases, such as $G$ being the free group on two or more generators. 

 The lattice of space-time subcones of $\Cal C$ is also interesting. If $G$ is abelian, then this lattice satisfies the ascending chain condition (which amounts to finitely many \quotes{peak} points, or generators in a suitable sense). This can be deduced from the Hilbert basis theorem. On the other hand, for the free group, unsurprisingly, the lattice of space-times cones in $\Cal C$ is about as far as possible from the ascending chain condition, and is in natural bijection with the paths of the Cayley tree. One would expect that for a nilpotent group, e.g., the discrete Heisenberg group, $H_3$, that the lattice of subcones would be similar to that for an abelian group. Rather surprisingly, this turns out to be false---the lattice is not noetherian.

 Because $\Cal C $  is finitely generated in an appropriate sense, maximal subcones exist. For $G = \Z^d$, there are only finitely many (but at least $d+1$), depending on the choice of $f$, and are easy to describe. For $G =F_2$ with the natural choice for $f$, they correspond to the path space. But for groups with intermediate properties (e.g., nilpotent), there are only countably many, and in the case of $H_3$ with the special choice for $f$, they are easy to describe. 

 In general, we can say something about the maximal subcones, but it is easier to explain if we switch to the language of dimension groups and Bratteli diagrams. There is a lexicon which more or less does the translation (which is not entirely bijective, but is close enough). 

Let $A$ be one of $\R$ or $\Z$. Form the direct limit of iterated left multiplications by $f$, $A_f = \lim \Arrow f\times ; AG. AG$, as a partially ordered vector space (abelian group). We obtain what amounts to an order ideal, $R_f = \lim \Arrow f\times; A^{S^n} . A^{S^{n+1}}$; this corresponds to the cone $\Cal C$. There is a natural shift $\Arrow \SS; R_f.R_f$ which acts as a positive endomorphism (with additional properties); it amounts to increasing the time index (the second coordinate of $G \times \Z$) by $1$. The perfidious (unfaithful) pure traces on $R_f$ are precisely the pure traces that kill the image of $\SS$ on $R_f$, and we are led to study the quotient partially ordered abelian group, $R_f/\SS R_f$ (corresponding to the complement of the obvious space-time subcone in $\Cal C$). 

 This can also be expressed naturally as a direct limit, $\lim \Arrow \overline f \times ; A^{\Gamma_k}. A^{\Gamma_{k+1}}$ where $\Gamma_k \subset S^k$ and $\overline f \times$ is restriction and compression of $f \times$. Determining the $\Gamma_k$ is interesting. There are several choices (leading to the same direct limit), and  the simplest one is given by 
 $$
 \Gamma_k = S^k \setminus S^{k-1}.  
 $$
 This is somewhat awkward to deal with and is larger than it need be. 

 To refine it, say the pair $(G,S)$ (where $S$ is  admissible) {\it satisfies WC\/} if for all positive integers $k$, 
$$
\Set{g \in G}{\exists\, m \in \N \text{ \st } 
S^m g \subseteq S^{m+k}}\quad \text{is finite.}
$$
 The reader can experiment with this concept by looking at $G = \Z$ and the holey $f = x^n + x + 1 + x^{-1} + x^{-n}$ for $n > 2$. For all $k > n-1$ but for no smaller $k$, $\supp f^k$ is the interval of integers $\brcs{-k,-k+1, \dots, k-1,k}$. 
 
 This condition, WC, amounts to the function 
$\Arrow \tilde l_S; G. \Z^+$ 
given by $$g \mapsto \inf \Set{k} {\exists\, m \in \N, S^m g \subset S^{m+k}}
$$ 
being a weight function. The weight function corresponding to $S$, $l_S$, given by $g \mapsto \inf \Set{k}{g \in S^k}$ is obviously related to this, but in general is different from $\tilde l_S$ (as in our holey example, and also for any finite group). When $(G,S)$ satisfies WC,  we can replace $\Gamma_k$ by $\Gamma_k':= \tilde l_S^{-1}(k)$, and this is   easier to deal with. 

 It is not clear whether WC really depends on the support $S$, or holds independently of the choice (that is, if $(G,S)$ satisfies WC for the  one admissible $S$, then $(G,S')$ satisfies WC for all admissible $S'$). 
 
 However, it is true for abelian groups, and a stronger property called SWC, which depends only on $G$ is true for nontrivial amalgamated products and hyperbolic groups (the last is due to Elisabeth Fink). Rather surprisingly, the WC property fails for our standard $H_3$ example, and probably for all torsion-free nonabelian nilpotent groups, and this is closely related to the structure of $R_f/\SS R_f$.  In contrast,  the  stronger form, SWC, does hold for the infinite dihedral group---the usual bad boy in this context---because it is a nontrivial free product. Among abelian by finite groups, those that satisfy SWC are characterized; it turns out only to depend on the natural rational action of the finite quotient group, and is independent of (some) extension data. 
 
 Even if $(G,S)$ does not satisfy WC, we can still excise $\Gamma_k$ to obtain a useful realization for the quotient $R_f/\SS R_f$. Then we have a rather startling result on the maximal order ideals (equivalently, maximal space-time subcones of $\Cal C$), namely that they all arise from eventually stationary systems of finite width, and in particular, have unique trace. 
 
 Returning to the traces on $R_f/\SS R_f$ (in fact, we never really left them), we analyze them in full detail in our $H_3$ example. In addition to the faithful pure traces, which corresponding to the real characters on $G/G' \iso \Z$, we have four families of suitably multiplicative traces corresponding the the standard random walk on $\N$ (move forward or  one step or stay put with equal probability), and eight families of  discrete traces, arising from the  sequences of elements having unique predecessor in  a suitable subcone. The natural action of $D_4$ is transitive on the two sets of pure traces; moreover, the set of  limit points of each family of discrete traces consists of one of the four families of suitably multiplicative traces.

 \SecT (Some) definitions and notation, and statement of results

Let $G$ be a (countable) discrete group, and form the group rings $\R G$ and  $\Z G$---the integral group ring. These are respectively a real algebra and a ring, and we denote them $AG$, where $A$ is either $\R$ or $\Z$ (or occasionally, $\Q$, the rationals). The elements of $AG$ can be written in the form $x = \sum (x,g)g$, using inner product notation to denote the coefficient of the basis element $g$ corresponding to $g \in G$; sometimes we write $x(g)$ for $(x,g)$. This admits natural positive cone, $(AG)^+ = \Set{x \in AG}{(x,g) \geq 0 \text{ for all $g \in G$}}$; that is, the positive cone consists of elements of the group algebra/ring all of whose coefficients are nonnegative. 

Suppose $f \in (AG)^+$. 
Define the big direct limit (as ordered vector spaces), $A_f := \lim \Arrow f\times; A.A$, given by multiplicaton on the left by $f$. This preserves the positive cone. If $G$ is abelian, the direct limit (as a vector space, and as an $A$-module) is simply the real algebra/ring $AG[f^{-1}]$, shorthand for $AG[X]/(Xf -1)$, formally inverting $f$. (If $A$ has no zero divisors, e.g., if $G$ is torsion-free abelian, then the kernel of the map $AG \to AG[f^{-1}]$ is zero.)

We also note that since $AG$ is lattice-ordered, $A_f$ is a dimension group (and in fact, a dimension space as defined in [H4]), \wrt the direct limit ordering. As usual for direct limit constructions, the typical element of the direct limit is denoted $[h,k]$ with $h \in A$ and $k \in \Z^+$, the latter identifying the time (or the copy of $A$ to which $h$ belongs), starting from $k = 0$. Of course, $[h,k] = [fh,k+1]$. And the positive cone obtained from the direct limit ordering consists of elements $[h,k]$ for which there  exists a positive integer $n$ \st the group algebra/ring element $f^n h$ has only nonnegative coefficients. 

There is a shift function, $\Arrow \SS ;A_f.A_f$  given by $\SS\([h,k]\) = [h,k+1]$, as well as a function induced by left multiplication, $\Arrow f\times; A_f.A_f$ given by $[h,k] \mapsto [fh,k]$. Since $[fh,k+1] = [h,k]$, it follows that $f\times$ and $\SS$ are mutually inverse; both are clearly order preserving, so they induce inverse order automorphisms of $A_f$ as partially ordered vector spaces/abelian groups.  

Define $\Bf 1 = [1,0]$ where the   $1$ in the first coordinate represents the identity element of $G$. Let $R_f$ denote the order ideal generated by $\Bf 1$, that is,
$$
R_f = \Set{[h,k]}{\exists\ K \in \N \text{ \st } -K \Bf 1 \leq [h,k] \leq K \Bf 1 }.
$$
In particular, $[h,k] \in R_f$ iff there exists $m$ \st $- K f^{m+k}\leq f^{m} h \leq K f^{m+k}$, where of course, the ordering is coefficientwise (that is, the usual ordering on $AG$). Some observations:

\item{(i)} If $G$ is abelian, then $R_f$ is a partially ordered (real) algebra, having $\Bf 1$ as an order unit; the multiplication operation is described by $[a,k][b,l] = [ab, k+l]$. If $G$ is not abelian, then  $R_f$ is not generally a ring.
\item{(ii)} As it is an order ideal in a dimension group, $R_f$ is a dimension group.

\vskip 5pt
If $h \in AG$ and $p \in (AG)^+$, we write $h \prec p$ to mean that there exists $K \in \N$ \st for all $g \in G$, $|(h,g)| \leq K (p,g)$; equivalently, $\supp h \subset \supp p$. 

If $D$ is a partially ordered abelian group and $d,d'$ are elements of $D$ with $d' \in D^+$, then we also write $d \prec d'$ to denote that there exists a positive integer $N$ \st $-Nd' \leq d \leq N d'$; equivalently, that $d$ belongs to the order ideal generated by $d'$, that is, $d \in \langle d'\rangle$. 

We are interested in a couple of things. One is $R_f^+$; determining this is equivalent to determining $\Set{h \in AG}{\exists\ n, m \in \N \text{ \st $h \prec f^n$ and $f^m h \in (AG)^+$}}$. The second is to try to describe the extremal (or pure, ergodic, irreducible, \dots---they all mean the same thing) traces (harmonic functions) on $R_f$. In fact, the latter yields information on the former.

By (i) and [H2, Proposition 1.1], when $G$ is abelian, every pure trace on $R_f$ is multiplicative; in this case, every multiplicative trace on $R_f$ is pure but  not every pure trace on $R_f$ can be extended to a multiplicative trace on $A_f$.

\def\LLL{{\Cal L}}
\def\RR{{\Cal R}}

\noindent{\bf Results\ }
Our first result concerns realizing $R_f/\SS R_f$ efficiently; that is, as a direct limit of the form $A^{\Gamma_m} \to A^{\Gamma_{m+1}}$ where $\brcs{\Gamma_i}$  is a disjoint family of finite subsets of $G$, approximately the group elements at distance exactly $m$ from the identity (this is not precise), but we want to do it avoiding dead ends ([GrH]) and other phenomena. This leads to the definition of $\Gamma_m'$, the Goldilocks of choices, in section 1, and it also leads to a property of $(G,S)$  (where $S$ is an admissible subset of $G$), known as WC, and its strengthening SWC. 

These allow the limit realization for $R_f/\SS R_f$ to be relatively pleasant. We see quickly  (section 2) that a lot of big and small groups, such as abelian, free,  most amalgamations, satisfy WC, and in fact the latter two classes satisfy the much stronger SWC, which abelian groups do not. However, some abelian by finite groups do satisfy SWC, and we characterize those that do (Theorem \Xone), in terms of the rational representation of the quotient. The methods involved in dealing with this are extended to very detailed structural results. 

Sections 4 and 5 discuss endomorphisms of $R_f$, and various properties related to bounded endomorphisms (a restrictive class).  For example, there is a natural notion of local order boundedness (for positive endomorphisms), but for groups satisfying WC, this implies much stronger properties (for example, Lemma \Bthr). We introduce  a few properties suggested by those of group rings, but put in the natural ordered setting. For example, in section 6, a property motivated by Jacobson's conjecture (for commutative rings, subsequently proved, then posed in the noncommutative case, partially solved), that $\cap \SS^n R_f$ be zero, is established at  least for left orderable groups (it can fail for general groups). 

Sections 7 and 8 deal with traces, and maximal order ideals. Traces are the translation of harmonic functions, but here restricted to relatively small space-time cones. The pure (or extreme, ergodic, indecomposable, minimal) traces of interest are not the globally defined ones, but those that are not faithful (and cannot be defined globally); they admit a partial action (that is, by a subsemigroup of the group), and this is used in dealing with maximal order ideals. 

The main result of section 9, is that in this generality, the quotients $R_f/M$ for $M$ a maximal order ideal are what is known as stationary dimension groups, those arising from iterating a square matrix with nonnegative real (if $A = \R$) or integer coefficients (if $A = \Z$) (Theorem \twoone). This is surprising because $R_f$ arises from what amounts to a very complicated infinite stationary system (repeatedly applying an infinite matrix), yet the quotients  by maximal order ideals are fairly simple. Section 10 concerns special maximal order ideals (that will be of importance in our case study, that of the Heisenberg group). And section 11 shows that an arbitrary square primitive $0-1$  matrix can be realized as such a quotient for suitable admissible $f$ when $G$ is the free group on two generators (Example \onetty). For general groups, the size of the matrices is usually severely limited (an extreme example, for abelian groups, size one). 

Sections 12--14 constitute a case study of the Heisenberg group: generators $g,h$ and relations $h g = z gh$ with $z$ central, with admissible $f = 1 + g + h +  g^{-1} + h^{-1}$. In spite of expectations, WC fails drastically, the set of pure faithful traces of $R_f$ fails to be dense in the set of pure traces of $R_f$, the set of space-time cones contained in the standard one (starting at element $1$ of the group at time zero) fails to be noetherian (despite the corresponding group ring being right and left noetherian as a ring), \dots. All of these properties are in start contrast to those of the abelian case. 

We describe completely the pure traces on $R_f$: (a) the pure faithful traces (which come from characters of the group; (b) their limits,  which are multiplicative traces on $R_f$ (and cannot be extended to characters); and (c) eight families of discrete traces, which consist of isolated points in the pure trace space. The first two types correspond to the pure traces in the abelianized version, and can be identified with the points of the Newton polytope; the points in the eight families (there is natural $D_4$-action on everything here). 

A crucial result in  the determination of the pure traces is calculating the limit points of the families of discrete traces (it sufices, by the $D_4$-action to deal with only  one family). We use Szerkeres' asymptotic formula for restricted partitions (of integers, with a bounded number of parts) to deduce (Theorem \onetth) that the nontrivial limit points are precisely those in (b). Most of the work   involves, not $R_f$ itself, but one of its quotients by an order ideal, denoted $\overline R_f$ implemented by $\overline f = g +  h$ (which is {\it not\/} admissible). 

\SecT Space-time cones, direct limits, order ideals, \dots

{\it General considerations.} Let $S$ be a countable set, form $S \times \Z^+$, and let $P = (P_{s,n})$
be an array of nonnegative integers with the property that for each $n$, the set $\Set{s \in S}{P_{s,n}> 0}$ is finite (even this condition can be relaxed). Form the (countable-dimensional) partially ordered vector space (or abelian group) with basis $\brcs{{}s}$, $V$; it has positive cone $V^+ :=\Set{v}{(v,{}s) \geq 0 \text{ for all $s \in S$}}$, where $(v,{}s)$ denotes the coefficient of ${}s$ in the decomposition of $v$ as a linear combination of the basis vectors. Then $\Arrow P; V.V$ is a partially ordered vector space (abelian group), and we can take the direct limit $V_P := \lim \Arrow P; V.V$; its elements are the equivalence classes $[v,n] $ arising from the equivalence relation generated by $[Pv,n+1] = [v,n]$.

Then $V_P$ admits a positive cone arising from the direct limit construction, $$V_P^+ = \Set{[v,n]}{\exists N \in \N \text{ \st $P^N v \in V^+$}}.$$

The infinite matrix $P$ (which need not have its column sums equalling one, nor even bounded) implements a directed graph structure on $S \times \Z$, so that it becomes a Bratteli diagram, where the arrows are weighted by nonnegative real numbers. Specifically, we say there is an edge from $(s,m) \to (t,m+1)$ if $P_{ts} > 0$. Then the notion of path (meaning directed path, of course) makes sense; we say there is a path from $(s,m) \to (s',m')$ if $m' > m$ and there is a path (of edges), $(s,m)\to (s_1,m+1) \to \dots \to s_{m'-m -1},m'-1) \to (s',m')$. Of course, this is equivalent to $(P^{m'-m})_{s',s} \neq 0$.

In this context, define a {\it space-time cone\/} to be a nonempty subset $U $ of $ S \times \Z^+$ with the following properties:
\item{(a)} if $(s,n) \in U$ and there is a path from $(s,n) \to (s',n')$, then $(s',n') \in U$;
\item{(b)} if $(s,n) $ has the property that every path emanating from it eventually hits $U$, then $(s,n) \in U$.

\noindent The obvious example, the improper space-time cone, is $S \times \Z^+$ itself. Space-time cones are closed under finite intersections and arbitrary joins (the join of a two cones is the smallest one containing both). A cone $U$ is {\it finitely generated\/} if there exists a finite set of points in $S \times \Z^+$ \st the smallest cone containing the finite set is $U$.

For $P$ fixed, there is a natural bijection between the set of space-time cones (associated with $P$) and the order ideals of $V_P$. Given an order ideal $I$ of $V_P$, define $U = \Set{(s,n)}{[{}s,n] \in I}$; it is routine to check that this is a space-time cone. Conversely, given $U$, let $I$ be the order ideal of $V_P$ generated by $\Set{[{}s,n]}{(s,n) \in U}$. The finitely-generated cones correspond exactly to the order ideals with order units (if $U$ is generated by the finite set $\brcs{s(i),m(i))}$, then $\sum_i [{}{s(i)},m(i)]$ is an order unit for $I$, etc).

In our situation, $S = G$ and $P = f\times$, where $f \in \R G^+$ (or $\Z G^+$), and $V_P = A_f$. The order ideal generated by $\pmb 1 = [1,0] \ in A_f$ is $R_f$, and its corresponding space-time cone is the set $U_f:= \Set{(g,m)}{\exists N \text{ \st $P^N g \prec P^{N+m}$}}$ (in many cases, we can assume that $N = 0$, that is, there are no holes).

The collections of order ideals (or space-time cones) of $A_f$ and by restriction of $R_f$ have a natural structure of a lattice: we can sum any family of order ideals, and the result will be an order ideal, and the intersection of a finite family of order ideals is an order ideal (this is true for all dimension groups). In general, the intersection of a countable descending chain of order ideals need not be an order ideal.

\Lem Lemma {0.1}. Suppose $f \sim f_0$ and $f$ is admissible. Then there is a natural lattice isomorphism between the order ideal lattices $\Cal I(A_f)$ and $\Cal I(A_{f_0})$; this induces a lattice isomorphism between $\Cal I(R_f)$ and $\Cal I(R_{f_0})$.

This is obvious, and is a special case of the following. Let $\Gamma_n$ be countable (possibly infinite) sets, and let $\brcs{P_{n,\gamma,\gamma'}}$ be an array nonnegative real numbers indexed by $(\gamma,\gamma') \in \Gamma_n \times \Gamma_{n+1}$ \st for each $\gamma \in \Gamma_n$, there are only finitely many $\gamma' \in \Gamma_{n+1}$ is not zero. Suppose $\brcs{P'_{n,\gamma,\gamma'}}$ is a similar array, with the property that $P'_{n,\gamma,\gamma'} \neq 0$ iff $P_{n,\gamma,\gamma'} \neq 0$. Then the collection of space-time cones in $\cup \Gamma_n$
determined by $P$ and $P'$ have naturally isomorphic lattices of order ideals, sending the finitely generated ones to the finitely generated ones.

A reasonable conjecture is that if $U$ is a space-time cone contained in $U_f$, then $U$ should be finitely generated; equivalently, every order ideal of $R_f$ has an order unit (that is, a relative order unit, not an order unit of $R_f$). Alternatively, every increasing chain of order ideals of $R_f$ is eventually stationary.
This is true when $G = \Z^d$ (because $R_f$ is a commutative noetherian ring), and also for the infinite dihedral group (because it comes from a matrix-valued random walk in only one real variable). It obviously fails when the group $AG$ is not noetherian (as a ring; that is, right and left noetherian). But the property fails  even for abelian by finite groups. 

It fails generically (that is, generically \wrt the choice of $f \in \R G^+$ \st $1 \in \supp f$ and $\supp f$ generates $G$ as a semigroup) for the  group, $\Z^2 \times _{\theta} C_2$ ($C_2$ denotes the cyclic group of order two; $\Z_2$ would be a natural notation, but there are too many $\Z$s running around) with the obvious action, $(1,0) \mapsto (0,1)$ and $(0,1) \mapsto (1,0)$. Viewing $\Z^2$ multiplicatively, set $f_0 = x + c$ (where $c$ represents the transposition of the generators). It is easy to check that $f:= x^{-1}y^{-1}f_0^4$ is admissible, but it is easier to work with $f_0$ itself (the constructions of $A_f$, $R_f$, etc are the same).

The corresponding matrix-valued random walk is given by $\(\smallmatrix x & 1 \\ 1 & y \\\endsmallmatrix\)$. This is the matrix appearing in [H4,  section 7], where it is observed that the corresponding bounded subring has order ideals without order units. This in turn means that $R_f$ (which is a faithful ordered module over the bounded subring) contains order ideals without order units, equivalently, not all space-time subcones of $U_f$ are finitely generated. Moreover, as a consequence of [op\.cit\.], this property is generic (we would still have to show that it is generic for the matrix-valued random walks that come out of $f\in \R G^+$, but this is straightforward) for $G = \Z^d \times_{\theta} K$ with $K$ a finite group when $d \geq 2$ (even with trivial action of $K$). 

A torsion-free   group for which the group ring $AG$ is noetherian (that is, right and left noetherian), but for which the lattice of space-time cones in section \One.

\SecT 1 Weight functions and Goldilocks realizations of $R_f/\SS R_f$

Let $G$ be a finitely generated  group. We say a finite subset $S$ of $G$ is {\it
admissible\/} if $1 \in S$ and $\cup S^n = G$. If $f \in AG^+$ and $\supp f$ is admissible, then we say $f$ is {\it admissible.}

We define a {\it weight function\/} (there are several different definitions in the literature) on a group $G$ to be a function $\Arrow l; G.\Z^{+}$ with the following properties. 
 \item{(0)} $l(1) = 0$; 
  \item{(i)} for all $g,h \in G$, $l(gh) \leq l(g) + l(h)$; 
  \item{(ii)}for all nonnegative integers $k$, the set $l^{-1}(\brcs{k})$ is finite.

  If $l$ is a weight function, we can also define $l$ on finite subsets of $G$ via $l(T) = \max_{g \in T}{l(g)}$; then it is elementary that if $T$ and $U$ are finite sets, $l(T\cdot U) \leq l(T) + l(U)$.  We abbreviate $l^{-1} (\brcs{k})$ to $l^{-1}(k)$, and $l^{-1}\(\Set{r}{0 \leq r \leq s}\)$ to $l^{-1}(\leq s)$ (this will especially apply when we discuss real weight functions, in which case $r,s$ are nonnegative real numbers).

  If $S$ is an admissible set for infinite $G$,   we  define the weight function $l_S$ via $l_S (g) = k$ if $g \in S^k \setminus S^{k-1}$, and these weight functions are the most significant. They also have two additional properties, not generally satisfied by weight functions:
 \item{(iii)} for every (other) weight function $l'$ on $G$, there exists $c_0 \in \Z^+$ \st for all nonnegative  integers $a, m$, 
  $$ 
  l'(l^{-1}(\leq m+a)) \leq l' (l^{-1}(\leq m)) + a c_0.
  $$
  \item{(iv)} there exists $d \in \N$ \st for all $k \geq d$, the set $\cup_{i=0}^{d} l^{-1}(k-i)$ is nonempty. 
 
 Here $l^{-1}(\leq t)$ is an abuse of notation for $\Set{g \in G}{l(g) \leq t}$. To see that $l = l_S$ satisfies the additional property (iii), just note that $l^{-1}{(\leq n)} = S^n$, and 
 $$
  l'(S^{m+a}) \leq l'(S^m) + l'(S^a) \leq  l'(S^m) + a l' (S).
 $$ 
  And when $l = l_S$, we can choose $d = 0$ (since $G$ is infinite).

  We call a weight function satisfying (iii) and (iv) an {\it elaborated weight function.}
 
  If $l$ is any weight function (not necessarily elaborated), then we can define $\Arrow \tilde l; G.\Z \cup \brcs{-\infty}$ (unlike the original function $l$, it could conceivably have negative values) via 
 $$ 
 \tilde l(g) = \inf\Set{k}{\text{for almost all $m \in \N$, $l^{-1}(\leq m)\cdot g \subset l^{-1}(m+k)$}}.
$$
 In the case that $l = l_S$ for an admissible set $S$, the existence of one  $m$ \st $l^{-1}(\leq m)\cdot  g \subset l^{-1}(m+k)$ implies  it holds for all larger values of $m$ (this is not true for more general  weight functions). In particular, 
$$\eqalign{
\tilde l_S (g) &= \inf\Set{k}{\text{there exists $m$ \st $S^m g \subset S^{m+k} $}}.\cr
}$$
  It is easy to check that  $\tilde l$ satisfies (0), and if $l = l_S$, then $\tilde l_S$ is subadditive, that is, it satisfies (i). If $G$ is finite, then $\tilde l_S$ is identically $-\infty$. However, if $G$ is infinite, then $\tilde l_S$ is nonnegative-valued: $S^m g \subset S^n$ entails $|S^m| \leq |S^n|$; since $G$ is infinite and $S$ is admissible, the sequence $\dots \subset S^m \subset S^{m+1} \subset \dots$ is strictly increasing, and therefore $|S^m| \leq |S^n|$ entails $m \leq n$. Thus the $k$ in the definition of $\tilde l_S$ is nonnegative (and could be zero). 

In general, it is not true that $\tilde l_S$ is a weight function, that is, that $\tilde l_S^{-1}(k)$ is finite for all $k$. We say that $(G,S)$ satisfies the {\it WC property\/} if $\tilde l_S^{-1}(k)$ is finite for all $k$---in other words, that $\tilde l_S$ is a weight function. We will discuss WC  in more detail in the next few sections. 

The main reason for introducing $\tilde l_S$ is that for each $g\in G$, the element $x_g := [g, t]$ (in $A_f$) belongs to $R_f$ if and only if $t \geq \tilde l_S (g)$; so the definition is motivated by the direct limit structures of $A_f$ and $R_f$.

  If $H$ is a subgroup of $G$ and $S$ is an admissible set for the latter, then $l:= l_S|H$ is a weight function on $H$, and may under some circumstances be   elaborated.

\Lem Lemma \firstone. For any admissible set $S$ of the infinite group $G$, the set $G_0:= \tilde l_S^{-1}(0)= \Set{g \in G}{[ g,0]\in R_f}$ is a torsion subgroup of $G$.  

\Pf The element $g$ belongs to $\tilde l^{-1}(0)$ iff there exists $m$ \st $S^m g \subseteq S^m$. Finite cardinality of the sets implies $S^m g = S^m$; as $1 \in S \subset S^m$, we have $g^k \in S^m$ for all $k$, but finiteness implies $g$ has finite order.

If $g,h \in \tilde l_S^{-1}(0)$, then $\tilde l_S (gh) \leq \tilde l_S (g) + \tilde l_S (h) = 0$; so $gh \in \tilde l_S^{-1}(0)$. As $\tilde l_S^{-1}(0)$ is thus a submonoid of a group and all of its elements have finite order, it is also a group.
\qed

If WC holds, then obviously $\tilde l_S^{-1}(0)$ is a finite group (it need not be normal).  This suggests that infinite finitely generated   torsion groups might fail WC for every choice of admissible $S$.

Given admissible $f \in (AG)^+$ and $S = \supp f$, define $\Gamma_k = S^k \setminus S^{k-1}$ and 
$$
\Gamma_k'' =\Set{\gamma \in G}{[{}{\gamma},k] \in R_f, \text{ and $[{}{\gamma},k-1] \not\in R_f$} }.
$$
Alternatively, $\Gamma_k'' =  \tilde l_S^{-1}(k)$. Then $\supp f^k \subseteq \supp \cup_{i \leq k} \Gamma_i''$, but the reverse inequality fails in general. For example, by Lemma \firstone, $\Gamma_0''$ can contain more than just $1_G$, so $\Gamma_0$ need not equal $\Gamma_0''$. If $(G,S)$ satisfies WC, then each $\Gamma_k''$ is finite.
 Also define $\Gamma_k' = S^k \cap \tilde l_S^{-1}(k)$. This has the advantage over $\Gamma_k''$ that it is necessarily finite; it is better in another way as well, in terms of representing $R_f/\SS R_f$  as a direct limit of ordered groups, $\lim A \Gamma_{n}' \to A \Gamma_{n+1}'$. 

 Here is a simple example of admissible  $S$ inside  $G = \Z$; as usual, we regard the elements of $A \Z$ as Laurent polynomials in the variable $x$. Set $f = x^4 + x^3 + x + 1 + x^{-1}$ (we can make it symmetric by adding $x^{-4} + x^{-3}$ if desired, but this does not change the essential properties).
Then $\supp (x^2 f) = \brcs{6,5,3,2,1} \subset \brcs{8,7,6,5,4,3,2,1,0,-1,-2} = \supp f^2$. Hence $[x^2,1] \in R_f$, but $x^2 \not \in \Log f$. So $2 \not\in \Gamma_1$, but $2 \in \Gamma_1''$.

Recall $A = \Z$ or $\R$, and $f \in (AG)^+$ is admissible.

Following [H4], let $f \in A G$ be admissible, and partition $G$ as follows. Define as we have, $\Gamma_0 = \brcs{1_G}$, $\Gamma_1 = \supp f \setminus \Gamma_0$, and more generally, $\Gamma_n = \supp f^n \setminus \cup_{i < n} \supp f^i$. This means that $\Gamma_n$ consists of the states (group elements) that can be reached for the first time by $n$ iterations of the random walk coming from left multiplication by $f$. Since we assume that $1_G \in \supp f$, we have $\supp f^{n-1} \subset \supp f^n$, and since we also assume that $\supp f$ generates $G$ as a semigroup, we have $G = \dot\cup \Gamma_n$, and $\supp f^n = \dot\cup _{i\leq n} \Gamma_i$.

We then have $A G = \oplus A \Gamma_n$, and left multiplication by $f$ has what amounts to matrix representation \wrt this decomposition. Of importance for use in examining $R_f/\SS R_f$ is that we can throw away most of the big matrix.

For each $n$, let $\Arrow f_n; A \Gamma_n. A \Gamma_{n+1}$ by the map obtained by restricting and compressing left multiplication by $f$; that is, if $\gamma \in \Gamma_n$, then $[{}{\gamma}, n] \in R_f$, and we consider $f{}{\gamma}$; this has support in $\supp f^{n+1}$, and we discard all the elements of $\supp f^n$, yielding a member of $A \Gamma_{n+1}^+$. That the map is well defined is routine, and we then observe that there is a natural map $R_f/\SS R_f \to \lim \Arrow f_n; A \Gamma_n. A \Gamma_{n+1}$ which is an order isomorphism. This is explained in [H4] under more general circumstances; the resulting dimension group is called the {\it future dimension group.} (A minor difference---in the reference cited, the emphasis is on real rather than integer coefficients.) We are interested in the simple quotients of the future dimension group (by order ideals).

\noindent{\it Successors and predecessors.} Using $\Gamma_k = S^k \setminus S^{-1}$  in the direct limit to describe $R_f/\SS R_f$ is not optimal, because of what are known as {\it dead ends\/} ([GrH]); that is, given $g \in \Gamma_k$, there need not exist $h \in \Gamma_1$ \st $hg \in \Gamma_{k+1}$, i.e., $g$ has no successors. On the other hand, every element of $\Gamma_k$ has a predecessor, almost by definition, that is, there exists $j \in \Gamma_{k-1}$ and $\gamma \in \Gamma_1$ \st $\gamma j = g$. So the transition matrices have no zero rows, but may have zero columns. 

 If instead,  with  $\Gamma_k'' = \tilde l_S^{-1}(k)$, then it is also true that $R_f/ \SS R_f$ can be obtained as an order direct limit, $\lim A\Gamma_k'' \to A\Gamma_{k+1}''$, but this time, the sets need not be finite (finiteness of the sets is precisely the condition WC), but even if they are finite, predecessors need not exist (although successors do). Fortunately, with  $\Gamma_k'$, we obtain both predecessors and successors (hence the Goldilocks reference in the title of this section), and $R_f/\SS R_f$ is naturally order isomorphic to the corresponding limit of the (semi-obvious) maps $A\Gamma_k' \to A \Gamma_{k+1}'$.

Just as we have natural maps $A \Gamma_k \to A \Gamma_{k+1}$ induced by $f$, we obtain corresponding maps $A\Gamma_k' \to A \Gamma_{k+1}'$ (and also with double primes); we will show that in fact the two direct limits are naturally isomorphic (and to $R_f/{\SS R_f}$). First, we recall the map $A \Gamma_k \to A \Gamma_{k+1}$. We   write $\supp f^k = \dot\cup_{i=0}^k \Gamma_{i}$, and thus $A [\supp f^k]$ (the set of $A$-valued functions on $\supp f^k$) is $\oplus_{i=0}^k A \Gamma_i$. Let $\pi_k$ be the projection from this onto $A \Gamma_k$. Then we can define the map $\Arrow F_k; A \Gamma_k . A \Gamma_{k+1}$ via $F_k = \pi_{k+1} (f \times) |A \Gamma_k$, where $f \times$ of course represents left multiplication by $f$ on $A G$.

We can also define the corresponding map $\Arrow F_k'; A \Gamma_k' . A \Gamma_{k+1}'$, first by sending ${}{\gamma} $ to $f {}{\gamma}$, and then projecting
onto $A \Gamma_{k+1}'$, that is, removing from the support of $f{}{\gamma}$ all group elements that do not belong to $\Gamma_{k+1}'$. We have that $[{}{\gamma}, k] \in R_f$; then consider $[f{}{\gamma},k+1] = \sum_{h \in \supp f} (f,{}h){}{hg}$. Each $[{}{hg},k+1] \in R_f$, and if $[{}{\gamma},j] \in R_f$ for some $j < k$ (that is, ${}{\gamma}\not \in \Gamma_{k}'$), then for every $h$ in $\supp f$, $[{}{hg} ,j+1] \in R_f$, and thus none of the $hg$ belong to $\Gamma_{k+1}'$. Hence the only elements $\gamma $ \st ${}{\gamma} \in \oplus_{j \leq k+1} A \Gamma_{j}'$ is sent to a nonzero element of $A \Gamma_{k+1}'$ by left multiplication by $f$ and subsequent projection must belong to $\Gamma_k'$.

It can however, happen that $F_k'({}{\gamma}) = 0$ for $\gamma \in \Gamma_k'$, and indeed this is what happens in the example over $\Z$  given above.

\Lem Lemma \firsttwo. The maps $F'_k$ yield a natural isomorphism $\lim_k \Arrow F'_k;A\Gamma_k'. {A\Gamma'_{k+1}} \iso R_f/\SS R_f $.

\Pf 
Given $\gamma \in \Gamma_k$, we send ${}{\gamma}$ (the characteristic function of $\brcs{\gamma}$) to zero in $A \Gamma_k'$ if $\gamma \not\in \Gamma_k'$, and to itself (that is, the characteristic function of the singleton set) viewed as an element of $A \Gamma_k'$ if $\gamma \in \Gamma_k'$. The kernel of this map $\Arrow \phi_k ; A \Gamma_k. A \Gamma_k'$ is spanned by elements ${}{\gamma}$ \st $\gamma \in \Gamma_k$ but $[{}{\gamma},j] \in R_f$ for some $j < k$. The latter entails that $(\supp f^M)\gamma \subset \supp f^{M+j} $ for all sufficiently large $M$. For each $h \in \supp f^m$, we have that $hg \in \supp f^{M+j}$. It follows that $F_{M+k-1}\cdot F_{M+ k-2}\cdot \cdots \cdot F_k ({}{\gamma}) = 0$. Hence $\ker \phi_k \subset \ker ( A \Gamma_k \to \lim_j A \Gamma_j)$.

We check that the diagram
$$ \CD 
A\Gamma_k  @>{F_k}>>  A{\Gamma_{k+1}}@>{F_{k+1}}>> \dots \\   
\phi_k  @VVV \phi_{k+1}   @V VV \\  
A\Gamma'_k  @>F'_k >>  A{\Gamma'_{k+1}} @>F'_{k+1} >> \dots	 \\ 
 \endCD   
$$ 
\comment
$$
\diagram
A\Gamma_k&\rTo^{F_k} & A{\Gamma_{k+1}} & \rTo^{F_{k+1}} &\dots  \\
\dTo^{\phi_k} &&\dTo^{\phi_{k+1}}&  \\
A\Gamma'_k&\rTo^{F'_k} & A{\Gamma'_{k+1}} & \rTo^{F'_{k+1}} &\dots  \\
\enddiagram
$$
\endcomment
commutes. For $\gamma \in \Gamma_k$, $F_k' \phi_k ({}{\gamma})$ is either zero (when $\gamma \not \in \Gamma_k'$) or $F_k' ({}{\gamma})$. In the former case, since $\gamma \in \Gamma_k$, we have $[{}{\gamma},k ] \in R_f$, so $[{}{\gamma},j] \in R_f$ for some $j < k$. Then $[f{}{\gamma}, j+1] \in R_f$, and thus $h\gamma \not \in \Gamma_{k+1}'$ for all $h \in \supp f$, and it follows that $\phi_{k+1}({}{h\gamma}) = 0$, so $\phi_{k+1}F_{k} ({}{\gamma}) = 0 = F_k' \phi_k({}{\gamma})$.

If $\gamma \in \Gamma_k'$, then $F_k' \phi_k ({}{\gamma})$ is just $F_k'({}{\gamma}) = \sum_{h \in T'} (f,{}h) {}{h\gamma}$ where $T' = \Set{h \in \supp f}{h\gamma \in \Gamma_{k+1}'}$, On the other hand, $F_k({}{\gamma}) = \sum_{h \in T} (f,{}h) {}{h\gamma} \in \Gamma_{k+1}$ (where $T= \Set{h \in \supp f}{h\gamma \in \Gamma_{k+1}}$), and this is mapped by $\phi_{k+1} $ to $ \sum_{h \in T \cap T'} (f,{}h) {}{h\gamma} \in \Gamma_{k+1}$. For $h \in \supp f$ and $g \in \Gamma_k$, we have $g \in \supp f^k$, so $h\gamma \in \supp f^{k+1}$; hence if in addition, $h \notin T$, then $h\gamma \not \in \Gamma_{k+1}$, so there must exist $j < k+1$ \st $h\gamma \in \supp f^j$. But this entails $[{}{h\gamma},j] \in R_f$, so $h\gamma \not\in \Gamma_{k+1}$, and thus $T \cap T' = T'$. So $\phi_{k+1}F_{k} ({}{\gamma}) = F_k' \phi_k({}{\gamma})$ in this case as well.

This yields an order-preserving map $\Arrow \Phi; \lim A \Gamma_k . \lim A \Gamma_k'$. We have seen that the kernel of the map is zero. We now show it is onto, and an order-isomorphism. To do both, it suffices to show that given $\gamma \in \Gamma_k'$, there exists an integer $M> 0$ and an element $p \in (A \Gamma_{M+k})^+$ \st $\phi_{M+k}(p) = F_{M+k-1}' \cdot F_{M+k-2}'\cdot \dots \cdot F_{k}' ({}{\gamma})$.

From $\gamma \in \Gamma_k'$, we have $[{}{\gamma}, k ] \in R_f$, so there exists $M$ \st $(\supp f^M)\gamma \subset \supp f^{M+k}$, and thus $h\gamma \in \supp f^{M+k}$ for all $h \in \supp f^M$. Define, as in the previous argument, $T = \Set{h \in \supp f^M}{h\gamma \in \Gamma_{M+k}}$. If $h \in \supp f^M \setminus T$, then $h\gamma \in \oplus_{i < M+k} \Gamma_i$, and thus $h\gamma \not\in \Gamma_{M+k}'$. It follows that if we set $p = \sum_T (f^M,{}h){}{h\gamma} \in A \Gamma_{M+k}$, then $\phi_{M+k}(p) = \sum_{T' \cap T} (f^M,{}h){}{h\gamma}$. As in the preceding argument, if $h \in T' \setminus T$, then $h\gamma \not \in T'$; thus $T' \subseteq T$, and so $\phi_{M+k}(p) = \sum_{T'} (f^M,{}h){}{h\gamma}$. On the other hand, it is easy to check that $F_{M+k-1}' \cdot F_{M+k-2}'\cdot \dots \cdot F_{k}' ({}{\gamma}) = \sum_{T'} (f^M,{}h){}{h\gamma}$, and we are done. \qed

A similar result holds with $\Gamma_k'' = \tilde l_S^{-1}(k)$ replacing $\Gamma_k'$. However, the former need not be finite, and moreover, there need not be predecessors. Using the original $\Gamma_k$, there are always predecessors, but there need not be successors ({\it dead ends\/} [GrH]). 
 The Goldilocks situation occurs with $\Gamma'_k: = S^k \cap \tilde l_S^{-1}(k)$. Predecessors and successors exist, as we now show. 

\Lem Lemma \firstthr. (Predecessors and successors)
\item{(i)} Given $\gamma \in \Gamma_k'$, there exists $h \in \supp f$ \st $h\gamma \in \Gamma_{k+1}'$.
\item{(ii)} If $\gamma \in \Gamma_k \cap (\cup_{i < k}\Gamma_i'')$, then for all $h \in \supp f$,
$h \gamma \in (\cup_{i < k+1}\Gamma_i')$.
\item{(iii)} For every $\delta \in \Gamma_{k+1}'$, there exists $\gamma \in \Gamma_k''$ and an $h \in \supp f \cap \Gamma_1'$ \st $\delta = h\gamma$.

\Rmk
Part (i) says that for the transition matrix of $F_k$ \wrt the obvious basis of $A\Gamma_k$, the columns corresponding to elements of $\Gamma_k'$ contain at least one nonzero entry; part (iii) says that the rows corresponding to elements of $\Gamma_{k+1}'$ contain at least one nonzero entry. Part (ii) says that the columns corresponding to elements of $\Gamma_k \setminus \Gamma_k''$ are identically zero. This means that we can discard the latter columns, and also discard the rows corresponding to the elements of $\Gamma_{k+1}\setminus \Gamma_{k+1}''$, and so obtain induced maps $A \Gamma_k' \to A \Gamma_{k+1}'$; these have the property that every column and every row contains a positive entry.

\Pf (i) For $\gamma \in \Gamma_k''$, we have $\gamma \in \supp f^k$ but $[{}{\gamma},j] \not\in R_f$ for all $0 \leq j < k$. For all $h \in \supp f$, we have $hg \in \supp f^{k+1}$. If for $h \in \supp f$, $hg \not\in \Gamma_k''$, then either $h\gamma \not\in \Gamma_{k+1}$ or $h\gamma \in \cup_{j < k} \Gamma_{j}'$. The former entails (since $h\gamma \in \supp f^{k+1} = \dot\cup_{j\leq k+1} \Gamma_j$) that $h\gamma \in \Gamma_j$, and thus $[{}{h\gamma},j] \in R_f$, whence $h\gamma \in \cup_{j < k} \Gamma_{j}'$ in any event.

Assume $h \gamma \in \cup_{j < k+1} \Gamma_j'$ for all $h \in \supp f$. Then $[{}{h\gamma},k] \in R_f$ for all $h$ (since for an arbitrary $h$, $[{}{hg},j] \in R_f$ for $j < k+1$ entails $[{}{hg},j'] \in R_f$ for all $j' \geq j$). Hence $[{}{\gamma,k-1}]= [f{}{\gamma},k] \in R_f$, so that $\gamma \in \Gamma_{k-1}'$, a contradiction.

\noindent (ii) We have $[{}{\gamma},j] \in R_f$ for some $j < k$, so that $[f {}{\gamma},j+1] \in R_f$, and thus $[{}{h\gamma}, j+1] \in R_f$ for all $h \in \supp f$, so $h \gamma \in (\cup_{i < k+1}\Gamma_i'') $.

\noindent (iii) Given $\delta \in \Gamma_{k+1} \supp f^{k+1}$, there exists $\gamma \in \supp f^k$ and $h \in \supp f$ \st $\delta = h\gamma$. If $\gamma \not\in \Gamma_k''$, there exists $j < k$ \st $[{}{\gamma},j] \in R_f$, and so $[{}{h\gamma},j+1] \in R_f$, contradicting $\delta = h\gamma \not\in \cup_{j< k}\Gamma_{j+1}'$. Hence $\gamma \in \Gamma_k'$. If $h \not\in \Gamma_1'$, then as $h \in \supp f = \Gamma_1 \cup \Gamma_0$, we must have $h \in \Gamma_0'$, and again we easily derive a contradiction.
\qed

Actually, (iii) proves a bit more: if a group element $\delta$ appears in $\Gamma_{k+1}'$, then whenever $h \in \supp f^j$ and $\gamma \in \supp f^{k+1-j}$ and $\delta = h\gamma$, then $h \in \Gamma_{j}'$ and $\gamma \in \Gamma_{k+1-j}'$. 

\Lem Lemma \firstfou. Assume $f$ is an admissible element of the infinite group $G$; set $G_0$ to be the torsion subgroup,  $ \tilde l_{S}^{-1}(0)$. Each $\tilde l_S^{-1}(k)$ is a right and left $G_0$-set, and  each is faithful. In particular, if $G_0$ is infinite, then so is  $\tilde l_S^{-1}(k)$ for all positive integers $k$.

\Pf Pick $g \in G$ and $h \in G_0$; then $h^{-1} \in G_0$ and $\tilde l(gh), \tilde l(hg) \leq \tilde l(g) + \tilde l(h) = \tilde l(g)$. Also, $\tilde l(g) \leq \tilde l(gh) + \tilde l(h^{-1}), \tilde l(h^{-1})+\tilde l(hg)$, yielding the reverse inequalities. So $\tilde l(gh) = \tilde l(hg) = \tilde l(g)$.

If $G$ is infinite, then $\Gamma_k''$ is nonempty for all $k$ (as $\Gamma_k'' \cap S^k$ is), and the action is clearly faithful. \qed

\SecT \WC\ Property WC

Let $G$ be a finitely generated  group and $S$ an admissible subset $G$. We recall that the pair $(G,S)$
 has the {\it WC property\/} if for
every nonnegative integer $l$, the set
$$
\Set{g \in G}{\exists \ m \in \Z^+ \text{ \st } S^m g \subset S^{m+l}}
$$
is finite. Alternatively, there exists an integer $m_0$ (depending on $l$) \st
if for some $m$, $S^m g \subset S^{m+k}$, then $S^{m_0}g \subset S^{m_0 + k}$. 
Equivalently, if $f $ is any element of $\R G$ with support equalling $S$, and
all coefficients nonnegative, then for each $k$, the set $\Set{g \in G}{[g,k]
\in R_f}$ is finite.   For an admissible set $S$, whether  $\tilde l_S$ is a weight function is of course equivalent to $(G,S)$ satisfying WC. 

If $(G,S)$  has the WC property for every admissible subset $S$ of $G$, then we say that $G$ {\it satisfies WC.}

Given  admissible $S$, we define the function $\Arrow \tilde l_S;G.\Z^+$ as
the smallest integer $k$ for which there exists an integer $m$
\st $S^m g \subset S^{m+k}$. It is easy to check that this notation is consistent, in the sense that applying the definition of $\tilde l$ above to $l = l_S$ yields $\tilde l_S$. 
So yet another characterization of WC is that
$\tilde l^{-1}(l)$ is finite for all $l \in \Z^+$. We note that for $g,h \in
\cup S^n$, $\tilde l(gh) \leq \tilde l(g) + \tilde l (h)$. Moreover, if $g \in
S^n$, then $\tilde l(g) \leq n$ (but the inequality is almost always strict).

\def\cvx{\text{cvx\,}}

For example, if $G$ is finite, then it is clear that $\tilde l(g) = -\infty$ for all $g \in G$, no matter what the choice of admissible $S$. If $G = \Z = \langle x\rangle$, then with $f = x^3 + x^2 + 1 + x^{-1}$ and $S = \supp f$, it follows  that $\tilde l(x) = 1$, even though $x \not\in \supp f$; more generally, if $f = (x^n + x^{-n})(1+x + x^{-1})$, then $\tilde l(x^i) = 1$ for all $|i| \leq n$ (look at the $n$th power of $f$; its string of coefficients has no gaps).

It is tempting to conjecture that if $j$ is a torsion element of $G$ and $G$ is infinite, then $\tilde l(j) = 0$. This does happen frequently, but not always; here is a simple class of examples.

\Lem Example \wcone. Let $G = \Z \times C_n$ (the direct product), where $\Z = \langle x\rangle$ and the cyclic group $C_n = \langle g \rangle$. Set $f = 1 + x + x^{-1} + g$ and $S = \supp f$. Then $S$ is admissible, and $\tilde l_S (g^s) = s$ for $1 \leq s \leq n-1$.

\Pf Regard $A G$ as the free $A[x,x^{-1}]$ module with basis $\brcs{1,g,g^2, \dots, g^{n-1}}$. Form $f^m$, and consider the monomials in $x^{\pm1}$ that appear in the coefficients of $g^t$ in the decomposition of $f^m$. When $t = 0$, we obtain $\Set{x^i}{|i| \leq m}$, when $t = 1$, the corresponding monomials are $\Set{x^i}{|i| \leq m-1}$, and in general we see that for $t \leq n-1$, we have the coefficient of $g^t$ contains precisely the monomials $\Set{x^i}{|i| \leq m-t }$.

Since $g^s \in \supp f^s$, we have $\tilde l(g^s) \leq s$. If $\tilde l(g^s) < s$, then there would exist $m$ \st $\supp f^m g^s \subset \supp f^{m+s-1}$. Look at the coefficient of $g^s$ in both products; for the left side, it is the coefficient of $g^0 = 1$ in $f^m$, which has monomials $\Set{x^i}{|i|\leq m}$; for the right side, by the preceding paragraph (with $m$ replaced by $m+s-1$ and $s = t$), the corresponding set of monomials is $\Set{x^i}{|i| \leq m+s-1 - s}$. Hence the terms $x^{\pm m}g^s$ appear in $f^m g^s$ but not in $f^{m+s-1}$, yielding a contradiction.
\qed

We will show that abelian by finite groups and nontrivial amalgamated products satisfy WC, and Elisabeth Fink has shown  that non-elementary hyperbolic groups satisfy an even stronger property. However, for at least one choice of (standard) admissible $S$, the discrete Heisenberg group $(H_3,S)$ fails to satisfy WC. 

Let $f$ be an admissible element of $AG$, and set $S = \supp f$. The motivation for considering $\tilde l_S$ lies in the fact that for $g \in G$, we have $[g, k] \in R_f$ if and only if $k \geq \tilde l_S$. (Of course,  $l_S$ and $\tilde l_S$  do not depend on the choice of nonzero coefficients appearing in  $f$.)

\Lem (feeble) Lemma \wctwo. Let $K$ be a finite normal subgroup of $G$ \st $G/K$
satisfies WC. Then $G$ satisfies WC.

\Pf Let $S$ be an admissible set in $G$, and let $\pi$ be the quotient map. Fix $k$.
Then there are finitely many $g \in K$ \st $\pi(S)^m g \subset \pi(S)^{m+k}$ for some
$m$, and $\pi^{-1}\brcs{g}$ is finite.
\qed

\Lem Lemma \wcthr. If $G = \Z^d$ and $S$ is an admissible subset, then $|\tilde l_S^{-1}(\leq t)| \leq |tK \cap \Z^d|$ where $K = \cvx S$ inside $\R^d$, and this is sharp.

\Pf  Here we use additive notation.
 Let
$S$ be an admissible subset of $G$. Regard $\Z^d $ as a sublattice of $\R^d$.
Form the compact convex polytope $K = \cvx S$. There exist finitely many linear
functionals $\alpha_i $ on $\R^d$ \st $K = \cap \alpha_i^{-1}(\geq \beta_i)$ for a
corresponding set of real numbers $\beta_i$.

Fix $l \in \Z^+$ and suppose that $mS + v\subset (m+l)\subset (m+l)S$ (as usual,
$tS$ means the set of sums of $t$ elements of $S$). Then $mK + v \subset
(m+l)K$. Pick $\alpha_i$, and $z$ (depending on $i$) in $K$ \st $\alpha_i(z) =
\beta_i$ (the minimum possible; the set of such $z$ is a face of $K$). Then $mz
+ v \in (m+l)K$ entails $\alpha_i(v) \geq (m+l)\beta_i - m \beta_i$, so
$\alpha_i(v) \geq l \beta_i$. Applying this with every $\alpha_i$, we obtain $v
\in lK$, so $v \in lK \cap \Z^d$, which of course is finite. In particular,
$|\tilde l^{-1}(l)| \leq |lK \cap \Z^d|$, and this is easily
shown to be sharp whenever $K = sK'$ with $s \geq d-1$, where $K'$ is another
lattice polytope.
\qed

 Thus for torsion-free abelian groups,  $|\tilde l_S^{-1}(\leq t)| \sim |l_S^{-1}(\leq t)| = c t^d + \Oh{t^{d-1}}$.

\Lem Corollary \wcfou. Finitely generated abelian groups satisfy WC.

We have a bit more. 
\Lem Proposition \wcfiv. Abelian by finite groups satisfy WC. 

The proof will require slightly more convex geometry, and will not yield such a precise estimate for $|\tilde l_S^{-1} (\leq t)|$ as in the statement of Lemma \wcthr.  We begin with the following well known observation. 

\comment
 \Lem Proposition. Let $(G = \Z^d,S)$ be a finitely generated free abelian  group with admissible set, and let $C = \cvx S \subset \R^d$. Then $\tild{}S^{-1}(\leq k) \subseteq kC \cap \Z^d$. 

 In fact, $|\tilde l_S^{-1}(\leq k)| = c k^d$ for some constant $c$, but this takes a little more work.
\endcomment

 \Lem Lemma \wcsix. Let $A,B$ be compact convex polytopes in $\R^d$. For $v \in \R^d$, if $A +v \subseteq A + B$, then $v \in B$. 

 \Pf By using normals to the facets of $A$ and $B$, we can find disjoint index sets $\brcs{\alpha}$ and $\brcs{\beta}$, linear $f_{\alpha}$, $g_{\beta}$, and real numbers $r_{\alpha}$ and $s_{\beta}$ \st $A = \cap f_{\alpha}^{-1}(\leq r_{\alpha})$ and $B = \cap g_{\beta} ^{-1} (\leq s_{\beta})$. Dropping the requirement that the linear functions expose facets, we can combine the two sets, creating $\brcs{f_{\gamma}}$ \st $A = \cap f_{\gamma}^{-1 }(\leq a_{\gamma})$ and $B = \cap f_{\gamma}^{-1 }(\leq b_{\gamma})$ for some real numbers $\brcs{a_{\gamma}, b_{\gamma}}$; we can also assume that $\sup f_{\gamma}|A = a_{\alpha}$ and $\sup f_{\gamma}|B = b_{\gamma}$. 

 Suppose $v \notin B$. Then there exists $\gamma$ \st $f_{\gamma}(v) > b_{\gamma}$.  Choose $w \in A$  \st $f_{\gamma} (w) = a_{\gamma}$. Thus $f_{\gamma} (w + v) > a_{\gamma} + b_{\gamma}$, contradicting $\sup f_{\gamma} |(A+ B) \leq \sup f_{\gamma}|A + \sup f_{\gamma}|B = a_{\gamma} + b_{\gamma}$. \qed 

 \Lem Lemma \wcsev. Suppose that $H \triangleleft G$ and $K:= G/H$ is finite. Let $S$ be an admissible subset of $G$. Suppose that for all nonnegative integers $k$, we have $|\tilde l_S^{-1} (\leq k) \cap H|<\infty$. Then $(G,S)$ satisfies WC, and $|\tilde l_S^{-1}(\leq t)| \leq |K|\cdot |\tilde l_S^{-1} (\leq t) \cap H|$. 

 \Pf Let $\Cal K$ be a (fixed) transversal of $K =G/H$ containing $1$ (that is, the representative of $1 \in K$ is $1 \in G$). Then every element of $G$ is uniquely expressible in the form $\kappa h$ where $\kappa \in \Cal K$ and $h \in H$. 

 Suppose that $\tilde l_S (\kappa h) \leq t$. Then there exists $N$ \st for all $n \geq N$, we have $S^n \kappa h \subset S^{n+t}$. There exists an integer $m$ \st $\Cal K^{-1} \subset S^m$. Pick $n > \max\brcs{m, N}$. Then
 $$
 S^{n-m} \subset S^{n-m}S^{m}\kappa = S^n \kappa. 
$$
 Hence $S^{n-m}h \subset S^n \kappa h \subset S^{n+t}$; thus, $\tilde l_S (h) \leq t$. By hypothesis, the set of such $h$ is finite; since there are only finitely many choices for $\kappa \in \Cal K$, $\tilde l_S^{-1} (\leq t)$ is finite, and in fact, $|\tilde l_S^{-1}(\leq t)| \leq |K|\cdot |\tilde l_S^{-1} (\leq t) \cap H|$. \qed 

The following are standard definitions of characteristic subgroups of any group $G$ (see [P; pp~ 115 \& 117])
$$\eqalign{
 \Delta \equiv \Delta(G) & = \Set{g \in G}{g \text{ has only finitely many conjugates}}\cr
 \Delta^+ & = \Set{g \in \Delta}{g \text{ has finite order.}} 
}$$

 \noindent {\it Proof of Proposition \wcfiv.} Let $H$ be an abelian normal subgroup of finite index in $G$; since $G$ is finitely generated and $H$ is of finite index, $H$ is finitely generated. Since $K= G/H$ is finite, $H\subset \Delta (G)$. Then $\Delta^+(G)$ is a finite normal subgroup of $G$, so by Lemma \wctwo, we can factor it out, and thus assume that $\Delta^+(G)$ is trivial. Hence $\Delta$ is a torsion-free abelian group of finite index. Thus we may assume $H = \Delta $. 

 Let $S $ be an admissible subset of $G$. Identifying $H = \Z^d$, we regard elements of $H$ as elements of $\R^d$ and use additive notation for elements thereof. We will verify the finiteness criterion of the previous lemma. 

 Fix a transversal $\brcs{\kappa_k} \subset G$ of $K$ with $\kappa_{1} = 1_G$ (the notation is designed so that $\pi(\kappa_k) = k$).   Define  $\Arrow c; K \times K. H$ via $\kappa_k \kappa_{k'} =  \kappa_{kk'}c(k,k')$.
 
 For each integer $n$, chop $S^n$  into pieces via $S^n = \dot\cup_{k\in K} \kappa_{k} S_{n,k}$ where $S_{n,k} \subset H$ (this is well-defined). It is routine to check that 
 $$
 S_{n+1,k} = \bigcup_{jj' = k} \(S_{1,j}^{j'} +c(j,j')+ S_{n,j'}\).
$$
 Here superscript $j'$ indicates the effect of the corresponding automorphism of $H$. Iterating this, we see that $S_{n+t,k}$ can be expressed as a union of finite sets, each of the form $T_{j,j',k,t} + S_{n,j'}$, the $T_{j,j',k,t}$ not depending on $n$. Define $C_t$ to be convex hull of the union of $T_{j',k,t}$ over all possible $j,j',k$ without the constraint that $jj' =k$. Then $C_t$ is compact and convex, and for all $k$, $\cup_k S_{n+t,k} \subset \cup_k S_{n,k} + C_t$.

 Pick $h \in H$ with $\tilde l_S(h) \leq t$. Then $S^n h \subset S^{n+t}$ for all sufficiently large $n$. This translates to (changing to additive notation) to $S_{n,k} + h \subset S_{n+t,k}$ for all $k \in K$.  Set $A = \cvx \(\cup_k S_{n,k}\)$, so that $A + h \subset \cvx \(\cup_k S_{n+t,k} \) $. The latter is contained in $\cvx \( \cup_k S_{n,k} + C_t\)$, which is $A + C_t$. Thus $h \in C_t$, so that $|\tilde l_S^{-1}(\leq t) \cap H| \leq |C_t \cap H|$, and the latter is finite, since $C_t$ is compact. \qed

For $g \in G$ (and the admissible set $S$ fixed), recall the notation, $x_g = [g, \tilde l(g)] \in R_f$; angle brackets  indicate order ideals; thus if $T \subset A_f^+$, then $\langle T \rangle $ denotes the smallest order ideal containing $T$, that is, 

$$\Set{x \in A_f}{\exists N \in \N \text{ and a finite subset $T'$ \st } -N\sum_{t \in T'} t \leq x \leq N\sum_{t\in T'}t }.
$$
If $T$ does not consist of positive elements in $A_f$, there is no guarantee that there is a smallest order ideal containing $T$. 
In case $T$ is a singleton, $\brcs{t}$, we write $\langle t \rangle$. 
As $A_f$ is a dimension group, finite sums and intersections of order ideals are order ideals; also, $R_f = \langle \1\rangle$, and any order ideal of $R_f$ is also an order ideal of $A_f$. 

\Lem Lemma \wceig. Let $S$ be an admissible set in the group $G$. Let $g,h \in G$, and set $l = l_S$.
\item{(i)} There is an inclusion of order ideals, $\langle x_g\rangle \subset \langle x_h\rangle$ iff $\tilde l (gh^{-1}) = \tilde l(g) - \tilde l(h)$;
\item{(ii)} $\langle x_g\rangle =\langle x_h\rangle$ iff $\tilde l(gh^{-1})= 0$.

\Rmk As a consequence of (i), the condition implies $ \tilde l(g) \geq \tilde l(h)$, and in (ii), the conclusion implies $\tilde l(g) = \tilde l(h)$.

\Pf We can assume that $G$ is infinite.
\noindent (i) If the inclusion holds, then there exists an integer $m$ \st $S^{m+ \tilde l(h)}g \subset S^{m+ \tilde l(g)}h $, whence $S^{m+ \tilde l(h)}gh^{-1} \subset S^{m+ \tilde l(g)} $. Since $G$ is infinite, the function $n \mapsto |S^n|$ is strictly increasing, hence $ \tilde l(g) \geq \tilde l(h)$. Thus $S^{m+ \tilde l(h)}gh^{-1} \subset S^{m+ \tilde l(h) + (\tilde l(g) - \tilde l(h))}$. Therefore $\tilde l(gh^{-1}) \leq \tilde l(g) - \tilde l(h)$. Since $\tilde l(g) \leq \tilde l(gh^{-1}) + \tilde l (h)$, we deduce the reverse inequality.

Conversely, suppose $\tilde l (gh^{-1}) = \tilde l(g) - \tilde l(h)$. For all sufficiently large integers $m$, $S^m gh^{-1} \subset S^{m+\tilde l(g) - \tilde l(h)}$; hence $S^m g \subset S^{m+\tilde l(g) - \tilde l(h)}h$. Increasing $m$, we deduce $S^{m+ \tilde l(h)} g \subset S^{m+\tilde l(g) }h$. Thus $[g,-\tilde l(h)] \prec [{}h, -\tilde l(g)]$ in $A_f$, so that (on applying $\SS^{\tilde l(g) + \tilde l(h)}$) $[g, \tilde l(g)] \prec [{}h, \tilde l(h)]$.

\noindent (ii) If equality of the order ideals hold, then (i) yields both $\tilde l (gh^{-1}) = \tilde l (g) - \tilde l(h)$ and $\tilde l(hg^{-1}) = \tilde l (h) - \tilde l (g) = - \tilde l (gh^{-1})$. Since the values of $\tilde l$ are nonnegative, this forces both ends to be zero.

The converse is straightforward.
\qed

\Lem Corollary \wcnin. Assume $G$ is infinite. For all $g \in G$, the element $x_g = [g, \tilde l_S (g)] \in R_f$ does not belong to $\SS R_f = \langle [1,1]\rangle$.

\Pf If $x_g \in \SS R_f$, then we could write $[g, \tilde l_S(g)] = [a,k+1]$ where $a \in \R G$ and $[a,k] \in R_f$. Since $[a,k+1]$ is in $R_f^+$, we can assume that $a$ is already in $\R G^+$ (increasing $k$ as necessary). Then there exists $m$ \st $f^{m+k+1} g \prec f^{m+\tilde l (g)} a \prec f^{m+\tilde l(g) + k}$, yielding $\tilde l(g) \leq \tilde l(g) -1$, a contradiction.
\qed

 \comment
  \Lem Lemma. Suppose that the infinite finitely generated group $G$ satisfies SWC.  
  Let $l$ be an elaborated weight function on $G$. Then $\tilde l \geq l'/c_0 - C'$ for some positive constants $c_0$ and $C'$, and  $\tilde l^{-1} (k)$ is finite for all integers $k$. In particular, $G$ satisfies WC. 
 
  \Pf Fix $k$,  and let  $g$ be any element of $G$  \st $\tilde l(g) = k$. Then for all sufficiently large $M$, there exists $s \in l^{-1} (\leq M)$ \st $l(sg) = k + M$. Set $r = l(W)$ (that is, $\max_{w \in W} l(w)$), and let $d$ be the integer arising in property (iv). Let $M$ be any integer exceeding $r+ d$. We may  find an element $h \in l^{-1}(\leq M-r-d)$ \st $l'(h) = l'(l^{-1}(\leq M-r-d))$.  There exists $w \in W$ \st $l'(hwg) \geq l'(h)+l'(g) - C$. Since $l(hw) \leq  l(h)+ l(w) \leq M$,   for sufficiently large $M$, we have $l(hwg) \leq M+k$. Thus $hw \in l^{-1}(\leq M)$ and $hwg \in l^{-1} (M + k)$. 
 
  We have 
 $$\eqalign{
  l'(l^{-1}(\leq M + k)) & \geq l'(hwg) \geq l'(h) + l'(g)- C \cr
 & = l'(l^{-1}(\leq M-r-d)) + l'(g) - C; \quad\text{hence, by property (iii),} \cr
 l'(g) -C &\leq l'(l^{-1}(\leq M+k)) - l'(l^{-1}(\leq M-r-d))  \leq  (d+r + k )c_0\quad \text{so}\cr
  k c_0 &\geq l'(g) - C - (r+d)c_0, \quad\text{and thus} \cr
\tilde l &\geq \frac {l'}{c_0}- \frac C{c_0} -(r+d).\cr
 }$$

 Thus $l^{-1}(k) \subset \cup_{t = 0}^{c_0 (r+d) + c} (l')^{-1}(t)$, so is finite. In particular, if $l = l_S$, then $\tilde l_S$ is a weight function.  \qed
\endcomment

We wish to provide a mass of examples of groups satisfying WC  and stronger properties. Let $\alpha $ be a nonnegative real number. We say that the countable discrete group $G$ {\it satisfies\/} SWC($\alpha$)  if there exists a weight function $l$ on $G$, a finite subset $W \subset G$, and a nonnegative real number $C$ \st for all $g,g' \in G$, 
$$
l(gWg') \geq l(g) + \alpha l(g') - C
$$
(recall that if $T$ is a finite subset, then $l(T)$ is defined to be $\max_{t \in T} l(t)$). When this is the case, we say that $(l,W,C)$ {\it implements\/} SWC($\alpha$). If additionally, we can choose $l$ to be of the form $l_S$ for some admissible subset $S$, then $G$ satisfies SSWC($\alpha$), and this is implemented by $(S,W,C)$. If $\alpha =1$ (the largest possible value), then we use the notation SWC and SSWC. These are the strongest conditions in this family, although I do not know of an example which satisfies SWC($\alpha$) for some $\alpha > 0$, but not SWC. 

We can also allow $\alpha = 0$; then the semidirect product $\Z^2 \times_{\theta}\Z_2$, where $\theta$ is multiplication by $-1$ satisfies SWC(0), but not SWC($\alpha$) for any $\alpha > 0$. (In contrast, if we take $\theta$ to be either of the other two (up to conjugacy) nontrivial representations of $\Z_2$ on $\Z^2$, then the crossed product does not even satisfy SWC(0), for trivial reasons.) Even the weakest condition has a consequence: if $G$ satisfies SWC(0), then the centre of $G$ is finite. [If $z$ is in the centre, then $l(W) = l(zWz^{-1}) \geq l(z) -C $, hence $l$ is bounded on the centre.]

\Lem Lemma \wcten. Suppose that $G$ satisfies SWC($\alpha$) for some $\alpha > 0$, and this is implemented by $(l,W,C)$.
\item{(a)} Let $S$ be an admissible subset of $G$. Then there exists $C' \geq 0$  \st $\tilde l_S \geq \alpha l/l(S) - C'$; in particular, $G$ satisfies WC.
\item{(b)} Suppose that  $K$ is a finite group and $\Arrow 
\theta; K. \Aut G$ is a group homomorphism. Then there exists $\alpha' > 0$ \st the semidirect product $G \times_{\theta} K$  satisfies SWC($\alpha'$).  

\Rmk If we replace $l_S$ by an elaborated weight function $l'$, then the corresponding result in (a) is still true, $\tilde l' \geq \alpha l/l((l')^{-1}(\leq 1)) - C'$, with the proof modified in an elementary manner.  

\Rmk Part (b) will be improved to finite extensions in Proposition \Xthr; it is proved there for SWC = SWC(1), but the proof works for SWC($\alpha$) for any $\alpha > 0$. 

\Pf If $G$ is finite, there is nothing to do, so we may assume $G$ is infinite. We adopt the convention that for $v \in G$ and $k \in K$, $k v =v^k k$ where $v^k = \theta(k) v \theta(k^{-1})$.

\noindent (a) Suppose that $S$ is an admissible subset of $G$, and $\tilde l_S (g) = k$ for some $g \in G$; then $S^m g \subset S^{m+k}$ for some positive integer $m$; this remains true on increasing $m$ ad lib. There exists a positive integer $a$ \st $W \subset S^a$, as $S$ is admissible. We may assume that $m $ is larger than $a$. Find $h\in S^{m-a}$ \st $l(h) = l(S^{m-a})$. There exists $w \in W$ \st $l(hwg) \geq l(h) + \alpha l(g) - C$. Thus 
$$
\eqalign{
\alpha l(g) & \leq l(hWg) - l(h)   + C \cr
& \leq l(S^{m+k}) - l(S^{m-a})   + C; \text{ since $l(S^{m+k}) \leq l(S^{m-a}) + l(S^{a+k})$,} \cr
& \leq l(S^{a+ k}) + C \leq (a+ k) l(S) + C.\cr
}$$
Hence $l(g) \leq (a+k)l(S)/\alpha + C/\alpha$; since $\tilde l_S (g)  = k$, we obtain $l/l(S) \leq \tilde l_S/\alpha +  (a+C)/\alpha$, and thus $$\tilde l_S \geq  \frac {\alpha l}{l(S)} - C'$$
where $C' = a + C$. Since $l$ is a weight function, for any nonnegative integer $b$, $l^{-1}(\leq bl(S)/\alpha  + C')$ is finite, and thus so is $\tilde l_S ^{-1} (b) $. Hence $(G,S)$ satisfies WC for all admissible $S$, and so  $G$ satisfies WC. 

\noindent (b) Set $W_K := \cup_{k,k' \in K} W^k k'$ inside $G_0 = G \times_{\theta} K$. Define $\Arrow L;G.\Z$ via $L(v ) = \max_{k \in K} l(v^k)$; then $L$ is a maximum  of finitely many weight functions, and is thus itself a weight function. Moreover, $L$ is invariant under right or left action by $K$. It follows that the function on $G_0$, $l_K$, defined by 
$$
l_K (vk) = \cases   L(v) & \text{if $k = 1$}\\
 L(v) + 1 & \text{if $k \neq 1$}\\
\endcases$$
for $v \in G$ and $k \in K$ is a weight function (on $G_0$); this is well-defined, since the representation as $vk$ is unique. We claim that $(l_K, W_K, C +4)$ implements SWC($\alpha$). 
Pick $g = vk, g' = v'k' \in G_0$. 

 Pick $j \in K $ \st $l(v^j) = \max_{k \in K} l(v^k)$, and $j' \in K$ \st $l((v')^{j'}) = \max_{k \in K} l((v')^k)$.  Find $w_0 \in W$ \st $l(v^j w_0 (v')^{j'}) \geq l(v^j) + \alpha l((v')^{j'}) - C$; this last is just $L(v) + \alpha L(v') - C$.

 Now consider the product $g w_1 g' := vk \cdot (w_0^{k^{-1}j^{-1}} k^{-1}j^{-1} j')\cdot v'k'$. The parenthesized term, $w_1 = w_0^{k^{-1}j^{-1}}k^{-1}j^{-1} j'$, belongs to $W_K$, and the product simplifies 
$$\eqalign{g w_1 g' & =
 vk \cdot (w_0^{{jk}^{-1}}(jk)^{-1} j')\cdot v'k' \cr & = vk \cdot k^{-1}j^{-1} w_0 \cdot j' v' k\cr
 & = vj^{-1} w_0 (v')^{j'}j'k \cr
 &= j^{-1}\( v^j w_0 (v')^{j'}\) j'k.\cr
}$$
 Thus 
$$\eqalign{
l_K(g w_1 g') + 2 &\geq l_K (v^j w_0 (v')^{j'}) \cr
& \geq L(v) + \alpha L(v') - C\cr
 & \geq l_K(vk) + \alpha l_K(v'k') - C - 1- \alpha\cr
& = l_K(g) + \alpha l_K(g') - C - 1 - \alpha.\cr
}$$ So $l_K (gw_1 g') \geq l_K (g) + \alpha l_K(g') - C-4$.\qed

\Lem Proposition \wcele.   Free products   of nontrivial 
finitely generated groups satisfy SSWC.

\Rmk This includes $\Z_2 * \Z_2$, the infinite dihedral group. 

\Rmk  Elisabeth Fink has shown  that   nonelementary hyperbolic groups satisfy SWC, and with $|W| = 3$. 

\Pf Let $H$ and $K$ be nontrivial groups. We find a finite set $W$ and an admissible set $T \subset G = H*K$ \st $(T,l_T,W)$ implements SSWC. There
is a normal form for elements of the  free product, and we use this
to construct a weight function with the required properties.

Pick  admissible subsets $S_H$ and $S_K$ for $H$ and $K$ respectively; then $T:= S_H 
\cup S_K$ is an admissible subset for $G$. On a product with more than one letter,  $g = h_1 k_1 h_2 k_2 \dots$ (using the normal form),
define $l'(g) = \sum l_H(h_i) + \sum l_K(k_i)$. It is easy to check that $l' = l_T$. 

Select elements $h_{(1)} \in H$ and $k_{(1)}\in K$. Set $W$ to be
$\brcs{h_{(1)},k_{(1)} ,1}$. There are four possibilities for the ordered pairs consisting of the
terminal letter of $h$ and the initial letter of $g$; we can always pick an
element $w$  of $W$ so that there is no cancellation in $h w g$.
\qed

This also includes $D_{\infty} = \Z_2 * \Z_2$, for which WC is otherwise awkward
to prove directly.

Amalgamated free products apparently often result in SWC: if $G \cap H = L$, and $L$ satisfies SWC, then it seems likely that $G *_L H$ also does. (Note a necessary condition: if $L$ contains an infinite subgroup of the centre, then the amalgamated free product cannot be SWC, as the centre of an SWC group must be finite.)

\comment

\Lem Proposition \wctwe. If $G$ is free, then every finitely generated subgroup of $G$ satisfies WC.

\Pf Every finitely generated subgroup is free; if on one variable, it is
abelian, and we know that WC holds; if on more than one variable, it is a
nontrivial free product, hence the preceding applies.
\qed
\endcomment  

There are weakenings of SWC that are still sufficient to show the group satisfies WC.  We say a group $G$ satisfies {\it ssWC\/} if there exists a weight function $l$ on $G$ and a nonnegative constant $C$ \st for all $g$, there exists a finite subset $W(g) \subset G$ \st for all $h \in G$, 
$$
l(hW(g)g) \geq l(g) + l(h) - C.
$$
If additionally, we can choose $W(g)$ so that $\sup_{g \in G} |W(g)| < \infty$, then we say that $G$ satisfies {\it sWC}. 

Trivially, for a group $G$, $\text{SWC} \implies \text{sWC} \implies \text{ssWC}$.  Moreover, the argument of Lemma \wcten(a) (slightly modified) yields that $\text{ssWC} \implies \text{WC}$. 

Torsion-free finitely generated abelian groups satisfy sWC, with $|\sup W(g)| = 2$. To see this, here $ G =  \Z^n$; set $l(x) = \| x\|_{\infty}$, so that $l = l_S$ with $S$ consisting of all $n$-tuples each of whose entries belong to $\brcs{0,1}$. Let $g = v = (v(i))$, set $m =\sup_i |v(i)|$, and define $W(g) = \brcs{w, -w} $ where $w = 2m (\sgn v(i))$ (here the sign function takes values $\pm 1$). Set $C = 0$. To check that the definitions apply, let $h = (u(i))$, and suppose $j$ is such that $|u(j)| = \| h \| = \max |u(i)|$. If $\sgn u(j) = \sgn v(j)$, then the $j$th coordinate of $v + w + h$ is $v(i) + 2m \,\sgn v(i)  + u(i)$, and because of the sign condition, the absolute value of this is  $2m + |v(i)| + |u(i)| \geq l(v) + l(h)$. 

If instead, $\sgn u(i) = - \sgn v(i)$, then the $j$th coordinate of $v -w +h$
is $v(i) - 2m \, \sgn v(i) + u(i)$, so has absolute value at least $|u(i)| + m$, and this is $l(v) + l(h)$. 

This provides another way of proving finitely generated abelian groups satisfy WC, but does not give the quantitative estimate of Lemma \wcthr. Many of the results for SWC also apply to ssWC and sWC. But it is unclear how useful the latter two are. 

The property ssWC for groups is roughly analogous (or at least superficially similar) to that of strong primeness for rings  [HL], sWC is analogous to bounded strongly prime, and SWC is analogous to uniformly strongly prime.

\SecT \X\  SWC for abelian by finite groups

Since the infinite dihedral group $D_{\infty}$ is a   free product of nontrivial groups ($\Z_2 * \Z_2$), it thus satisfies SWC, and the proof of Proposition \wcele\  shows that the $W$ can be chosen to have three elements. But $D_{\infty}$ is also a semidirect product $\Z \times_{\theta} C_2$, where $\theta \in \Aut \Z$ is multiplication by $-1$; we see fairly quickly that there exists a two-element choice for $W$. 

This suggests the problem of determining when an   abelian by finite group satisfies SWC, which we will answer in this section. We also show that any $N$ is a normal (not necessarily abelian) subgroup of finite index in $G$, and $N$ satisfies SWC, then so does $G$. The converse of course fails, as shown by the infinite dihedral group, $D_{\infty}= \Z_2 * \Z_2$: it contains a copy of $\Z$ as a subgroup of index two. 

Suppose $H$ is a finitely generated  abelian group that is a normal subgroup of a group $G$ \st $K:= G/H$ is finite. We wish to decide when $G$ satisfies SWC. There is an immediate reduction to torsion-free $H$ (by Lemma \wctwo), so $H \iso \Z^d$ for some positive integer $d$.   Since $H$ is abelian, there is a natural group homomorphism $\Arrow \theta; K. \Aut (H) = \gl(d,\Z)$, that is, an integral representation of $K$. We may construct this map first by taking a cross-section $\Cal K = \brcs{E_k}_{k\in K}$ of $K$; that is, $\brcs{E_k H}_{k \in K}$ is a complete set of cosets of $H$ in $G$ \st the quotient map $G\to K$ sends $E_k \mapsto k$; we may also assume that $E_1$ is the identity of $G$. Then $\theta$ is given by $\theta (k)(x) = E_k x E_k^{-1}$; since $H$ is abelian, this is independent of the choice of representative $E_k$ of the coset $E_k H$, and it is easy to see that $\theta $ is a homomorphism. We have a few different notations for $\theta(k)(x)$, e.g., $x^k$ (where $\theta$ is understood). Sometimes we use Greek letters for elements of $K$, e.g., $\phi$.

We may tensor this with the rationals, creating $\Arrow \Theta:= \theta \otimes 1_{\Q}; K.\gl(d,\Q)$. We will also have to work with the corresponding real representations. 

 We say a  finite-dimensional rational  representation of a finite group is {\it trivplicity-free\/} if it is multiplicity-free (that is, no irreducible appears with multiplicity exceeding one) and the trivial representation does not appear. We show that trivplicity-freeness  of $\theta \otimes 1_{\Q}$ is necessary and sufficient  for finite extensions of $\Z^d$ to satisfy SWC. 

 \Lem Theorem \Xone. Let $H$ be a finitely generated torsion-free abelian group that is a normal subgroup of a group $G$, and suppose that $K:= G/H$ is finite. Let $\Arrow \theta; K.\gl(d,\Z)$ be the group homomorphism induced by $K$. Then $G$ satisfies SWC iff $\theta \otimes 1_{\Q}$ is trivplicity-free.

A consequence is that much of the extension data for $H \to G \to K$ is irrelevant (for example, whether the map splits); another consequence is that the criterion, instead of involving the  relatively subtle integral representation $\Arrow \theta; K .\gl(d,Z)$ only requires dealing with    the coarser rational representation $\Theta$. 

The proof actually yields a lot of structural information. If $\Theta$ is trivplicity free, we construct  implementations of SWC, $(l,\Cal K)$, where the $l$ is obtained geometrically from the dual action(s) of $\theta$ on $\Z^{1\times d}$ and $\R^{1\times d}$; the weight functions $l$ come from the normal vectors to the various possible polytopes obtained as the convex hull of orbits of $\theta$ on $\Z^{d\times 1}$ (these are closely related to the corresponding gauge functions of the convex hulls; see Appendix \Y).  The resulting weight functions satisfy the additional property, $l(h^n) = nl(h)$ for all positive integers $n$ and $h \in H$ (later in this section, the notation changes to additive when dealing with weight functions on abelian groups; thus, it will appear as $l(nh) = nl(h)$). 

In the converse direction, we also  show that if $l$ satisfies  $l(h^n) = nl(h)$ (when restricted to the normal abelian subgroup), then it arose from one of the constructions in the first part of the argument. 

In this section, $H$ will be a finitely generated torsion-free abelian group normal and of finite index in a group $G$; $N$ will be used for normal (but not necessarily abelian) subgroups of general groups. 

We now begin the proof that $\Theta$ trivplicity-free entails $G$ satisfies SWC, one direction of Theorem \Xone.

Let $\Arrow \lambda; G.\Z^+$ be a function. We say it is a {\it semi-weight function\/} [GrH] if 
\item{(a)} $\lambda^{-1}(t)$ is finite for every $t \in \Z^+$;
 \item{(b)} there exists $D\geq 0$ \st for all $g,g' \in G$, we have $\lambda(gg') \leq \lambda(g) + \lambda(g') + D$. 

\noindent There is no requirement that $\lambda (1) = 0$. 

For functions $\Arrow \lambda, \lambda'; G.\Z^+$, we define the usual equivalence relation, $\lambda \sim \lambda'$, if $\sup_{g \in G} |\lambda(g)- \lambda'(g)| < \infty$ (often written as $|| \lambda - \lambda' || < \infty$). Obviously, properties (a) and (b)  are preserved by this equivalence relation (with a possibly different choice of $D$), as is the property that is part of determining SWC,
\itemitem{($*$)} There exist a finite subset $W$ of $G$, and a positive real number $C$ \st for all $g,h \in G$, we have $\sup_{w \in W}\lambda (gwh) \geq \lambda (g) + \lambda (h) - C$,

\noindent although the constant $C$ might change (the $W$ remains the same). 

In dealing with group extensions $H \triangleleft G$ with $G/H$ finite, semi-weight functions arise naturally when we try to extend  weight functions on $H$ to $G$. Fortunately, there is a very simple result showing that any semi-weight function is equivalent to a weight function. 

\Lem Lemma \Xtwo. Let $\lambda$ be a semi-weight function on $G$. Then there is a weight function $\Arrow l;G.\Z^+$ \st $\lambda \sim l$. 

\Pf Let $D$ be the constant for $\lambda $ arising in the definition of semi-weight function. Pick a nonnegative integer $M$, and define $\Arrow \lambda_M;G.\Z^+$ via
$$
\lambda_M (g) = \cases 0 &\text{if $g = 1$}\\
\lambda(g) + M & \text{else.}
\endcases
$$
First, $\lambda_M^{-1}(\leq t) \subseteq \lambda^{-1}(\leq t-M) \cup \brcs{1}$, so the left side is finite. Now we verify that $\lambda_M$ is subadditive if $M \geq D$. Let $g,g' $ be elements of $G$. 

If $gg' = 1$, then $\lambda_M(gg') = 0 \leq \lambda_M(g) + \lambda_M(g')$ trivially. If $g = 1 \neq g'$, then $\lambda_M (gg') = \lambda_M (g') \leq \lambda_M (g)+\lambda_M (g')$, and similarly, if $g\neq 1 = g'$, subadditivity occurs.

So we may assume $g,g',gg'$ are all not the identity. Then 
$$\eqalign{
\lambda_M(gg') & = \lambda(gg') + M\cr
& \leq \lambda(g) + \lambda(g') + D + M \cr
& \leq \lambda(g) + M + \lambda(g')  + M = \lambda_M (g) + \lambda_M(g').
}$$
Thus $\lambda_M$ is a weight function. It is clear that $||\lambda_M - \lambda|| = \max \brcs{\lambda(1),M}$. Thus if we set $l = \lambda_D$, then $l$ is a weight function equivalent to $\lambda$.\qed

\Lem Proposition \Xthr. Let $N$ be a normal subgroup of $G$ \st $G/N$ is finite. If $N$ satisfies SWC, then so does $G$. 

\Pf Suppose that $(l_0, W_0, C)$ implements SWC for $N$. 
Let $\Cal K = \brcs{E_k}_{k \in K}$ be a set of coset representatives of $G $ modulo $N$; we may assume that $E_1 = 1$. For each $k \in K$ and $h \in N$, define $l_k (h) = l_0(E_k h E_{k}^{-1})$, and $l(h) = \sup_{k \in K} l_k (h)$. Then $l$ is a weight function on $N$. There exists a two-cocyle $\Arrow c; K \times K.N $ \st $E_k E_{k'} = c(k,k')E_{kk'}$ for all $k,k' \in K$. 

Now define $\Arrow\lambda;G.\Z^+$  (depending on the specific set of class representatives)
$$
\lambda (hE_k) = l(h)
$$
for $h\in N$. This is well-defined, and we show it is a semi-weight function on $G$. Obviously, $\lambda^{-1}(\leq t) \subset \cup_{k \in K} l^{-1}(\leq t)\cdot E_k$, so is finite. Now we verify $\lambda(gg') \leq \lambda(g) + \lambda(g') + D$ for suitable $D\geq  0$ and all $g,g' \in G$.

Set $X = \max_{j,j' \in K}\brcs{l  (c(j,j'),  l_0 (c(j,j')^{-1})}$. Write $g = hE_k$ and $g' = h'E_{k'}$ with $h,h' \in N$. Then $gg' = h \cdot (E_k h' E_k^{-1})\cdot E_k E_{k'}$, and this in turn expands as $h \cdot (E_k h' E_k^{-1})\cdot c(k,k') E_{kk'}$. Thus 
$$\eqalign{\lambda(gg') & = l(h \cdot (E_k h' E_k^{-1})\cdot c(k,k'))\cr
&\leq l(h) + l(E_k h' E_k^{-1}) + l(c(k,k'))  \cr
& \leq \lambda(hE_k) + l(E_k h' E_k^{-1})  + X. \cr
 }$$

There exists $j \in K$ \st $l(E_k h' E_k^{-1})) = l_0 (E_j E_k h' E_k^{-1}E_j^{-1})$, and we rewrite the last term in parentheses as $c(j,k)E_{jk}h' E_{jk}^{-1}c(j,k)^{-1}$. Thus 
$$\eqalign{
l(E_k h' E_k^{-1})& = l_0 (c(j,k)\cdot E_{jk}h' E_{jk}^{-1}\cdot c(j,k)^{-1}) \cr
&\leq l_0(E_{jk}h' E_{jk}^{-1}) + 2X \cr
& \leq l(h') + 2X = \lambda(h'E_{k'}) + 2X. 
}$$
Set $D = 3X$; we see that $\lambda$ is a semi-weight function.

\comment

If $k = 1$, then $\hat l (gg') = \hat l(hh' E_{k'})$. If additionally, $k' =1$, $\hat l(gg') = l(hh') \leq l(h) + l(h') = \hat l(g) + \hat l(g')$. If instead $k = 1 \neq k'$, then $\hat l(gg') = l(hh') + M \leq l(h) + l(h') + M = \hat l (g) + \hat l (g')$. Hence we can assume that $k \neq 1$. 

If $k' = 1 \neq k$, $gg' = hE_k h' = h \cdot E_k h' E_{k}^{-1} \cdot E_k$, so that $\hat l(gg') = l (h \cdot E_k h' E_k^{-1}) +M$, and this is then l

$$\eqalign{
 gg' &= h \cdot  E_{k}h' E_{k}^{-1} \cdot E_{k} E_{k'} \cr
&= h \cdot  E_{k}h' E_{k}^{-1} \cdot  c(k,k')E_{kk'}; \quad \text{thus,}\cr
\hat l (gg') & = \cases l(h\cdot E_{k}h' E_{k}^{-1} \cdot  c(k,k')) + M & \text{if $kk' \neq 1$} \\ l(h\cdot E_{k}h' E_{k}^{-1} \cdot  c(k,k^{-1})) & \text{if $kk'= 1$} \\
\endcases \cr
& \leq \cases l(h)+ l(E_{k}h' E_{k}^{-1}) +  l(c(k,k')) + M & \text{if $kk' \neq 1$} \\ l(h)+ l(E_{k}h' E_{k}^{-1}) +l  (c(k,k^{-1})) & \text{if $kk'= 1$}. \\
\endcases \quad  \text{We also have}\cr
l(E_{k}h' E_{k}^{-1}) &= \sup_{j \in K} l_0 (E_j E_k h' E_k^{-1}E_{j}^{-1}) \cr
& = \sup_{j \in K} l_0 (c(j,k)E_{jk} h' E_{jk}^{-1}c(j,k)^{-1}) \cr
& \leq \sup_{j \in K} l_0 (c(j,k))   +  l_0 (E_{jk} h' E_{jk}^{-1}) + l_0 (c(j,k)^{-1})\cr
& \leq l(h') + 2X. \quad \text{Hence }\cr
\hat l (gg') & \leq \cases l(h)+ l(h') +  3X + M & \text{if $kk' \neq 1$} \\ l(h)+ l(h') +3X & \text{if $kk'= 1$}. \\
\endcases 
}$$

\endcomment

 Define $W = \Cal K^{-1}\Cal K^{-1}W_0 \Cal K$; that is, 
$$
 W = \Set{E_{k(1)}^{-1} E_{k(2)}^{-1} w_0 E_{k(3)}}{k(i) \in K; w_0 \in W_0}.
 $$

 We will show that $\lambda$ satisfies ($*$) \wrt this choice of $W$. To this end, let $g = hE_k$ and $g' = h'E_{k'}$ be arbitrary elements of $G$. There exist $j,j' \in K$ \st both $l(h) = l_j (h)$ and $l(h') = l_{j'}(h')$ with $h,h' \in N$. For $v \in G$, define 
$$\eqalign{
 x &= E_j g v g' \cr & = E_j h E_j^{-1} \cdot E_j E_k v E_{j'}^{-1}\cdot E_{j'}h' E_{j'}^{-1}\cdot E_{j'}E_{k'}.\cr
}$$

 We have $\lambda (x) \leq \lambda  (E_j) + \lambda (gvg') + D$, and thus 
$$
\lambda (gvg')\geq \lambda(x) - \lambda (E_j) - D.\tag 1
$$

 We may find $w_0 \in W$ \st 
$$\eqalign{
 l_0 (E_j h E_{j}^{-1}\cdot  w_0 \cdot E_{j'}h' E_{j'}^{-1}) & \geq l_0 (E_j h E_{j}^{-1}) + l_0 (E_{j'} h E_{j}^{-1}) - C \cr 
 & = l(h) + l(h')- C.
}\tag 2$$ Solve for $w_0 = E_j E_k v E_{j'}^{-1}$, that is, set $v = E_k^{-1}E_j^{-1}w_0 E_{j'} \in W_0$. Let $L = \lambda(\Cal K^{-1})$; obviously $\lambda(\Cal K) = 0$. 
$$
\lambda(x) \geq \lambda (x E_{k'}^{-1}E_{l'}^{-1}) - \lambda (x E_{k'}^{-1}) - \lambda(x E_{j'}^{-1}) - 2D. \tag 3
 $$
 Finally, 
 $$\eqalign{
\lambda (g v g')& \geq  \lambda(x) - \lambda(E_j) - D \quad\text{by (1)}\cr 
 & \geq \lambda(x E_k^{-1}E_{j'}^{-1}) - \lambda(E_j)- \lambda(E_{k'}^{-1})- \lambda(E_{j'}^{-1}) - 3D \quad \text{by (3)}\cr 
 & = \lambda(E_j h E_j^{-1} \cdot w_0\cdot E_{j'}h' E_{j'}^{-1}) - \lambda(E_j)- \lambda(E_{k'}^{-1})- \lambda(E_{j'}^{-1}) - 3D\cr 
 & \geq  l(E_j h E_j^{-1} \cdot w_0\cdot E_{j'}h' E_{j'}^{-1}) - 2L - 3D\cr 
 & \geq  l_0(E_j h E_j^{-1} \cdot w_0\cdot E_{j'}h' E_{j'}^{-1}) - 2L - 3D\cr 
 & \geq l(h) + l(h')- C - 2L - 3D \quad \text{by (2)}\cr
 & = \lambda (hE_k) + \lambda (h'E_{k'})- C - 2L - 3D. \cr 
 }$$
Hence with $C' = C + 2L + 3D$, $\lambda $ satisfies ($*$). By Lemma \Xtwo, there exists a weight function $l'$ on $G$ with $l' \sim \lambda$, and thus $(l',W)$ implements SWC on $G$.\qed

In the course of the proof, we constructed $W$ so that $|W| \leq |K|^3 \cdot |W_0|$. There should exist a construction so that $|W| \leq |K|\cdot |W_0|$. If the extension splits, that is, $G$ is a semidirect product of $N$ by $K$, then $\lambda$ is already a weight function. 

The converse fails (the infinite dihedral group is a counter-example); but it is still plausible that for $N \triangleleft G$ and of finite index and $G$ satisfying (S)WC, then $N$ satisfies WC.

\Lem Lemma \Xfou. Suppose that  $N$ is a subgroup of finite index in the finitely generated group $G$ and $N$ satisfies WC. Let $S$ be an admissible subset of  $G$ \st for $l_0:= l_S|N$, $\tilde l_0^{-1} (\leq k)$ is finite for all $k$. Then $\tilde l_S^{-1} (k)$ is finite for all $k$. 

\Pf A subgroup of finite index in a finitely generated group is itself finitely generated ([P; Lemma 1.7]), so $l_0$ is defined. There exists an integer $r$ \st $S^r$ contains a complete set of coset representatives and their inverses, that is, $\cup_{i=1}^n g_i N = G$ and $\brcs{g_i^{\pm1}} \subset S^r$. Suppose that $\tilde l_S^{-1}(k)$ is infinite for some nonnegative integer $k$. Then there exists an integer $j \in \brcs{1,2,\dots, n}$ \st $U:=  \Set{g \in \tilde l_S^{-1}(k)}{gN=  g_jN }$ is infinite. For each $g \in U$, define $h_g = g_j^{-1}g \in N$. Then $\tilde l(h_g) \leq r + k$.

Suppose that $h \in N$ and $\max_{l_0(s) \leq M} l(sh) - M = t$. Then $\max_{l_S(s) \leq M+ r} l(sh) - M \geq \max_{l_0(s) \leq M}l(g_i s h) - (M +r) \geq \max_{l_0(s) \leq M}l( s h) - (M +2r) $. Hence $\tilde l(h) \geq \tilde l_0(h) -2r $. Letting $g$ vary over $U$, there exists $h_g \in N$ \st $\tilde l_0(h_g)$ is arbitrarily large (since $U$ is infinite), but $\tilde l_0(h_g) \leq \tilde l(h_g)  +2r$ is bounded above, a contradiction. \qed

 The following are all either well-known or easy; see also Appendix  \Y.

 \Lem Lemma \Xfiv. Let $v_i$ ($i=1,2,\dots, m$) be nonzero elements of $\R^{1\times d}$. Define 
$$Y = \Set{w\in \R^{d\times 1}}{v_i w \leq 1 \text{ for all $i$}}.$$
 \item{(a)} Then $Y $ is a closed convex subset of $\R^{d\times 1}$ containing an open ball centred at the origin.
 \item{(b)} If $\brcs{v_i}$ contains a basis for $\R^{1 \times d}$, and there is a relation of the form $\sum_{i=1}^m \alpha_i v_i = 0 $ with all $\alpha_i > 0$, then $Y$ is compact. {\par }
\noindent From now on, assume the hypothesis of (b) holds.
 \item{(c)} If $v_i \in \Q^{1\times d}$, then the extreme points of $Y$ have only rational entries. 
 \item{(d)} Given a facet $F$ of $Y$, there exists $i$ \st $F $ is exposed by $v_i$ (that is, $F = Y \cap v_i^{-1}(1)$). 
\item{(e)} For the facet $F$ exposed by $v_i \equiv v_F$, define $C_F $ to be the convex hull of $\brcs{0,F}$. Then  $\cup_F C_F = Y$ and for any $x \in C_F$, $\max_j v_j x =v_F x$. 

 \Pf (a) is trivial. (b) It suffices to show $Y$ is bounded. Let $E_j$ be the standard basis elements of $\R^{1\times d}$. By relabelling, we may assume $\brcs{v_i}_{i=1}^d$ is a basis of $\R^{1\times d}$. For each $j$, we can write $E_j = \sum_i \alpha_{ij}v_i$ with $\alpha_{ij}$. For each $\alpha_{ij} < 0$, we may replace $v_i$ by $-\sum_{k\neq i}^m \alpha_k v_k$; we thus obtain $E_j$ is a nonnegative linear combination of all the $\brcs{v_i}$. We can also obtain $E_j$ as a non-positive linear combination of $\brcs{v_i}$, by replacing those $v_i$ for which $\alpha_{ij} > 0$ by $-\sum_{k\neq i}^m \alpha_k v_k$. Thus we can write $E_j = \sum_{i=1}^m \beta_{ij} v_i = -\sum_{i=1}^m \gamma_{ij} v_i$, where $\beta_{ij}, \gamma_{ij}$ are all nonnegative. Suppose $x \in Y$; then $E_j x \leq \sum_{i=1}^m \beta_{ij} $, and $-E_j x \leq \sum_{i=1}^m \gamma_{ij}$. Thus $|E_j x|$ is bounded above by the maximum of the two. Hence $\| x\|_{\infty}  = \max \brcs {|E_jx|}$ is bounded above by the maximum over $j$ and thus $Y$ is bounded. 

 \noindent (c) The extreme points are determined uniquely by equations of the form $v_j x = 1$ for some collection of $v_j$s inside $\brcs{v_i}$, and every affine linear system with rational entries that has a real solution has a rational solution---so uniqueness guarantees that the one solution must have rational coordinates. 

 \noindent (d) Really well-known. 

 \noindent (e) Pick nonzero $x \in Y$; let $R$ be the ray  $tx$ (with $t \geq 0$). Since $0$ belongs to the interior of $Y$, there exists a smallest $t \equiv t_0$ \st $t_0 x $ belongs to the boundary of $Y$; necessarily, $t_0 \geq 1$. The boundary of $Y$ is the union of its facets, hence there exists a facet  $F$ \st $t_0 x \in F$. Then $x = t_0^{-1}(t_0)x + (1- t_0^{-1}0) \in \cvx \brcs{F,0} =C_F$. 
 Now suppose $x \in C_F$; then we can write (since $F$ is convex), $x =\lambda f + (1-f)0$ for some $f \in F$ and $0\leq \lambda \leq 1$, that is, $x = \lambda f$. Then $v_F x = \lambda v_F f = \lambda$, while for all the other $v_i$s, $v_i x = \lambda v_i f \leq \lambda v_i f$ (since $f\in Y$). \qed

\noindent {\it Proof that trivplicity-freeness entails SWC.} Let $\Arrow \Theta; K.\gl(d,\Q)$ denote any rational representation of the finite group $K$; this is an action of $K$ on $\Q^{d\times 1}$ (the last is usually denoted $\Q^d$). We obtain an action of $K$ on the dual space, $\Q^{1\times d}$ via $\Theta^*(g)(v) = v\theta(g^{-1})^T$; of course, $(\Theta^*)^*$ is equivalent to $\Theta$. Now we claim that if $\Theta$ is trivplicity-free, then so is $\Theta^*$. It is sufficient to show that if $\Theta$ is irreducible, then so is $\Theta^*$ (since we can apply this to $\Theta^*$), but this is elementary.  

We will construct a weight function on $H$ that extends to $G$, \st if $\Cal K$ is a cross-section of $K$ in $G$ (that is, $\Cal K$ consists of exactly one representative of the equivalence classes modulo $H$), then $(l,\Cal K)$ implements SWC. For each $k\in K$, select $E_k \in G$ whose image in $K$ is $k$; without loss of generality, we may assume $E_1 = 1 $. Let $\Cal K = \brcs{E_k}_{k \in K}$.

 We have $\theta$, the action on $H = \Z^d$, and $\theta^*$ acting on the dual space, in this case, $\Z^{1\times d}$. Obviously $\theta$ acts on $\Q^d $ and $\R^d$, the second space yielding $\Theta$, and $\theta^*$ on $\Q^{1\times d}$ (it is immediate that $\theta^* \otimes \Q  = \Theta^*$), and on $\R^{1\times d}$. Now $\Theta$  is trivplicity-free, so that $\Theta^*$ is as well. This means that the action of $K$ on $\Q^{1\times d}$ is cyclic (this part only requires the representation to be multiplicity-free). Let $v$ be a cyclic vector. Then the orbit of $v$ under the action of $K$ spans the $K$-module $\Q^{1\times d}$ as a rational vector space, and it is immediate that it spans $\R^{1\times d}$ as a real vector space. Moreover, the orbital sum  of $v$ (summing over its orbit) is invariant; since the trivial representation does not appear in $\Theta^*$,  it follows that $\sum \theta^*(k) (v) = 0$, so the hypotheses of Lemma \Xfiv(b)  are satisfied. 

 Write the orbit as $Z:= \brcs{\Theta^* (k)v}_{k \in K}$, and define $Y = \Set{w\in \R^{d\times 1}}{zw \leq 1 \text{ for all $z \in Z$}}$. Then $Y$ is a compact convex body (by Lemma \Xfiv(a,b)), a polyhedron, and there is an obvious action of $K$ on $Y$. Because $Y$ is defined by inequalities defined by linear functions, there must exist a facet $F$ exposed by one of the linear functionals, i.e., $F = Y \cap z_0^{-1}(\leq 1)$ for some $z_0 \in F$. Since $K$ acts transitively on $Z$, it follows that every one of the $z$s in $Z$ exposes a facet. Moreover, if $F_0$ is an arbitrary facet, it must be exposed by one of the $z$s. It follows that $K$ acts (the action obtained from $\Theta$) transitively on the set of facets of $Y$. We also see that the extreme points of $Y$ have only rational coordinates ((c) above). 

 Now pick a positive integer $M$ so that $Mv$ and $M$ times every extreme point of $Y$ has only integer coordinates. The outcome is that $MY$ is a lattice polytope with interior, and $\Theta(K)$ acts transitively on its facets; moreover, restricting to lattice points yields the action of $\theta$ (and in the dual space, $\theta^*$). 

 Define a weight function $\Arrow l_0; \Z^d . \Z^+$ via $l_0 (w) = \sup \Set{Mz w}{z = \Theta^*(k)v, k \in K}$ (using additive notation for elements of $H = \Z^d$). We check that this is a weight function, among other things. Index the facets of $MY$, $F_z$ (they are in bijection with the points in the orbit of $v$), and for each one, define the lattice cone, 
 $$
 \Cal C_z = \(\bigcup_{n=1}^{\infty}\cvx \brcs{nF, 0}\) \cap \Z^d 
 $$
 Each $\Cal C_z = \cup n Y_{F_z}$; moreover, $\Cal C_z$ is a subsemigroup of $\Z^d$, and $\Z^d = \cup_{z} \Cal C_z$. The latter follows since $MY = \cup_z \cvx\brcs{F_z,0}$ by  Lemma \Xfiv(e). 

 Now we claim that if  for some  $w \in \cvx \brcs{F_z,0}$ (not necessarily a lattice point) for some $z$, then $Mz w \geq Mz' w$ for all $z' \in Z$. But this is an immediate consequence of the second part of Lemma \Xfiv(e). In particular, $l_0| \Cal C_z $ is a semigroup homomorphism (additive). 

 Now  for any $t$, $l_0^{-1}(\leq t)$ consists of the lattice points in a multiple of $Y$, so is finite. Moreover, $l_0 (w_1 + w_2) \leq l_0 (w_1) +l_0 (w_2)$ is obvious, as is the fact that the range of $l_0$ is nonnegative (since every $w$ belongs to at least one $\Cal C_z$). So $l_0$  is a weight function on $H$, and it is invariant under the action of $K$.

 It thus can be extended to a weight function on all of $G$ as in the proof of Lemma \Xthr.

 Now set $W = \brcs{E_k}_{k \in K}$. For $w,w'$, there exists $k \in K$ \st $\theta(k)(w')$ belongs to the same $\Cal C_z$ as $w$ (since $\theta^*$ acts transitively on the faces). Then $l (w E_k w') = l(w E_k w' E_k^{-1}E_k)$, and this is $l_0 (w ) + l_0 (E_k w' E_k^{-1}) = l_0 (w) + l_0(w')$. It follows from Lemma \Xfou\ 
that $(l,W)$ implements SWC on $G$.\qed

 As to the converse, one part is easy: if $\Theta = \theta \otimes 1_{\Q}$ contains the trivial representation, then $Z(G) \cap H$  is not trivial, hence infinite; but SWC implies the centre is finite. The necessity of multiplicity-freeness is  more difficult, as we do not have much control on the choice of weight function $l$; on the other hand, by the next result, we can assume that $W$ is any cross-section of $K$ in $G$. 

 \Lem Lemma \Xsix. Let $H$ be a torsion-free abelian group that is a normal subgroup of $G$, with $K:= G/H$ finite. Suppose $(l,W)$ implements SWC for $G$. Let $\Cal K = \brcs{E_k}_{k \in K}$ be a cross-section of $K$ inside $G$.  Then $(l, \Cal K)$ implements SWC for $G$.

\Pf Set $g = hE_j$ and $g' = h'E_{j'}$ (with $h,h' \in H$). There exists $w \in W$ \st $l(gwg') \geq l(g) + l(g') - C$. Write $w = h'' E_k$. Then $gwg' = h E_j h'' E_k h'' E_{j'}$, which we rewrite as $h \cdot  E_j h'' E_j^{-1}\cdot E_j E_k h'' E_{j'}$. Thus,
$$\eqalign{
 l (h E_{jk} h') & \geq l(h \cdot  E_j h'' E_j^{-1}\cdot E_j E_k h'' E_{j'}) - l(\theta(j)(h'')^{-1})) - l(c(j,j')^{-1}) - l(E_{j'}^{-1})\cr
 & \geq  l(g) + l(g') - C - \sup_{j, h''} l(\theta(j)(h'')^{-1})) - \sup_{j,j'}  l(c(j,j')^{-1}) - \sup_{j}  l(E_{j}^{-1}) 
}$$
 All three of the suprema are over finite sets independent of the choice of $h$ and $h'$, so simply contribute to the constant. Finally, Lemma \X.2 applies.\qed

 In the case that $H \to G$ splits, i.e., $G = H \times_{\theta} K$, we can just take $K = W$. It can certainly happen that a smaller subset of $K$ can be chosen for $W$; this already happens for the non-trivial representation of $C_6$ on $\Z^2$ (there is only one, up to conjugacy), wherein $W$ can be chosen to be the copy of $C_3$ inside $C_6$.

\Lem Corollary \Xsev. Let $H$ be a torsion-free abelian normal subgroup of finite index in the finitely generated group $G$. Let $\Arrow \pi;G . G/H:= K$ be the quotient map. Suppose that $G_0 $ is a subgroup of $G$ \st $\pi(G_0) = K$.  If $G$ satisfies SWC, then so does $G_0$. 

\Pf Let $\Cal K$ be a cross-section of $K$ in $G_0$. By Lemma \X.6, there exists a weight function $l$ on $G$ \st $(l,\Cal K)$ implements SWC for $G$. Since $\Cal K \subset G_0$, it follows immediately that $(l|G_0,\Cal K)$ implements SWC for $G_0$. \qed

 Define a {\it  real weight function,} also denoted {\it $\R$-weight,} on a group $G$ to be a function $\Arrow\LL; G.\R^+$ \st  
 \item{(i)} $\LL(1) = 0$;
 \item{(ii)} for all $g,g' \in G$, $\LL (gg') \leq \LL (g) + \LL (g')$; 
 \item{(iii)} for all $t\in \R^+$, $\LL^{-1}(\leq t)$ is finite.

We continue to use additive notation when the group is abelian.

 \Lem Lemma \Xeig.  Let $H$ be a torsion-free abelian normal subgroup of finite index in a group $G$. Let $K = G/H$, let $\Arrow \theta; K. \Aut H$ be the corresponding integral representation, and let $\Cal K = \brcs{E_k}_{k \in K}$ be a cross-section in $G$ of $K$. Let $(l, \Cal K,C)$ implement SWC on $G$. Define $\Arrow \LL;G.\R^+$ via $\LL (g) = \lim_{n\to \infty} l(g^n)/n$, and denote the restriction of $\LL$ to $H$ by $\Arrow \LL_0;H.\R^+$. The following hold.
 \item{(0)} For all $k \in K$ and $h \in H$, $\LL (h^{\theta(k)}) = \LL (h)$;
 \item{(a)} for all $h \in H$, $l(h)\geq \LL_0 (h) \geq l(h) - |K|(C + l(\Cal K))$; 
 \item{(b)} $ \LL_0$ is a real weight function on $H$;
 \item{(c)}  for all $h',h \in H$, there exists $k \in K$ \st $\LL(h' + h^{\theta(k)}) =\LL(h') + \LL (h) $.

\Rmk We will show that $\LL_0$ is the gauge of a compact convex polytope having the origin in its interior. When $l|\Z^d= l_S$ for some admissible $S \subset \Z^d$, this follows from Lemma \Yfiv. However, here there is no reason to think that the restriction to $\Z^d$ is of this form. 

\Rmk Almost never is $\hat l$ a real weight function on $G$ (Lemma \Ytwo). 
 \Rmk Part (0) does not depend on SWC. It can be rephrased as $\LL \circ \theta = \LL$.

\Rmk As it stands, $\LL$ is real-valued. It is conceivable that $\LL$ in this context is integer-valued; this would have considerably simplified the proof of the converse. 
  
 \Pf That $\LL_0$ is well-defined is sometimes called Fekete's (Black's) theorem; it also obviously satisfies (i) and (ii) of the definition of real weight function.
  
  \noindent (0) Follows from Lemma \Yone\ (Appendix).
\comment  
Set $D = l(\Cal K \cup \Cal K^{-1})$. For $k \in K$ and $h$ in $H$, by definition, $h^{\theta(k)} = E_k h E_k^{-1}$, and this also applies to powers of $h$, that is, $nh^{\theta(k)} = E_k (nh)E_k^{-1}$ for $n \geq 1$. Thus $l(nh^{\theta(k)}) \geq l(nh) - 2D$. Dividing by $n$ and taking limits, we see that $\LL(h^{\theta(k)}) \geq \LL(h)$. From $nh = E_k^{-1}(E_k (nh)E_k^{-1})E_k$, we see that $l(nh) \geq l(nh^{\theta^{k}}) - 2D$, so $\LL(h^{\theta(k)}) \leq \LL(h)$. 
\endcomment  
 
 \noindent (a,b) Fix an element $h$ of $H$. Let  $M$ be a positive integer and consider (with $M$ $h$s)
$$\eqalign{
Ml(h) - (M-1)C  &\leq \sup l \(h \Cal K h \Cal K h  \dots     \Cal K h    \)  \cr 
 & \leq  l\( h + h^{\theta(k(1))}+ \dots + h^{\theta(k(M-1))} \) + l(\Cal K^M)\cr
 & \leq l\(\sum_{k \in K} f_M (k) h^{\theta (k)}\)+ Ml(\Cal K);\cr
}$$
 here $\LL(\Cal K^M)$ is the supremum  of the $\LL$ values over all products of $M$ elements of $\Cal K$, the  $k(i)$ are various elements of $K$, and $f_M(k)$ are nonnegative integers \st $\sum_k f_M(k) = M$. For each $M$, there exists $k_M$ \st $f_M (k_M) \geq M/|K|$. There thus exists an element $k' \in K$ and  infinite subset $S$ of $\N$ \st $k_M = k'$ for all $M \in S$. By subadditivity, for $M \in S$,
 $$\eqalign{
 l\(\sum_{k \in K} f_M (k) h^{\theta (k)}\) &\leq l(f_M (k')h) + \sum_{k \neq k'}f_M(k)l (h^{\theta (k)})\cr & = l(f_M (k')h) + (M- f_M(k'))l(h), \quad \text{so}\cr 
l(f_M (k')h) & \geq  Ml(h) - (M-1)C - Ml (\Cal K) - (M- f_M(k'))l(h)\cr 
& = f_M(k') l(h) - M(C + l (\Cal K)); \quad\text{thus}\cr
  \frac{l(f_M(k')h)}{f_M(k')} & \geq  l(h) - \frac{M(C + l (\Cal K))}{f_M(k')} \cr 
 & \geq l(h) - |K|(C+ l(\Cal K)).\cr
}$$

 Since $f_M(k') \geq M/|K|$ and $S$ is infinite, $\LL (h) = \lim_{M \in S} l(f_M(k')h)/f_M(k')$, and thus $\LL (h) \geq l(h) - |K|(C + l(\Cal K))$. This proves (a), and now property (iii) in the definition of real weight function follows from $\LL^{-1}(\leq t) \subseteq l^{-1}(\leq t + |K|(C + l(\Cal K)))$; the former is thus finite.  This finishes the proof of (b).

 \noindent (c) For any positive integer $M$, there exists $k^{(M)}$  \st $l((Mh') + (Mh)^{\theta(k^{(M)})}) \geq l(Mh') + l(Mh') - C$. The   former is $l(M(h' + h^{\theta(k^{(M)})}))$. There exists $k \in K$ \st $k^{(M)} = k$ for infinitely many $M$. For such $M$, divide the inequality by $M$ and take the limit; we obtain $\LL(h' + h^k) \geq \LL(h') + \LL(h)$. The reverse inequality is trivial, since $\LL$ is $\theta$-invariant (part (0)).\qed          

 \comment 
 \item{(c)} for all $h',h$ in $H$, there exists $k \in K$ \st $\LL(h' + h^{\theta(k)}) \geq \LL(h') + \LL (h) - |K|C$;
\endcomment

\comment
 Assume that $\Theta:= \theta \otimes 1_{\Q}$ is not multiplicity-free, but does not contain the trivial representation. If also $G$ satisfies SWC, we will obtain a contradiction.

 Assume that $G$ satisfies SWC; by Lemma \Xsix, we may assume it is implemented by $(l,\Cal K, C)$, where $\Cal K$ is a cross-section. Now $\Theta$ is not multiplicity-free; it is routine to check that there exists a $\theta$-invariant subgroup $J$ of $H$ \st $J = J_1 \oplus J_2$ and there is an isomorphism from $J_1$ to $J_2$ commuting with the action of $\theta$. For convenience, we rewrite this as $J = J_1 \oplus J_1$, where, with  $H_0:= J_1 \oplus 0$ and $H_1:= 0 \oplus J_1$, each of $H_0$ and $H_1$ is $\theta$-invariant, and the action of $\theta$ on each of the two copies is identical. The last means $\theta (k) (j,j') = (j^k,( j')^{k}$ (for $j,j' \in J_1$) the action being obtained from (a copy of) $J_1$, that is, the action is the same on each coordinate. 

 Now construct $\LL$ as previously. 
\endcomment
 \def\AA{{\Cal A}}

 The following subsum principle is obvious, but very useful. 

 \Lem Lemma \Xthy\ (Subsum principle). Suppose that $H$ is an abelian group, and $\Arrow l; H.\R^+$ is subadditive; also suppose that $A $ is a finite subset of $H$, and $l\(\sum_{a \in A} a\) = \sum_A l(a)$. 
\item{(i)}For all nonempty $B\subset A$, we have $l\(\sum_{a \in B} a\) = \sum_B l(a)$.
\item{(ii)} If additionally, $l(nx) = nl(x)$ for infinitely many positive integers $n$, then for all functions $\Arrow f; A. \Z^+$, $l\(\sum_{a \in A} f(a)a\) = \sum_A f(a) l(a)$.

 \Pf (i) is trivial; as for (ii), $l(nx) = n l(x)$ for all nonnegative integers $n$ follows from (i), and now (i) applies to $(\max f(a))\cdot \sum_{A}a$. \qed

We have a batch of definitions. For $x$ in $H$, and $\LL$ as in Lemma \Xeig, 
$$\eqalign{
B(x) & = \Set{\phi \in K}{\LL(x + x^{\phi}) = 2 \LL (x)}\cr
 \AA_0 (x) & = \Set{x' \in H}{ \LL(x + x') = \LL(x) + \LL (x')}\cr
 \AA(x)& = \Set{x' \in H}{ \LL(x + {x'}^{\phi}) = \LL(x) + \LL (x')\text{ for all $\phi \in B(h)$}}\cr
}$$
 From the construction of $\LL$, $B(x)$ contains $1$---but is not generally a group (in all examples I've examined, it is a union of groups). Also, $\cap_{\phi\in K} \AA_0(x^{\phi}) = \AA(h)$. If $B(h)$ is a minimal set in the collection of sets, $\brcs{B(x)}$ as $x$ varies over $H$, we say that $B(h)$ is {\it minimal.} Minimal $B(h)$ are often the one-element group, but need not be.  Now we have a number of elementary results. 

 \Lem Lemma \Xnin. Suppose $x,h \in H$, and $B(h)$ is minimal. 
\item{(a)} If $\phi \in B(x)$, then $B(x + x^{\phi}) \subseteq B(x)$ and  $\phi \cdot B(x + x^{\phi}) \subseteq B(x)$. 
 \item{(b)} $B(h)$ is a group. 
 \item{(c)} If $x_1 \in \AA_0(x)$, then $B(x+x_1) \subseteq B(x) \cap B(x_1)$.
 \item{(d)} If for some $h_1$ in $H$, $B(h_1)$ is  minimal, then there exists $\alpha \in K$ \st $B(h_1) = \alpha B(h)\alpha^{-1}$. 
\item{(e)} $\AA(h)$ is a subsemigroup of $H$, stabilized by $B(h)$; we also have that $\AA(h)$ is relatively convex (no holes), and generates $H$  as a group.
\item{(f)} For all $x \in H$, the set $\Set{\phi \in K}{x^{\psi} \in \AA(h)}$ is nonempty and is a union of cosets of $B(h)$ in $K$.

 \Pf (a) For $\psi \in B(x + x^{\phi})$, we have $\LL(x + x^{\phi} + x^{\psi} + x^{\phi\psi}) = 2\LL (x + x^{\phi}) = 4\LL(x)$. By the subsum principle (applied to $x + x^{\psi}$), $\psi \in B(h)$; the subsum principle applied to $x + x^{\phi\psi}$ yields $\phi\psi \in B(x)$.
 
 \noindent (b) By minimality, $B(h + h^{\phi}) = B(h)$ for all $\phi \in B(h)$; by (a), $\phi B(h) = B(h)$, that is, $B(h)$ is closed under the operation. Since the identity belongs to $B(h)$ and $K$ is finite, $B(h)$ is a group. 

 \noindent (c) For $\psi \in B(x+x_1)$, we have $\LL (x+x_1 + x^{\psi} + x_1^{\psi}) = 2\LL(x + x_1) = 2\LL (x) + 2\LL (x_1)$. By the subsum principle applied to $x + x^{\psi}$ and $x_1 + x_1^{\psi}$, we see that $\psi \in B(x)$ and $\psi \in B(x_1)$ respectively. 

 \noindent (d) There exists, by Lemma \Xeig(a,c), $\rho \in K$ \st $\LL(h + h_1^{\rho}) = \LL(h) + \LL(h_1^{\rho})$. Applying (c) with $x_1 = h_1^{\rho}$, and using minimality of $B(h)$, we have $B(h + h_{1}^{\rho}) = B(h) \subset B(h_1^{\rho})$. But $B(h_1^{\rho}) = \rho B(h_1) \rho^{-1}$. Minimality of $B(h_1)$ implies minimality of any of its conjugates, so $B(h) =\rho B(h_1) \rho^{-1}$; set $\alpha = \rho^{-1}$. 

\noindent (e) Pick $h_1 \in \AA(h)$. We first show that $h + \sum_{B(h)} h_1^{\phi} \in \AA(h)$, by induction on the combined statement, $h + \sum_{T} h_1^{\phi} \in \AA_0(h) $ for subsets $T \subseteq B(h)$ and $B(h + \sum_{T} h_1^{\phi} \in \AA(h)) = B(h)$. If $|T| = 1$, the result is true by hypothesis and (c). For $|T| > 1$, write $T = T_0 \dot\cup \brcs{\psi}$, and set $j = h + B(h + \sum_{T_0} h_1^{\phi}$; then $B(h) =  B(j)$ by the induction hypothesis, so $\LL(j + j^{\psi}) = 2\LL(j) = 2\LL(h) + 2(|T|-1)\LL(h_1)$. By the subsum principle, $\LL\( h + \sum_{T} h_1^{\phi}\) = \LL(h) + |T|\LL(h_1)$, so that $\sum_{T} h_1^{\phi} \in \AA_0(h)$; now by (c) with $x = h$ and $x_1 =\sum_{T} h_1^{\phi}$, we have that $B(h + \sum_{T} h_1^{\phi} \in \AA(h) = B(h)$, completing the induction. 
 
 Since $e:= \sum_{B(h)} h_1^{\phi} \in \AA_0(h)$ and is invariant under every element of $B(h)$, it follows that $e \in \AA(h)$. 

 Now pick $h_2 \in \AA(h)$. Then ${}2:= \sum_{B(h)} h_2^{\phi}$ is $B(h)$-invariant, and belongs to $\AA(h)$. There exists $\psi \in K$ \st $\LL (h + e + (h + {}2)^{\psi}) = \LL(h+e ) + \LL(h+{}2)^{\psi}$; but this is just $2\LL (h) + |B(h)|\LL(h_1) + |B(h)|\LL(h_2)$. By the subsum principle applied to $h + h^{\psi}$, we see that $\psi \in B(h)$. By invariance, the sum is just$\LL(h + h^{\psi} + e + {}2) $. By the subsum principle, $\LL(h + e + {}2) = \LL (h) + |B|\LL(h_1) + |B|\LL({}2)$. Again by the subsum principle, $\LL(h + (h_1 + h_2)^{\phi}) \in \AA_0(h)$ for all $\phi \in B(h)$, and thus $h_1 + h_2 \in \AA(h)$. It obviously is stabilized by $\AA(h)$. 

 Suppose that $a \in H$ and for some positive integer $M$, $Ma \in \AA(h)$. Then the subsum principle once more shows that $a \in \AA(h)$. So $\AA(h)$ is relatively convex. 

 Let $J_h = \AA(h) - \AA(h)$, the group generated by $\AA(h)$. Then each of $J_{h^\phi}$ is conjugate to $J_h$ within $H$, and so if one is of infinite index in $H$, then they all are. But $H = \cup J_{h^{\phi}}$, so by  [NN], so at least one of them is of finite index in $H$. But conjugacy within $H$ yields that they all are of finite index within $H$. Finally,  relative convexity shows that $J_h = H$. 

\noindent (f) By numerous applications of the subsum principle, $\AA(\sum_{B(h)} h^{\phi}) = \AA(h)$; set $p =  \sum_{B(h)} h^{\phi}$. There exists $\psi \in K$  \st $\LL(p + x^{\psi} ) = \LL (p) + \LL (x)$, that is, $x^{\psi} \in \AA_0(p)$. Since $p$ is $B(h)$-invariant, $\AA_0(p) = \AA(p) = \AA(h)$. Since $\AA(h)$ is stabilized by $B(h)$, $x^{\psi\phi}\in \AA(h)$ for all $\phi \in B(h)$, that is, $\psi B(h)$ is in the set. \qed

If $B(h)$ is minimal, then as we will see later, $B(h)/\ker \theta$ is a malnormal subgroup of $K/\ker \theta$ (a subgroup of a finite group is {\it malnormal\/} if the intersection of its conjugates consists of the identity element). 

 If (minimal) $B(h)$ is trivial, then the statements are even easier to prove (and some become trivial, e.g., (d) and the stabilization part of (e)). But minimal $B(h)$ need not be trivial, even for the examples constructed in the course of proving the SWC property; whether this occurs depends on the choice of the cyclic vector in the dual space. Minimal $B(h)$ is trivial iff the  orbit of the cyclic vector  under $K$ is full (that is, of cardinality equalling that of $K$) and $\ker \theta$ is trivial; but it can easily happen (with $\ker \theta$ trivial) that the orbit is not full (e.g., for $K = S_3$, the orbits can be of size 3 or 6). On the other hand, it follows from (d) that $B(h) = \ker \theta$ if $K/\ker \theta$ is abelian, or marginally more generally, if every subgroup of prime order is normal (this includes Hamiltonian groups).

 \noindent {\it Conclusion of the proof of the converse.} Pick $h$ \st $B(h)$ is minimal. Then $\AA(h)$ is subsemigroup of $H$, and if $C$ is a complete set of coset representatives of $B(h)$ in $K$ containing $1$, then $H = \cup_{\psi \in C}\AA(h^{\psi})$ (notice that $\AA(h) = \AA(h^{\phi})$ for $\phi \in B(h)$, and we can even replace $h$ by $\sum_{B(h)} h^{\phi}$, that is, assume invariance).  It is tempting to consider what would amount to the boundary $\cup_{\psi \in C\setminus\brcs{1}} (\AA(h) \cap \AA(h^{\psi})) \neq \AA(h)$ (some of the intersections can consist only of  $0$), and a fundamental domain for the $K$-action, etc. But we don't really need all that to prove the converse. 

For a fixed $h$ \st $B(h)$ is minimal, let $C(h)$ denote the convex hull of $\AA(h)$ in $\R^{d\times 1}$. We will extend $\LL$ to a convex-linear function on $C(h) \to \R^+$ whose only zero is at the origin. First,  $\LL$ extends to an additive function $H = \AA(h) - \AA(h) \to \R$ in the obvious way: if $x= x_1 - x_2 = x_3 - x_4$ with $x_i \in H$, we have $x_1 + x_4 = x_2 + x_3$, so that $\LL (x_1) + \LL(x_4) = \LL (x_2) + \LL(x_3)$, hence $x_1 - x_2 \mapsto \LL (x_1) - \LL (x_2) $ is well-defined, and obviously additive. Hence there exists unique $v \in \R^{d\times 1}$ \st $x \mapsto vx$ implements this extension; in particular, for $x \in \AA(h)$, we have $\LL(x) = vx$. This yields a map $\Arrow L; \R^{d\times 1}.R$ given by left matrix multiplication by $v$. 

It suffices to show that \wrt  dual action of $K$ (on $\R^{1\times d}$) is cyclic, in particular that the span of $\brcs{v^{\phi}}_{\phi \in K}$ is all of $\R^{1\times d}$. If not,  there exists nonzero  $y\in \R^{d\times 1}$ \st $v^{\phi}y = 0$ for
all $\phi \in K$. This amounts to $vy^{\phi} = 0$  for all $\phi \in K$. For each
positive integer $N$, there exists $x_N \in \Z^{d\times 1} = H$ \st $\| Ny - x_N
\|_{\infty} \leq 1$; of course, the set $\brcs{x_N}$ is not finite. Since $K$
acts  as a finite set of linear transformations, there exists $A
\in \R^{+}$ \st for all $N$  and $\phi$, $\| Ny - x_N^{\phi} \|_{\infty} \leq
A$. Applying $v$, we obtain $|vx_N^{\phi}| < A\| v\|$ (the latter norm is as a linear functional on $\R^{d\times 1}$, the latter equipped with the $\infty$-norm). But this means
$\brcs{\LL (x_N)}$ is bounded, contradicting finiteness of $\LL^{-1}(\leq t)$ for
all $t\geq 0$. \qed

 \comment 
A standard argument now yields that if $x = \sum \lambda_i x_i$ with $x_i \in \AA(h)$ and $\lambda_i$ nonnegative rationals, then $\sum \lambda_i \LL( x_i) = vx$: multiply by a positive integer to clear denominators, and use additivity of $\LL$. Finally if $y = \sum \sigma_i x_i$ with $\sigma_i$ nonnegative real numbers and $x_i \in \AA(h)$, then again $ \sum \sigma_i \LL(x_i) = vy$, since we may approximate the $\sigma_i$ both above and below by nonnegative rational numbers. 
 
 Now we observe that $\LL(x) = 0$  for some $x \in H$ entails that $x = 0$. This follows immediately from $\LL(nx) = n \LL(x)$, and $\LL^{-1}(0)$ is finite. By the preceding, for nonzero $y = \sum \sigma_i x_i\in C(h)$ (with $\sigma_i > 0$, $x_i \in \AA(h) \setminus \brcs{0}$), $vy = \sum \sigma_i vx_i > 0$. 

 Now we claim that the action of $K$ on $\R^{1\times d}$ (the dual  of the action on $\R^{d\times 1}$ obtained from the original action on $H$) is cyclic with generator $v$. If not, the subspace spanned by the orbit of $v$, $\brcs{v^{\phi}}$ is not full, and so there would exist a nonzero vector $y \in \R^{d\times 1}$ \st $v\phi y = 0$ for all $\phi \in K$. This amounts to $v w^{\phi^{-1}} = 0$ for all $\phi \in K$. Since $H = \cup \AA(h^\psi) = \cup \AA(h)^{\psi}$, it follows that $\cup C(h)^{\psi} = \R^{d\times 1}$. Hence there exists $\psi \in K$ \st $y^{\psi} \in C(h)$. Thus $vy^{\psi} >0$, contradiction. Hence the action of $K$ on $\R^{1\times d}$ is cyclic, and it easily follows from the reduction that $\Theta$ is cyclic. \qed 
\endcomment

\noindent {\it Faithfulness---or not---of $\theta$.} Here we elaborate a little on the structure of minimal $B(h)$. It is clear that $\ker \theta \subset B(x)$ for any $ x \in H$, and when $B(x)$ is a group, $\ker \theta$ is normal in it. Let $B(h)$ be minimal. We show that $B(h)/\ker \theta$ is malnormal in $K/\ker \theta$, that is, the intersection of its conjugates is trivial. 

 \Lem Lemma \Xten. Given $k \in K\setminus \ker\theta$, there exists $j_k \in H$ \st $\LL(j_k + j_k^k) < 2\LL (j_k)$. 

\Pf Assume to the contrary  that for all $j \in H$, we have $\LL(j_k + j_k^k) = 2\LL (j_k)$ (the left side cannot be larger than the right, because of subadditivity and $\theta$-invariance of $\LL$). Let $S$ be a collection of coset representatives of $\langle k\cdot \ker \theta \rangle$ in $K/\ker \theta$, and for any $j \in H$, define $x \equiv x(j,S)= \sum_{a \in S} j^a$. 

 Since the trivial representation does not appear in $\Theta$, it follows that $\theta$ has no nontrivial fixed points, and thus for any $x \in H$, $\sum_{b \in K} x^b = 0$. Let $s$ be the order of $k$ modulo $\ker \theta$ (that is, the smallest positive integer \st $k^s \in \ker \theta$. We then have $\sum_{i=0}^{s-1} x^{k^i} = 0$. We claim we can choose $j$ so that the corresponding $x(j,S)$ is not zero for some choice of the coset representatives, $S$. 

 To this end, pick $j \in J$ and $t \in S$ \st $j^t \not\equiv j^{tk}\pmod {\ker  \theta}$ (that this exists follows from the definitions). If for some $S$, $x(j,S) = 0$, then replace one element $t$ of $S$ by $t' = tk$, creating $S'$. If also $x(j,S') = 0$, then subtracting one expression from the other, we obtain $j^{tk} - j^t$, which is nonzero by hypothesis. Hence at least one of $x(j,S)$ and $x(j,S')$ is not zero. Call whichever is nonzero,  $x$. Then $x \neq 0$, but $\sum x^{k^i}= 0$. 

 Now iterate the map $j \mapsto j + j^k$, starting with $x$, that is, we obtain a sequence defined by  $x_1 = x + x^k$, $x_n = x_{n-1} + x_{n-1}^k$. By the very first assumption, $\LL (x_1) = 2\LL (x)$, $\LL(x_2) = 2 \LL(x_1) = 4\LL(x)$, and in general, $\LL(x_n) = 2^n \LL(x)$. Write $x_n = \sum_{i=0}^{s-1} f_n(i)x^{k^i}$; here $\sum f_n(i) = 2^n$. It is easy to check that the  $s$-tuple of coefficients  $(f_n(0), f_n(1), \dots , f_{n-1})$ is obtained by applying the matrix $(\I_s + P)^n$ to ${}1 = (1,0,0,\dots)$, where $P$ is the cyclic permutation matrix given by $0 \to 1 \to 2 \dots \to s-1 \to 0$. In particular, ${}1(\I_s + P)^n = 2^n (1,1,\dots,1)/n + v$  where $\| v \|_1 = \Oh {\lambda^n}$ for some $\lambda < 2$ ($\lambda$ is the second largest absolute value and is easily computed---we just don't have to). Since $\sum x^{k^i} = 0$, this yields $x_n = \sum a_i x^{k^i}$ where $\sum a_i = \oh{2^n}$ and we can always ensure that the $a_i$ are nonnegative by subtracting $\min_i \brcs{f_n(i)}$ from each term. But this gives $\LL (x_n)\leq \sum a_i \LL (x)$, contradiction. 

 Thus there exists $j\equiv j_k \in H$ \st $\LL(j+j^k) <2 \LL (j)$. \qed

A subgroup $J$ of a group $L$ is {\it malnormal\/} if it contains no nontrivial normal subgroup of $L$. 

\Lem Corollary \Xele. If $B(h)$ is minimal, then $B(h)/\ker \theta$ is malnormal in $K/\ker\theta $.

\Pf By Lemma \Xten, $\ker \theta = \cap_{x \in H} B(x)$; this remains true if we consider the intersection over only minimal $B(x)$. But all minimal $B(x)$ groups are conjugate in $K$ to $B(h)$, the intersection on the right is $\cap_{\phi \in K} \phi B(h)\phi^{-1}$. In particular, every  subgroup of $B(h)$ that is normal in $K$ is contained in $\ker \theta$. \qed

In all of this, there does not appear to be a reduction to the case that $\theta$ is faithful (its kernel is trivial). We pose this as a question. 

\Lem Question. Suppose $G$ is an SWC group that is an extension of an abelian group by a finite group. Does there exist an abelian  normal subgroup $H$ of $G$ \st the map $\Arrow \theta; K := G/H. \Aut (H)$ is one to one?

This question is stated for abelian groups, not just torsion-free abelian groups, but there is an immediate reduction to the latter. A counter-example that is minimal in $|K|$ would have to have various properties, e.g., no nontrivial abelian subgroup of $\ker \theta$ is  normal in $K$,  and no nontrivial subgroup of $\ker \theta$ can split $K$. \vskip 5pt
 
\noindent {\it Symmetry of $\LL$.}
 With $G$ a group, let $\Arrow \xi; G.G$ be the antiautomorphism given by $x \mapsto x^{-1}$. For functions $\Arrow l, \LL; G.\R^+$, define $l_-$, $\LL_-$ to be the composition with $\xi$. Now suppose that $G$ satisfies SWC, and this is implemented by $(l, W)$. An obvious question is whether we can arrange (by modifying $l$) so that $l \sim l_-$, or $\LL = \LL_-$, that is $\LL(x^{-1}) = \LL(x)$ for all $x \in G$. When this occurs, computations tend to be simpler. If $l$ is of the form $l_S$ for some admissible subset $S$, then sufficient for this to occur is that $S = S^{-1}$, i.e., $S$ is symmetric. 

 We restrict to our current situation, wherein $ H$ is torsion-free abelian and of finite index in $G$, with $\Arrow \theta; K=G/H. \Aut (H) $; we may as well assume that $\theta \otimes 1_{\Q}$ is trivplicity-free, and SWC is implemented by $(l,\Cal K)$. As before, we construct $\LL$ restricted to $H$. Using additive notation, we   ask when $\LL(x) = \LL(-x)$. 

The answer is that {\it it depends.} It obviously happens when $\xi \in \theta(K)$, but it also happens when the convex polytope constructed there is symmetric under $x \mapsto -x$. For example, if $K = S_3$ (the permutation group on three elements), and the convex hull of the orbit in the dual space is a hexagon, the corresponding $\LL$ will be invariant under $\rho$, whereas if the convex hull is a triangle, it is not (these correspond respectively to the cases that minimal $B(h)$ be trivial or of order two; in the latter case, there exists  $x \in H$ \st $B(x)$ is minimal, but $B(-x)$ is not.) 
 
 There is a very easy argument to show that for general $K$ of odd order, there is no weight function $l$ with $(l,W)$ implementing SWC \st even merely $\LL = \LL_-$. Thus if $l$ is determined by an admissible set (as is the case with the examples constructed in the course of the proof), the latter cannot be symmetric. We obtain a slightly stronger result. 

 Continuing the notation used in the proof of the converse (that is, $(l,\Cal K)$ implements SWC, $\LL$ is defined via $\lim l(nx)/n$). Define, for nonzero $x \in H$, the subset of $K$ given by
$$
N(x) = \Set{\psi \in K}{\LL(x -x^{\psi}) = \LL(x) + \LL(-x)}.
$$
 By Lemma \Xeig(c), this is nonempty, and obviously misses the identity.

 \Lem Lemma \Xtwe. Let  $x$ be a nonzero element of $H$, and let $\psi \in N(x)$. 
\item{(i)} $N(-x) = N(x)^{-1}$ (the latter is the set of inverses of elements of $N(x)$). 
 \item{(ii)} $B(x - x^{\psi}) \subseteq B(x)$ and  $N(x)\cdot B(x-x^{\psi}) = N(x)$; $\psi B(x)\psi^{-1} \subseteq B(-x)$. 
\item{(iii)}  If additionally, $\psi \in N(-x)$,  then $\(N(x) \cap N(-x)\)\cdot N(x-x^{\psi}) \subseteq B(x)$ and  $N(x- x^{\psi}) \subset N(x)$.

\Pf (i) Trivial.

\noindent (ii) Pick $\phi \in B(x - x^{\psi})$. Then 
$$\eqalign{
 \LL(x + x^{\phi} - x^{\psi} - x^{\psi \phi})& = 2 \LL(x - x^{\psi})\cr
 & = 2 \LL(x) + 2\LL (-x). \cr
}$$
 By the subsum principle, we deduce $\phi\in B(x)$, $\psi\phi \in N(x)$, and $\psi \phi \psi^{-1} \in B(-x)$. The first and third yield the corresponding statements. The second says $N(x)\cdot B(x-x^{\psi}) \subseteq N(x)$, but the identity belongs to $B(x -x^{\psi})$, so cardinality gives equality. 
\noindent (iii) Select $\rho \in N(x-x^{\psi})$. 
 Then 
$$\eqalign{
 \LL\(x - x^{\psi} - (x - x^{\psi })^{\rho}\)& =  \LL(x - x^{\psi}) + L((-x) - (-x)^{\rho})\cr
 & = 2 \LL(x) + 2\LL (-x). \cr
}$$
 By the subsum principle, $\psi \rho \in B(x)$ and $\rho \in N(x)$, yielding the results.\qed

  Now we assume that $N(x) = N(-x)$  and $B(x) = B(-x)$ (for example, these are consequences of $\LL = \LL_-$).

  \Lem Lemma \Xthi. Suppose that $N(x) = N(-x)$ and $B(x) = B(-x)$ for all $x \in H$. Then there exists $j \in H$ \st $B(j)$ is minimal, and $B_0(j) := B(j) \cup N(j)$ is a group, in which $B(j)$ is normal and of index exactly two. 

 \Pf Find $h$ \st $B(h)$ is minimal, and choose $\psi \in N(h)$. Then $B(h - h^{\psi}) = B(h)$ and $N(h - h^{\psi}) \subseteq N(h)$. If the latter inclusion is strict, set $h_1 = h - h^{\psi}$ and choose $\psi_2 \in N(h_1)$; then with $h_2 = h_1- h_1^{\psi_2}$, we have $B(h_2) = B(h_1) = B(h)$ and $N(h_2) \subseteq N(h_1) \subseteq N(h)$. This process may be iterated. It must  eventually terminate, yielding $j:= h_n$ \st $B(j) = B(h)$ (so is minimal) and $N(j-j^{\psi}) = N(j)$ for all $\psi \in N(j)$. 

 Then  $N(j)^2 \subseteq B(j)$ and $N(j)B(j)= N(j)$. The latter says $N(j)$ is a union of right $B(j)$-cosets. If it contained two or more cosets, then $N(j)^2$ (all products of two elements of $N(j)$) would have strictly larger cardinality than $B(j)$. Hence $N(j)$ consists of a single coset. It cannot contain the identity (since $j \neq 0$). Writing $N(j) = \tau B(j)$, we have $\tau \phi \tau \in B(j)$ for all $\phi \in B(j)$. 

 Using $B(x) = B(-x)$, we see that $\tau B(j) \tau^{-1} \subset B(j)$, so equality holds, and thus $\tau B(j) = B(j)\tau$. It follows easily $\tau B(j) \dot\cup B(j)$ is a group with $B(j)$ of index two therein. \qed

 In particular, if $\LL = \LL_-$, then $K$ has to have even order. 
\vskip10pt 

In contrast, the corresponding result for SWC(0) is much easier, at least for finite  extensions of abelian groups. In general, if $G$ is an infinite group satisfying SWC(0), then its centre is finite; this is an obvious consequence of the definitions. In the case that $G$ is abelian by finite, the converse holds. The only problem  is that it's not really clear what the advantage is of $G$ satisfying SWC(0).

\Lem Proposition \Xftn. Let $G$ be an abelian by finite group. The  following are equivalent. 
 \item{(a)} $G$ satisfies SWC(0); 
\item{(b)} the centre of $G$ is finite; 
  \item{(c)} any (or one) extension $H \triangleleft   G \to G/H:= K$ with abelian $H$ and finite $K$, the corresponding representation $\Arrow\theta; K. \Aut (H)$ does not contain the trivial representation.

 \Pf That (a) implies (b) is straightforward, and (b) implies (c) is trivial. Assume (c) for one choice of normal abelian group of finite index. This easily implies that the centre of $G$ is finite, and we thus have (b) equivalent to (c) (the latter with {\it any\/}). Assume (c), so there exist $H$, $\theta$, etc, with the indicated properties. Let $\Cal K$ be a transversal of $K$ in $H$. Let $l$ be a weight function on $H$ \st for all  integers $n \geq 0$, we have $l(nh) = n l(h)$ (such exist; for example, we may take the $l^1$-norm on $H$ modulo its torsion subgroup). By replacing $l $ by $\sum_{k \in K} l \circ \theta(k)$, we may assume that $l $ is $K$-invariant as well. This permits us to define a weight function on $G= \dot\cup \kappa H$ as in the proof of \Xthr, via $l_1 (\kappa h) = l(h) + C$ for non-identity  $\kappa$, extending $l$. 

  Without loss of generality, we can assume $H$ is torsion-free. Hence for any $h\in H$, the orbital sum, $\sum h^k$ (where $h^k = \kappa h \kappa^{-1}$ with $\kappa \mapsto k \in K$) is zero. We have $n l(h)= l(nh) = l\(\sum_{k \in K} (h + (h')^k)\)$.  There thus exists $k' \in K$ \st $l(h + (h')^{k'}) \geq l(h)$. It follows as in \Xthr, that $l_1$ satisfies SWC(0) with $W = \Cal K$. \qed

\SecT \M\ Right multiplications

There are a lot of endomorphisms of $A_f$ (and to a lesser extent, of $R_f$) arising from right multiplication. If $(a,k) \in AG \times \Z^+$, then we define the right multiplication operator $\Arrow \RR_{a,k}; A_f.A_f$ given by $[b,m] \mapsto [ba,m+k]$.  It is easy to check that this is well-defined (since right multiplication  on $AG$ commutes with left multiplication by $f$), and commutes with $\SS$ and any  well-defined left multiplication operator. 

 The brackets $[ \cdot,\cdot]$ are missing from the subscript, because $\RR_{a,k}$ and $\RR_{fa,k+1}$ are almost always different if $af \neq fa$. If $a \in (AG)^+$, then $\RR_{a,k}$ is a positive endomorphism---but if merely $[a,k] \in A_f^+$, this is not so clear. 

It is of interest to determine precisely when $\RR_{a,k}$ is an endomorphism of
$R_f$ or when it is positive. Obviously necessary for $\RR_{a,k}$ to be an
endomorphism of $R_f$ is that $[a,k] \in R_f$. This is sufficient if either
$\Supp a \subseteq S^k$ or if $a \in (AG)^+$. We define $\End A_f$ to be the ring of group endomorphisms of $A_f$ that commute with the shift $\SS$. Right multiplications belong to this ring. We define a positive cone of $\End A_f$, $\End^+ A_f$, to consist of the elements $\theta$ of $\End A_f$ \st $\theta (A_f^+) \subset A_f^+$ (warning: it is not clear that $\End^+ A_f$ generates $\End A_f$ as a ring, equivalently, as an abelian group, so we normally would reduce  to the ring it generates).  Similar definitions and comments apply to $\End R_f$. It is trivial that $\RR_{a,k}$ is a difference of elements of $ \End^+ A_f$ for any $(a,k) \in AG \times \Z^+$.

\Lem Lemma \Mone. Suppose that $[a,k] \in R_f$. 
\item{(a)} If  $a$  is in the centre of $\R G$,
then $\RR_{a,k} \in \End R_f$.
\item{(b)} If $\Supp a \subseteq S^k$, then $\RR_{a,k} \in \End R_f$.
\item{(c)} If $a \in (A G)^+$, then $\RR_{a,k} \in \End^+
R_f$.

\Rmk In (b), the assumption that $[a,k]$ belong to $R_f$ is of course redundant. 
\Rmk Neither of the following is likely to be true: (i) if $[a,k] \in R_f$,
then $\RR_{a,k} \in \End R_f$; (ii) if $[a,k] \in R_f^+$ and $\RR_{a,k} \in \End
R_f$, then $\RR_{a,k} \in \End^+ R_f$.

\Pf To show that $\RR_{a,k}\in \End R_f$, it is sufficient to show that if $g
\in G$, then $\RR_{a,k}([g,\tilde l(g)]) \prec [1,0]$, that is, there exists
$m$ \st $f^m ga \prec f^{m+k + \tilde l(g)}$. For all sufficiently large $n'$,
$\Supp f^{n'}g \subset \Supp f^{n' + \tilde l(g)}$.

\noindent  (a) For all sufficiently large $n$, $\Supp f^n a \subseteq \Supp
f^{n+k}$, whence $\Supp f^n a g \subset \Supp f^{n+k} g$.  Hence for all
sufficiently large $n$, $\Supp f^{n}ag \subset \Supp f^{n + k + \tilde l(g)}$. 
Thus if $a$ is central, then $\Supp f^n ga \subset \Supp f^{n+k + \tilde l(g)}$,
whence $\RR_{a,k} ([g,\tilde l(g)]) \prec [1,0]$.

\noindent (b) By  hypothesis, $\Supp f^{n'+k}a \subset \Supp f^{n'+k+\tilde l}$,
so $\Supp f^{n'}ga \subset \Supp f^{n' +k + l}$.

\noindent (c) For all sufficiently large $m$,  $f^m a \prec f^{m+k}$.  For
sufficiently large $m$, $f^m g \prec f^{m+\tilde l(g)}$. As $a \geq 0$, it
follows that $f^m ga \prec f^{m+\tilde l(g)} a \prec f^{m+ \tilde l(g) + k}$.
\qed

A special case of part (c) arises if $a = g$, an element of the group; then $\RR_{g,k}$ is an endomorphism of $R_f$ iff $k \geq \tilde l(g)$. It has the additional property (which it shares with $\SS$) that $\RR_{g,k} ([b,m]) \in A_f^+$ entails $[b,m] \in A_f^+$.

\noindent  {\it Positivity of $\RR_{a,l}$}. To test whether $\RR_{a,l}$ (with $a \in A G$) is a positive endomorphism of $A_f$, we may assume $l = 0$. The condition reduces to:
\item{}for all $g \in G$, there exists a nonnegative integer $m(g)$ \st $f^{m(g)}a \in (A G)^+$.

 \noindent If the $m(g)$ can be chosen to be bounded (unlikely, and difficult to verify), then it can be shown that if $G$ is ICC, then $\RR_{a,0}$ is positive implies $a \in (A G)^+$.

 A little more interesting is the question of characterizing $(a,l) \in A G \times \Z^+$ \st $\RR_{a,l}$ is an endomorphism of $R_f$.

 First, obviously sufficient is that $\tilde l_S(\Supp a) \leq l$ (since $\RR_{g,\tilde l_S (g)}$ is an endomorphism of $R_f$). If $a \in (A G)^+$, then it is necessary as well: just observe that if there exists $g \in \Supp a$  \st $\tilde l(g) > l$, then $\RR_{g,l}$ is not an endomorphism of $R_f$, and since it is  bounded above by $\RR_{a,l}$, the latter cannot be an endomorphism either. 

 \Lem Lemma \Mtwo. Suppose that $(G,f)$ is a finitely generated group with an admissible element. Then sufficient for $\RR_{a,l}$ be be an endomorphism of $R_f$ is that $\tilde l_S(\Supp a) \leq l$. If $a \in (A G)^+$, this is necessary as well.

\def\End{\text{End}\,}

 \SecT \B\ Bounded endomorphisms

Let $M$ be a partially ordered abelian group, and form $\End^+ M$, the set of
group endomorphisms of $M $ sending $M^+$ to $M^+$, that is, the positive
endomorphisms of $M$. Say $\phi \in \End^+ M$ is {\it locally order-bounded\/}
if for all $x \in M^+$, there exists an integer $N_x$   \st $\phi(x) \leq N_x
x$; this is equivalent to the positive endomorphism $\phi$ sending every order
ideal to itself. If $\sup_{x \in M^+} N_x  < \infty$, then we say that $\phi $
is {\it order-bounded} (sometimes simply called {\it bounded,\/} as in
[H1; section I]).  

Suppose that $g \in G$; we have defined the (positive) automorphism $\Arrow \RR_{g,k}; A_f. A_f$ given by $\RR_{g,k}([h,l]) = [hg,l+k]$; if $g \in \supp f^k$, then $\RR_{g,k}$ induces an order embedding (not generally onto) $R_f \to R_f$.

Let $H$ be a partially ordered group, together with an order embedding $\Arrow \phi; H.H$. Pick an element $v \in H^+$. Define $\Cal E(H,\phi)^+$ to consist of those positive endomorphisms $\Arrow \psi; H.H$ \st $\psi \phi = \phi\psi$. We define $\Cal E(H,\phi)$ to be the group generated by $\Cal E(H,\phi)^+$, that is, $\Cal E(H,\phi)^+ - \Cal E(H,\phi)^+$. it is easy to check that this is an ordered ring with positive cone $\Cal E(H,\phi)^+$; for  $\rho \in \Cal E(H,\phi)$, then $\rho \in \Cal E(H,\phi)^+$ iff $\rho(H^+) \subset H^+$.

Obviously, the identity $\I$ belongs to $\Cal E(H,\phi)$. We define the {\it bounded subring\/}
$$
\Cal E_b \equiv \Cal E_b(H,\phi)= \Set{\rho \in \Cal E(H,\phi)}{\exists\ N \in \N \text{ \st $-N\I \leq \rho \leq N\I$}}.
$$
It is easy to check that this is the order ideal of $\Cal E(H,\phi)$ generated by $\I$ (the identity element of $\Cal E(H,\phi)$, and then it follows easily that $\Cal E_b (H,\phi)$ is an ordered ring with $\I$ as order unit. If $\rho \in \Cal E(H,\phi)^+$, belonging to $\Cal E_b (H,\phi)$ is a very restrictive condition: it amounts to showing that there exists a positive integer $N$ \st for every $h \in H^+$, $\rho(h) \leq Nh$. The following properties are elementary [H2; 1.1].  

Let $(D,u)$ be an  unperforated partially ordered abelian group with order unit (this includes all dimension groups with order unit). The subset $\Inf D$, the {\it infinitesimal subgroup of $D$}, can be defined as $$\Set{d \in D}{d  + D^{++} \subset D^{++}},
$$
where $D^{++}$ is the set of order units of $D$. Equivalently, $\Inf D = \cap  \,\ker \tau $, $\tau$ varying over all the traces of $D$ (or we can restrict to pure normalized traces; the intersection of the kernels is the same), and this is the more usual definition. Alternatively, the elements of the infinitesimal subgroup are characterized by the property that $-u \leq Nd \leq u$ for all integers $N$. 

\Lem Lemma \Bone. Let $(H,\phi)$ be a partially ordered abelian group with self-order embedding.
\item{(a)} For any order ideal $J$ of $H$, $\Cal E_b(H,\phi)(J) \subset J$.
\item{(b)} Every pure trace on $(\Cal E_b(H,\phi), \I)$ is multiplicative, $\Inf \Cal E_b(H,\phi) = \cap_{L \in \partial_e S(\Cal E_b(H,\phi), \I)} \ker L$ and is a convex two-sided (ring) ideal of $\Cal E_b(H,\phi)$, and $\Cal E_b(H,\phi)/\Inf \Cal E_b(H,\phi)$ is commutative.
\item{(c)} If $H$ admits an order unit $v$, then for every pure trace $\tau \in \partial_e S(H,v)$ and $\psi \in \Cal E_b(H,\phi)^+$ there exists unique $L \in \partial_e S(\Cal E_b(H,\phi),\I)$ \st $\tau\circ \psi = L(\psi)\cdot \tau$.

Now replace $(H,\phi)$ by $(A_f,\SS)$ and form the corresponding bounded subring, which we will temporarily denote $C_b (A_f)$ (the $C$ refers to centralizer of the shift  $\SS$---recall that $\SS([a,k]) = [a,k+1]$); it is easy to check that this is exactly the same as $\Cal E_b(R_f,\SS)$, so we will denote it $C_b(R_f)$. Since $R_f$ has an order unit, (c) of the preceding can be applied. In general, $\Cal E_b (R_f,\SS)$ is small, in that it is commutative modulo an ideal generated by   infinitesimals.

\Lem Example \Btwo. Let $G = \Z \times S_3$; let $\Z = \langle x\rangle$, $\sigma = (1,(123))$ and $g = (1,(12))$, and set $f = 1 + x +x^{-1} + \sigma + g$. Then $f$ is admissible, $\RR_{\sigma,2}$ is a bounded endomorphism, but $\tilde l_S(\sigma) = 1$ and $\RR_{\sigma,1}$ is not bounded.

\Pf Clearly $\sigma \in \supp f$, so $\tilde l(\sigma) \leq 1$. We first show $\tilde l(\sigma^2) > 1$ (which forces $\tilde l (\sigma^2) = 2$). We show that for all $m$, $f^m \sigma^2 \not \prec f^{m+1}$. As usual, look at the monomials (in $x$) appearing the coefficient of $\sigma^2$ on left and right. On the left, these are the monomials appearing in the$1$ term of $f^m$, that is, $\Set{x^i}{|i| \leq m}$. On the right, these are the monomials appearing in the coefficient of $\sigma^2$ in $f^{m+1}$---since at least two multiplications that don't involve $x^{\pm1}$ are necessary, the monomials consist of $\Set{x^i}{|i| \leq m+1-2 = m-1}$.

Hence $f^m \sigma^2 \not\prec f^{m+1}$, and in particular, $[{}{\sigma^2},1] \not\in R_f$. If $\RR_{\sigma,1}$ were bounded, then on applying it to $[g,1]$, we would obtain $\RR_{\sigma,1}([g,1]) \leq N [g,1]$, that is, $f^m g\sigma \prec f^{m+1}g$, which amounts to $f_m g\sigma g^{-1} \prec f^{m+1}$. Since $g \sigma g^{-1} = \sigma^2$, the last paragraph proved this is impossible. Thus $\RR_{\sigma,1}$ is not bounded.

Since $\max \brcs{\tilde l(\sigma),\tilde l (\sigma^2)} \leq 2$ (and thus equals $2$), it follows that $\RR_{\sigma,2}$ is bounded. It is easy to check that $\tilde l (\sigma) \neq 0$.
\qed

\Lem Lemma \Bthr. Let $(g,k) \in  G \times \Z^+$. 
\item{(a)} As a positive
endomorphism of $A_f$, $\RR_{g,k}$ is locally order-bounded iff 
$$
k \geq \sup
\Set{\tilde l_S (hgh^{-1})}{h \in G}.
$$
\item{(b)}Local order-boundedness of $\RR_{g,k}$ as an endomorphism of $A_f$ implies it is a
positive endomorphism of $R_f$ and locally order-bounded as an endomorphism
thereof.
\item{(c)} If $\Arrow \phi; R_f. R_f$ is a   positive endomorphism commuting
with $\SS$ (as an endomorphism of $R_f$)
, then it can be extended uniquely to a positive endomorphism $\tilde \phi$ of
$A_f$ commuting with $\SS$, and if $\phi$ were locally order-bounded, then so
would be $\tilde \phi$.
\item{(d)} If $g$ has only finitely many conjugates in $G$, then $\RR_{g,k}$ is order bounded if and only if $k \geq \sup
\Set{\tilde l_S (hgh^{-1})}{h \in G}.$
\item{(e)} If $(G,S)$  satisfies WC, then local order-boundedness of $\RR _{g,k}$
implies order-boundedness, and also implies that $g$ has only finitely many conjugates.

\Rmk So if the supremum in (a) is infinite, then $\RR_{g,k}$ is not locally
order bounded for any choice of $l$. On the other hand, if $g$ has only finitely many conjugates, then it follows from (a,d) that $\RR_{g,k}$ is order-bounded, not just locally order-bounded.

\Rmk Concerning (d), for the discrete Heisenberg group (discussed in gruesome detail later in this article), and $f$ a specific admissible element, WC fails---but it is still true that local order boundedness implies order boundedness.

\Pf (a) Assume $\RR_{g,k}$ is locally order-bounded. Applied to $[{}h,0]$,
there exists a positive integer$N_h$  \st $[{}{hg},k] \leq N_h [{}h,0]$. Hence
there exists an integer $m(h)$ \st $S^{m(h)}gh \subseteq S^{m(h) + k} h$; thus,
$S^{m(h)}ghg^{-1} \subseteq S^{m(h) + k}$. This entails $\tilde l(ghg^{-1}) \leq
k$.

On the other hand, if for some $h$, $\tilde l(ghg^{-1}) \leq k$, then
$S^{m}ghg^{-1} \subseteq S^{m + k}$ for some $m$ depending on $h$, and the
argument goes in reverse.

\noindent (b) Since $R_f$ is the order ideal generated by $[1,0]$ any locally
order-bounded endomorphism of $A_f$ is automatically an endomorphism of $R_f$,
and  it is trivial to check that it is locally order-bounded as an endomorphism
thereof.

\noindent (c) Pick $[a,l] \in A_f$; then for all sufficiently large $k$,
$[a,l+k] \in R_f$ (for example, if $\Supp a \subset S^{l+k}$). So we attempt to
define  $\tilde \phi([a,l]) = \SS^{-k}\phi([a,l+k]$ whenever $[a,l+k] \in R_f$.
To see that this is well-defined, suppose to begin with that $[a,l] \in R_f$;
then $\phi([a,l+k]) = (\phi\circ\SS^k)[a,l]$, and since $\phi$ commutes with
$\SS$, we have $\phi([a,l+k]) = (\SS^k \circ \phi)([a,l])$. Since $\SS$ is an
automorphism of $A_f$, $(\SS^{-k}\circ\phi)([a,l+k]) = \phi([a,l])$. For general
$[a,l] \in A_f$, it follows that if $[a,l+k] \in R_f$, and $k' > k$, then
$\SS^{-(k'-k)}\circ \phi([a,l+k']) = \phi(a,l+k)$.

That $\tilde \phi$ is positive is routine, and similarly, it is straightforward
that if $\phi$ is locally order-bounded, then so is $\tilde \phi$.

\noindent (d) Set $z = \sum
g'$ where $g'$ runs over the finite set of conjugates of $g$. Then $z$ is in the
centre of the group ring, and moreover, $[z,k] \in R_f$ (since $\tilde l(g')
\leq l$).

We observe that $\RR_{z,k}$ is an order-bounded endomorphism of $A_f$: pick
$[a,l] \in R_f^+$; then there exists $N$ \st $f^N a \in (\R G)^+$. Now
$\RR_{z,k}[a,l] = [f^N az,k+l+N] = [zf^N a,k+l+N]$. Since $f^m z \prec f^{m+k}$
for sufficiently large $m$, say $f^m z \leq K f^{m+l}$, we have $f^m z(f^N a)
\leq  K f^{m+l}(f^N a)$. Hence $[za, l + k] \leq K[a,l]$, so that $\RR_{z,k}$ is
order-bounded.

Since $\RR_{g,k}$ is positive and order-bounded above (as an endomorphism of
$R_f$) by $\RR_{z,k}$, it follows that $\RR_{g,k}$ is also order-bounded.

\noindent (e) From (a) and WC, $g$ has only finitely many conjugates, and $k$
is at least as large the maximum of $\tilde l_S$ on its conjugates.  \qed

\Lem Corollary \Bfou. Suppose that for some choice of admissible $f \in (A G)^+$,
there exists $b \in (A G)^+$ and a positive integer $k$ \st $\RR_{b,k}$ is
locally order-bounded. If in addition, $\Supp b$ generates $G$ as a group and
$(G,S)$ satisfies WC, then $G$ is central by finite.

\Pf If $g \in \Supp b$, then $\RR_{g,k}$ is bounded above by a multiple of
$\RR_{b,k}$, hence is itself locally order-bounded. By the preceding, $g$ has
finitely many conjugates. Hence $g^{-1}$ has only finitely many conjugates, and
this property is preserved by products; thus $G$ is   an FC-group. A finitely
generated FC-group is central by finite. \qed

There are significant positive endomorphisms that do not commute with $\SS$.
For example, if $A = \R$ and $\tau \in F_0$, then $\tau$ composed with the
inclusion $\R \to \R \cdot 1 \subset \R G$ is a positive endomorphism with $\tau
\circ \SS = 0$.

In general, not every positive endomorphism of $R_f$ or $A_f$ is of the form
$\RR_{a,k}$---for many groups, any nontrivial linear combination of nonzero endomorphisms in
$\brcs{\RR_{b(i),k(i)}}$ is nonzero when the $k(i)$ are distinct, and is not
of the form $\RR_{a, k}$. This leads to the next pair of definitions.

\def\rr{{\rho}}

We say the pair $(G,f)$ (where $f$ is admissible) {\it satisfies EP\/} if for all $(a,k), (a',k') \in \R G
\times \Z^+$,
$$
 \RR_{a,k} = \RR_{a',k'} \text{ implies either }\cases a= a' = 0 &
\text{or}\\ a = a' \text{ and }k = k'.\endcases
$$ 
And {\it $G$ satisfies EP\/} if  $(G,f)$ does for every admissible $f \in AG$. 

We say that $(G,f)$  {\it satisfies EEP\/} if for any finite subset $J \subset \Z^+$, and every finite set of nonzero elements $\brcs{a_j}_{j \in J} $ with $a_j \in AG$, the set of right multiplication operators $\brcs{\RR_{a_j,j}}_{j \in J}$ is linearly independent as a subset of $\End A_f$. And $G$ satisfies EEP if $(G,f)$ does for every admissible $f$. 

First note that $[a,k] = 0$ (in $R_f$) means there exists $m$ \st $f^m a = 0$. This can happen (with admissible $f$) whenever there is a non-identity torsion element in $G$. For example, let $f_0$ be any admissible element, and let $h$ be an element of $G$ with $h^n = 1$ for some $n > 0$. Set $f = f_0\cdot \sum_0^{n-1} h^j$; as both terms in the product have only positive coefficients, $f$ is admissible. But $a = 1-h$ is killed by $f$. (It is also possible to have $f$ being a right zero divisor but not a left zero divisor, e.g., if $G$ is a nontrivial free product of two groups not both $\Z_2$, and $G$ has an element of finite order.) In particular, a necessary condition for no admissible $f$ to be a zero divisor in $A G$ is that $G$ have no elements of finite order.

We observe that failure of EP for $(G,f)$ is equivalent to the following property:

\itemitem{}\noindent There exist distinct $a,b \in A G$, not both zero, and nonnegative integers $l > k$  \st for all $g \in G$ there exists a nonnegative integer $m(g)$ \st 
$$
f^{l+m(g)} g a = f^{k+ m(g)} g b. \tag 1
$$

The failure of EEP is equivalent to the following:   
 
\itemitem{}\noindent There exist  $a(i) \in A G$ ($i =1,2,\dots , n$), not all zero,  \st for all $g \in G$ there exists  a nonnegative integer $m (g)$ \st 
$$
\sum_{i=0}^n f^{i+m(g)} g a(i) = 0.  \tag2
$$
 
 For the EP property, we first observe that if there is an equation of the form $f^m z = f^k z'$ where $f,z,z' \in AG$, $k < m$,  and $z,z'$ are central, then $\RR_{a,k} = \RR_{b,m}$, which violates EP if not both $a,b$ are zero. (The converse holds, and we will deal with it soon.)
 
We follow the now-classical development of properties of group rings ([P; chapter 4]). Recall the two characteristic subgroups of $G$ defined earlier, 
$$\eqalign{
 \Delta \equiv \Delta(G) & = \Set{g \in G}{g \text{ has only finitely many conjugates}}\cr
 \Delta^+ & = \Set{g \in \Delta}{g \text{ has finite order.}} 
}$$
 A first observation concerning EP, analogous to that of determining when $[a,l ] = 0$ (and almost equally trivial), is that if $\Delta^+$ is not trivial, then there exists an admissible $f$  and a nonzero central element $z$ \st $fz = 0$; so in particular, $\RR_{z,k} = 0$, violating EP. 

 To see this, assume $ \Delta^+\neq \brcs{1}$. By [P, p 118, 1.8], there exists a finite subgroup, $H$, of $\Delta^+$  \st $H$ is normal in $G$. Set $x = \sum_{g \in H} h \in AG$. It is easy to check that  $x$ is central in $AG$ (since $H$ is normal in $G$) and $x^2 = |H|x$. With $m =|H|$, we thus have $x^2 = mx$. Let $f_0$ be any admissible element of $AG$, and set $f = f_0 x$; again, as $x$ has only positive coefficients, $\supp f_0x = (\supp f_0)(\supp x)$, so $f$ is admissible. But $f\cdot (m-x) = 0$ and $m-x$ is central. 

 Generically, admissible $f$ will not be zero divisors, but to obtain results, we usually exclude nontriviality of $\Delta^+$. 
The condition that $\Delta^+$ be trivial is equivalent to primeness of the group ring $AG$ (due to Connell [C]). Here is an elementary sufficient condition for $(G,f)$ to satisfy EP. Its proof is postponed until we recall the usual instruments used in the study of  group rings. 

 \Lem Proposition \Bfiv. Let $(G,f)$ be a group with an admissible element, \st $f$ is not a right zero divisor in $AG$. If $\Delta^+(G)$ is trivial and $G$ is not central by finite, then $(G,f)$ satisfies EP. 

\Lem Corollary \Bsix.   Suppose that $\R G  $ has no
zero divisors and $G$ is not abelian. Then $G$ satisfies EP. 

\Pf Let $f \in \R G$ be admissible. Since there are no zero divisors, $\Delta^+$ is trivial and $G$ is torsion-free. A torsion-free central by finite group is abelian, so $G$ is not central by finite, and now Proposition \Bfiv\  applies.
 \qed 

For this and other, related results, we remind the reader of the central expectation developped by M~ Smith [S] in characteristic zero, reminiscent of the corresponding central expectation  in von~Neumann algebras. This is presented in [P; pp 124--128]. 

 For each $x \in \Delta$, let $c_x $ be the sum over the conjugates of $x$  (in [P], an unintelligible Fraktur character was used), and let $|c_x|$ be the number of conjugates, that is, the index of the centralizer of $x$ in $G$. The collection $\brcs{c_x}_{x \in \Delta}$ is a basis (when $A = \R$) for $Z(\R G)$, the centre of $\R G$. Define a map $\Arrow \sharp; \R G. Z(\R G)$ (a huge sharp-like character was used in op\.\,cit, and moreover it acted on the right) via  
$$
 \sharp \(\sum r_x x\) = \sum_{x \in \Delta} \frac {r_x}{|c_x|} c_x.  
 $$
 Crucial  is the following result given in [P], restricted to the field $\R$.

 \Lem Proposition \Bsev. [P; Lemma 2.5, p\,124] Suppose that  $a_i, b_i$ ($i=1,2,\dots, n$) are elements of $\R G$ \st for all $g \in G$, $\sum_i a_i g b_i = 0$. Then each of the following hold. 
\item{(a)} $\sum \sharp (a_i) b_i = 0$; 
 \item{(b)} $\sum_i a_i \sharp (b_i) = 0$; 
 \item{(c)}$\sum_i \sharp (a_i) \sharp (b_i) = 0$.

 In particular, if $\Delta (G)$ is trivial, $\sharp$ takes values in the scalars; groups with $\Delta = \brcs{1}$ are known as ICC groups (every non-identity element has infinitely many conjugates). In this latter case, it is not difficult to reconstruct $\sharp$, since one can show that if $Q $ is a finite subset of $G$, there exists an increasing sequence of finite subsets $T_1 \subset \dots \subset T_n \subset T_{n+1} \subset \dots$ \st for all $q \in Q\setminus \brcs{1}$, the sequence $\( \sum_{t \in T_n} t q t^{-1}/|T_n|\) \to 0$ weakly in W$^*(G)$, which is a factor). 

 \Lem Lemma \Beig. Suppose that $(G,f)$ is a group with an admissible element. 
\item{(i)} If $f \in Z(\R G)$, then $G$ is central by finite. 
 \item{(ii)} Suppose there exist  $z(i) \in Z(\R G)$ ($i=0,1,\dots,n$) \st 
$$
\sum_{i= 0}^n f^i z(i) = 0
$$ and $z(n)$ is not zero and either  has only nonnegative coefficients, or is not a zero divisor in $Z(\R G)$. Then $G$ is abelian by finite. 

 \Rmk The result says that if merely {\it one\/} admissible $f$ satisfies the properties, then $G$ is respectively central by finite, abelian by finite. We will improve these (somewhat). In (ii), when $z(n) = 1$, $f$ is {\it integral\/} over the centre. 

 \Pf (i) Since $f$ is central, it follows that $\supp f \subset \Delta$. As $G = \cup \supp f^n$, we have that $G = \Delta$. Hence $G$ is an FC group, and a finitely generated FC group is central by finite. 
 
 \noindent (ii) First, we consider the case $z(n)  = 1$. From the equation $f^n = - \sum_{i=0^{n-1}} z_i f^i$ where $z_i$ are central, we have that for all $m$, we can write $f^m = \sum_{i=0}^{n-1} A(i,m)f^i$, where $A(i,m)$ are central, in fact, polynomials in $\brcs{z_i}$. Since $f$ is admissible, taking supports, we see that the number of cosets of $\Delta$ in $G$ is bounded above by $\left|\cup_{i=0}^{n-1} \supp f^i\right|$ (since $\supp A_{i,m} \subset \Delta$). This says that $G/\Delta$ is finite. As $G$ is finitely generated and $\Delta$ is of finite index therein, $\Delta$ is itself finitely generated. Being a finitely generated  FC group, $\Delta$ is central by finite. Hence $G$ is (central by finite) by finite, and so is abelian by finite. 

If $z(n)$ has only nonnegative coefficients, then $\supp fz(n) = (\supp f)(\supp z(n))$ and $fz(n)$ has only nonnegative coefficients, and thus is admissible. Multiplying the displayed equation by $z(n)^{n-1}$ exhibits $fz(n)$  as integral over the centre, and so the previous paragraph applies.

Assume that $z = z(n)$ is not a zero divisor in the centre. Then it is not a zero divisor in $A G$ [P]. We have $f^n z = - 
\sum_{i=0}^{n-1} f^i z(i)$. For any  $m > 0$, (left) multiply by $(fz)^m$, and substitute in the obvious way; we obtain $A(i,m) \in Z(A G)$ \st 
$$
f^{m+n}z^{m+1} = -\sum_{i=0}^{n-1} f^i A(i,m).
$$
Since $\supp A(i,m) \subset \Delta$, the only cosets of $\Delta$ that can possibly be represented on the right side are those of $\cup_{i=0}^{n-1}\supp f^i$, and this for every $m$. 

Now pick $g \in G$; there exists $m$ \st $g \in \supp f^{m+n}$. Decompose $f^{m+n} = \sum_{t \in \supp f^{m+n}} \lambda_t t$ according to the $\Delta$ cosets; that is, for each coset of $\Delta$ in $G$, pick a representative $c$, and form $q_c = \sum_{t \in c\Delta} \lambda_t t$. Then $f^{m+n} = \sum_c q_c$. Since $\supp z^{m+1} \subset \Delta$, for $c \neq c'$ (representatives of distinct cosets of $\Delta$, we have $\supp q_{c'} z^{m+n} \cap \supp q_c z^{m+n} $ is empty. 

Pick $c_0$ \st $g \in c_0 \Delta$. Then $q_{c_0} z^{m+1} \neq 0$ (since $z$, and therefore all of its powers, are not zero divisors in $\R G$), and moreover, $\supp q_{c_0} z^{m+1} \subset c_0 \Delta = g \Delta$. So there exists $g' \in \supp q_{c_0} z^{m+1} $ \st $g' \Delta = g\Delta$. 

For $c \neq c_0$, $\supp q_{c_0} z^{m+n} \cap \supp q_c z^{m+n} = \emptyset$, and thus $g' \in \supp  f^{m+n} z^{n+1}$. Hence $g' \in \cup_{i=0}^{n-1}\supp f^i\cdot \Delta$. But this says that the coset $g\Delta = g'\Delta$ is represented by one of the (uniformly in $m$) finite set of cosets on the (with $z(n)= 1$) that this forces $\Delta$ to be central by finite, hence $G$ is abelian by finite. \qed

 \Pf (Proposition \Bfiv) Suppose that $\RR_{a,k} = \RR_{b,l}$ with $k < l$ and  both $a,b \in A G$ and not both $a$ and $b$ are zero. Applied to $[g,0]$  for $g \in G$, we obtain $[gb,l] = [ga,k] = [f^{l-k}ga,l]$. Hence for each $g$, there exists a nonnegative integer  $m(g)$ \st $f^{m(g)}gb = f^{m+l-k}ga$. Since $f$ is not a right zero divisor, $gb = f^{l-k}ga$ for all $g \in G$. This rewrites as $f^{l-k}ga - gb =0$  for all $g \in G$. Pick $h \in \supp a$, so that $1 \in \supp ah^{-1}$. Apply \Beig\ to the equations $f^{l-k}gah^{-1} -gbh^{-1} = 0$ for all $g$; this yields central $z_1 = \sharp (ah^{-1})$ and $z_2 = \sharp (bh^{-1})$ \st $f^{l-k}z_1 = z_2$. Moreover, since $1 \in \supp ah^{-1}$, it follows that $z_1 \neq 0$. 

 Since $\Delta^+$ is trivial, $z_1$ is not a zero divisor in $Z(A H)$, and  this implies $z_1$ is not a zero divisor in $A G$. Since $z_2$ is central, for all $w \in A G$, we have $wf^{l-k}z_1 = w z_2 = z_2 w = f^{l-k} z_1 w = f^{l-k}w z_1$; thus $(wf^{l-k}- f^{l-k}w)z_1 = 0$. Since $z_1$ is not a zero divisor, $f^{l-k}$ is thus in the centre. As $f$ is admissible, so is $f^{l-k}$, and thus $G$ is central by finite, contradicting the hypotheses. \qed

 If we weaken the nonzero divisor hypothesis, but strengthen considerably the group properties, we obtain a corresponding EEP result for ICC groups ($\Delta = \brcs{1}$). For an element $s$ of a ring, we denote the right annihilator of $s$ (the set of elements in the ring right multiplication by which kill $s$) by $\rr (s)$.  

 \Lem Proposition \Bnin. Suppose that $(G,f)$ is a group with an admissible element \st $\Delta^+(G)$ is trivial, and suppose there exists a nonnegative integer $t$ \st $\rr(f^{t}) = \rr(f^{t+1})$. Then either $(G,f)$ satisfies EEP or $G$ is abelian by finite. 
 
 \Rmk The annihilator condition holds if, for example, $f$ is not a right zero divisor ($t = 0$), or $f = f^*$ (that is, $f$ is symmetric under $g \mapsto g^{-1}$) with  $t =1$ (more generally, this holds if $f^*f = ff^*$), or if $\R G$ is a Goldie ring (this implies that all chains of annihilators are of uniformly bounded length in the presence of $\Delta^+ = \brcs{1}$). Sufficient conditions for $A G$ to be Goldie are given in [P; pp 609, 611--612, and other places]; I do not know how up-to-date these are. I suppose it is plausible that if $G$ is torsion-free, then at the very least the annihilator condition holds, but this probably leads to the zero divisor conjecture. 

\Rmk If $G$ is abelian by finite and with trivial $\Delta^+$ (and sometimes  with nontrivial $\Delta^+$), then EEP fails for every admissible $f$. 

 \Pf Suppose $\sum_{i=0}^n \RR_{a(i),i} = 0$, where not all  $a(i) $ in $AG$ are zero. Let $k$ be the smallest $i$ \st $a(k) \neq 0$. Applied to $[g,0] $, we have $\sum_k^n [ga(i),i] = 0$; this  translates to $\sum_k^n [f^{n-i}g a(i), n] = 0 $. Hence there exists a nonnegative integer  $m(g)$ \st $f^{m(g)}\sum_k^n f^{n-i}g a(i) = 0$. The annihilator chain condition implies that for all $g \in G$, $f^t\sum_k^n f^{n-i}g a(i) = 0$. Pick $h \in \supp a(k)$, and multiply the last equation on the right by $h^{-1}$; this rewrites to 
 $$
 f^{t+n - k}g a(k)h^{-1} + f^{t+n - k-1}ga(k+1)h^{-1} + \dots + f^tg a(n)h^{-1} = 0 \qquad\text{for all $g \in G$.}
 $$
Applying $\sharp $  and Proposition \Bsev, we obtain the equation, $f^m z = - \sum_{i < m} f^i z_i$, where $z = \sharp (a(k)h)$ and $z_i$ are central in $\R G$. Since $1 \in \supp ah^{-1}$, it follows directly from the definition of $\sharp$ that $z \neq 0$. 

Since $\Delta^{+}$ is trivial, $\Delta$ is torsion-free abelian; hence $z$ is not a zero divisor in $\R \Delta$, and thus not a zero divisor in $Z(\R G)$. Now Proposition \Beig\ applies. 
\qed

 \Lem Lemma \Bten. Suppose that $G$ has an admissible $f$ satisfying an equation of the form, 
$$
 f^n g a(n) + \sum_{i=0}^{n-1} f^i g a(i) = 0 \quad \text{for all $g \in G$,}
 $$ 
with $a(i)\in \R G$.
 If the two-sided ideal of $\R G$ generated by $a(n)$ contains a nonzero element with no negative coefficients, then $G$ is abelian by finite. 

 \Pf The hypothesis says that there exist a finite collection $\brcs{c(i), d(i)} \subset \R G$ \st $p = \sum c(i)a(n)d(i)$, where $p = \sum_{h \in \supp p} r_h h$  and $r_h$ are positive real numbers. By expanding the $c(i)$ and $d(i)$, we obtain an extravagantly large sum, $p = \sum r_{j}g_j a(n) h_j$ where $g_j, h_j \in G$ and $r_j \in \R$. We have $\sum_0^n f^i gr_j g_j a(i)h_j = 0$ for all $g$, taking linear combinations of these over $j$, we obtain $f^n g p + \sum_0^{n-1} f^i g b(i) = 0$ for some $b(i) \in \R G$ and all $g \in G$. Select $h_0 \in \supp p$; as before,  we have $f^n g p h_0^{-1} + \sum_0^{n-1} f^i g b(i)h_0^{-1} = 0$. 

 Applying \Bsev(b), we have $f^n \sharp (ph_0^{-1}) + \sum_{0}^{n-1} f^i z_i = 0$ where $z_i$ are central. Since $1 \in \supp ph_0^{-1}$, $z := \sharp (ph_0^{-1})$ is a nonzero central element; it also has no negative coefficients. Now Lemma \Beig(ii)  applies, so that 
 $G$ is abelian by finite. \qed  

 The displayed equation is precisely what we obtain from an equation of the form $\sum_{0}^m \RR_{a(i),m} = 0$ for some $m$ (possibly differing from $n$, as extra powers of $f$ may be necessary) when the chain condition  on annihilators of powers of $f$ holds.

The conditions on annihilators in the EP and EEP results are annoying, because they are not generally easy to verify. 
We also obtain similar results for the simpler problem,  when does $\RR_{a,l} = 0$ imply that $a = 0$? Necessary is obviously $\(\cup_n \rho (f^n)\) \cap Z(A G) = \brcs{0}$, and this is sufficient if there is an integer $n$ \st $\rho(f^n) = \rho(f^{n+1})$.

\SecT  \IP\ Intersection property

 The pair $(G,f)$, where $G$ is a group and $f$  is an admissible element of $AG$, satisfies IP ({\it intersection property\/}) if 
$$
 \bigcap_{n=1}^{\infty} \SS^n R_f = \brcs{0} .
 $$
If for all admissible $f$, $(G,f)$ satisfies IP, then we say $G$ satisfies IP. For groups with torsion elements, e.g., $D_{\infty}$, there typically exist admissible $f_0$, $f_1$ \st $(G,f_0)$ satisfies IP, while $(G,f_1)$ does not. 
 All  finite groups (including the one-element group) fail to satisfy IP for drastic reasons, since if $G$ is finite, then $A_f = R_f$ and so $\SS$ acts as an automorphism on $R_f$. On the other hand, all  torsion-free abelian groups satisfy IP. It is conceivable that all infinite torsion-free groups satisfy it; but I could only prove IP when $G$ is left orderable, a fairly strong property. Even then, the proof is tedious.

 The notion of IP is suggested by Jacobson's conjecture in commutative ring theory (intersection of powers of the Jacobson radical), which was eventually proved, and even generalized  to some right and left noncommutative noetherian rings. 

A group is {\it left orderable\/} [P; p\,586---there called right orderable] if there exists a total ordering on $G$ \st for all $x,y,z \in G$ with $y < z$, we have $xy < xz$. 

\Lem Theorem \IPone. If $G$  is left orderable, then for all admissible $f$, $\cap \SS^n R_f = \brcs{0}$, that is, $(G,f)$ satisfies IP. 

 \Rmk In the course of the proof, there is a possibility of inadvertently assuming the ordering is bi-invariant (that is, $G$ is orderable), but I'm fairly sure that I managed to avoid this. 

 \Pf Suppose $[a,k] \in \cap \SS^n R_f$ for some $a \in AG$ and nonnegative integer $k$. Then $[f^m a,0] \in \SS^n R_f$ for all nonnegative integers $m$, since $\SS^{-1}$ is an order automorphism of $A_f$. Conversely, $[f^m a, 0] \in \SS^n R_f$ for all $n$ iff $[a,k]$ does. 

 So assume $[a,0] \in \cap \SS^n R_f$ for all nonnegative $m$. Let $x$ be the maximum (\wrt the hypothesized left invariant total ordering on $G$) of the inverses of the elements in $\supp a$. That is, $x \geq z^{-1}$ for all $z \in \supp a$, and $x^{-1}\in \supp a$. On left multiplication by any $z \in \supp a$, we obtain $zx \geq 1$. Set $a' = ax$, so that all elements in the support of $a'$ are at least as large as $1$. We claim that $[a',0] \in  \SS^n R_f$ for all $n$. 

 Set $m = \tilde l(x)$.  As $\RR_{x, m}$ is an order preserving endomorphism of $R_f$ commuting with $\SS$, it follows that it is an endomorphism of  $\SS^n R_f$. Hence $[a',m] \in   \SS_n R_f$. Thus $[a',0] \in  \cap \SS_n R_f$. 

 Now set $z_0$ to be the maximum of the elements of $\supp a'$. Then $z_0 \geq z$ for all $z \in \supp a'$. For each positive integer $k$, $[f^k a',0] \in R_f$ (since it belongs to $\SS^n R_f$, and the latter contained in $R_f$). Hence for each $k$, there exists a positive integer $m(k)$ \st $\supp f^{m(k) + k }a' \subset \supp f^{m(k)}$. Fix $k$, and pick $h \in \supp f^{m(k) + k}$. Then $hz_0 \geq hz$ for all $z \in \supp a$.  Among the elements of the form $\brcs{hz_0}_{h \in \supp f^{m(k)} + k}$, there is a maximal element, say given by $h_0 z_0$ (we are not assuming that $h_0 = \max \brcs{h}$; this would not be useful unless the ordering was also right invariant, which it need not be). Then the coefficient of $h_0 z_0$ in the product $f^{m(k)+ k}a'$ is not zero: if it were zero, there must be another pair $(h',z') \in \supp f^{m(k)+ k} \times \supp a'$ with $h_0 z_0 = h'z'$. We have $h'z_0 \geq h'z'$ (by left invariance, as above), so that $h_0 z_0 \leq h'z_0$; but the definition of $h_0$ ensures equality, and thus $h' = h$, and therefore $z_0 = z'$, a contradiction. 

 Hence $h_0 z_0 \in \supp f^{m(k)}$; since $\supp f^{m(k)} \subset f^{m(k)+ k}$, it follows that $h_0 z_0 \geq (h_0 z_0)z_0 = h_0 z_0^2$ (using the definition of $h_0$). However, $z_0 \geq 1$ (since all the elements of $\supp a'$ are at least as large as $1$), so left invariance implies $h_0 z_0^2 \geq h_0 z_0$. We thus have $h_0 z_0^2 = h_0 z_0$, whence $z_0 = 1$. But this forces either $\supp a' = 1$ or $a' = 0$; if the former, then    $\supp f^{m(k) + k}a = \supp f^{m(k)+ k}$, and this is contained in $\supp f^{m(k)}$. Since $G$ is infinite, this is impossible. So $a'=0$, and thus $a = 0$. \qed 

 (Left orderability implies $AG$ has no zero divisors, and thus $[a,k] $ being zero entails  $a = 0$, explaining the drastic conclusion.) 

When $G$ admits torsion elements, there exist admissible $f$ \st $(G,f)$ does not have IP; this is essentially for trivial reasons. 

\Lem Lemma \IPtwo. Let $G$ be a group. 
\item{(a)} If $G$ has torsion, then there exists admissible $f_0$ that is a right zero divisor in $AG$.
\item{(b)} Suppose  $f_0$ is an admissible element for $G$ and is a right zero divisor in $AG$; if $f = 1+ f_0$, then    $\cap \SS^n R_f $ is not zero, that is, $(G,f)$ does not satisfy IP. 

\Pf (a) Suppose $f_1$ is any admissible element of $AG$, and let $\theta$ be an element of order $n > 1$. Set $f_0 = f_1 \cdot (1 + \theta + \dots + \theta^{n-1})$; this is admissible, and $f_0 \cdot (1 -\theta) = 0 $.  

\noindent (b) Suppose $f_0\cdot m = 0$ for some nonzero $m \in AG$. Then $fm = m$, and thus $f^n m =m$ for all nonnegative integers $n$. This entails $[f^n m,0] \in R_f$ for all $n \geq 0$: there exists $k$ \st $\supp m \subset f^k$, so $\supp f^n m = m \subset f^{k}$ and thus $[f^{n-k}m,0] \in R_f$ for all $n \geq k$. Hence $[m,0] \in \cap \SS^n R_f$. Moreover, $[m,0]$ is not zero in $R_f$, since $f^n m \neq 0$.
\qed

If $G$ is also nilpotent and $\theta$ is a commutator (this might not be necessary), the element $[1-{}{\theta},0]$ constructed in the proof is also a nonzero element of $\Inf R_f$. 

 \Lem Example \IPthr. A group $G$ with   admissible elements, $f, f_1$, \st $(G,f)$ satisfies IP, but $(G,f_1)$ does not. 

\noindent Set $G = D_{\infty}$ with generators $\brcs{x,\theta}$ subject to $\theta^2 = 1$ and $\theta x \theta = x^{-1}$. With $f = (1 + x + x^{-1})(1+ \theta)$, we see that there is an order isomorphism $A_f  \to A_F$ intertwining the shifts sending $[1,0] \in A_f$ to $[1,0] \in  A_F$, where $F = 2(1 + X + X^{-1})$ with $X$ being the generator of $\Z$.  In particular, $\cap \SS^n R_f = \brcs{0}$. On the other hand, whenever a group $G$ has torsion, there exists an admissible $f$ \st $\cap \SS^n R_f$ is nonzero. Lemma \Beig(b) applies to give the existence of $f_1$. 
\qed

There may be a connection between $\Inf R_f \neq \brcs{0}$ and $\cap \SS^n R_f \neq \brcs{0}$; however, for the discrete Heisenberg group $H_3$ and the standard choice of admissible $f$, $\Inf R_f$ is not zero (as we will see in section \One) but $\cap \SS^n R_f$ is zero (since $H_3$ is orderable).

\SecT \TR\  Traces and harmonic functions

Let $U$ be a partially ordered abelian group with positive cone $U^+$; a {\it trace\/} on or of $U$ 
 is a nonzero positive group homomorphism $\Arrow \tau; U. \R$; {\it positivity\/} means $\tau(U^+) \subset \R^+$. If $U$ is a partially ordered real vector space, then $\tau$ is automatically real linear. If $U$ has an order unit $u$ (that is, a positive element \st for all $x \in U$, there exists a positive integer $N$ \st $-N u \leq x \leq Nu$), then we say a trace $\tau$ is {\it normalized at $u$} if $\tau (u) = 1$. The set of  traces normalized at the order unit $u$ is denoted $S(U,u)$, and using the natural embedding $S(U,u) \subseteq \R^{U}$, \wrt the weak topology, $S(U,u)$ is a compact convex set. The set of extreme points of $S(U,u)$ is denoted $\partial_e S(G,u)$, and inherits its topology from $S(G,u)$. Terminology varies: extreme points may also be called any of pure, ergodic, irreducible, indecomposable, minimal, \dots. We use both extreme and pure.  

When $U$ is a dimension group (as will be the case for all of $R_f$ and its order ideals), then $S(U,u)$ is a Choquet simplex. 

When $U$ does not have an order unit, we can still topologize the cone of all traces, and the extreme points correspond to extreme rays of points. 

A trace on $U$ is {\it faithful\/} if $\ker \tau \cap U^+ = \brcs{0}$; we call it  {\it unfaithful\/} or {\it perfidious\/} otherwise. For a trace $\tau$,  the sum of all the order ideals $\ker \tau$ contains is itself an order ideal (true in any dimension  group), so  it contains a largest order ideal, denoted $\ker^+ \tau$; this is generated as an abelian group by the set of positive elements killed by $\tau$. The trace  $\tau$   is unfaithful precisely when  $\ker^+ \tau$ is nonzero.  

In our situation, $U = A_f$ or $R_f$, and in the latter case, we usually take as its order unit $u = [1,0]$. The question is then to describe the (pure) traces. On $A_f$, the traces correspond exactly to harmonic functions on the improper space-time cone $G \times \Z^+$, while those of $R_f$ correspond to the harmonic functions defined on the subcone generated by  $(1,0)$. It is well known (and we will provide yet another proof) that the pure faithful traces of $(R_f,u)$ are precisely the pure  traces that extend to traces $A_f$ (and are parts of extremal rays of traces thereon), and moreover, the extension is unique, and this yields all the pure rays of traces on $A_f$. It is also well known (in the more general context of Pruitt's theorem) that the pure traces on $A_f$ correspond to left positive eigenvectors of the operator, right multiplication by $f$.

When $G = \Z^d$, i.e., torsion-free abelian (all groups are finitely generated in this article), the pure rays of traces on $A_f$ are easy to describe: they are point evaluations at points of the positive orthant, $(\R^d)^{++} = \Set{x= (r_i) \in \R^d}{x_i > 0}$. The extremal normalized traces of $R_f$ can be described as the closure of the image of this set under the natural moment map corresponding to $f$; this describes the pure trace space ($\partial_e S(R_f, [1,0])$) as the Newton polytope of $f$, with the faithful ones in the interior of this polytope, the unfaithful ones on the boundary. For a lot more details, see [H1, H2]. 

In this particular case, the set of faithful pure (normalized) traces of $R_f$ was dense in the pure trace space. This phenomenon hardly holds at all in the nonabelian case. In the case of $G = H_3$ (the simplest discrete Heisenberg  group), in sections \One--\FTN, we will completely describe $\partial_e S(R_f, [1,0])$, including its topology, when $f$ is the most obvious choice of admissible element. It turns out to be surprisingly complicated, and an early observation is the the set of faithful pure traces is  not dense in the set of pure traces. 

An alternative point of view is to consider, for each $\lambda \in \R^+$ (including $\lambda = 0$), the set of traces $F_{\lambda} := \Set{\tau \in S(R_f, \1)}{\tau \comp \SS = \lambda \tau}$ (recall that $\1 = [f,1] = [1,0] \in R_f$); these are each closed faces of $\partial_e S(R_f,\1)$, and $\lambda > 0$ is the eigenvalue corresponding to the (pure)  faithful trace $\tau \in F_{\lambda}$. Then $\partial_e S(R_f,\1) = \cup_{\lambda \in \R^+} F_{\lambda}$. All the unfaithful pure traces reside in $F_0$, and these are the (more) interesting ones. 

Recall from section \B, the notion of order-bounded endomorphism; these are endomorphisms of $A_f$ (and $R_f$ in this case) that commute with the shift $\SS$ and are bounded by a multiple of the identity as operators on $A_f$; they form a ring, denoted $\Cal E_b$. There is another notion of bounded endomorphism of $A_f$, where we do not insist on commuting with the shift, merely 
$$
\End _b (A_f) = \Set{\phi \in \End A_f}{\exists \in  N \in \N \text{ \st} -NI \leq \phi \leq N \I}. 
$$
Of course, the condition means that for $\phi$ to belong to $\End_b (A_f)$, there exists a positive integer $N$ \st  all $x \in A_f^+$,  we have $ \phi(x) \leq Nx$. For many choices of $(G,f)$, it happens that members of $\End_b (A_f)$ automatically commute with the shift, so that in these cases, the ring equals $\Cal E_b$. This is again a partially ordered ring and admits the identity as order unit. (It thus is almost commutative, in the sense that the commutator ideal consists of infinitesimals.)
The following is elementary and well-known in other contexts. Recall that $\I$ is the identity operator on $A_f$ (and by restriction, on $R_f$) and should be distinguished from $\1$, the distinguished element of $R_f$. 

\Lem Lemma \TRone. Let $\tau$ be a pure trace on $(R_f,u)$. There exists a multiplicative (hence pure) normalized trace $ \tau_b \in S( \End_b (A_f), \I) $ \st for all $\phi \in \End_b (A_f)$, we have $\tau \circ \phi =\tau_b(\phi)\cdot \tau$, and $\tau_b (\phi) = \tau(\phi(\1))$. 

\Pf First, assume $\phi \in E_b:=\End_b(A_f)^+$, in particular, $\phi(R_f^+) \subset R_f^+$. Then for all $x \in R_f^+$, $\phi (x) \leq N \phi(x)$. Hence $\tau \circ \phi \leq N\tau$. Purity of $\tau$ yields a nonnegative real number $\lambda_{\phi}$ \st $\tau \circ \phi  = \lambda_{\phi} \cdot \tau$. Evaluating at $\1$, we obtain $\lambda_{\phi} = \tau(\phi(\1))$. Now let $\psi$ be another element of $E^+$; then so is the product, $\phi \psi  = \phi \circ \psi$. We thus have $\tau (\phi \psi) = \lambda_{\phi \psi}\cdot  \tau$. On the other hand, this equals $(\tau \circ \phi) \circ \psi = (\lambda_{\phi}\tau) \circ \psi$, and the latter is of course, $\lambda_{\phi} \lambda_{\psi} \tau$. Since $\tau$ is nonzero, we have $\lambda_{\phi \psi} =\lambda_{\phi} \lambda_{\psi}$. Even easier is $\lambda_{\phi + \psi} = \lambda_{\phi} + \lambda_{\psi}$. 

Since the identity is automatically an order unit of $E$, every element of $E$ is a difference of two positive elements. It follows immediately that $\phi \mapsto \lambda_{\phi}$ (defined on $E^+$) extends uniquely to a positive, multiplicative map $E \to \R$, in particular, a trace. Multiplicativity yields purity (an easy exercise). \qed

Obviously $\SS$ is a bounded endomorphism. The result  almost yields that $\SS$ is in the centre of $\End_b (A_f)$, i.e., that a bounded endomorphism automatically commutes with the shift. We note that $(\phi \circ \SS - \SS \circ \phi)(R_f) $ is killed by every pure trace (as is true for every additive commutator), so if $R_f$ has no infinitesimals, then $\phi$ and $\SS$ commute (or even better, $\End_b (A_f)$ is commutative). For some groups, this is indeed true: obviously for abelian groups, but also for large groups (for which there are no nontrivial bounded endomorphisms other than scalars). However,  for the  discrete Heisenberg group and the standard admissible $f$, not only are there nonzero infinitesimals, but in fact $\End_b (A_f)$ is not commutative. 

Let $\tau$ be a pure trace of $R_f$. Set $\lambda = \lambda_{\SS}$, i.e., $\tau \circ \SS = \lambda \tau$. If $\tau$ is faithful, then $\lambda > 0$. The converse is (not surprisingly) true, and moreover, in that case, $\tau$ extends uniquely to a trace on $A_f$, $\tau_A$, \st $\tau_A \circ \SS = \lambda\tau_A$; moreover, all pure rays of traces of $A_f$ contain such a trace. This is the content of the next (elementary) result. Because we are working with $\SS$ (the inverse of multiplication by $f$, acting on $A_f$), the $\lambda$s go to zero, rather than the more usual infinity. 

\Lem Lemma \TRtwo. Let $\tau$ be a pure trace of $(R_f,\1)$ \st $\lambda \equiv \lambda_{\SS} \neq 0$. 
\item{(a)} Then $\tau$ extends uniquely to a trace, $\tau_A$, on $A_f$ \st $\tau_A \circ \SS = \lambda \tau_A$. 
\item{(b)} $\tau_A$ is a faithful trace on $A_f$ (and thus $\tau$ is a faithful trace on $R_f$). 
\item{(c)}$\tau_A$ is a member of a pure ray of traces on $A_f$. 
\item{(d)} Every pure ray of traces of $A_f$ consists of traces $\tilde\tau$  \st there exists $\lambda > 0$ with $\tilde\tau \circ \SS^k = \lambda^k $ for {\it all\/} $k \in \Z$; up to normalization, these are of the form $\tau_A$ for some faithful pure trace $\tau$ of $R_f$. 
\item{(e)} If $G$ is infinite, then for every $\lambda > 0$, there exists a pure ray of traces satisfying the condition in (d).
 \item{(f)} Let $\sigma$ be a pure trace of $(R_f, \I)$. If $\sigma$ is not faithful, then $\sigma \circ \SS = 0$.

\Rmk In (d), negative powers of $\SS$ also occur. The fact that $\tilde\tau \circ \SS = \lambda \SS$ is equivalent to $\tau$ being a left eigenvector (with no negative coefficients) for the operator $f \times$ (multiplication from the left) on $\R G$, if we view $\tau$ as an element of $(\R^+)^G$ (the set of nonnegative real-valued functions on $G$, equipped with the pointwise topology).

\Pf   We have observed that $\SS$ is an order-automorphism of $A_f$ (in particular, $\SS^{-1}$ is defined and order-preserving),  $A_f = \cup_{n\in \Z^+} \SS^{-n}R_f$, and $A_f^+ = \cup_{n\in \Z^+} \SS^{-n}R_f^+$. 

\noindent (a) We check that the assignment given by $[a,l] \mapsto \lambda^{l-k}\tau([a,k])$ whenever $[a,k] \in R_f$,  is well-defined. If $[a,k], [a,k'] \in R_f$, say with $k' > k$, then $[a,k'] = \SS^{k'-k}([a,k])$, and thus $\tau ([a,k'] = \lambda^{k'-k}\tau([a,k]))$. Hence $\lambda^{l-k}\tau([a,k]) = \lambda^{l-k'}\tau([a,k'])$. We denote the resulting function $\tau_A$. It is clearly a trace on $A_f$ satisfying $\tau_A \circ \SS = \lambda \tau_A$, and uniqueness \wrt this property follows from $A_f = \cup_{n\in \Z^+} \SS^{-n}R_f$. 

\noindent (b) It suffices to show that $\tau$ is faithful. If not, there exists $g \in G$ \st $\tau([g, l]) = 0$, for some $l \geq \tilde l(g)$. By applying $\SS$, we obtain $\tau([g, l]) = 0$ for {\it all\/} $l \geq \tilde l(g)$.   If $h \in \supp f^k$ for some $k$, then $[{}{hg},l+k] \in R_f^+$, and $[{}{hg}, l+k] \prec [g, l]$. Hence $\tau{[{}{hg}, l+k]} = 0$. Setting $h = g^{-1}$ and observing that   $\cup \supp f^k  = G$, we deduce $\tau([1,t]) = 0$ for some nonnegative integer $t$. Since $\tilde l(1) = 0$, this forces $\tau(\1) = 0$, which in turn forces $\tau$ to be identically zero on $R_f$, a contradiction. 

\noindent (c) Now we show that $\tau_A$ is part of an extremal ray of traces. Suppose that $\rho$ is a trace on $A_f$ \st for some positive real number $N$, we have $\rho \leq N \tau_A$ (that is, for all $[a,l] \in A_f^+$, $\rho([a,l]) \leq N\tau_A ([a,l])$). If $\sigma:= \rho|R_f$ is not zero, then $\sigma$ is a trace, and we have $\sigma \leq N \tau$ (we do not have to normalize $\sigma$ in what follows; what is important is that it can be normalized). As $\tau$ is pure, there exists a positive real $\alpha$ \st $\sigma = \alpha \tau$. Evaluating at $1$ yields $\alpha =  \sigma(\1) = \rho(\1)$.

Applying this to $(\rho \circ \SS^s)|R_f \leq N(\tau_A \circ \SS^{s})|R_f = N\tau \circ \SS^{s} = N \lambda^s \tau$ for each integer $s$, we obtain  $\alpha_s $ \st $(\rho \circ \SS^s)|R_f = \alpha_s \tau$; evaluation at $\1$ yields $\alpha_s = \rho(\SS^s(\1))$. Thus
$$
(\rho \circ \SS^s)|R_f =  \rho (\SS^s (\1))\cdot \tau. 
$$
Thus for $[a,l] \in A_f$ \st $[a,k] \in R_f$, 
$$\eqalign{
\rho ([a,l]) & = \rho\circ \SS^{l-k}([a,k])\cr
& =  \rho(\SS^{l-k}(\1))\tau([a,k])\cr
& = \lambda^{k-l}\rho(\SS^{l-k}(\1)) \tau([a,l]). 
}$$
This is true for all $k' \geq k$ (since $[a,k'] \in R_f$). We claim this implies $\rho (\SS^{m}(\1)) = \lambda^m$  for all $m \geq 0$, yielding that $\rho$ is a multiple of $\tau$. 

By (b), $\tau$ is faithful. For each $l \geq 0$, set $a= {}1$; then $\tau([{}1,l]) \neq 0$, and we can take any $k \leq l$. The left side is then $(\rho \circ \SS^l)(\1)$, and the right side is $\lambda^{k-l}\rho(\SS^{l-k}(\1)) \lambda^{l}$. This yields $(\rho \circ \SS^l)(\1) = \lambda^k \rho(\SS^{l-k}(\1)) $ for all $0 \leq k \leq l$. Setting $k =l$ results in $(\rho \circ \SS^l)(\1) = \lambda^l \rho(\1)$. Thus $\lambda^{k-l}\rho(\SS^{l-k}(\1)) = \rho(1)$ for all eligible choices, and in particular, $\rho = \rho(\1) \cdot \tau$. 

Now suppose $\rho(\1) = 0$. Then $\rho|R_f$ is zero. Now consider $\rho \circ \SS^{-l} \leq N \tau \circ \SS^{-l} = N \lambda^{-l} \tau$. Since $\rho$ is nonzero on $A_f$, there must exist an integer $m > 0 $ \st $\rho ([f^m,0]) > 0$ (since the order ideal generated by  $\brcs{[f^m,0]}$ is all of $A_f$). Thus $\rho \circ \SS^{-m}|R_f $ is not zero (its value at $\1$ is $\rho ([f^m,0])$). We deduce $\rho \circ \SS^{-m}|R_f \leq \lambda^{-m} N \tau$, and the result of previous paragraph says that $\rho \circ \SS^{-m}$ is a scalar multiple of $\tau$. Composing with $\SS^m$, we deduce that $\rho$ is a scalar multiple of $\lambda^m \tau$, hence a scalar multiple of $\tau$. 

\noindent (d) Let $\tilde \tau$ be part of a pure ray of traces of $A_f$. If $r = (f,1) > 0$ (the coefficient of the identity of the group element in $f$), then $\SS \leq r \I$ (as endomorphisms of $A_f$. Hence $\tilde \tau \circ \SS \leq r\tilde \tau$; purity entails there exist $\lambda \geq 0$ \st $\tilde \tau \circ \SS = \lambda \tilde\tau$.  We observe that $\tilde \tau$ must be faithful by the same argument (applied to $\tau_A$) in (b). Evaluating at $[1,0] = \1$, we obtain $\lambda = \tilde\tau (\1) \neq 0$. It follows that for all positive $k$, we have $\tilde \tau \circ \SS^k = \lambda^k \tilde \tau$. 

Now consider negative powers of $\SS$. Since $\SS$ is an order automorphism of $A_f$, $\tau_l:=\tilde\tau \circ\SS^l$ is part of a pure ray of traces of $A_f$ for every integer value of $l$, in particular for negative integer values. Applying the preceding to $\tau_{-1}$, we have, for all integers $k \geq 0$, that for some $\lambda_{-1}> 0$,  $\tau_{-1} \circ \SS^k = (\lambda^{-1})^k \tau_{-1}$. This yields $\tilde\tau \circ \SS^{k-1} = \lambda_{-1}^k \tau_{-1}$ for all $k \geq 1$. We $k = 1$, we obtain $\tau_{-1} = \lambda^{-1}\tilde \tau$. In general, $\tau_{-2} = (\tau_{-1})_{-1} $, and we see by induction that $\tau_{l} = \lambda^l\tilde \tau$ for all negative integers $l$.

\noindent (e) This is a simple consequence of Pruitt's theorem (see [H4; Theorem 2.4].

\noindent (f) Since $\SS$ is a bounded endomorphism of $R_f$, $\sigma \circ \SS \leq N \sigma$ for some positive integer $N$; by purity of $\sigma$, either $\sigma \circ \SS = 0$, or there exists $\lambda > 0$ \st $\sigma \circ \SS = \lambda \sigma$. The latter contradicts (b) as $\sigma$ was not faithful. 
\qed


In section \B, we gave necessary conditions for $\RR_{g, l}$ to be a bounded endomorphism (these obviously commute with the shift), and when $(G,f)$ satisfies WC, the conditions are sufficient. We also saw that there are usually {\it lots\/} of endomorphisms, as most reasonable groups satisfy EEP. There is another source of bounded endomorphisms, which we now investigate. 

It is tempting to try to define left multiplication (rather than right multiplication) by elements of $AG$ on $A_f$. In general, these are not well-defined. Let $\LLL_{a,l}$ denote the {\it candidate\/} function ostensibly given by $[b,k] \mapsto [ab,k+l]$. We have an elementary result on when it is defined; essentially, $af = fa$. 

Before dealing with left multiplications, let us revert to the general situation. Let $D$ be a partially ordered abelian group, and $\Arrow \phi; D.D$ a positive endomorphism. Then $\phi$ is order-bounded if $\phi \leq N \I$ as endomorphisms of $D$. We can define the operator norm on the semiring of (positive) bounded endomorphims of $D$, and this is  $\| \phi \| = \inf \Set{\alpha \in \R^+}{\phi \leq \alpha \I}$ (in general, the infimum is not achieved, that is,  the norm itself need not belong to the set on the right). Now assume $D$ has an order unit, $u$. Then we have a double dual representation $D \to \Aff S(D,u)$ given by $d \mapsto \hat d$, where $\hat d(\tau) = \tau(d)$ (for $\tau$ a normalized trace). This is a positive homomorphism of partially ordered abelian groups (with the  obvious ordering on $\Aff S(D,u)$), and we have a norm on the latter, yielding a pseudo-norm on $D$, $\| \widehat d\|  = \sup_{\tau \in S(D,u)} |\tau(d)|$. We can restrict the traces to extremal ones without changing the supremum, and we also have that the supremum is achieved by an extremal trace (this is standard Choquet theory). 

\Lem Lemma \TRthr. Let $(a,l) \in AG \times \Z^+$. 
\item{(a)} The candidate function $\LLL_{a,l}$ is defined if and only if there exists a nonnegative integer $m$ \st $f^{m}\cdot (af-fa) = 0$. 
\item{(b)} When $\LLL_{a,l}$ is defined, it is an endomorphism of $A_l$ commuting with the shift, and $\LLL_{fa,l+1} = \LLL_{a,l}$; if additionally, $a \in AG^+$, then $\LLL_{a,l}$ is a positive endomorphism.
\item{(c)} When $\LLL_{a,l}$ is defined and $[a,l] \in A_f^+$, it is a bounded endomorphism iff $[a,l] \in R_f$, and moreover, $\left\|{\LLL_{a,l}} \right\| = \left\|\widehat{[a,l]} \right\|$, the latter having norm determined by $u = \1$.

\Pf \comment
For (a) and (b), we may assume $l =0$, since $\SS$ is an order automorphism of $A_f$.
\endcomment
\noindent (a) To check definedness, it is necessary and sufficient to show that  $[ab,k+l ] = [afb,k+ l+1]$ for all $(b,k) \in AG \times \Z^+$. If this equality holds, then for there exists $m (b)$ (depending on  $b$) \st $f^{m+1}ab = f^m afb$. Setting $b = {}1$ (the identity element of the group), we deduce $f^{m+1}a = f^m af$.

Conversely, if $f^{m}\cdot (af-fa) = 0$, then $f^{m+1}ab = f^m afb$ holds for all $b$, and thus $[ab,k+l] = [afb,k+l+1]$. 

\noindent (b) The first is a straightforward application of (a); the rest are obvious. 

\noindent (c) Because $\LLL_{a,l}$ commutes with the shift, to show it is   bounded (when defined), it suffices to show there exists a positive integer $N$ \st for all $g \in G$,  $-N[g,0] \leq [ag,l] \leq N [g,0]$. By assumption, there exists $N$ \st $-N[1,0] \leq [a,l] \leq N[1,0]$. Hence there exists an integer $j$ \st in $AG$, we have $-Nf^{j+l} \leq f^j a \leq N f^{j+l}$; right multiplication by $g$ is order preserving, so $-Nf^{j+l}g \leq f^j ag \leq N f^{j+l}g$, and thus $N[g,0] \leq [ag,l] \leq N[g,0]$, as desired.

Conversely, $\LLL_{a,l}([1,0]) = [a,l]$, and thus the latter belongs to the order ideal generated by $\1$, that is, to $R_f$. 

Now if $\LLL_{a,l} \leq \alpha \I$, applied to $\1$, we obtain $[a,l] \leq \alpha [1,0]$, so for any normalized trace, $\tau ([a,1]) \leq \alpha$, and thus $\left\| \widehat {[a,l]}\right\| \leq \left\|{\LLL_{a,l}} \right\|$. 

Suppose that $\beta = \left\| \widehat {[a,l]}\right\|$. Then for all $\epsilon > 0$, for all normalized traces  $\tau$, we have that $\tau ((\beta + \epsilon)[1,0] - [a,l]) > 0$. By [GooH; 4.1], $x_{\epsilon}: = (\beta + \epsilon)[1,0] -[a,l] \in R_f^+$. Hence there exists $m$ \st $f^m ((\beta+ \epsilon)f^l - a) \in (AG)^+$. It follows that $\LLL_{(\beta+ \epsilon)f^{m+l},0} \geq \LLL_{f^m a,0}$ is a positive endomorphism; applying the shift, $(\beta + \epsilon)\LLL_{f^{m+l},m+l} \geq \LLL_{f^m a, m+l}$. By repeated applications of part (b), we have $(\beta + \epsilon)\I \geq \LLL_{a,l}$, and thus $\left\|{\LLL_{a,l}} \right\| \leq \beta + \epsilon$ for all $\epsilon$. 
\qed

In particular, if $\LLL_{a,l}$ is defined and bounded, then for all pure traces $\tau$ of $R_f$, we have $\tau \circ \LLL_{a,l} = \tau([a,l])\cdot\tau$. 
Of course, $\LLL_{f,1}$ is always defined, and is the identity. On the other hand, $\RR_{f,1}$  is usually quite different from the identity, which is $\RR_{1,0}$ (for example, $\LLL_{f,1} \neq \RR_{f,1}$ if $f$ is not a right zero divisor and $G$ is not central by finite). 

Now consider the effect of various endomorphisms on the trace space of $R_f$. If $\tau$ is a normalized trace of $R_f$, and $0 \leq h \prec f^l$, then we may form $\tau\comp \RR_{h,l}$. If $\tau([h,l]) \neq 0$, then $\tau_{h,l} := (\tau([h,l]))^{-1} \tau \comp \RR_{(h,l)}$ is also a normalized trace, and we have a map $\tau \mapsto \tau_{h,l}$ between normalized traces.

From L \TRtwo, the unfaithful pure traces of $R_f$ are precisely the pure traces that do not extend to $A_f$, or equivalently, kill $\SS(\1) = [1,1]$; thus they constitute the extremal boundary of $F_0$ (corresponding to eigenvalue $\lambda = 0$; there is an argument that the $\lambda$s in the preceding should be replaced by their reciprocals, and then $F_0$ would be relabelled $F_{\infty}$). Obviously $G$ acts on each $F_{\lambda}$ with $\lambda > 0$; multiplication by $[g,0]$ is an order automorphism of $A_f$, so the action is simply $\tau \mapsto \tau \circ \RR_{g,0}$. The action may be trivial (as occurs if $G$ is abelian) or faithful (as occurs if $G = F_2$ and we take the standard admissible element $f = 1 + g + g^{-1} + h + h^{-1}$). 

However, there generally is no action of $G$ on $F_0$, only a partial action, that we now describe. 

\Lem Lemma \TRfou. Let $\tau \in \partial_e S(R_f,[1,0])$, and suppose that $g$ is an element of $G$ and    $k \geq \tilde l_S(g)$. If $\tau \circ \RR_{(g,k)}$ is not zero, then it is (after renormalizing) a pure trace of $R_f$.

\Pf First, $\RR_{g,k}$ is a positive endomorphism of $R_f$, so that  $\phi:= \tau \circ \RR_{g,k}$ is either zero or a trace; in the latter case, we can renormalize by dividing by $\tau([g,k])$.

Assuming $\phi$ is nonzero, we apply the purity criterion  of [GoH; Theorem 3.1]. Since $R_f$ is a dimension group, given two positive elements $x,y \in R_f$ and $\epsilon > 0$, there exist $z \geq x,y$ \st $\tau(z) \leq \max\brcs{\tau(x),\tau(y)} + \epsilon$. Pick $[a,l], [b,m] \in R_f$, and consider the elements $x = [ag,k+l]$, $y = [bg,k+m]$ in $R_f$. If we assume that $x,y \geq 0$, then it easily follows that $[a,l], [b,m] \in R_f^+$; we can thus assume (by replacing, for example, $[a,l]$ by $f^r a, r+l$  for suitably large $r$ and relabelling, that $a, b \in (AG)^+$. Applying the purity criterion to $\tau$, given $\epsilon$, there exists $z_0 \geq x,y$ \st $\tau(z_0) \leq \tau (x)+ \epsilon, \tau(y) + \epsilon$. We can write $z_0 = [c,n]$ with $c \in R_f^+$ \st $f^p c \geq f^{p+n-l-k} ag, f^{p+n-m-k} bg$ coordinatewise for all sufficiently large $p$.

For all sufficiently large  $t$, $f^t g \prec f^{t +k}$. If $p$ is chosen large enough, then $f^{p+n-l-k}a g\prec f^{p+n-k}g \prec f^{p+n}$; similarly, $f^{p+n-m-k}bg \prec f^{p+n}$.
Let $c_1$ be the element of the group ring obtained by taking the coordinatewise maximum, that is $(c_1,{}h) = \max\brcs{(f^{p+n-l-k} ag,{}h),(f^{p+n-m-k} bg,{}h)}$ for all $h\in \supp f^{p+n}$; then $f^p c \geq c_1$ coordinatewise, so that $z_0 \geq [c_1,p+n]$. Also $f_p c \geq c_1 \geq f^{p+n-l-k} ag, f^{p+m-m-k} bg$ coordinatewise. Moreover, $c_1 g^{-1} \geq f^{p+n-l-k} a, f^{p+n-m-k} b$, and has the same support as the sum. That is, $\supp c_1 g^{-1} \prec f^{p+n-l}$.

Set $z = [c_1 g^{-1},p+n-l] \in R_f^+$. Obviously $\RR_{g,k}(z) \geq \RR_{g,k}([a,l]), \RR_{(g,k)}([b,m])$. Moreover, $\tau \circ \RR_{(g,k)}([a,l]) , \tau \circ \RR_{(g,k)}([b,m]) \geq \tau(z_0) -\epsilon \geq \tau\circ \RR_{(g,k)} (z) -\epsilon$. This verifies the purity criterion.
\qed

\Lem Corollary \TRfiv. Let $\tau$ be a pure normalized unfaithful trace of $R_f$, and let $g \in G$. 
\item{(a)} If  $\tau(x_g) \neq  0$, then $(\tau(x_g))^{-1}\tau \circ \RR_{g,\tilde l(g)}$ is a pure normalized trace of $R_f$. 
\item{(b)} Suppose that $[g,k] \in R_f$ for some $g \in G$. If either $\tau(x_g) = 0$ or $k > \tilde l_S(g)$, then $\tau \circ \RR_{g,k} = 0$. 

This give a partial  action on $\partial_e F_0$, the set of pure unfaithful traces of $R_f$. Let $\tau$ be one such; define $\Cal H_{\tau} = \Set{g \in G}{\tau(x_g) \neq 0}$ and $\Cal G_{\tau} = \Set{g \in \Cal H_{\tau}}{\tau \circ \RR_{g,\tilde l(g)} = \tau(x_g) \tau}$. Then $\Cal G_{\tau}$ is a subsemigroup of $G$, and if $g, h \in \Cal G_{\tau}$, then $x_g x_h = x_{gh}$, that is, $\tilde l(gh) = \tilde l(g) + \tilde l(h)$. In addition, $\Cal H_{\tau}$ admits a groupoid-like structure (if $g$ and $h$ belong to $\Cal H_{\tau}$, then $gh$ also does provided $\tilde l(gh) = \tilde l(g) + \tilde l(h)$ and $\tau(x_{gh}) \neq 0$). We also see that  $\Cal H_{\tau} \Cal G_{\tau} = \Cal H_{\tau}$. 

The subsemigroup (actually a monoid) $\Cal G_{\tau}$ may contain units (invertible elements), but not very many aside from the identity. If $g,g^{-1} \in \Cal G_{\tau}$, then $\tilde l(g) = \tilde l(g^{-1}) = 0$. Hence the only invertible elements are torsion and belong to $\tilde l^{-1}(0)$. So if $(G,f)$ satisfies WC, there are only finitely many units in $\Cal G_{\tau}$; if instead, $G$ is torsion-free, there are no units except the identity therein. 

We will come back to this construction when we deal with maximal order ideals.

\comment
Let $\tau$ be a trace of $R_f$. We define $\ker^+\tau$ to be the largest order ideal contained in $\ker \tau$. As sums of order ideals are order ideals (for dimension groups), we can equally define $\ker^+ \tau$ to be the sum of all the order ideals that $\ker \tau$ contains. Alternatively, $\ker^+ \tau$ is an order ideal and $\ker^+ \tau \cap R_f^+ = \ker \tau \cap R_f^+$. Obviously, $\tau$ is faithful iff $\ker^+\tau$ is zero.  
\endcomment


\SecT \MI\ Maximal order ideals

We prepare for a surprising structure theorem (Theorem \twoone) for maximal order ideals of $R_f$. An order ideal $J$ of a dimension group is {\it maximal\/} if it is not the whole group and not contained in any other proper order ideal. When the dimension group has an order unit, maximal order ideals exist. In the case that the dimension group is $R_f$, maximal order ideals correspond to maximal space-time subcones of the cone generated by $(1,0)$.

For a partially ordered abelian group $D$ with positive cone $D^+$, we denote the set of order units $D^{++}$  (assuming it has order units). 

\Lem Lemma \MIone. Let $(D,u)$ be a dimension group with order unit. Suppose that $d \in D^+$. Then $d \in D^{++}$ iff for all maximal order ideals $M$, $d \not\in M$.

\Pf As $d$ is in $D^+$, $d$ is an order unit iff the order ideal generated by $d$, $I(d):= \Set{a \in D}{\pm a \prec d}$ is all of $D$. One direction of the lemma is trivial. If $d$ is not an order unit, then $I(d)$ is a proper order ideal; since $D$ admits an order unit, there exists a maximal order ideal $M$ containing $I(d)$, so $d \in M$.
\qed

\Lem Lemma \MItwo. Let $M$ be a maximal order ideal of $R_f$. Then $\SS R_f \subset M$.

\Pf As $R_f/M$ is a simple dimension group, it thus  admits a pure trace; this induces a pure trace on $R_f$, $\tau$, \st $\ker^+ \tau = M$; necessarily,  this $\tau$ admits a positive element in the kernel, and so $\tau$ is not a faithful trace. Since every such pure trace satisfies $\tau \circ \SS = 0$, we are done.
\qed 

Thus $\cap M$ (where $M$ varies over all maximal order ideals of $R_f$) contains ${\SS R_f}$. The reverse inclusion fails in general, as we will see later.

\Lem Lemma \MIthr. Let $M$ be a maximal order ideal in $R_f$.
\item{(i)} Suppose that $g,h \in G$ and $[g, m] \in R_f$ and $[{}h,n] \in M$ for some nonnegative integers $m$ and $n$. Then $[{}{gh},m+n]\in M$.
\item{(ii)} For all $g \in G$ and $m \in \N$ \st $[g,m] \in R_f$, there exists a maximal order ideal $M$ \st $[g,m] \in M$.
\item{(iii)} Let $p \in AG$, $q \in A^+G$ and suppose that  $[p,m] $ and $[q,n]$ are order units in $R_f$. Then $[pq,m+n]$ is an order unit in $R_f$.

\Rmk In (iii), we really require that $q \in A^+$ or something similar, not merely that $[q,n]$ be an order unit. The same proof works if we only assume for all (pure) traces $\tau$, that $\tau \circ \RR_{q,n}$ be traces (or zero).

\Pf (i) There exists $k$ \st $f^k g \prec f^{m+k}$; hence $f^k {}{gh} \prec f^{m+k} {}h$. Thus $[{}{gh},j ] \prec [f^m {}h, j]$ (as elements of $A_f$) for all nonnegative integers $j$. Setting $j = m+n$, we infer $[{}{gh},j] \prec [{}h, n]$; as the latter belongs to $M$, so does the former.

\noindent (ii) Since $z:= [g,m] \in R_f^+$, it suffices to show that $z$ is not an order unit in $R_f$ by Lemma \MIone. If $z$ were an order unit, then $\Bf 1 \prec [g,m]$, whence there exists $k$ \st $f^{m+k} \prec f^k g$. This entails $\supp f^{m+k} \subset (\supp f^k)g$, so that $|\supp f^{m+k}| \leq |\supp f^k|$; as $\supp f^j \subset \supp f^{j+1}$ (strictly), we deduce $m + k = k$, so $m=0$, contradicting $m \in \N$.

\noindent (iii) We first note that $[pq,m+n] \in R_f$; this requires only $[p,m], [q,n] \in R_f$ and $q \in A^+$. By hypothesis, there exist $k,l \in \N$ \st $f^{m+k} \prec f^k p$, $f^{n+l} \prec f^l q$, and $m+k \leq l$; right multiplying the first by $q \in A^+$, we have $f^{m+k}q \prec f^k pq$ (this requires $q$ to be in $A^+$; otherwise, we would obtain an intractable expression such as $f^k p f^{k'}q$), so that $f^{m+k-l} (f^l q) \prec f^k pq$. Then $f^{m+k -l +n+l} \prec f^k pq$, and thus $\Bf 1 \prec [pq,m+n]$. \qed

 In the following,   (iii) implies (i) and (ii), and is implied by a weaker condition.

\Lem Lemma \MIfou. Assume $G$ is infinite. For $g \in G$, the element $[g,n]$ is an order unit of $R_f$ iff the following conditions all hold.
\item{(i)} $n = 0$
\item{(ii)} $g$ has finite order
\item{(iii)} there exists an integer $m$ \st $(\supp f^m)\cdot g = \supp f^m$.

\Pf Suppose $[g,n]$ is an order unit. Then there exists positive integer $K $ \st $\Bf 1 = [{}{1_G},0] \leq K [g,n]$. This implies there exists $m$ \st $f^{m+n} \leq f^mg$, which in turn implies $\supp f^{m+n} \subset (\supp f^m)g$. This entails $|\supp f^{m+n}| = |\supp f^m|$; however, as $G$ is infinite, the sequence $(\supp f^j)$ must increase strictly. Thus $n = 0$ (i), and $\supp f^m \subseteq (\supp f^m)\cdot g$; since the two sets have the same finite cardinality, they must be equal, yielding (iii).

By induction, for any $h \in \supp f^m$, we have that $hg^k \in \supp f^m$ for all $k$. As the latter is a finite set, we must have $hg^k = hg^l$ for some $k < l$, so $g^{l-k}$ is the identity, yielding (ii).

Conversely, suppose (iii) holds. Then there exist positive integers $K$ and $L$ \st $ [f^m g,m] \leq K[f^m,m] \leq L[f^m g, m]$, whence $[g,0] \leq K \Bf 1 \leq L [g, 0]$ yielding both that $[g,0] \in R_f$ and is an order unit therein.
\qed

In general, (i) and (ii) together do not imply (iii) (examples are ubiquitous, e.g., if $G$ is the infinite dihedral group).

If $[g,0]$ belongs to $R_f$ (without assuming it is an order unit), then we obtain $f^k g \leq K f^k$ (in $A^+$), and therefore $(\supp f^k)g \subset \supp f^k$; but again by cardinality of the finite sets, we must have equality. Hence $[g,0] \in R_f$ is necessary and sufficient for $[g, 0]$ to be an order unit in $R_f$.

\comment

Assume the order of $g$ is infinite, and $g \in \Gamma_n$. Let $M$ be a maximal order ideal of $R_f$ not containing $[g,n]$, and let $N = \phi_g^{-1}(M) $; the latter can be rewritten as $\phi_g^{-1}(M \cap \langle [g,n]\rangle )$. Since $[g, n] \not\in M$, it follows that $\Bf 1 \not\in N$. We wish to conclude that $N$ is a maximal order ideal. We observe that $R_f$ is a dimension group, so that $M \cap \langle [g,n]\rangle$, being an intersection of two order ideals, is itself and order ideal. Since $\Arrow \phi^{-1}; \langle [g,n]\rangle . R_f$ exists and is an order isomorphism, we obtain that $N$ is an order ideal, and $R_f/N \iso \langle [g,n]\rangle/(M \cap \langle [g,n]\rangle)$.

Now the map $ \langle [g,n]\rangle \to R_f \to R_f/M$ is positive, has kernel $M \cap \langle [g,n]\rangle$, the latter being an order ideal in $M \cap \langle [g,n]\rangle$. Hence we obtain an embedding $ \langle [g,n]\rangle/(M \cap \langle [g,n]\rangle) \to R_f/M$; moreover the range is an order ideal. As $R_f/M$ is simple, $ \langle [g,n]\rangle/(M \cap \langle [g,n]\rangle) \to R_f/M$ is onto.

So there is a $G$ action on the set of maximal order ideals? $M \mapsto \RR_{g,\tilde l(g)}^{-1} (M)$? If this is the case and the action is nontrivial (it is trivial for abelian $G$), having finitely many order ideals would imply a proper subgroup of finite index.
\endcomment

\Lem Lemma \MIfiv. Suppose $g,h \in G$, and $[g,k] \prec [{}h,k]$. Then $[{}h,k]
\prec [g,k]$ and $gh^{-1}$ has finite order.

\Pf Follows from Lemmas  \wceig\ and \firsteig.\qed

In $R_f$, and for fixed $k$, what is the smallest integer $t$ for which there
exist a $t$-element subset $\brcs{g_i} $ of $G$ with the property that
$[{}{g_i},k] $ all belong to $R_f$ and the order ideal generated by
$\brcs{[{}{g_i},k] }$ is $R_f$? To this end, define $\Arrow t_S \equiv t; \Z^+.\N$
(where $S = \supp f$) via
$$
t(k) = \min\Set{t \in \N}{\exists \text{$\brcs{g_i}_{i=1}^t \subset S^k$ and $N \in \N$ \st
$[{}{g_i},k] \in R_f$ and $[1,0] \leq N\sum [{}{g_i},k]
$}}.
$$

Obviously, $t(k) \leq |S^k|$. For $G = \Z^d$ and $S = \supp f$, we will
verify that $t$ is constant, and equal to the number of extreme points of $\cvx
S$ in $\R^d$. However, when $G$ has exponential growth, $t$ is practically
asymptotic with the growth.

\Lem Lemma \MIsix. $|\tilde l_S^{-1}(k)| \geq t(k) \geq \limsup_m |S^{m+k}|/|S^{m}|$.

\Rmk The upper bound is slightly sharper than $|S^k \setminus S^{k-1}|$.

\Pf Suppose that $x = \sum [{}{g_i},k]$ is an order unit for $R_f$, and suppose that $\tilde l_S(g_1) < k$. Then $y: = [{}{g_1},k] = \SS[{}{g_1},k-1] \in \SS R_f$. Hence $y$ belongs to every maximal order ideal of $R_f$. As a consequence $x-y$ is a positive element of $R_f$ that belongs to no maximal order ideal; hence it must be an order unit. This means that we can throw away all the $g_i$ \st $\tilde l_S(g_i) < k$ and still have an order unit, whence the upper bound.

Suppose that $[{}{g_i},k] \in R_f$ for $i =1,2,\dots,t$, and their sum is
an order unit. Then for all sufficiently large $m$, we have $S^m g_i \subset
S^{m+k}$ and $\cup S^{m+k} = \cup_i S^m g_i$. Hence $|S^{m+k}| \leq \sum |S^m
g_i| = t |S^m|$, and therefore $t \geq |S^{m+k}|/|S^{m}|$.
\qed

The lower bound is interesting only when $G$ has
exponential growth.

In the case of $G = \Z^d$, all $R_f$ have only finitely many maximal order
ideals, and these correspond to the extreme points of the convex hull of $S$.
This is closely related to constant $t(k)$.

\Lem Lemma \MIeig. Suppose $R_f$ has $N$ maximal order ideals. Then $t(k) \leq N$ for
all $k$.

\Pf Since $[1,0] = [f^k,k] \prec \sum_{g \in S^k} [g,k] $, given a maximal
order ideal $M$, there exists $g_M \in S^k$ \st $[{}{g_M} ,k] \not\in M$. Form
$x =\sum_M [g,k] \in R_f$. Then $x \in R_f^+$, and since $x$ does not belong
to any maximal order ideal, it must be an order unit.
\qed

If $G = \Z^d$, then the maximal order ideals are in bijection with the extreme
points of $K :=\cvx S$. For each maximal order ideal $M_v$, $\ker M \cap
\Set{x^w/f^k}{w \in \Log f^k}$ (adopting the notation of [H1, H2]) consists of
everything except $x^{kv}/P^k$. Hence the only choice (for $k$ fixed) for the
$x^w/P^k \not\in M_v$ is $x^{kv}/P^k$. Each of these belongs to all the other
maximal order ideals. Hence $t(k) $
is at least as large as the number of extreme points of $K$, yielding a constant
$t$.

On the other hand, we will see in section \FTN\ that for the Heisenberg group $H_3$, with the natural admissible $f$, $t(m) = 4m^2 + \Oh{m}$, while the corresponding growth is of order $m^4$. 

\comment
\Lem Lemma. Suppose that $t$ is bounded by $N$. Then $R_f$ has at most $N$
maximal order ideals.

\Pf Let $\brcs{M_i}$ vary over $N+1$ order ideals.
\endcomment

\SecT \Two\ More on maximal order ideals

As before,   $\Gamma'_k = S^k \setminus \brcs{ \tilde l_S^{-1}(\leq k-1)}$, so that $R_f/{\SS R_f}$ is isomorphic to the direct limit,
$\lim \Arrow\psi_k; A\Gamma'_k . A\Gamma'_{k+1}$, where the map $\psi_k$ is given by left multiplication by $f$, compressed and restricted; that is, $\psi_k (g) = \sum (f,{}h)hg$ restricted to those $h \in \supp f$ \st $hg \in \Gamma_{k+1}$. For $g \in \Gamma'_k$, we may identify its image in $R_f/{\SS R_f} $ with the image of $[g, \tilde l_S(g)] =: x_g$.

Let $M$ be an order ideal of $R_f$ containing $\SS R_f$; this includes all maximal order ideals, by Lemma \MItwo. We define $\Gamma'_{k,M} = \Set{g \in \Gamma'_k}{[{}{g},\tilde l_S(g)] \not\in M}$. We can realize the dimension group $R_f/M$ as the limit restricted and compressed maps $\lim \Arrow\psi_{k,M}; A\Gamma'_{k,M} . A\Gamma'_{k+1,M}$. We will show the following.

\Lem Theorem \twoone. Suppose $f$ is an admissible element of $A G$, and let $M$ be a maximal order ideal of $R_f$. Then 
\item{(a)} $\sup_k|\Gamma'_{k,M}| < \infty$;
\item{(b)}  $R_f/M$ has unique trace; 
\item{(c)} $R_f/M$ is isomorphic to a stationary dimension group, via a periodic telescoping of $\psi_{k,M}$.

Part (c) means that there exist positive integers $L, k$ \st for each $n \geq L$, the composite maps $\Arrow \psi_{nk+k-1,M} \circ \psi_{nk + k-2,M}\circ \cdots \circ \psi_{nk,M}; A\Gamma_{kn,M}'. A \Gamma_{(n+1)k,M}'$ can be identified with each other. Parts (a) and (b) are  consequences of (c), but we have to prove them first. Examples exist where both the $L$ and $k$ cannot be chosen to be $1$.

 If $G =\Z^d$, the maximal ideals are in bijection with the extreme points of $\cvx S$, and in this case $|\Gamma_{k,M}'| = 1$. If $G$ is free on two or more letters and $f =1+ \sum g_i^{\pm 1} $, then again $|\Gamma_{k,M}'| = 1$ for all $k$ and $M$, and in both cases, if $A = \Z$, the direct limits are of the form $\lim \Arrow n \times; \Z. \Z \iso \Z[1/n]$ as ordered groups, for $n$ being the coefficient of the one element of $g \in \Gamma'_{1,M}$ in $f$. The more complicated construction of section \IMI\ realizes an arbitrary primitive stationary system as a quotient by maximal order ideal in $R_f$ for suitable choices of $f \in (\Z G)^+$.

Towards the proof of the theorem, consider a maximal order ideal $M$. Then $\Gamma'_{k,M}$ is not empty, since $[1,0] = [f^k, k]$ is an order unit, hence not in $M$. This means that for at least one $g \in S^k$, $[g,k]$ does not belong to $M$. Since $\SS R_f \subset M$, it follows that $\tilde l(g) = k$, hence $g \in \Gamma'_k$. Since $[g,k] + M$ is an order unit of $R_f/M$, it follows that there exists an integer $l$ \st there are paths from $g \in \Gamma'_{k,M}$ to every $\gamma \in \Gamma_{k+l,M}'$. That is, for each such $\gamma$, there exists $h \in \Gamma'_{l}$ \st $\gamma = hg$; moreover, this is true for all $l' \geq l$.

We have a (limited) action of some elements of $G$ on  the set of  order ideals. Suppose $M$ is an order ideal and $x_g:= [g, \tilde l_S(g)] \not\in M$ for some $g \in G$ (that is, $g \in \Gamma_{\tilde l_S(g), M}'$). Define a new set, denoted $g\circ M$, via $g \circ M = \RR_{g,\tilde l_S(g)}^{-1}\(M \cap \langle [g,\tilde l_S(g)]\rangle\)$. 

 \Lem Lemma \twosix. Suppose that $M$ is a proper order ideal of $R_f$ and $g$ is a element of $G$ \st $x_g \not\in M$. Then the set
$$
g\circ M:= \RR_{g,\tilde l_S(g)}^{-1} \( M \cap \langle x_g\rangle\)
$$
is a proper order ideal. \item{(a)} If additionally, $\langle x_g \rangle + M = R_f$, then $R_f/M \iso R_f/g \circ M$ as partially ordered abelian groups. 
\item{(b)} If $\tau$ is a pure unfaithful trace and $\tau(x_g) \neq 0$, then $\ker ^+(\tau \circ \RR_{g,\tilde l_S(g)}) = g \circ \ker^+\tau$.  

\Rmk Here $\RR_{g,\tilde l(g)}^{-1} (\cdot )$ means the inverse image for  the map $\Arrow \RR_{g,\tilde l(g)}; R_f.R_f$. If $M$ is maximal, then the comaximality hypothesis of part (a) is redundant, as sums of order ideals in dimension groups are order ideals.

\Pf We can assume that $G$ is infinite (else no proper nonzero order ideals exist). From the definitions, $y:= [c,n]$ (for $c \in A G$ and we may assume that $c \prec f^n$) belongs to $g \circ M$ iff $[cg,n+\tilde l(g)] \in M \cap \langle x_g \rangle$. We note that $[cg, n + \tilde g] \prec [g, \tilde l(g)]$ in any case (since $S^n g \subset S^{n+\tilde l(g)}$). Hence for $y\in R_f$, $y \in g \circ M$ iff $\RR_{g,\tilde l(g)} (y) \in M$.
This entails that $[1,0] \not\in g\circ M$, so that the latter is proper if $M $ is.

It is trivial that $g\circ M$ is a subgroup of $R_f$, and if $A = \R$, it is a vector subspace thereof.

\noindent{\it $g \circ M$ is convex.} Suppose that $[y_1, n(1)] \leq [z,n] \leq [y_2,n(2)]$ with $[y_i,n(i)] \in M$ and $[z,n] \in R_f$. Since $\RR_{g,\tilde l(g)}$ is order preserving, it follows that $[y_1g, n(1)+ \tilde l(g)] \leq [zg, n + \tilde l(g)] \leq [y_2g , n(2)+ \tilde l(g)]$. Since both ends belong to $M$, the middle term is thus in $M$, and so $[zg_2, n(2)+ \tilde l(g)] \in M$. Hence $[z,n] \in g \circ M$.

\noindent {\it $g \circ M$ is generated as an abelian group by its positive elements---that is, it is directed.} If $y = [c,n] \in g \circ M$, then $z = [cg,n+ \tilde l(g)] \in M \cap \langle x_g\rangle$. As the latter is an order ideal of $R_f$ (any intersection of finitely many order ideals in a dimension group is itself an order ideal), we can write $z = [p_1, k(1)] - [p_2, k(2)]$ where $[p_i , k(i)] \in M^+ \cap\langle x_g \rangle^+ $. Since $[p_i,k(i)] \prec x_g$, there exist $m_0$ \st for all $m \geq m_0$ both $f^{m+ \tilde l(g)}p_i $ have only nonnegative coefficients and $f^{m+ \tilde l(g)}p_i \prec f^{m+k(i)}g$. By increasing $m$ if necessary, we may assume that $p_i$ each have no negative coefficients.

Thus $f^{m+\tilde l(g)}p_ig^{-1} \prec f^{m+k(i)}$. Therefore $k(i) \geq \tilde l(g) $ and $\tilde l(hg) \leq k(i) - \tilde l(g)$ for every $h \in \supp p_i$. Thus $y_i:= [p_i g^{-1}, k(i)- l(i)] \in R_f$, and since the $y_i$ have no negative coefficients, it follows that $y_i \in R_f^+$.

Moreover, $\RR_{g, \tilde l(g)}(y_i) = [p_i, k(i)] \in M \cap \langle x_g \rangle$, and thus $y_i \in g\circ M$. Finally, $y = \RR_{g,\tilde l(g)}^{-1} = y_1 - y_2$, a difference of positive elements of $g \circ M$.

Thus $g \circ M$ is a proper order ideal.

\noindent (a) Now suppose that $M + \langle x_g \rangle = R_f$. Then the second isomorphism theorem applies (for order ideals in dimension groups), $R_f/M \iso \langle x_g\rangle /(M \cap \langle x_g\rangle)$ as pointed partially ordered abelian groups (and vector spaces, if they are real vector spaces); under this map, $[1,0] + M \mapsto x_g + M \cap \langle x_g\rangle$.

We have the map $R_f \to \langle x_g\rangle$ given by $\RR_{g,\tilde l(g)}$; it is a consequence of the definitions that the kernel of the induced map $R_f \to \langle x_g\rangle/(M \cap \langle x_g \rangle)$ is exactly $g \circ M$, and this yields an order isomorphism $R_f/g \circ M \to \langle x_g\rangle/(M \cap \langle x_g \rangle) $.

\noindent (b) This is tautological.
\qed

In particular, if $M$ is a maximal order ideal and $x_g \notin M$, then $g \circ M $ is also a maximal order ideal. Sometimes this action can be inverted, but not in an obvious way (since it can easily happen that $x_{g^{-1}} \in M$). This gives an action of a part of $G$ on $\Cal M = \brcs{\text{maximal order ideals of $R_f$}}$. If $x_g \not\in M$ (there are plenty of such elements if $G$ is infinite), then we form $g \circ M$; then we can find $g_2 $ \st $x_{g_2} \not\in g\circ M$, and form $g_2 \circ (g\circ M)$. It is routine to verify that $\tilde l_S(g_2 g) = \tilde l_S(g_2) + \tilde l_S(g)$, so that $x_{g_2 g} = x_{g_2} x_g$ (the latter has the obvious interpretation), and then, as we will see below, $g_2\circ (g \circ M) =( g_2 g) \circ M$. We can obviously continue this process indefinitely. However,  it can happen that $g\circ M = M$ for all eligible $g$, or at the other extreme, there does not exist $g $ with $\tilde l_S(g) > 0$ \st $g\circ M = M$, or even that for no  $N \in \Cal M$, is $M = g\circ N$ for  any non-identity $g$.

In any event, we will show that the quotients $R_f/M$ are stationary systems in a fairly strong way, and in particular, there is only one trace that kills $M$ and $M$ is of cofinite rank in $R_f$.

Assume $G$ is infinite. Let $M$ be a maximal order ideal of $R_f$. Define for $k > 0$, $\Gamma'_{k,M} = \Set{g \in \Gamma'_k}{x_g \not\in M}$. If $g \in \Gamma'_k$, then by definition, $x_g = [g,k]$, that is, $\tilde l(g) = k$. If necessary, we formally define $\Gamma'_{0,M} = \brcs{1}$. There are a few obvious remarks, in analogy with what happens for the larger sets, $\Gamma'_k$.

A simple non-membership criterion for elements of the form $x_g$ is available. 

\Lem Lemma \twothr. Let $M$ be a maximal order ideal of $R_f $, and $b$ an element of $G$ \st $x_b \notin M$. Then for $a \in G$, we have $x_a \notin b \circ M$ if and only if $x_{a} x_{b} \notin M$. When the latter occurs, $\tilde l_S(ab) = \tilde l_S(a) + \tilde l_S(b)$ and $x_{ab} = x_a x_b$. 

\Pf We have that  $x_a \notin b \circ M$ iff  $z:= [{}{ab}, \tilde l(a) + \tilde l(b)] \notin M \cap \langle x_b\rangle$. We claim $z \in \langle x_b\rangle$: there exists $m$ \st $S^m a \subset S^{m+ \tilde l(a)}$ (where $S = \supp f$, as usual). Then $S^m ab \subset S^{m+\tilde l(a)}b$, and thus $f^m ab \prec f^{m+ \tilde l(a)}b$, so that there exists a positive integer $N$ \st for all $l$, $[{}{ab}, l+ \tilde l(a)] \leq N [{}b, l ]$; set $l = \tilde l(b)$. 

Hence $x_a \notin b \circ M$ iff $z \notin M$; of course, $z = x_a x_b$. Since every maximal order ideal contains $\SS R_f$, we have that if $x_a x_b \notin M$, then  $\tilde l(ab) \geq \tilde l(a) + \tilde l(b)$, yielding equality, and also $x_{ab} = x_a x_b$. \qed

\Lem Lemma \twofou. Suppose that $h \in S^k$, $j \in S^l$ and $g=hj \in \Gamma'_{k+l}$.
Then $h \in \Gamma'_k$ and $j \in \Gamma'_l$.

\Pf If either $\tilde l (h) < k$ or $\tilde l(j) < l$, then $\tilde l(g) < k+l$, contradicting $\tilde l(g) = k+l$.\qed

\Lem Lemma \twofiv. (a) $x_g \not \in \SS R_f$ for all $g \in G$.
\item{(b)} If $g \in \Gamma'_k$, then there exists $h \in \Gamma'_1$ \st $hg \in \Gamma'_{k+1}$.
\item{(c)} If $g \in \Gamma'_k$ and $k > 0$, then there exist $h \in \Gamma'_1$ and $j \in \Gamma'_{k-1}$ \st $g = hj$.

\Pf (a) If $x_g:= [g, \tilde l(g)] = \SS [h,k] =[h,k+1]$ for some $h \in AG$ with $h \prec f^k$, then there exists $m$ \st $f^{m+k+1} g = f^{m+\tilde l(g)}h \prec f^{m+\tilde l(g) + k}$, which would yield $\tilde l(g) + 1 \leq \tilde l(g)$, a contradiction.

\noindent (b) We have
$$\eqalign{
[g, \tilde l(g)] & = [fg, \tilde l(g) + 1] \cr
&= \sum_{h\in S} (f,g) [{}{hg}, \tilde l(g)+1].\cr
}$$
Since all the coefficients, $(f,g)$, are positive, there exists $h \in S$ \st $[{}{hg},\tilde l(g) + 1] \not \in \SS R_f$. This forces $\tilde l(hg) = \tilde l(g) +1$, so $h g \in \Gamma'_{k+1}$. By Lemma \twofou, $h \in \Gamma'_1$.

\noindent (c) Since $g \in S^k$, there exist $h \in S$ and $j \in S^{k-1}$ \st $g = hj$. The lemma above now applies.
\qed

\Lem Lemma \twosev. (a) Every maximal order ideal contains $\SS R_f = \langle [1,1]\rangle$.
\item{(b)} For all $k$, the set $\Gamma'_{k,M}$ is not empty.
\item{(c)} If $g \in S^k$ and $x_g \not\in M$, then $\tilde l_S(g) = k$ and $g \in \Gamma'_{k,M}$.
\item{(d)} If $g = hj$ with $h \in S^a$ and $j \in S^b$, then $g \in \Gamma'_{a+b,M}$ entails $j \in \Gamma_{b,M}'$ and $h \in \Gamma_{a,j\circ M}'$, and in particular, $\tilde l_S(hj) = \tilde l_S(h) + \tilde l_S(j)$.
\item{(e)} For all $k > 0$ and all $g \in \Gamma'_{k,M}$, there exist $j \in \Gamma'_{k-1,M}$ and $h \in \Gamma_{1,j\circ M}'$ and \st $g = hj$.

\Pf (a) A consequence of Lemma \MItwo.

\noindent (b--e) are proved similarly to the methods used for the sets $\Gamma'_k$.
\qed

\Lem Corollary \twoeig. If $g \in \Gamma'_{k,M}$, then for all positive integers $\lambda$, $\Gamma_{\lambda,g\circ M}\cdot g \subset \Gamma'_{\lambda + k,M}$. Conversely, if $h \in S^{\lambda}$ and $hg \in \Gamma'_{k,M}$, then $h \in \Gamma_{\lambda,g \circ M}$.

\Lem Lemma \twontn. Let $M$ be a maximal order ideal of $R_f$. Suppose that $a,b \in G$, and both $b \circ M$  and $a\circ (b \circ M)$ are  defined. Then $(ab)\circ M$ is defined, and $a\circ (b \circ M) = (ab)\circ M$; moreover, $\tilde l(ab) = \tilde l(a) + \tilde l(b)$.  

\Rmk The hypotheses are that $x_b \notin M$, and that $x_a \notin b \circ M$; the conclusion includes $x_{ab} \notin M$. 

\Pf Set $z:= [{}{ab}, \tilde l(a) + \tilde l(b)]$.  By Lemma \twothr, $z \notin M$  and $z = x_{ab} = x_a x_b$; moreover, $\tilde l(ab) = \tilde l(a) + \tilde l(b)$. We  can form $(ab)\circ M$. For $g \in G$, we have that $x_g \in (ab)\circ M$  entails $[{}{gab}, \tilde l(g) + \tilde l(ab)] \in M$, and  this rewrites as 
$[{}{gab}, \tilde l(g) + \tilde l(a) + \tilde l(b)] \in M$. Thus $\RR_{b,\tilde l(b)}\([{}{ga}\tilde l(g) + \tilde l(a)]\) \in M$, so that $[{}{ga},\tilde l(g) + \tilde l(a)] \in b \circ M$, and again applying the definitions, we deduce $x_g = [g, \tilde l(g)] \in M$. Since any order ideal is the integer (real) span of elements of the form $x_g$, it follows that $(ab) \circ M \subseteq M$; since both are maximal order ideals, equality holds. \qed 

\Lem Corollary \twonin. Suppose $M$ is a maximal order ideal. Define $\Cal G_M: = \Set{g \in G}{g \circ M = M}$. Then $\Cal G_M$ is a subsemigroup of $G$, and the restriction $\tilde l_S|\Cal G_M$ is additive. Moreover if both $g$ and $g^{-1}$ belong to $\Cal G_M$, then $\tilde l_S(g) = \tilde l_S (g^{-1}) = 0$.

\Pf The identity belongs to $\Cal G_M$, so the latter is nonempty. By the preceding, $g \circ (h \circ M)$ is defined for all $g,h$ in $\Cal G_M$, and equals $(gh)\circ M$, verifying that $gh \in \Cal G_M$. Additivity of the restriction of $\tilde l$ is an immediate consequence of the preceding, as is the final statement.\qed

As $\tilde l_S^{-1}(0)$ is   a torsion group (Lemma \firstone),  if $G$ is torsion-free, no nonidentity element of $\Cal G_M$ is invertible therein. If $(G,S)$ satisfies WC, then $\tilde l_S^{-1}(0)$ is a finite group, so only finitely many elements of $\Cal G_M$ are invertible therein. Unfortunately, these subsemigroups are typically quite small.

\Lem Lemma \twoten. Suppose that $M$ is a maximal order ideal of $R_f$, and $g$ is an element of $G$ \st $\tilde l_S (g) = k$ and $x_g \notin M$. Then $\Gamma'_{k+l,M} \subseteq \Gamma'_{l,g\circ M}\cdot g$ for all nonnegative integers $l$. 

\Pf Pick $\gamma \in \Gamma'_{k+l,M}$; by Lemma \twofiv, there exists $h \in \Gamma'_{l}$ \st $\gamma = hg$. 
We claim that that $h$  belongs to $\Gamma'_{l,g\circ M}$. First, $\tilde l (h) = l$, since $ h \in \Gamma'_l$; then if $[{}h, l] \in g\circ M$, we would have $[{}{hg}, \tilde l(g) + \tilde l (h)] \in M$. But $\gamma = hg \in \Gamma'_{k+l, M}$ so that both $\tilde l(gh) = \tilde l (\gamma) = k+ l = \tilde l (g) + \tilde l (h)$ and $[{}{gh},k+l] \not\in M$. \qed

In particular, for all sufficiently large $l$, 
$$
|\Gamma'_{k+l,M}| \leq |\Gamma'_{l,g\circ M}|.
$$
This can be arranged for every $g' \in \Gamma'_{k,M}$ simultaneously at the cost of increasing the minimal $l$ for which it holds (since $\Gamma'_k$ is finite). We then see that each of the $\Gamma'_{l,g_i\circ M}g_i$ are equal (as $g_i$ varies over $\Gamma'_{k,M}$).

Hence we have the following preliminary result.

\Lem Lemma \twoele. Let $M$ be a maximal order ideal of $R_f$. Either of the following conditions is sufficient for $\sup_k |\Gamma'_{k,M}| < \infty$.
\item{(i)}There exists $g \in \cup_k \Gamma'_{k,M}$ \st $g\circ M = M$, or
\item{(ii)} The number of maximal order ideals of $R_f$ is finite.

\Pf (i) There exists $k$ \st $g \in \Gamma'_{k,M}$ and by the preceding, for all $l \geq l'$, we have $|\Gamma'_{k+l,M}| \leq |\Gamma'_{l,g\cdot M}| = |\Gamma'_{l,M}| $. Then $\max\brcs{|\Gamma'_{j,M}|}_{j \leq l'+k-1}$ is an upper bound for $|\Gamma'_{j,M}|$.

\noindent (ii) Given a maximal order ideal, $M$, there always exists $g \in \Gamma'_{1,M}$. Hence we can construct an infinite sequence of maximal order ideals, $(M, g_1\circ M, g_2 \circ (g_1\circ M), \dots$.   By Lemma \twontn, when it is defined, $g\circ (g' \circ M) =$ equals $(gg')\circ M$, so we can write the elements of the sequence $h_n\circ M$, where $h_n \in S^n$; it easily follows from $[{}{h_n} , n] \not\in M$, that $\hat l(h_n) = n$. Since there are only finitely many maximal order ideals, there exist $h_n \in \Gamma'_{n,M}$ and $h_m \in \Gamma'_{m,M}$ \st $h_n \circ M = h_m \circ M$. Writing $h_m = jh_n $, we see that $j\circ (h_n\circ M) = h_n \circ M$. By (i), we have $\sup_k |\Gamma'_{k,h_n \circ M}| < \infty$.

We also have, for all sufficiently large $l$, that $|\Gamma'_{n+l,M}| \leq |\Gamma'_{n+l,h_n\circ M}| $. Thus $\limsup_k |\Gamma'_{k,M}| \leq \limsup |\Gamma'_{k,h_n\circ M}| < \infty $.
\qed

\Lem Lemma \twotwe. If $\brcs{|\Gamma'_{k,M}|}_k$ is bounded, then $R_f/M$ has unique trace.

\Pf The maps $ \Arrow\psi_k; A\Gamma'_k . A\Gamma'_{k+1}$ are restrictions of multiplication by the same $f$; hence the nonzero coefficients of the corresponding matrices are bounded below and above. Since the limit is simple, and there is a bound on the width of the matrices, it easily follows that there is a unique trace (there is a uniform bound on how many terms are necessary to obtain a strictly positive matrix, and then Birkhoff's criterion can be used).
\qed

  The next is easy, but the reverse inclusion with bounded index of primitivity is trickier.

\Lem Lemma \twothi. For all $k, \lambda \in \N$ and $g \in G$ with $\tilde l_S(g) = k$ and $g \in S^k$ with $x_g \not\in M$,
$$
\Gamma'_{\lambda, g\circ M}\cdot g \subseteq \Gamma'_{k+\lambda,M}.
$$

\Pf If $h $ belongs to the left side, then $hg^{-1} \in \Gamma'_{\lambda,g \circ M}$, and thus $\tilde l(hg^{-1}) = \lambda$ and $[{}h, \tilde l (g) + \lambda] \not\in M$. We obtain $\tilde l (h) = \lambda +k$, so $h \in \Gamma'_{\lambda + k, M}$. \qed

We can write $R_f/M$ as the direct limit, $\Arrow \phi_k; A\Gamma'_{k,M} . A\Gamma'_{k+1,M}$ where $\phi_n$ is the compression and restriction of $f \times$ (this is a general property of order ideals occurring in dimension groups). The previous results say that the only $j$ in $\supp f $ that contribute to $\phi_k$ and yield a particular element $g \in \Gamma'_{k+1,M}$ in the support of the their image are those $j \in \Gamma'_{1,g\circ M}$.
In particular, we can write $\phi_k$ in terms of natural bases (up to the choice of ordering on the bases) as a matrix. Obviously, we can do the same with $1$ replaced by any positive integer $\lambda$, by replacing $f$ by $f^{\lambda}$, compressed and restricted to a positive map $A\Gamma'_{k,M} \to A\Gamma'_{k+\lambda,M}$. The matrices representing these are much easier to describe when they are strictly positive (that is, every entry is positive).

The following will be improved to $\lambda_0$ {\it independent of $k$.}

\Lem Lemma \twoftn. (Here $M$ is a maximal order ideal.) For $g \in \Gamma'_{k,M}$, there exists $\lambda_0 \equiv \lambda_0 (k)$ \st for all $\lambda \geq \lambda_0$,
$$
\Gamma'_{\lambda,g\circ M}\cdot g = \Gamma'_{k+\lambda,M}.
$$

The proof requires yet another version of a predecessor/successor result.

\Lem Lemma \twoffn. Suppose $\rho \in S^k$, $\sigma \in S^{\lambda}$, and $\sigma\rho \in \Gamma'_{k+ \lambda,M}$. Then $\rho \in \Gamma'_{k,M}$ and $\sigma \in \Gamma'_{\lambda, \rho \circ M}$.

\Pf We have $[{}\sigma \rho, k+ \lambda] \prec [{}{\sigma,k}]$; if the latter were in $M$, then so would be in the former, contradicting $\sigma\rho \in \Gamma'_{k + \lambda, M}$. Thus $[{}{\rho, k}] \not\in M$. This entails $\tilde l(\rho) = k$, so $\rho \in \Gamma'_{k,M}$. If $[{}\sigma,\lambda] \in \rho \circ M$, then $[{}{\rho \sigma, \lambda + k}] \in M$, again a contradiction. As $\tilde l(\sigma) \leq \lambda$, it follows that $\tilde l(\sigma) = \lambda$. Thus $\sigma \in \Gamma'_{\lambda, \rho\circ M}$.
\qed

\Pf (of Lemma \twoftn) Since $R_f/M \iso \lim A\Gamma'_{k,M} \to A\Gamma'_{k+1,M}$, and the limit dimension group is simple, provided there are no zero rows or columns in the matrices, given $k$, there exists $\lambda_0$ (depending on $k$) \st the matrix representing $A\Gamma'_{k,M} \to A\Gamma'_{k+\lambda_0, M}$ is strictly positive, and this holds for all larger $\lambda$.

In particular, given $\gamma \in \Gamma'_{\lambda+k,M}$, there exists $h \in S^{\lambda}$ \st $\gamma = hg$. By the preceding lemmas, $h \in \Gamma'_{\lambda,g\circ M}$, so that $\gamma \in \Gamma'_{\lambda,g\circ M}\cdot g$. Hence $\Gamma'_{k+ \lambda,M} \subseteq \Gamma'_{\lambda,g \circ M}\cdot g$.
\qed

Now we can describe (at least for $\lambda \geq \lambda_0(k)$), the matrices given by the map $A\Gamma_{k,M}' \to A\Gamma_{k+ \lambda,M}'$.

The columns correspond to the elements, $\gamma \in \Gamma'_{k+\lambda,M}$. From $\Gamma'_{\lambda,g\circ M}\cdot g = \Gamma_{k+\lambda,M}'$, we see that if $g,g' \in \Gamma'_{\lambda,M}$, then $\Gamma'_{\lambda,g\circ M} = \Gamma'_{\lambda,g'\circ M}g'g^{-1}$; in particular, $|\Gamma'_{\lambda,g\circ M} | = |\Gamma'_{\lambda,g'\circ M}g'g^{-1}| = |\Gamma'_{\lambda+k,M}|$.

The $g$-indexed column has as its $j$th entry, $(f^{\lambda},{}j)$ where $j$ runs over $\Gamma'_{\lambda,g\circ M}$. The cardinalities match, and the resulting matrix is square. We still have to pick orderings on the elements of
the pairs of sets $\Gamma'_{k+\lambda,M}$ and $\Gamma'_{k,g\circ M}$ for each $g$, in order to obtain actual (ordered) bases.

Next, we show that $\lambda_0$ does not depend on $k$.

\Lem Lemma \twosxn. We can choose $\lambda_0(k)$ to be independent of $k$.

\Pf Fix $k>1$. Pick $\gamma \in \Gamma'_{k+1,M}$; by Lemma \twoffn, we can factor $\gamma = h_1h$ where $h_1 \in \Gamma'_{1,h\circ M}$ and $h\in \Gamma'_{k,M}$. Assuming $\lambda_0(k)$ is monotone increasing in $k$, pick $\lambda \geq (\lambda_0(k)-1 )\vee \lambda_0(1)$; we then have
$$\eqalign{
\Gamma'_{\lambda,\gamma \circ M}\cdot \gamma
&= \(\Gamma'_{\lambda, h_1 \circ (h\circ M)}\cdot h_1\) \cdot h\cr
&= \Gamma'_{\lambda+1,h\circ M}\cdot h\cr
& = \Gamma'_{\lambda + 1 + k,M}.\cr
}$$
Hence $\lambda_0(k+1)$ can be chosen to be at most $\lambda_0(1)$.
\qed

So we can simply write $\lambda_0$ in place of $\lambda_0(k)$.

\Lem Corollary \twosvn. Assume that $\lambda \geq \lambda_0$.
\item{(i)} For all $h \in \cup_{j > 0}\Gamma''_{j,M}$ ($\supset \Gamma'_{j,M}$), we have $\Gamma'_{\lambda,h\circ M}\cdot h = \Gamma'_{\lambda + \tilde l (h),M}$;
\item{(ii)} $|\Gamma'_{\lambda,M}| \leq \max_{N \in \Cal M} |\Gamma'_{\lambda_0,M}| < |\Gamma'_{\lambda_0}|$.

\Pf (i) is restatement of the previous result. (ii) Set $k = \lambda - \lambda_0$ and pick $h \in \Gamma'_{k,M}$. Then $\Gamma'_{k+\lambda,M} = \Gamma'_{\lambda,h \circ M}\cdot h$, so $|\Gamma'_{\lambda,M} | = |\Gamma'_{\lambda,h \circ M}| \leq \max_{N \in \Cal M} |\Gamma'_{\lambda_0,M}| < |\Gamma'_{\lambda_0}|$. \qed

\noindent {\it Conclusion of proof of Theorem \twoone.} At this point, we can conclude that $R_f/M$ has a number of properties. First, it has a unique trace, since it is given as the direct limit of square matrices with strictly positive and uniformly bounded entries (bounded by $\max_{g \in S^{\lambda_0}}(f^{\lambda_0},g)$)---Birkhoff's criterion applies to yield uniqueness of the trace.
Next, $R_f/M$ is of finite rank (explicitly, at most $\liminf_{k\to \infty} |\Gamma'_{k,M}|$; the sizes of the $\Gamma'_{k,M}$ can cycle around).

Up to this stage, we did not care about the ordering (as in {\it ordered basis\/}) on the rows and columns. To obtain stationarity, we now have to deal with it.

A {\it multiset\/} is an unordered collection of objects (in our case, sets) permitting multiplicities; we use the notation $\mset{\cdot}$. Of interest are the multisets
$$
\Cal T_{k,i} = \Cal T_{i,k} (M):= \mset{
\Gamma'_{i,h\circ M}}_{h \in \Gamma'_{k,m}
},
$$
where $M$ is a maximal order ideal. These represent the columns of the transition matrix $\Arrow \phi_{k,i}; A \Gamma'_{k,M}. A \Gamma'_{k+i,M}$, at least if $i \geq \lambda_0$. Fix $i \geq \lambda_0$. Since the set of entries of the matrices form a subset of $\brcs{(f^i,{}j)}_{j \in S^i}$, and this is finite, for each $k$, there exists $l\equiv l(i)$ \st $\Cal T_{k,i} = \Cal T_{k+l,i}$ (equality as multisets).

We will show that $\Cal T_{k,i} = \Cal T_{k+l,i}$ implies $\Cal T_{k+1,i} = \Cal T_{k+l+1,i}$. This will be enough to prove that after a uniform telescoping (replacing $f$ by $f^m$ for some $m$ that is divisible by $\lambda_0 \cdot l$), the system representing $R_f/M$, $\lim A\Gamma'_{tm,M}\to A\Gamma'_{(t+1)m,M}$ is stationary.

We first observe that $\Cal T_{k,i} = \Cal T_{k+l,i}$ iff there exists a bijection $\Arrow \alpha; \Gamma'_{k,M} . \Gamma'_{k+l,M}$ \st $\Gamma'_{i,h\circ M} = \Gamma'_{i,\alpha(h) \circ M}$ for all $h\in \Gamma'_{k,M}$. We want to construct a bijection $\Arrow \beta; \Gamma'_{k+1,M} . \Gamma'_{k+l+1,M}$ \st $\Gamma'_{i,h\circ M} = \Gamma'_{i,\beta(h) \circ M}$.

For each $j \in \Gamma'_{k+1,M}$, we can pick $p(j) \in \Gamma'_{k,M}$ \st there exists $\gamma_j \in \Gamma'_{i,p(j)\circ M}$ with $j = \gamma_j \cdot p(j)$; we can even do this so that $p$ is one to one (if $i$ is divisible by $\lambda_0$, which we have assumed). Now set $\beta(j) = \gamma \cdot \alpha(p(j))$.

We verify easily that $\beta$ is well-defined, that is, $\gamma_j \in \Gamma_{i,\alpha(p(j))\circ M}$;  then $\beta(j) \in \Gamma_{k+l+1}$, and that $\beta$ is a bijection, and $\Gamma'_{i,j\circ M} = \Gamma_{i,\beta(j)\circ M}'$.

Thus we have a sequence of equalities $\Cal T_{k,i} = \Cal T_{k+l,i}$, $\Cal T_{k+1,i} = \Cal T_{k+l+1,i}$, $\dots $ yields periodicity, that is $\Cal T_{k+s,i} = \Cal T_{k+s + ml,i} $. The composition of $l$ consecutive maps (beginning with $A \Gamma_{k,M} \to A\Gamma_{k + li,M}$) are thus identical, and so (after a telescoping) the system is stationary.
\qed

\SecT \IMI\  Invariant maximal order ideals

If $M$ is a maximal order ideal of $R_f$, we have seen that $R_f/M$ has unique trace; this yields a corresponding trace, denoted $\tau_M$,  on $R_f$, \st $\tau_M(M) = (0)$. 

\Lem Lemma \IMIone. Suppose that $M$ is a maximal order ideal, and $x_g := [g,\tilde l_S(g)]$. Then
$$
\tau_M \circ \RR_{g,\tilde l_S(g)} = \cases \tau_M(x_g)\cdot \tau_{g\circ M} & \text{if $x_g \not\in M$}\\
0 & \text{if $x_g \in M$}.\\
\endcases
$$

\Pf Suppose that $x_g \not\in M$, so that $g\circ M$ is defined, and is a maximal order ideal. Evaluating at $[1,0]$, we have that $\tau\circ \RR_{g,\tilde l(g)}$ is not zero. If $a = [c,k] \in g\circ M$, then we can assume that $\supp c \subset S^k$ (by replacing $[c,k]$ by $[f^l c, k+l]$ and $k$ by $k+l$ for sufficiently large $l$), and $[cg, k+ \tilde l(g)] \in M$; but $[cg, k+ \tilde l(g)] = \RR_{g,\tilde l(g)}(a)$. Hence $g \circ M \subset \ker \(\tau_M \circ \RR_{g,\tilde l(g)}\)$; since $g\circ M$ is an order ideal, we have $g \circ M \subset \ker^+ (\tau_M \circ \RR_{g,\tilde l(g)})$. The latter is (by definition) an order ideal. Since $g\circ M$ is a maximal order ideal, we have $g \circ M = \ker^+ (\tau_M \circ \RR_{g,\tilde l(g)})$, and thus $\tau_M \circ \RR_{g,\tilde l(g)}$ induces the unique trace on $R_f/g\circ M$. Hence $\tau_M \circ \RR_{g,\tilde l(g)}$ is a scalar multiple of $\tau_{g\circ M}$, and evaluating at $[1,0]$ yields the scalar.

If $x_g \in M$, then for all $a = [c,k] \in R_f$, we have $\RR_{g,\tilde l(g)}(a) = [cg,k+\tilde l(g)] \prec x_g \in M$. Hence $\RR_{g,\tilde l(g)}(R_f) \subset M$, and thus $\tau_M \circ \RR_{g,\tilde l(g)}= 0$.
\qed

\Lem Lemma \IMItwo. Suppose that $M$ is a maximal order ideal in $R_f$, $g,g' \in \Gamma'_{k,M}$, $h \in \Gamma'_{k',M}$, and $h \circ M = M$.
\item{(a)} For all sufficiently large $k$ (depending only on $S$), $a := g'hg^{-1} \in \Gamma'_{k',M}$ and $a\circ (g\circ M) = g'\circ M$; in addition, $\tilde l_S (g'hg^{-1}) = \tilde l_S(h)$. In particular, if $g' =g$, then $a = ghg^{-1}\in \Gamma'_{k',M}$ and $a \circ (g\circ M) = g\circ M$.
\item{(b)} For all sufficiently large $k$, $gh \in \Gamma'_{k+k',M}$ and thus $\tilde l_S(gh) = \tilde l_S(g) + \tilde l_S (h)$. In addition, $(gh)\circ M$ is defined and equals $g\circ M$.
\item{(c)} If $k$ is sufficiently large, then $a = g h^{-1} \in \Gamma'_{k-k', M}$ and $a \circ M = g\circ M$; moreover, $\tilde l_S(gh^{-1}) = k-k'$.
\item{(d)} If $n\tilde l_S(h) -\tilde l_S(g)$ is sufficiently large, then $b = h^n g^{-1} \in \Gamma'_{nk' - k, g\circ M}$ and $b \circ (g \circ M) = M$.

\Pf For all sufficiently large $l,l'$ with $l+k = l'+k'$, we have
$$
\Gamma'_{l,g\circ M} g = \Gamma_{l', M} h = \Gamma_{l+k,M},
$$
using $h\circ M = M$.

\noindent (a) If $l' = k$, then $g'h$ belongs to the middle term, so there exists $a \in \Gamma'_{l,g\circ M}$ \st $ag = g'h$, and so $a = g'hg^{-1}$. Since $l' = k$, it follows that $l = k'$. Since $a \in \Gamma'_{k',g\circ M}$, $a \circ (g\circ M)$ is defined, and thus is equal to $(ag)\circ M = (g'h)\circ M = g'\circ (h \circ M) = g'\circ M$.

\noindent (b) Using $l'+k' = l+k$ and setting $l' = k$, we see that $gh$ belongs to the middle term, so belongs to the right term, which is $\Gamma_{k'+k,M}$. Hence $\tilde l(gh) = k+k' = \tilde l(g) + \tilde l(h)$. The last statement is straightforward.

\noindent (c) Set $l' = \tilde l(g) - \tilde l(h) = k-k'$, so that $l +k = l'+k' = k$. Then $h$ belongs to the right side, so $gh^{-1} \in \Gamma'_{k-k',M}$. The rest is routine.

\noindent (d) Set $l = nk' - k$. By (b) applied inductively to $g = h$, $g = h^2$, etc, we obtain $h^n \in \Gamma_{nk',M}$. Hence $h^n $ belongs to the right side, and thus $h^n g^{-1} \in \Gamma'_{nk'-k, g \circ M}$, and the rest is straightforward.
\qed

 A normalized trace $\tau$ on $R_f$ is {\it  multiplicative\/} if $\tau ([aa',k+k'] = \tau([a,k])\cdot \tau([a',k'])$ whenever $x = [a,k]$ and $y= [a',k']$ belong to  $R_f$ (the product $aa'$ is that of the group ring). When the group is abelian, the multiplicative traces are precisely the pure ones, but for nonabelian groups, multiplicativity is a relatively rare occurrence. Examples are those obtained from characters of the group, but there are others.

\Lem Lemma \IMIthr. For a maximal order ideal $M$, of $R_f$, the following are equivalent.
\item{(a)} For all $g \in G$ \st $x_g \not\in M$, we have $g\circ M = M$;
\item{(b)} the trace $\tau_M$ is multiplicative.

\Pf (a) implies (b) 
Suppose that for all $g$ \st $x_g \notin M$, we have $g \circ M = M$. We verify multiplicativity. It suffices to prove it in the case that $x,y \in A_f^+$, and this reduces to the case that $x = [{}{j}, k]$ and $y = [{}{j'}, k']$ with $k \geq \tilde l(j)$ and $k' \geq \tilde l(j')$. If either $k > \tilde l (j)$ or $k' > \tilde l(j')$, then $\tau_M (x)\cdot \tau_M(y) =0$, and since $k + k'> \tilde l(j) + \tilde l(j') \geq \tilde l(jj')$, we also have that $\tau_M([{}{jj'}, k+k']) = 0$.

Hence we reduce to the case that $x = x_j$ and $y = x_{j'}$. If both do not belong to $M$, then we have already seen that multiplicativity. If $x_{j'} \in M$, then $[{}{jj'},k+k'] \prec [{}{j'},k'] $, so $\tau_M ([{}{jj'},k+k']) = 0$, so multiplicativity holds in this case as well.

Finally, suppose that $[{}j,k] \in M$ and $[{}{j'},k'] \not\in M$. If $\tilde l(jj') > \tilde l(j) + \tilde l(j')$, then $[{}{jj'}, k+k'] \in \SS R_f \subset M$, so that $\tau_M([{}{jj'}, k+k'])$. This leaves the case that $\tilde l(jj') =\tilde l(j) + \tilde l(j')$, so that $[{}{jj'}, k+k'] = x_{jj'}$. If $x_{jj'} \not\in M$, then $(jj')\circ M$ is defined, as is $j'\circ M$. If $x_{j} \in j'\circ M$, then $[{}{jj'},k+k'] \in M$ (from the definition of $j'\circ M$). Hence $x_{jj'} \in M$, a contradiction.

\noindent (b) implies (a). Pick $g \in G$ \st $x_g \not\in M$. As $\tau_M $ is multiplicative, for any $h \in G$, $\tau_M ([{}{hg}, \tilde l(g) + \tilde l(h)]) = \tau(x_g)\cdot \tau_M(x_h)$. If $x_h \not\in g\circ M$, then $\tau_M ([{}{hg}, \tilde l(g) + \tilde l(h)]) \neq 0$, so that both factors are nonzero, and $\tilde l(g) + \tilde l(h) = \tilde l(hg)$, and we deduce $\tau_M (x_h) \neq 0$; since $x_h$ is in $R_f^+$, this implies $x_h \not\in M$.

It is easy to verify that any maximal order ideal $J$ is spanned (additively) by $\SS R_f$ together with $\Set{rx_j}{r \in A; \text{$j \in G$ \st $x_j \not\in J$}}$. Thus $M \subset g\circ M$; since both are maximal order ideals, they must be equal.
\qed

\Lem Lemma \IMIfou. Suppose $g,g' \in \Gamma_{k,M}$ and $g \circ M = g'\circ M$. Then $g'g^{-1}$ is torsion. If additionally, $G$ is torsion-free, then $g = g'$.

\Pf Select $\lambda \geq \lambda_0$; then $\Gamma_{\lambda, g\circ M} \cdot g = \Gamma_{\lambda +k,M} = \Gamma_{\lambda, g'\circ M} \cdot g'$. Since $g \circ M = g'\circ M$, we have $\Gamma_{\lambda, g\circ M} \cdot g = \Gamma_{\lambda, g\circ M} \cdot g'$, and therefore, $\Gamma_{\lambda, g\circ M} = \Gamma_{\lambda, g\circ M} \cdot (g'g^{-1})$. Since $\Gamma_{\lambda, g\circ M} $ is finite, some power of $g'g^{-1}$ must be the identity. \qed

\Lem Corollary \IMIfiv. Suppose that $G$ is torsion-free and $M$ is a maximal order ideal \st for all $g $ with $x_g \not\in M$, we have $g\circ M = M$. Then there exists $h \in G$ \st for all  $k$, $\Gamma_{k,M}' = \brcs{h^k}$. 

\Rmk This is a type of unique factorization property for powers of $h$, that is, if $h^k$ is a product of $k$ elements of $S$, then each of the factors must be $h$ itself.

\Pf By Lemma \IMIfou, $|\Gamma_{k,M}'| = 1$ for all. Then $h$ is defined by $\Gamma_{1,M}' = \brcs{h}$, and define $h_i$ via $\Gamma_{i,M} = \brcs{h_i}'$. By Lemma \IMIthr, $\tau_M ((x_h)^i) \neq 0$, and thus $[{}{h^i},i] \not\in M$, forcing $\tilde l(h^i) = i$ and $x_{h^i} \not\in M$; hence $h^i \in \Gamma_{i,M}'$, and uniqueness forces $h_i = h^i$. \qed 

One potential source of multiplicative pure traces, with nonzero positive kernel (not necessarily of the form $\tau_M$) is in $\partial_e F_0 \cap\(\overline{\cup_{\lambda > 0} \partial_e F_{\lambda}}\)$, if for all sufficiently small $\lambda$, every $F_{\lambda}$ contains a (normalized) character of $G$. This certainly happens when $G$ is nilpotent, and in that case, the density question (whether $\cup_{\lambda > 0} \partial_e F_{\lambda}$ is dense in $\partial_e S(R_f,[1,0])$) reduces to whether every $\tau \in \partial_e F_0$ is multiplicative. This happens rarely---with a modicum of noncommutativity, we can usually find a pair of elements $g,h$ of $G$, together with $\tau \in \partial_e F_0$ \st $\tau(x_{gh}) \neq \tau(x_{hg})$ (or even  one is nonzero, the other is zero).

This type of maximal order ideal appears almost ubiquitously, e.g., if $G$ is indicable, but can also appear in torsion-free non-indicable groups, such as $\Z \times_{\theta} \Z$ where $\theta (m) = -m$. This is a central extension of the dihedral group, and is not indicable. 


Let $G$ be finitely generated torsion-free nilpotent group of class two. Then $G/G' \iso \Z^d$ for some $d$, and $G' \subseteq Z(G)$, the centre. Let $\Arrow \pi; G. G/G'$ be the factor map. Let $\brcs{g_i}$ be a collection of representatives of elements of $G$ \st $\pi(g_i):= {}i$ is the standard basis for $\Z^d$. If $v = (v(1), v(2), \dots, v(d)) \in \Z^d$, let $g_v$ denote the word in $\brcs{g_i}$, $g_{1}^{v(1)}g_2^{v(2)}\dots g_{d}^{v(d)}$ (written in that order), so that $\pi(g_v) = v$.

If $W = \brcs{w^{j}}$ is a finite subset of $\Z^{d}$, let $K$ be the convex hull of $W$ inside $\R^d$. Let $\brcs{v_1, \dots , v_m}$ be the set of extreme points (often called vertices) of $K$; this is a subset of $W$. If $v$ is an extreme point of $K$, a {\it nearest neighbour\/} will denote a point in $W$ that is closest to $v$ along some edge in $K$ emanating from $v$---there is exactly one nearest neighbour to $v$ for each such edge.

Let $nW$ denote the set of sums of $n$ elements of $W$. Then it is easy to check that if $v$ an extreme point of $K$ and $v'$ is a nearest neighbour to $v$, then $(n-1)v + v'$ cannot be realized in any other way as a sum of $n$ elements of $W$.

Now let $f \in A G^+$ be  admissible. Then $\pi(f) \in (\Z^d)^+$ and $\pi(S)$ is admissible (\wrt $\Z^d$). Taking $K$ to be the convex hull of $\pi(S)$, we have that $K$ contains an open subset of $\R^d$, since $\pi(S)$ is admissible. Select a vertex $v$ of $K$, and let $p g_v$ with $p \in (A G')^+$ be the corresponding component of $f$, that is, $p = \sum_{z \in G'} (f,zg_v)z$, and let $v'$ be a nearest neighbour, with corresponding component $q g_{v'}$.

Then we can write $g_v g_{v'} = F(v,v') g_{v'}g_v$, with $F(v,v') \in G'$. At this point, we make the following assumption:

\noindent {\it Assumption.} $F(v,v') \neq 1$ (for at least one choice of extreme point $v$ and one of its nearest neighbours).

Now let $n$ be a positive integer; we consider the component of $g_{(n-1)v + v'}$ in $f^n$. Because of the uniqueness result for $(n-1)v + v'$, the only products of $n$ terms that will yield a contribution to $g_{(n-1)v + v'}$ are exactly the $n+1$ possibilities
$$
p^{n-1}q g_v \cdots g_v \cdot g_{v'} \cdot g_v \dots g_v = p^{n-1}q F(v,v')^{k-1} g_{v'}(g_v)^n
$$
the $v'$ term in position $k$. Since $F(v,v') \neq 1$ and $G'$ is torsion-free for each of $j= 0, 1,\dots, n-1$, there must be $z_j $ in the support of $F(v,v')^{j}$ that is not in the support of any of the other powers. This easily translates to the same property for the components of $g_{(n-1)v + v'}$, and means that in the Bratteli diagram at level $n$, there are $n$ points each with a unique predecessor, and they have at least one outgoing edge to one of the corresponding points at level $n+1$.

This yields infinitely many paths each corresponding to a maximal order ideal of the form (**). Moreover, since the corresponding traces depend on the order in which the terms are multiplied, they cannot be multiplicative.

We conclude that if the class two torsion-free nilpotent group and $S$ satisfy the assumption, then $R_f$ has infinitely many maximal order ideals, and has a non-multiplicative trace. The latter implies that the density condition, $\cup_{\lambda > 0} \partial_e F_{\lambda}$ is dense in $\partial_e S(R_f,[1,0])$, fails.

If $G$ is a nilpotent torsion-free properly class two (that is, $G' \neq 0$) group, then it is easy to find $S$ \st at least one extreme point satisfies the assumption, and in fact, practically all choices for $S$ will satisfy the assumption for most of their pairs of extreme point and nearest neighbours. I was not able to prove that {\it all\/} choices of $S$ will admit such a choice, but I imagine this is true. If it were, in the following, {\it for all choices\/} can be replaced by {\it for some choice.} As it stands, practically any $S$ will do.

\Lem Corollary \IMIsix. Suppose that $G$ is a finitely generated nilpotent torsion-free group \st for all choices of admissible $f$, either all $\tau \in \partial_e F_0$ are multiplicative, or $R_f$ has only finitely many maximal order ideals. Then $G$ is abelian.

\Pf If $G$ is not abelian, then $G' \neq 0$ and thus $G'' \neq G'$; so $G_0:= G/G''$ is nonabelian, torsion-free, and of class two. By the preceding, we can find admissible $f_0 \in (A G_0)^+$ \st both the density condition and the finitude of the set of maximal order ideals fail for $R_{f_0}$. We can replace $f_0$ by any power of itself (this does not change $R_{f_0}$). By raising it to a sufficiently large power, there is enough room in $\pi(f_0)$ so that a set of generators appear in the interior of the corresponding convex polytope $K$. It is thus easy to lift the power of $f_0$ to an admissible $f \in (A G)^+$ \st $\Arrow \Phi; G. G_0$ induces $\Phi(f) = f_0$.

This $\Phi$ also induces an onto homomorphism $R_f \to R_{f_0}$, so that the corresponding traces lift. The traces are discrete traces, so their kernels are maximal order ideals. It remains to solve the corresponding problem for the discrete Heisenberg group and arbitrary admissible $f$. 

We require one more result to finish the proof of Corollary \IMIsix.

\Lem Lemma \IMIsev. Let $G$ be a finitely generated torsion-free class two nilpotent group, and let $\Arrow \pi; G.G/G'$ be the quotient map. Let $H$ be a subgroup of $G$ \st $\pi(H)$ is of finite index in $\pi(G)$. Then {\par}
\item{(a)} $H$ is of finite index in G;
\item{(b)} $Z(H) \subset Z(G)$

\Pf For finitely generated torsion-free nilpotent groups, each of $G/G'$, $G'/G''$, etc, is finitely generated torsion-free. So we can write $G/G' \iso \Z^d$ for some $d >0$.

\noindent (a) We first show that if $\pi(H)$ is all of $\pi(G)$, then $H = G$.

Pick $h_i \in H$ \st $\pi(h_i)$ is the $i$th standard basis element of $\Z^d$. Let $\brcs{g_j}$ be a generating set for $G$, and write $z_{j,k} = g_j g_{k} g_{j}^{-1} g_{k}^{-1} $; then $\brcs{z_{j,k}}$ is a generating set for $G'$, and is contained in the centre of $G$ (since $G$ is of class two). Given $j$, there exists $h_{(j)} \in H$ \st $\pi(h_{(j)}) = \pi(g_{j})$. Hence $h_{(j)}g_j = z^{(j)}g_j h_{(j)}$ where $z^{(j)} \in G'$. It follows that $[h_{(j)},h_{(k)}] = [g_j,g_k] = z_{j,k}$. Hence $G'\subset H'$, and since $\pi(G) = \pi (H)$, it follows that $G = H$.

Now assume that $\pi(H)$ is of finite index. For each $i = 1,2, \dots,d$, there exists $g_i \in G$ \st $\pi(g_i) $ is the $i$th standard basis element of $\Z^d$. Then $\brcs{g_i}$ generates $G$, by the previous paragraph. For each $i$, there exists $n(i)> 0$ \st $\pi(g_i^{n(i)}) \in \pi(H)$. With $z_{i,j} = g_j g_{k} g_{j}^{-1} g_{k}^{-1} \in G'$ as before, the set of these generate $G'$, and an easy computation reveals that $[g_i^{n(i)},g_j^{n(j)}] = z_{i,j}^{n(i)\cdot n(j)}$. Hence $H'$ is of finite index in $G'$, and it is then immediate that $H$ is of finite index in $G$.

\noindent (b) There exists $N$ \st for all $g \in G$, $g^{N} \in H$. If $z \in Z(H)$, then $[z,g^N]=1$; but since $G$ is class two, $[z,g]^N = [z,g^N] = 1$. Thus $[z,g]$ is torsion; as $G$ is torsion-free, $[z,g] = 1$ for all $g \in G$.
\qed

\noindent {\it Conclusion of proof of Corollary \IMIsix.} Now we show no matter what the admissible $f$, there exists an extreme point of $K = \cvx S$, together with a nearest neighbour, that satisfies the assumption. Let $v$ be an extreme point, and let $w$ vary over its nearest neighbours. Let $J$ be the subgroup of $\pi(G)= \Z^d$ generated by $v$ and all its nearest neighbours. The convex hull of $v$ and all its nearest neighbours contains an open $d$-ball (true for any compact polytope with nonempty interior), hence the rank of $J$ is $d$, and thus $J$ is of finite index in $\Z^d$. Now let $H_v$ be the subgroup of $G$ generated by $\brcs{g_w} \cup \brcs{g_v}$.

By Lemma \IMIsev(b), the centre of $H_v$ is contained in the centre of $G$. However, if none of the $g_w$ satisfy the assumption, then $[g_v,g_w]=1$ for all nearest neighbours $w$. Hence
$g_v \in Z(H_v) \subset Z(G)$. Now do this for every extreme point $v$---failure of the assumption entails that all $g_v$ belong to the centre of $G$. However, $\cvx \brcs{v}_{v \in \partial_e K}$ contains an open ball, so that the subgroup of $\Z^d$ generated by the extreme points is of finite index. By Lemma \IMIsev(a), the subgroup of $G$ generated by all the $g_v$ is of finite index, and is contained in the centre of $G$. Hence $G$ is abelian (for a torsion-free nilpotent group $G/Z(G)$ is torsion-free).

Thus we have shown that if $S$ is admissible and $G$ is finitely generated torsion-free class two nilpotent, then there exists an extreme point together with a nearest neighbour that satisfy the assumption.
\qed

\SecT  \Sty\ Realizing stationary dimension groups

If $B$ is a nonnegative integer square matrix of size $n$, the {\it stationary dimension group\/} obtained from $B$ is the direct limit (as partially ordered abelian groups, each $\Z^n$ equipped with the coordinatewise ordering) with repeated multiplication by $B$, $\lim \Arrow B; \Z^n.\Z^n$. The limit dimension group is simple iff $B$ is primitive, and in that case, can be obtained from a primitive $0-1$ matrix (possibly required to be of larger size). Within the class of dimension groups, the simple stationary ones are relatively easy to characterize. Stationary simple dimension groups have unique trace, with values (after rescaling) in $\Q[\lambda]$, where $\lambda$ is the Perron eigenvalue of $B$. They can be classified by means of a combination of ideal classes in orders in number fields, and abelian extensions of torsion-free abelian groups.

Here we will show (slightly more than) every simple stationary dimension group can appear as an $R_f/M$ where $G = F_2$ (the free group on two generators), and some maximal order ideal $M$ of $R_f$ for suitable choice of admissible $f$, and $A = \Z$. In  fact we realize the matrix implementing the stationary dimension group, $B$, as the map $\Z \Gamma'_{n,M} \to \Z \Gamma'_{n+1,M}$ for all $n$. 

\Lem Example \onetty. Let $F_2$ be the free group  with generators $\brcs{g,h}$. Let $B$ be a $k \times k$ primitive $0-1$ matrix. Then there exists $f \in \Z F_2^+$ with the following properties.
\item{(0)} The coefficients of $f$ are all $0$ or $1$;
\item{(a)} $f$ is symmetric \wrt both $g \mapsto g^{-1} $ and $h \mapsto h^{-1}$;
\item{(b)} $f$ is admissible;
\item{(c)} there exists a quotient of $R_f/\SS R_f$ by an order ideal which is order isomorphic to the stationary dimension group, $\Arrow \lim B; \Z^k. \Z^k$;
\item{(d)} for every $n$, there exists a subset $\Gamma^B_n \subseteq \Gamma_n$ \st $|\Gamma^B_n| = k$, $f_n|\Z \Gamma^B_n$ has image in $\Z \Gamma^B_{n+1}$ and \wrt a natural ordered $\Z$-basis, its matrix is $B$, and if $\gamma \in \Gamma_n \setminus \Gamma^B_{n}$, then $\supp f {}{\gamma} \cap \Gamma^B_{n+1} = \emptyset$.

\Rmk Condition (d) is precisely what we need to prove (c). A consequence is that every stationary simple dimension group appears as a quotient by an order ideal of some $R_f$, if we restrict $G$ to be the free group on two generators. It is rather easy to see that if $G$ is abelian by finite, say with minimal index of the torsion-free abelian subgroups being $k$, then all simple quotients of $R_f$ have width at most $k$ (so all simple quotients have rank at most $k$).

\Pf For each $i = 1, 2, \dots, k$, define the element $v_i = g^{-(i-2)}h^{i-1}$, and set $w_{ij} = v_j g v_i^{-1}$. Let $f_0 = \sum_{B_{ij} = 1} {}{w_{ij}} $ and set $f = f_0 + \sum_{B_{ij} = 1} {}{w_{ij}^{-1}} + 1_G$. We will establish the various properties for $f$.

Explicitly, $w_{ij} = g^{-(j-2)} h^{j-1} g h^{-(i-1)}g^{i-2}$, and it is easy, although   tedious, to show that $w_{ij}$, $w^{-1}_{i',j'}$ are distinct from each other and from $1_G$ (it helps that the net degree of each $w_{ij}$ \wrt $g$ is one); so the nonzero coefficients of $f$ are all one.

Let $\Arrow d; G. \Z$ be the group homomorphism given by $g\mapsto 1$ and $h \mapsto 1$; then $d(w)$ is the total multiplicity of $g$ and $h$ in the word $w$, and we refer to it as simply the degree. We also have a similar group homomorphism for $g$ and $h$ separately, but we shall not use them as much. We note that $d(v_i) = 1= d(w_{ij})$, and $d_{w_{ij}^{-1}} = -1$.

We first note that $\Gamma_1$ is by definition $\supp f \setminus\brcs{1}$, so is $\brcs{w_{ij}} \cup \brcs{w_{ij}^{-1}}$, partitioning them into those of degree one and of degree $-1$. Then the elements of $\Gamma_n$ have degree between $-n$ and $n$, and those of degree $n$ are precisely the products of $n$ of the $w_{ij}$. Define the following subset of $\Gamma_n$, $\Gamma_n^B:= \brcs{v_i g^{n-1}}_{i=1}^k$.

First, $v_1 = g$, $v_2 = h$, and $w_{1i} = v_i g g^{-1} = v_i$. Thus $v_i \in \Gamma_1$, and we can write $v_i g^{n-1} = (v_i g v_j^{-1}) (v_j g^{n-2}) = w_{ji} v_{j}g^{n-2}$. This yields inductively that $v_i g^{n-1} \in \supp f^n$. Since $d(v_i g^{n-1}) = n$, it follows that $v_i g^{n-1} \not\in \supp f^l$ for any $l < n$. Hence $\Gamma_n^B \subset \Gamma_n$.

We also see that $g,h,g^{-1}, h^{-1} \in \supp f$, so that $\cup_{n\geq 0}\supp f^n = G$.

We consider $({}{v_1 g^{n-1}}, {}{v_2 g^{n-1}}, \dots, {}{v_k g^{n-1}})$ as an ordered basis for $\Z \Gamma_{n}^B$. The claim is that the restriction and compression of left convolution by $f$ yields sends $\Z \Gamma_{n}^B \to \Z \Gamma_{n+1}^B$, the matrix of the resulting transformation is just $B$, and most importantly, if there is a path of length $n-m$ from a point in $\Gamma_m$ to $\Gamma_n^B$ (permitted by the transition given by left multiplication by $f$), then it could only have arisen by transitions through $\Gamma_m^B$, $\Gamma_{m+1}^B$, \dots, $\Gamma_{n-1}^B$.

Under $f$, the permitted transitions from $v_i g^{n-1}$ to a point in $\Gamma_{n+1}$ are given by left multiplication by an element of the form $v_{s}gv_{t}^{-1}$ for some $s,t$ \st $B_{s,t} = 1$. If $v_s g v_{t}^{-1} v_i g^{n-1}\not\in \Gamma_{n+1}^B$, then we disregard it. If on the other hand, $v_s g v_t^{-1} v_i g^{n-1} \in \Gamma_{n+1}^B$, then we have an identity
$$\eqalign{
v_s g v_t^{-1} v_i g^{n-1} &= v_r g^n \quad \text{for some $r$; this is equivalent to } \cr
v_r^{-1} v_s g v_{t}^{-1}v_i &= g, \qquad \text{which reduces to }\cr
h^{1-r}g^{r-2}g^{2-s}h^{s-1}g h^{1-t}g^{t-2}g^{2-i}h^{i-1} & = g; \qquad\text{that is,}\cr
h^{1-r}g^{r-s}h^{s-1}g h^{1-t} g^{t-i}h^{i-1}g^{-1} & = 1_G. \cr
}$$
If $r=s$, we quickly see that $h^{1-t} g^{t-i}h^{i-1} = 1_G$, from which it follows that $i=t$. Conversely, if $i = t$, we deduce $h^{1-r}g^{r-s}h^{s-1} = 1_G$, whence $r=s$. So assume that $r \neq s$, and thus $i \neq t$. From the degree in $g$, we have $r-s + 1+ t-i -1 = 0$, so $(r-s) = (i-t)$.
If $s,t,i \neq 1$, the word is in reduced form, hence we reach a contradiction.

If $s=1$, the word reduces to $h^{1-r}g^{r-s+1}h^{1-t}g^{t-i}h^{i-1}g^{-1} = h^{1-r}g^{r}h^{1-t}g^{t-i}h^{i-1}g^{-1} $; if both $t \neq 1$ or $i \neq 1$, the word is again in reduced form, again a contradiction. If $t = 1$, we obtain a reduction to $h^{1-r} g^{r+t-i}h^{i-1} g^{-1}$, and if $i \neq 1$, obtain $h^{1-r} g^{r+1-i}h^{i-1}g^{-1}$; this reducing to the trivial word implies $r+ 1 =i$ and $r+1-i -1 = 0$, a contradiction. If $t = i = 1$, we have $h^{1-r} g^r g^{-1} = 1$, whence $r = 1 = s$ and $t = i$.

In all cases, $i = t$ and $r=s$, which is exactly what is needed to show that the matrix representation of $\Gamma^B_n \to \Gamma_{n+1}^B$ is just $B$ (\wrt the ordered bases $(v_i g^{n-1})$ and $(v_i g^{n})$).

Now we show that if $\gamma \in \Gamma_k$ and there exists $\beta = w_{ij} $ or $w_{ij}^{-1}$ or $\beta = 1$ \st $\beta\gamma \in \Gamma_{k+1}^B$, then $\gamma \in \Gamma^B_k$ and $\beta \in \brcs{w_{ij}}$. From degree of elements of $\Gamma^B_{k+1}$ being $k+1$, we must have $\deg \gamma = k$ and $\deg \beta = 1$. The latter forces $\beta = \brcs{w_{ij}}$ for some $(i,j)$ with $B_{i,j} = 1$. Write $\beta \gamma = v_t g^{k} \in \Gamma_{k+1}^B$. This yields
$$\eqalign{
g^{-(j-2)} h^{j-1} g h^{-(i-1)}g^{i-2} \gamma &= g^{-(t-2)}h^{t-1}g^k; \text{ that is,}\cr
\gamma &= g^{2-i}h^{i-1}g^{-1}h^{1-j}g^{j-t} h^{t-1}g^k.\cr
}$$
However, $\gamma \in \Gamma_k$ and of degree $k$ entails that $\gamma$ is a product of a string of $k$ of the $w_{rs}$, say $\gamma = w_{r_k,s_k} w_{r_{k-1},s_{k-1}}\cdots w_{r_1, s_1}$; we are required to show that $s_l = r_{l-1}$. We proceed by induction on $k$. From $w_{ij}w_{r_k,s_k} w_{r_{k-1},s_{k-1}}\cdots w_{r_1, s_1} = v_t g^k$, we have $\gamma = w_{ij}^{-1} v_t g^k$, so

This portion redone (previous paragraph). We show, by induction on $n$, that if
$$
v_i g^{-1}v_j^{-1} v_k g^{n-1} \text{ is a product of $n-1$ elements of $\Gamma_1^B$,}
$$
then $ j= k$. First, we consider the case $ n=1$. In that case, $v_i g^{-1}v_j^{-1}v_k = 1_G$ for some $l$. This leads to the equation,
$$\eqalign{
g^{2-i}h^{i-1}g^{-1}h^{1-j}g^{j-2}g^{2-k}h^{1-k} &= 1_G, \text{ that is,}\cr
g^{2-i}h^{i-1}g^{-1}h^{1-j}g^{j-k}h^{1-k} & = 1_G.\cr
}$$
Assume $j \neq k$. If both $i,j \neq 1$, then the word on the left is in reduced form (modulo the possibilities that $i=2$ or $k=1$), a contradiction. If $i=1$, we obtain $h^{1-j}g^{j-k} h^{1-k} = 1_G$, which is impossible. If $ j =1$ and $i \neq 1$, we obtain $g^{2-i}h^{i-1}g^{-k}h^{1-k} = 1_G$, which forces $k =1$, and thus $j = k$.

Now assume the result is true for all $t < n$, and suppose that $\gamma = v_{j_1}g v_{i_1}^{-1} \cdot v_{j_2}g v_{i_1}^{-2}\cdot \dots v_{j_{n-1}}g v_{i_{n-1}}^{-1}$ and $\gamma = v_i g^{-1}v_j^{-1}v_k g^{n-1}$ for some $i,j,k$, and we wish to conclude $j=k$. This yields the equation,
$$
(v_k^{-1}v_j)g(v_i^{-1}v_{j_1})g(v_{i_1}^{-1} v_{j_2})g \cdots g (v_{i_{n-2}^{-1}}v_{j_{n-1}})gv_{n-1}^{-1}g^{1-n} = 1_G.
$$

Now for any $a,b$, $v_a^{-1}v_b = h^{-(a-1)}g^{a-b} h^{b-1}$, so that $a\neq b$ implies $ v_a^{-1}v_b $ is in reduced form; that is, if $a\neq 1$ and $b \neq 1$, it is reduced with three monomials, if $a =1$ (and $b \neq 1$), it is $g^{1-b}h^{b-1}$ having two monomials, while if $b=1$ (and $a\neq 1$), it is $h^{1-a}g^{a-1}$. If any of the $i_l = j_{l+1}$, the expression for $\gamma$ becomes
$$
\gamma = v_{j_1}g v_{i_1}^{-1} \cdot v_{j_2}g v_{i_2}^{-1}\cdots g v_{i_{l-1}}^{-1} v_{j_l}g \cdot g (v_{i_{l+1}}^{-1}v_{j_{l+2}}) g \cdots v_{j_{n-1}}g v_{i_{n-1}}^{-1}.
$$

We write the equation
$$\eqalign{
(v_k^{-1}v_{j_1})g^{k_1} (v_{i_1}^{-1}v_{j_2}) g^{k_2}\cdot \dots \cdot (v_{i_{r-1}}^{-1}v_{j_r}) g^{k_r}v_{i_r}^{-1}g^{n-1} &= 1_G, \quad \text{ which rewrites to }\cr
(v_k^{-1}v_{j_1})g^{k_1} (v_{i_1}^{-1}v_{j_2}) g^{k_2}\cdot \dots \cdot (v_{i_{r-1}}^{-1}v_{j_r}) g^{k_r}(v_{i_r}^{-1}v_1)g^{-n} & = 1_G.
}$$
where none of $i_s = j_{s+1} $ or $k = j_1$ or $i_r = 1$. If none of the $j_1$ are one, then the word at left is in reduced form, hence the relation is impossible. Now we proceed by induction on the number of subwords of the form $(v_a^{-1}v_b) g^s$. Problems arise when some of the $i_s$, $j_t$ are equal to one, because of the following (put in reduced form):
$$
v_a^{-1} v_b = \cases h^{1-a}g^{a-b} h^{b-1} & \text{if $a,b \neq 1$}\\
h^{1-a}g^{a-1} & \text{if $a \neq 1$ and $b=1$}\\
g^{1-b}h^{b-1} & \text{if $a = 1$ and $ b\neq 1$.}\\
\endcases
$$
We do not consider the case $a = b =1$, since in our formulation, $a \neq b$ (since otherwise $r$ drops). In order for the relation to hold, we must have at least one subword of the form $h^v g^w g^x g^y h^z$, or $h^v g^x g^y h^z$ where $v, w \neq 0$ and $w + x + y$, respectively $x+y = 0$ (otherwise, there would be no further reduction in the big word). These situations can only arise under the following circumstances.

\noindent {(a)} $j_{s+1} = 1 = i_{s+1}$, arising from
$$g^{k_s}(v_{i_s}^{-1}v_{j_{s+1}})g^{k_{s+1}} (v_{i_{s+1}}^{-1}v_{j_{s+2}})g^{k_{s+2}} = g^{k_s}h^{1-i_s}g^{i_s -1} g^{k_{s+1}} g^{1-j_{s+2}} h^{ j_{s+2}-1}g^{k_{s+2}};$$
then the middle grouping of $g$s
being trivial implies $k_{i_{s+1}} = j_{s+2} - i_s$. The expression simplifies to $g^{k_s}h^{1 - i_s + j_{s+2} -1}g^{k_{s+2}} = g^{k_s} h^{j_{s+2} - i_s} g^{k_{s+2}}$, which gives us
$$
g^{k_s} h^{k_{s+1}} g^{k_{s+2}} = g^{k_s + k_{s+1}}(g^{1-(k_{s+1} +1)}h^{(k_{s+1}-1) + 1}g^{k_{s+2}}.
$$
This is the same as $g_{k_s + k_{s+1}} (v_1^{-1} v_{k_{s+1}})g^{k_{s+2}}$. The number of $(v_a^{-1}v_b) g^s$-words has just dropped by one.

\noindent (b) Consider the following situation where cancellation occurs (on the $g$s) in a subword of the form
$$
g^r (h^{1-u}g^{u-1})g^s (h^{1-z}g^{z-1})g^t,
$$
where $r,s,t > 0$ and $u,z > 1$. Then cancellation of the interior $g$s entails $s + u -1 = 0$, contradicting $s + u > 3$. Hence this cannot occur.

(c) Consider the following situation where cancellation occurs (on the $g$s) in a subword of the form
$$
g^r (h^{1-u}g^{u-z}h^{z-1})g^s (g^{1-w}h^{w-1})g^t,
$$
where $r, s, t > 0$, $u,z,w>1$, and $u \neq v$ (corresponding to $g^* (v_{i_a}))^{-1}v_{j_{a+1}}g^{r_{a+1}}(v_{i_{a+1}}^{-1}v_1)g^{r_{a+2}}$). For cancellation to occur, we would have to have $w = s+1$, whereupon the subword simplifies to
$$\eqalign{
g^r (h^{1-u}g^{u-z}h^{z-1})h^{w-1})g^t & = g^r h^{1-u}g^{u-z} h^{z+w-2}g^t\cr
& = g^r (h^{1-u}g^{u-1})g^{w-1}(g^{-{z+w-2}}h^{z+w-2})g^t,\cr
}$$
which is of the form $g^r (v_u^{-1}v_1)g^{w-1}(v_1^{-1}v_{z+w-1})g^t$, and there is no cancellation of the $g$ terms. We have thus replaced a subword with a subword of the same length, without interior cancellation of the $g$, and \st both the left and right parenthesized terms use a subscript $1$.

With the newly modified word, there must be cancellation of some consecutive $g$ terms, except now there are fewer places to look---either the new cancellation is of the type described above (with three consecutive $g$s), in which case the number of terms drops, or it is of type (c), in which case, it is eliminated by what we have just done. The process terminates when we run out of terms.
\qed



\def\Ud #1.{\lfloor #1 \rfloor \cdot \lceil #1 \rceil}

\def\flo #1{\lfloor  #1 \rfloor}
\def\cei #1{\lceil  #1 \rceil} 

\SecT \One\ A lot more than you ever wanted to know about the discrete Heisenberg group---a case study

Over the next two sections, we provide detailed information about the discrete  Heisenberg group $ G = H_3$ with  standard admissible  $f$. In this section, we  show that $(G,f)$  fails to satisfy WC, and does so in a spectacular way---$\tilde l_S $ is bounded (but not  constant, except trivially) on each  conjugacy class. This is probably what happens for any finitely generated nilpotent nonabelian torsion-free group and any admissible set. This is in sharp contrast to the abelian case.

We also show that the lattice of space-time cones contained within that generated by the group element $1$   at time zero (equivalently, the lattice of order ideals contained in $R_f$ for $f$ a specific admissible symmetric element) is {\it  not\/} noetherian, despite the fact that the group rings $\Z G$, $\R G$, $\C G$ are all right and left noetherian rings. This also contrasts sharply with the abelian case. 

The underlying Bratteli diagram for $R_f/\SS R_f$ is parabolic  in the sense that the cross-sections, that is, the levels, have the shape of a discretized parabola. 

 In the following two sections, 
we determine the pure trace space of $R_f$. It turns out that there are  four relatively discrete sets of  discrete pure traces, easily described, whose union is dense in the set of pure traces that kill an order ideal---again differing substantially from the torsion-free abelian case. But at least the trace space is a Bauer simplex.

The group has generators $\brcs{g,h,z}$ subject to the relations, $hg = zgh$ and $z$ is central (so that $z$ can be dropped from the set of generators). There is a normal form for elements, given by $z^r g^a h^b$ with $(r,a,b) \in  \Z^3$. Let $S = \brcs{g^{\pm1}, h^{\pm1}, 1}$. Then $S$ is an admissible set, and $l_S(z) = 4$ (as $z = hgh^{-1}g^{-1}$ and $ z \not\in S^3$).

Initially, we obtain an explicit formula for $l_S(z^r g^a h^b)$ (as a function of $(r,a,b)$), and use it to obtain $|a| + |b| \leq \tilde l_S(z^r g^a h^b) \leq |a| + |b| + 2$ for all $r$. A particular consequence is that if $u$ is any element of $G$, then $\RR_{g, \tilde l_S(g)+2}$ is an order-bounded endomorphism of $R_f$ (where $f$ is any admissible element of $\R G$ with $\supp f = S$). Lemma \Bthr(a,e) says this is what we would have expected---if $G$ satisfied WC, which it doesn't.

For $(a,b)$ in $\Z^2$ and $m$, define the finite (possibly empty) set of integers,
 $$
S(a,b,m) = \Set{r \in \Z}{z^r g^a h^b \in S^m}. 
$$
Thus $l_S(z^r g^a h^b)\leq m$  iff $r \in S(a,b,m)$. 

A set of integers $T \subset \Z$ will be called an {\it interval\/} of  integers if it consists of consecutive integers (alternatively, of the form $[k,k'] \cap \Z$). We will show that every nonempty $S(a,b,m)$ is an interval, and then determine its endpoints. They are determined by a piecewise quadratic formula. This will effectively determine $l_S(z^r g^a h^b)$, and then the results about $\tilde l_S$ easily follow. 

We will obtain an explicit and easily computed criterion for membership of $z^r g^a h^b$ in $S^m$, that is, determining $S(a,b,m)$.  First, we will show (provided $m \geq |a|+ |b| \geq 0$)   that  $S(a,b,m)$ is an interval of integers, that is, of the form $\brcs{-x, -x+1, \dots, 0, \dots, X}$. Thus $S(a,b,m)$ is determined by its maximum and minimum values. We can reduce (via Lemma \oneone) to the situation that $a,b \geq 0$. Then we give lower bounds for $\max S(a,b,m)$ and $-\min S(a,b,m)$ via explicit constructions (these obtained by optimization). The most arduous part is then to show these lower bounds are sharp. This is done by an unusual induction---but with a   large number of cases. 

 Every element of $S^m$ can be written (possibly in many different ways) in the form 
$$
h^{\gamma(1)}g^{\epsilon(1)}h^{\gamma(2)}\cdots h^{\gamma(k)} g^{\epsilon(k)} \tag*
$$ 
subject to the constraints $\gamma(i), \epsilon(i) \in\Z$;  if $i \neq 1$, then $\gamma(i) \neq 0$; if $i \neq k$, then $\epsilon(i) \neq 0$; and $\sum_i (|\gamma(i)| + |\epsilon(i)|) \leq m$.
The following is elementary.

\Lem Lemma \oneone. (a) $S(a,b,m) \subseteq S(a,b,m+1)$. 
\item{(b)} If $m < |a| + |b|$, then $S(a,b,m)$ is empty.
\item{(c)} If $m \geq |a| + |b|$, then $0 \in S(a,b,m)$. 
\item{(d)} $S(a,-b,m) = -S(a,b,m) = S(-a,b,m)$.
\item{(e)} If $m-a-b$ is even, then $S(a,b,m) = S(a,b,m+1)$.
\item{(f)} $S(a,b,m) = ab - S(b,a,m)$

\Pf (a) If an element of $G$ can be represented as a product of $m$ members of  $S$, then multiplication by $1$ makes it a product of $m+1$ members.  (b) Follows from the form (*), where $l_S(z^r g^a h^b) \geq |a| + |b|$. (c) $g^a h^b \in S^m$. (d) The map given by $g\mapsto g$, $h \mapsto h^{-1}$, $z \mapsto z^{-1}$ preserves the relation, so induces an automorphism of $G$. Under this automorphism, $S$ is left stable, and thus so is $S^m$ for every $m$. If $r \in S(a,b,m)$, then $z^r g^a h^b \in S^m$; applying the automorphism, we obtain $z^{-r} g^a h^{-b} \in S^m$. Thus $-r \in S(a,-b,m)$. The corresponding map with $g$ and $h$ interchanged yields the other equality of sets. 

\noindent (e) From the form (*), in addition to the constraints given there, we also have two more, $\sum \epsilon (i) = a$ and $\sum \gamma(i) = b$. Pick $r \in S(a,b,m+1) \setminus S(a,b,m)$, and a corresponding element $u = z^r g^a h^b \in S^{m+1} \setminus S^m$, so there exist  sequences of $\epsilon$s and $\gamma$s \st $\sum (|\epsilon(i)| + |\gamma(i)|) = m+1 $. Modulo $2$, we thus have 
$$ \eqalign{
m+1 &\equiv \sum \epsilon(i) + \sum \gamma(i) \cr
& \equiv a + b, \cr
}$$
contradicting the parity of $m-a-b$. 

\noindent (f) The assignment $g \mapsto h$, $h \mapsto g$, and $z \mapsto z^{-1}$ induces an automorphism that leaves each $S^m$ stable, and sends $z^r g^a h^b$ to $z^{-r} h^a g^b = z^{-r+ ab} g^b h^a$. Hence $r \in S(a,b,m)$ entails $ab -r \in S(b,a,m)$; conversely, if $s \in S(b,a,m)$, then $z^s g^b h^a \in S^m$; applying the automorphism yields $z^{-s} h^b g^a = z^{ab-s}g^a h^b$, so that $ab-s:= r \in S(a,b,m)$ and $ab-r = s$. 
\qed

Properties (b)--(f) of Lemma \oneone\  can be explained by a more general phenomenon. There is a natural representation of the dihedral group $D_4$ as a group of automorphisms of $H_3$. Specifically, we can construct four automorphisms uniquely determined by 
$$
g \mapsto g^{\varepsilon(1)}, \quad h \mapsto h^{\varepsilon(2)}, \quad z \mapsto z^{\varepsilon (1) \cdot \varepsilon(2)},
$$
where $\varepsilon (i) \in \brcs{\pm 1}$. Then the (uniquely determined) automorphism given by 
$$ 
g\mapsto h, \quad h \mapsto g, \quad z \mapsto z^{-1},
$$
together with the previous four generate $D_4$. The last displayed automorphism applied to  the generic word  sends $z^r g^a h^b \mapsto z^{ab-r}g^b h^a$. Of course, an alternative proof of Lemma \oneone\ can be obtained using this action of $D_4$. Of particular importance is the behaviour of the quadrants (restrict to the second and third coordinates) under these automorphisms.

Define the {\it defect\/} of the element $z^r g^a h^b$ or the triple $(a,b,m)$ to be $d:= m - |a| - |b|$. Now we want to show that $S(a,b,m)$ is an interval (if it is nonempty). By (d), we may assume $a,b \geq 0$, and by (e), we can assume that $m-a-b$ is even, and by (b), $d \geq 0$. In this case, $d = m-a-b$, and since $d$ is even,  we can write it  as $d = 2\partial$, where $\partial$ is a nonnegative integer.

First, we have a special case. It is convenient to remind the reader of the combinatorial properties of the coefficients appearing in the expansion of powers of $g+h$. We use inner product notation to denote the coefficient: thus if $j= \sum a_w w$ is an element of the group ring (with $w$ varying over a finite set of group elements), then the coefficient, $a_w$, of $w$ is denoted $(j, w)$.

 Let $p(r,a,m-a)$ denote the coefficient of $z^r g^a h^{m-a}$ appearing in $(g+h)^m$; we will shortly see that in order to be nonzero, necessary and sufficient is that  $0\leq r \leq (m-a)a$ and $0 \leq a \leq m$. As is well known, $p(r,a,m-a)$ has combinatorial significance; we sketch the correspondence. 

 Define (for $r,a,m-a \geq 0$),
 $$\eqalign{
 \PP(r,a,m-a) &= \Set{(s_1,s_2,\dots, s_a) \in \Z_+^m}{m-a \geq s_1 \geq s_2 \geq \cdots \geq s_a; r = \sum s_i}\cr 
 P(r,a,m-a) &= \left| \PP(r,a,m-a)\right|.
 }$$
 It is easy to verify that $\PP(r,a,m-a) \neq \emptyset$ provided that $0 \leq r \leq a(m-a)$. $P(r,a,b)$ counts the number of integer partitions of $r$ \st at most $a$ nonzero summands appear, and each is bounded by $b$ (which for much of our discussion is given by $m-a$. 

 To see that $P = p$, we consider a word in of length $m$ with $a$ $g$s and $m-a$ $h$s. This can be written in the form $h^{t(a)}gh^{t(a-1)}\dots gh^{t(1)}gh^{t(0)}$ where $t(j)$ are nonnegative integers adding to $a$.

  The exponent of $z$ will be $t(1) + 2t(2) + 3 t(3) + at(a)$ (obtained by starting from the right, and moving the accumulating $g$s to the left). Hence the multiplicity of $z^r g^a h^{m-a}$ in $(g+h)^m$ is the number of $a+1$-tuples $(t(0), t(1), \dots, t(a))$ of nonnegative integers \st $r = \sum_1^{a} jt(j) $ and $\sum_1^a t(j) \leq m-a$. There is an obvious bijection between the collection of these and $\PP(r,a,m-a)$, given by setting $s_1 = \sum_1^a t(j)$, $s_2 = \sum_2^a t(j)$, \dots, $s_a = t(a)$. In particular, $P(r,a,m-a) = p(r,a,m-a)$.

 The generic word  can also be written in the form $g^{t(0)}h g^{t(1)}\dots g^{t(m-a-1)}hg^{t(m-a)}$ where $t(j)$ are nonnegative integers adding to $a$.

  The exponent of $z$ will be $t(1) + 2t(2) + 3 t(3) + (m-a)t(m-a)$ (obtained by starting from the left). Hence the multiplicity of $z^r g^a h^{m-a}$ in $(g+h)^m$ is the number of $m-a+1$ tuples $(t(0), t(1), \dots, t(m-a))$ of nonnegative integers \st $r = \sum_1^{m-a} jt(j) $. This yields $p(r,a,m-a) = p(r,m-a,a)$. It is routine to verify another symmetry condition, $P(a(m-a)-r,a,m-a)) = P(r,a,m-a)$, so the corresponding result applies for $p$. In addition, it is well known that for $a$ fixed, the sequence  $\(P(r,a,m-a)\)_{r=0}^{a(m-a)}$ is unimodal. It follows that the maximum is $P(\flo{a(m-a)/2},a,m-a)$. 
 The following is elementary, but very useful. 

 \Lem Lemma \onetwo. Suppose that $0 \leq a \leq m$. 
If $1 \leq r \leq a(m-a)-1$, then $p(r\pm 1,a,m-a)\leq 2 p(r,a,m-a)$. 

 \Pf (i) Consider first the case of $p(r- 1,a,m-a)\leq 2 p(r,a,m-a)$. If $a= 0$ or $m$,  then the left side is $0$, so we are reduced to $1 \leq a \leq m-1$.  If $r \geq a(m-a)/2$, then $p(r+1,a,m-a) \leq p(r,a,m-a)$, since the sequence is unimodal and symmetric.  Hence we may assume that $r < a(m-a)/2$. We define two (set) maps, $\Arrow f_1; \PP(r,a,m-a). \PP(r+1,a,m-a)$ and $\Arrow f_2; D . \PP(r+1,a,m-a)$, where $D \subset \PP(r,a,m-a)$ and $\PP(r+1,a,m-a)$ is the union of the images of $f_1$ and $f_2$. 

 To define $f_1$, begin with an $a$-tuple $s_1 \geq s_2 \geq \dots \geq s_a$ \st $r = \sum s_j$ and $m-a \geq s_1$. Find the smallest $j$ \st $s_j \neq m-a$ ($j$ exists  since $r < a(m-a)/2$). Map $(s_i) \mapsto f_1((s_i))$ by replacing $s_j$ by $s_{j+1}$; the resulting tuple belongs to $\PP(r+1,a,m-a)$. 

 Now define, for $1 \leq t \leq a-1$, 
$$
 D = \Set{(s_i) \in \PP(r,a,m-a)}{m-a > s_1 = s_2 = \dots = s_t > s_{t+1}}, 
 $$
 and set $D = \cup D_t$. Obviously, this is a disjoint union.  Define $\Arrow f_2; D_t.\PP(r+1,a,m-a)$ by adding $1$ to $s_{t+1}$. This defines $\Arrow f_2; D. \PP(r+1,a,m-a)$. 

 Now we show that the union of the images of $f_1$ and $f_2$ is all of $\PP(r+1,a,m-a)$. The image of $f_1$ consists of the sequences $(s_i) \in \PP(a,r,m-a)$ \st either $s_1 > s_2$ or $m-a = s_1 = s_2 = \dots = s_{j-1} > s_{j+1}$ for some $j$. Set $Y $ to be the set of sequences in $\PP(r+1,a,m-a)$ \st $m-a > s_1 = s_2$. Then $\PP(r+1,a,m-a)\setminus \Im f_1 \subset Y$. Now pick $(s_j) \in Y$. If $s_1 = s_2 = \dots = s_a$, then $s_a > 0$, so $(s_1, \dots, s_{a-1},s_a-1) \in D_{a-1}$, and it is mapped to $(s_j)$ by $f_2$. 

 Otherwise, there exists $k < a$ \st $s_1 = s_2 =\cdots s_{k} > s_{k+1}$. Then $(s_1, s_2, \dots, s_{k-1}, s_{k}-1, s_{k+1}, \dots) \in D_{k-1}$, and its image under $f_2$ is $(s_j)$.  
 
 For the second part, apply the involution $r \mapsto a(m-a) -r$ and the identity $p(r,a,m-a) = p(a(m-a)-r,a,m-a)$ to the first part. \qed

\Lem Lemma \onethr.  For $0\leq a \leq m$, we have $S(a,m-a,m) =\brcs{0,1, \cdots, a(m-a)-1, a(m-a)} $.

\Pf Set $b = m-a$. Pick  $u = z^r g^a h^b \in S^m$. Since the defect is zero, all of its representations from (*) must have all the exponents positive (except possibly $\gamma(1)$ and $\epsilon(k)$ which could be zero). This means that $u$ appears (with nonzero coefficient) in $(g+ h)^m$, so that $p(r,a,b) \neq 0$. Conversely, if $p(r,a,b) \neq 0$, then $v$ appears in $(g+h)^m$, and obviously this implies $v \in S^m$. Since $p(r,a,b)\neq 0$ for $a,b,r \geq 0$ iff $r \leq ab$, we are done.
\qed

\Lem Lemma \onefou. If $S(a,b,m)$ is nonempty, then it  is an interval of integers containing zero. 

\Pf By Lemma \oneone(a--e), we may assume that  $a,b \geq 0$ and $m -a-b $ is a nonnegative even integer, 
Now for the rest of the cases. We assume inductively that $S(a,b,m)$ is an interval, and must prove that $S(a,b,m+2)$ is as well. That $S(a,b,a+b)$ (zero defect) is an interval follows from Lemma \onethr. 

If $a+b < m+2$, from the defect condition, we must have $a+b \leq m$. We will show that $S(a,b,m+2)$ is a union of intervals centred on elements of $S(a,b,m)$, together with the defect zero terms, $S(a,m+2-a,m+2)$---which we already know to be an interval. Since all these intervals contain $0$  (by (c)), it will follow that $S(a,b,m+2)$ will be an interval if $S(a,b,m)$ is; thus induction on $m$ applies. 

For $r \in S(a,b,m)$, pick a representation as in (*) for $u = z^r g^a h^b$. View the form as a string of $g$s, $h$s, $g^{-1}$s, and $h^{-1}$s ($1$s are not needed; their presence is reflected only in the difference $m - \sum_i (|\gamma(i)| + |\epsilon(i)|)$, an even integer; and $g g^{-1}$ cannot appear  consecutively, and neither can $h h^{-1}$). If $w$ is a substring, that is, $u$ is represented by the string $U = U_1 w U_2$, we can insert $g$ just before the first letter of $w$ and $g^{-1}$ after the last letter of $w$, creating $Ugwg^{-1}U_2 \in S^{m+2}$. Then   $gwg^{-1} = z^k w$, where $k$ is the (net) sum of all the exponents of $h$ that appear in $w$, and thus $Ugwg^{-1}U_2 = z^k u$. In particular, $r+k \in S(a,b,m+2)$. Similarly, $r-k \in S(a,b,m+ 2)$ (from $Ug^{-1}wg U_2$). 

Let $w$ be a substring for which $gwg^{-1} = z^k w$; temporarily, suppose that $k > 0$. Then we can find a sequence of initial substrings of $w$, $w_1 $, $w_2$, \dots, $w_{k-1}$ \st $gw_j g^{-1} = z^j w_j$. Hence all of $r, r+1, \dots, r +k \in S(a,b,m+2)$; using $g^{-1}w g^{-1}$, we obtain the interval $\brcs{r-k, r- k+1, \dots, r+k} \subset S(a,b,m+2)$. If instead $k < 0$, we can do the same, and obtain the same conclusion. 

Of course, we can also do this with $g$ replaced by $h$ (and by $h^{-1}$). Take the union of all the sets obtained by varying the group element $u$ that represent $r$,  the forms (*) that represent $u$, all substrings of each one, and conjugation of the substrings by $g^{\pm1}, h^{\pm1}$. The outcome is a union of intervals, each one centred at $r$; consequently the union is itself an interval. Call it $I(r)$.

Now let us make the inductive assumption that $S(a,b,m)$ is an interval, so that we obtain $\cup_{r \in S(a,b,m)} I(r)$ is a union of intervals for which it is easy to check that it is itself an interval. Call this monster interval $I \subset S(a,b,m+2)$. We claim that this covers all of $S(a,b,m+2)$. 

To see this, note that since we have assumed that $S(a,b,m)$ is not empty, the defect of $(a,b,m+2)$ is at least two. Now consider a group element $v = z^r g^a h^b$ in $S^{m+2}$. It has a representation of the form (*). Since we have assumed $a,b \geq 0$, if none of the $\epsilon(i)$ or $\gamma(i)$ are negative, then $\sum \epsilon(i) = b$ and $\sum \gamma(i) = a$, so that $v \in S(a,b,m)$, which of course is contained in $I$. If at least one of the exponents is negative, then we obtain a substring of one of the forms $g^{-1}w g$, $gwg^{-1}$, $h^{-1}w h$, or $hwh^{-1}$. Deleting the conjugating symbol yields an element of $S(a,b,m)$, and thus $S(a,b,m+2) \subset I$. \qed

 For $a,b,m$ nonnegative integers \st $m - a- b $ is an nonnegative even integer, define $\partial = (m-a-b)/2$ and  the  following, 
$$
F(a,b,m ) = \cases (b+ 
\partial)a & \text{if $\partial + b  \leq a$}\\ 
 (\partial + a)b & \text{if $\partial + a \leq b$}\\ \Ud \frac{
\partial + a + b}2.& \text{if $\partial + b  \leq a$ and $\partial + a \leq b$.}\\ \endcases 
$$
It is routine to check that this is well-defined. Our immediate goal is to prove the following. Its proof is unfortunately extremely tedious; a less cumbersome one would be desirable. 

\Lem Theorem \onefiv. Suppose that $a,b,m$ are nonnegative integers \st $m- a-b$ is a nonnegative even integer.  Then $\max S(a,b,m) = F(a,b,m)$. 

\Rmk From this, we obtain formulas for $\max S(a,b,m)$ and $\min S(a,b,m)$ for all $(a,b) \in \Z^2$ (\st $m - |a| - |b|$ is a nonnegative integer).

\Lem Lemma \onesix. Suppose the hypotheses of Theorem \onefiv\ are satisfied. Then $\max S(a,b,m) \geq F(a,b,m)$.

\Pf We exhibit elements of $S^m$ whose $r$ value is $F(a,b,m)$. 

Suppose that $\partial + b \leq a$. Then $u =  h^{b+\partial}g^a h^{-\partial}$ belongs to $S^m$ (since $b + \partial + a + \partial = m$), and equals $z^{a(b + \partial)}g^a h^b$. Thus $a(b + \partial) = F(a,b,m)$ belongs to $S(a,b,m)$, hence is less than or equal to $\max S(a,b,m)$.  Similarly, if $\partial + a \leq b$, set $u = g^{-\partial}h^b g^{a + \partial} = z^{b(a+\partial)}g^a h^b$. 

If $\partial  + a \geq b$ and $\partial + b \geq a$, set $u = g^{a - s} h^{b+t} g^s h^{-t}$  where $s = \flo{(\partial + a+ b)/2}$ and $t = \cei{(\partial +a-b  )/2}$. Thus $s,t \geq 0$ and $s+t = \partial + a = (m+a-b)/2$. Moreover, $s\geq a$ follows from $\partial + b \geq a$. Since  $ s-a + b+t  + s +t = 2(s+t) + b-a = m+a-b +b-a = m$, it follows that $u \in S^m$. Moreover, its $r$  value is $s\cdot (b+t) = \Ud (\partial + a+b)/2.$
\qed

\Pf (of  Theorem \onefiv) To prove the reverse inequality, we proceed by induction  on $(\partial,m)$ in that order (that is, we take the lexicographic ordering on a subset of $\Z_{+}^2$). The first case, that of $\partial = 0$,  has been done in Lemma \onethr.

Suppose $u \in S^m$ represents $\max S(a,b,m)$, that is, $u = z^r g^a h^b$ where $r = \max S(a,b,m)$. Write $u $ in the form $(*)$. We claim that we can assume that $\epsilon (1) \geq 0$ (that is, in the form $(*)$, the leading letter cannot be $h^{-1}$. If $u = h^{-m}$, then $S(a,b,m)=\brcs{ 0}$ and we can replace $u$ by $1$. Otherwise, we can write $u = h^{-s}w$ with $s> 0$ and where $w \in S^{m-s}$ and its leading letter is neither $h$ nor $h^{-1}$. Then form $u' = wh^{s}$; this clearly belongs to $S^m$, and $u' = h^{s}u h^{-s} = h^s z^r g^a h^b h^{-s} = z^{r + ab} g^a h^s = z^{ab}u$. Since $a,b \geq 0$, this has $r$-value at least as large. Thus  $ab = 0$, and we may replace $u$ by $u'$. 

Now we assume $\max S(a,b,n)= F(a,b,n)$ for all $n$ \st $n-a-b$ is a nonnegative even integer \st either $(n-a-b)/2 < \partial$ or both $(n-a-b)/2 = \partial  $ and $n < m$ (where $m-a-b$ is even). It suffices  to show $\max S(a,b,m) \leq F(a,b,m)$. We have numerous cases, most of which are straightforward. At a few points, the $a$s might have to drop below zero, but we can use  (d) to avoid this (the induction hypothesis requires the first and second coordinates to be nonnegative.

The case that $m=0$ is trivial. So we can assume $m > 0$. We can also assume that $\max S(a,b,m) > 0$. 
Here $\partial  = (m-a-b)/2$ and we will use the notation  $\partial'$ (and if necessary, $a'$ and $b'$) when we are dealing with other triples. 

Select nonnegative $u$ representing $\max  S(a,b,m)$, such that its leading letter is not $h^{-1}$.

 There are four possibilities for the leading symbol in one of the forms for $u$, that is, $1$, $g$, $g^{-1}$, and $h$. However, we can move all the leading $1$s to the other end without changing $u$, and since $u \neq 1$ nor $h^{-m}$ (else $\max S(a,b,m) = 0$), we will eventually find a form whose initial letter is one of $g^{\pm1}$ or $h$.  Thus we can write $u = gw$ or $u = g^{-1}w$, or $u = h^{-1}w$ arising from the form $(*)$, where $w \in S^{m-1}$.

\noindent {\bf Case 1:  $\partial + b \leq a$.} Thus $F(a,b,m) = a(b+\partial)$. 

{Case 1(a): $a,b > 0$}. 

\ \  \hskip 10pt {Case 1(a)(i): $u = gw$. } Thus $w = z^r g^{a-1}h^b$, so that $r \in S(a-1,b,m-1)$, and for $w$, we have $\partial' = \partial$, $m' = m-1$, $a' = a-1$, and $b' = b$. By the induction hypothesis, $\max S(a-1,b, m-1) = F(a-1,b,m-1)$. From $\partial + b \leq a$, we deduce that either $\partial' + b' \leq a'$, so that $\max S(a,b,m) = r \leq F(a-1,b,m-1) = a'(b' + \partial') = (a-1)(b + \partial-1) < F(a,b,m)$, which is impossible, or $\partial' + b' = a'+1$, that is, $\partial + b = a$, so that $F(a,b,m) = a^2$. On the other hand, $F(a-1,b,m-1) = \Ud (\partial' + a' + b)/2. = \Ud a- 1/2. < a^2$, again impossible.  

  \ \  \hskip 15pt {Case 1(a)(ii): $u = g^{-1}w$. }  Then $w = z^r g^{a+1} h^b$, and thus $r \in S(a+1,b,m-1)$. Hence $\partial' = \partial - 1$, so again the induction hypothesis applies (attempting to use  induction on $a$ would obviously run into trouble). Thus $r \leq \max S(a+1,b,m-1)$. Obviously, $\partial' + b' \leq a'$, so that $r \leq F(a+1,b,m-1) = (a+1) (b + \partial-1) = a(b + \partial)  + b + \partial - a -1 + a +1 = F(a,b,m) - (a - b -\partial) \leq F(a,b,m)$. 

  \ \  \hskip 15pt {Case 1(a)(iii): $u = hw$. } Here $w = h^{-1}u = h^{-1}z^r g^a h^b = z^{r-a} g^a h^{b-1}$. Thus $r-a \leq \max S(a,b-1,m-1)$. We have $\partial' = \partial$ and $m' = m-1$, so that the induction hypothesis applies, and thus $r -a \leq F(a,b-1,m-1)$. From $\partial + b \leq a$, we obviously have $\partial' + b' \leq a'$, and thus $r-a \leq a(b-1 + \partial)$, so that $r \leq a(b + \partial) = F(a,b,m)$. 

{Case 1(b)}: $b = 0$ and $a > 0$. Here $\partial = (m-a)/2$. Then $u = z^r g^a$.

\ \  \hskip 5pt {Case 1(b)(i): $u = gw$. } Thus $w = z^r g^{a-1}$, and $r \in S(a-1,0,m-1)$, for which we have $\partial' = \partial$ and $m' = m-1$. The same argument as in Case 1(a)(ii) now applies. 

  \ \  \hskip 15pt {Case 1(b)(ii): $u = g^{-1}w$. }  Then $w = z^r g^{a+1}$, and thus $r \in S(a+1,0,m-1)$. Hence $\partial' = \partial - 1$, and now the same argument as in Case 1(a)(ii) applies. 

  \ \  \hskip 15pt {Case 1(b)(iii): $u = hw$. } Here $w = h^{-1}u = h^{-1}z^r g^a h^{-1}  = z^{r-a} g^a h^{-1}$, and thus $r-a \in S(a,-1,m-1) = -S(a,1,m-1)$; hence $r \in S(1,a,m-1)$,  so $r \leq \max S(1,a,m-1)$. With $S(1,a,m-1)$, we have $\partial' = \partial -1 = (m-a-2)/2 $ (so the induction hypothesis applies, and thus $\max S(1,a,m-1) = F(1,a,m-1)$), $b' = a$, and $a' = 1$. Thus $\partial' + a' = \partial \leq a = b'$. Thus $ r \leq F(1,a,m-1) = b'(a' + \partial') = a\partial = F(a,0,m)$.  

{Case 1(c): $a = 0$.} Since $\partial + b \leq a$, this forces $\partial = b = 0 = a$, yielding $m = 0$, a contradiction.

\noindent {\bf Case 2:  $\partial + a \leq b$.} Thus $F(a,b,m) = b(a+\partial)$. We can obviously exclude $b = 0$ (as in the analogous situation, Case 1(c)).  

{Case 2(a): $a,b > 0$}. 

\ \  \hskip 15pt {Case 2(a)(i): $u = gw$. } Thus $w = z^r g^{a-1}h^b$, so that $r \in S(a-1,b,m-1)$, and for $w$, we have $\partial' = \partial$, $m' = m-1$, $a' = a-1$, and $b' = b$. By the induction hypothesis, $\max S(a-1,b, m-1) = F(a-1,b,m-1)$. From $\partial + a \leq b$, we deduce that  $\partial' + a' \leq b$, so that $\max S(a,b,m) = r \leq F(a-1,b,m-1) = b'(a' + \partial') = b(a + \partial-2) <  F(a,b,m)$, which is impossible.  

  \ \  \hskip 15pt {Case 2(a)(ii): $u = g^{-1}w$. }  Then $w = z^r g^{a+1} h^b$, and thus $r \in S(a+1,b,m-1)$. Hence $\partial' = \partial - 1$,  and the induction hypothesis applies . Thus $r \leq \max S(a+1,b,m-1)$. Now $\partial' + a' = \partial + a\leq b = b'$. Hence  $r \leq F(a+1,b,m-1) = b (a + \partial)  = F(a,b,m)$. 

  \ \  \hskip 15pt {Case 2(a)(iii): $u = hw$. } Here $w = h^{-1}u = h^{-1}z^r g^a h^b = z^{r-a} g^a h^{b-1}$. Thus $r-a \leq \max S(a,b-1,m-1)$. We have $\partial' = \partial$ and $m' = m-1$, so that the induction hypothesis applies, and thus $r -a \leq F(a,b-1,m-1)$. From $\partial + a \leq b$, we have either have $\partial' + a' \leq b'$ or $\partial' + a' = b'+1$. In the former case, $F(a,b-1,m-1) = b'(\partial' + a') = (b-1) (a+ \partial)$, so that $r \leq a + (b-1)(a+\partial ) = b(a+\partial) - \partial \leq F(a,b,m)$. 

If on the other hand, $\partial' + a' = b'+1$, then $b = \partial + a$, and thus $F(a,b,m) = b^2$, while $F(a,b-1,m-1) = \Ud (\partial' + a' + b')/2. = \Ud b - 1/2. = b(b-1)$. Hence $F(a,b-1,m-1) + a = b^2 + a-b \leq F(a,b,m)$.

{Case 2(b): $a= 0$}. 

\ \  \hskip 15pt {Case 2(b)(i): $u = gw$. } Thus $w = z^r g^{-1}h^b$, so that $r \in S(-1,b,m-1) = -S(1,b,m-1)$; thus $r + b \in  S(b,1,m-1)$. For $S(b,1,m-1)$, we have $\partial' = \partial-1$, $m' = m-1$, $a' = b$, and $b' = 1$. By the induction hypothesis, $\max S(b,1, m-1) = F(b,1,m-1)$. Now $\partial' + b' = \partial-1 + 1 = \partial \leq b = a'$. Hence $F(b,1,m-1) = a'(b'+\partial') = b\partial = F(0,b,m)$, and  $r \leq b(\partial-1)$.  

  \ \  \hskip 15pt {Case 2(b)(ii): $u = g^{-1}w$. }  Then $w = z^r g h^b$, and thus $r \in S(1,b,m-1)$. Hence $\partial' = \partial - 1$, and the same argument as in Case 2(a)(ii) applies. 

  \ \  \hskip 15pt {Case 2(b)(iii): $u = hw$. } Here $w = h^{-1}z^r  h^b = z^{r}  h^{b-1}$. Thus $r \leq \max S(0,b-1,m-1)$, and the same arguments as in Case 2(a)(iii) apply (there are two subcases therein). 

\noindent {\bf Case 3:  $\partial + a \geq b$ and $\partial + b \geq a$.} Thus $F(a,b,m) = \Ud (\partial + a  + b)/2.$. 

{Case 3(a): $a,b > 0$}. 

\ \  \hskip 15pt {Case 3(a)(i): $u = gw$. } Thus $w = z^r g^{a-1}h^b$, so that $r \in S(a-1,b,m-1)$, and for the latter, we have $\partial' = \partial$, $m' = m-1$, $a' = a-1$, and $b' = b$. Then $\partial' + b' = \partial + b - 1 \geq a-1 = a' $, while $\partial' + a' = \partial + a -2$. Hence we have three possibilities: $\partial' + b' \geq  a' $ and ($\alpha$) $\partial' + a' \geq b'$; and $(\beta$) $\partial' + a' = b'-1$; and $(\gamma)$ $\partial' + a' = b'-2$ (and in this case, $b \geq 2$).

\ \ \hskip 15pt ($\alpha$) Here $F(a-1,b,m-1) = \Ud (\partial + a-1 + b)/2. \leq \Ud(\partial + a + b)/2.$

 \ \ \hskip 15pt ($\beta$) This time, $F(a-1,b,m-1) = b'(\partial' + a') = b(b-1 )$, while $F(a,b,m) = b^2$. 

 \ \ \hskip 15pt ($\gamma$) This time, $F(a-1,b,m-1) = b(b-2)$, which is less than $F(a,b,m) = \Ud (2b-1)/2. = b(b-1)$. 

\ \  \hskip 15pt {Case 3(a)(ii): $u = g^{-1}w$.}  Then $w = z^r g^{a+1} h^b$, and thus $r \in S(a+1,b,m-1)$. Hence $\partial' = \partial - 1$,  and the induction hypothesis applies, so $r \leq \max S(a+1,b,m-1)$. Now $\partial' + a' = \partial + a \geq b = b'$, while $\partial' + b = \partial + b-1$; so there are again three possibilities: $\partial' + a' \geq b$ and ($\alpha$) $\partial' + b' \geq a'$; and ($\beta$) $\partial' + b' = a'-1 = a$; and ($\gamma$) $\partial' + b = a-1$.

\ \ \hskip 15pt ($\alpha$) Here $F(a+1,b,m-1) = \Ud (\partial + a + b)/2. = F(a,b,m)$.

 \ \ \hskip 15pt ($\beta$) This time, $F(a+1,b,m-1) = a'(\partial' + b') = (a+1)a$, while $F(a,b,m) = \Ud (\partial + a + b)/2 . = \Ud (2a+1)/2 .= a(a+1)$. 

 \ \ \hskip 15pt ($\gamma$) This time, $F(a+1,b,m-1) = a^2-1$, which is less than $F(a,b,m) = \Ud (\partial + a + b)/2 . = a^2$. 

\ \  \hskip 15pt {Case 3(a)(iii): $u = hw$.}  Here $w = h^{-1}u = h^{-1}z^r g^a h^b = z^{r-a} g^a h^{b-1}$. Thus $r-a \leq \max S(a,b-1,m-1)$. We have $\partial' = \partial$ and $m' = m-1$, so that the induction hypothesis applies, and thus $r -a \leq F(a,b-1,m-1)$. This time, $\partial' + a' = \partial + a \geq b > b'$ and $\partial' + b' = \partial + b-1$. Here there are only two subcases, $\partial' + a' \geq b'$ and ($\alpha$) $\partial' + b' \geq a'$; and ($\beta$) $\partial + b = a$.

\ \ \hskip 15pt ($\alpha$) Here $F(a,b-1,m-1) = \Ud (\partial + a + b-1)/2. \leq  \Ud (\partial + a + b)/2. = F(a,b,m)$.

 \ \ \hskip 15pt ($\beta$) This time, $F(a,b-1,m-1) = a'(\partial' + b') = a(\partial + b-1) = a(a-1)$, while $F(a,b,m) = \Ud (\partial + a + b)/2 . = a^2$.\par

{Case 3(b): $a= 0, b> 0$}. 

 \ \  \hskip 15pt {Case 3(b)(i): $u = gw$.} We have $w = z^r g^{-1} h^b$, so $r \in S(-1,b,m-1) = - S(1,b,m-1) = S(b,1,m-1)- b$, and $\partial' = \partial-1$, $a' = b$, $b' = 1$. We have $r \leq F(b,1,m-1)$. Now $\partial' + b' = \partial  \geq a' = b$  and $\partial' + a' = \partial  -1 = (m-b)/2$; either $(\alpha)$ $(m-b)/2 \geq 1$ or $(m-b)/2 = 0$. The latter is impossible, since $b \leq m-2$.This leaves  $F(b,1,m-1) = \Ud (\partial + b)/2.  = F(0,b,m)$.\par

  \ \  \hskip 15pt {Case 3(b)(ii): $u = g^{-1}w$. }  Then $w = z^r g h^b$, and thus $r \in S(1,b,m-1)$. Hence $\partial' = \partial - 1$, and the same argument as in Case 3(a)(ii) applies. 

 \ \  \hskip 15pt {Case 3(b)(iii): $u = hw$. } Here $w = h^{-1}z^r  h^b = z^{r}  h^{b-1}$. Thus $r \leq \max S(0,b-1,m-1)$, and the same argument as in Case 3(a)(iii) works.

{Case 3(c): $a > 0$, $b=0$. }

 \ \  \hskip 15pt {Case 3(c)(i): $u = gw$.}  Now $w = z^r g^{a-1}$, and the method of Case 3(a)(i) works.

 \ \  \hskip 15pt {Case 3(c)(ii): $u = g^{-1}w$.}  Then $w = z^r g^{a+1}$, and the method of Case 3(a)(ii) works. 

 \ \  \hskip 13pt {Case 3(c)(iii): $u = hw$.} Here $w = h^{-1}u = h^{-1}z^r g^a = z^{r-a} g^a h^{-1}$. Thus $r-a \leq \max S(a,-1,m-1) = - S(a,1,m-1) = S(1,a,m-1) -a$. Then $\partial' + b' \geq a' = 1$, while $\partial' + a' = \partial-1 + 1 = \partial \geq b' = a$. Thus $r \leq \Ud {(\partial + a)/2}. -a$, while $F(a,0,m) = \Ud (\partial + a)/2.$. 

{Case 3(d): $a = 0$, $b=0$. }  In this case, $F(0,0,m) = \Ud m/4. $). 

 \ \  \hskip 15pt {Case 3(d)(i): $u = gw$.} Here $w =z^r g^{-1}$, so $r \in S(-1,0,m-1) = -S(1,0,m-1) = S(0,1,m-1)$, so $r \leq F(0,1,m-1) = \Ud m/4. = F(0,0,m)$.  

  \ \  \hskip 15pt {Case 3(d)(ii): $u = g^{-1}w$.} This time, $w = z^r g$, and $r \leq \max S(1,0,m-1)$, whence $r \leq F(1,0,m-1) = \Ud {m/4}. = F(0,0,m)$. 

 \ \  \hskip 19pt {Case 3(d)(iii): $u = hw$.} Now $w = z^r h^{-1}$, and thus $r \in S(0,-1,m-1) = -S(0,1,m-1) = S(1,0,m-1)$; so $r \leq \Ud m/4 . = F(0,0,m)$. 
\qed

\comment
Finally, we notice that for nonnegative $a,b$, $F(a,b,m) = F(b,a,m)$. From $S(a,b,m) = ab - S(b,a,m)$, so $\inf S(a,b,m) = -(\max S(b,a,m) - ab) $, and thus $|\inf S(a,b,m)| = F(a,b,m ) - ab$. Hence $S(a,b,m)$ is the interval with least element $-(F(a,b,m) - ab)$ and maximum element $F(a,b,m)$. 
\endcomment

\Lem Theorem  \onesev. For $a,b\geq 0$ and $m-a-b$ an even nonnegative integer, we have the following.
\item{(a)} $S(a,b,m) = S(b,a,m)$; 
\item{(b)} $S(a,b,m)$ is the interval of integers with minimum $ab- F(a,b,m)$ and maximum $F(a,b,m)$. 

\Rmk For arbitrary integers $(a,b)$, the requirements are that $m - |a| - |b| $ is an even nonnegative integer. Then we can apply Lemma \oneone(a) to describe $S(a,b,m)$ exactly. 

\Rmk There probably is a way of proving $S(a,b,m) =S(b,a,m)$ directly, without going through the incredibly tedious computations above.

\Rmk This result has an immediate interpretation in terms of the Bratteli diagram for $R_f/\SS R_f$. At the $m$th level, restrict to the group elements $g = z^r g^a h^b$ (in $\Gamma_m$, so $m = |a|+|b|$), fix $a$; the possible choices for $r$ consist of the interval $0 \leq ab = a(m-a)$. As we let $a$ vary (from $0$ to $m$), the plot of possible $r$\,s is a discretized parabola. There are three other parts, corresponding to the signs of permitted $a$ and $b$, each with the same shape. 

\Pf  From the definition of $F$, we see that $F(a,b,m) = F(b,a,m)$. Hence $\max S(a,b,m) = \max S(b,a,m)$. Since $S(a,b,m) = -S(b,a,m) + ab$, we have that $\min S(a,b,m) = ab - \max S(b,a,m) = ab - F(b,a,m) = ab - F(a,b,m)$; this yields (b). But we also have $\min S(b,a,m) = ab - F(a,b,m) = \min S(b,a,m)$, and thus the endpoints of $S(a,b,m)$ and $S(b,a,m)$ are the same.
\qed

\Lem Lemma \oneeig. Assuming $a,b \geq 0$ and $m-a-b$ is an even nonnegative integer,
$$
F(a,b,m+2) - F(a,b,m)= \cases 
a & \text{if $\partial + b \leq a$}\\ 
b & \text{if $\partial + a \leq b$}\\
\cei{\frac{\partial + a + b}2} & \text{if $\partial + b \geq a$ and $\partial +a  \geq b$.} \endcases
$$
Moreover, $F(a,b,m+2) - F(a,b,m) \geq m/4$.

\Pf For the triple $(a,b,m+2)$, we have $\partial' = \partial +1$. Hence if $\partial + a \geq b$ and $\partial + b \geq a$, then the same inequalities hold for $(a,b, m+1)$, and thus the difference is $\Ud {(\partial + 1 + a +b)/2}. - \Ud (\partial + a + b)/2. = \cei{(\partial + a + b)/2}$. If $\partial +1 + b \leq a$, then the difference is $a(\partial + 1 + b) - a(\partial + b) = a$. If instead, $\partial + 1 + b = a + 1$ then $a = \partial + b$, and then $F(a,b,m+2) - F(a,b,m) = \Ud (\partial + 1 + b + a)/2. - a(b + \partial) = \Ud a + 1/2. - a^2 = a$. Finally, we can reverse the roles of $a$ and $b$ to obtain the middle equality. 

The inequality $\partial + b \leq a$ is equivalent to $(m-a-b)/2 + b \leq a$, that is, $3a \geq m+b$, so $a \geq m/3$. Similarly, $\partial + a \leq b $ is equivalent to $3 b \geq m+a$. Finally, $\partial + a + b = (m+a+b)/2$, so $(\partial + a + b)/2 \geq m/4$.
\qed 

Amusingly, this implies that if $m-a-b$ is divisible by four (it is already even), then $F(a,b,m+2) - F(a,b,m) = F(a,b,m+4) - F(a,b,m+2)$.

\Lem Corollary \onenin. Let $G$ be the discrete Heisenberg group, and let $S = \brcs{1, g^{\pm1} , h^{\pm 1}}$. Here $(r,a,b)$ run over $\Z^3$.
\item{(a)}  $$
\tilde l_S (z^r g^a h^b)  = \cases |a| + |b| & \text{if $|r| \leq |a|\cdot |b|$ and $\sign r = (\sign a)(\sign b)$ } \\ |a| + |b| + 2 & \text{else.}
\endcases
$$
In particular, $|a| + |b| \leq \tilde l_S(z^r g^a h^b) \leq |a| + |b| + 2$. 
\item{(b)} If $C$ is the conjugacy class of  $z^r g^a h^b$ in $G$, then $ \sup_{c \in C}\tilde l_S(c) \leq |a| + |b| + 2$; in particular, $\tilde l_S$ is bounded on each conjugacy class. 
\item{(c)} $G$ fails WC. 
\item{(d)}  If $w = z^r g^a h^b$ in $G$, then for all $k \geq |a| + |b| + 2$,   the endomorphism of $R_f$ (for any admissible $f$  with $\Supp f = S$) given by $\RR_{w,k}$ is locally order-bounded. 

\Rmk Part (d) leaves open the possibility that if $k \geq |a| + |b| + 2$, then $\RR_{u,k}$ is order-bounded. This is  true, and is proved in Corollary  \endsthr. 

\Pf (a) From factoring out the centre to obtain  $\Z^2$, it is obvious that  $\tilde l_S(z^r g^a h^b) \geq |a| + |b|$. Because $\tilde l_S$ is subadditive and $l_S (g^a h^b) = |a| + |b|$,   showing $\tilde l_S (z^r) \leq 2$ for all $r$ entails the last statement, so it becomes a matter of excluding $|a| + |b| + 1$.  

Pick $m > 4|r| +1$; we show that $S^m z^r \subset S^{m+2}$. For $u = z^s g^c h^d$ in $T:= (S^m\setminus S^{m-1}) \cup (S^{m-2}\setminus S^{m-3}) \cup \dots$ corresponding to the spheres of radius of the same parity as $m$ (so that $m \geq |c| + |d|$ and the difference is an even integer), first assume that the signs of $c$ and $d$ are equal to each other. Then $cd - F(|c|,|d|,m) \leq s \leq F(|c|, |d|, m)$, and so $cd - (F(|c|,|d|,m) - r) \leq s + r \leq F(|c|,|d|,m) + r$. If $r \geq 0$, $F(|c|,|d|,m+2) \geq F(|c|,|d|,m) + r$ (since $r < m/4$) and thus $cd - F(|c|,|d|,m+2)  \leq r + s \leq F(|c|,|d|,m + 2) $, and thus $u z^r \in S^{m+2}$. If $r \leq 0$, we work on the other side and obtain $uz^r \in S^{m+2}$. 

If $c$ and $d$ have opposite signs, then $-F(|c|,|d|,m)  \leq s \leq F(|c|,|d|,m) + cd$  (here $cd$ is negative), and the same arguments yield $u z^r \in S^{m+2}$ for all $u \in T$. For the remaining $T' := (S^{m-1}\setminus S^{m-2}) \cup (S^{m-3}\setminus S^{m-4}) \cup \dots$, we note that since $m-1> 4r$, we obtain for $u \in T'$, that $uz^r \in S^{m+1}\subset S^{m+2}$. Hence $S^m z^r \subset S^{m+2}$. 

Applying the action of $D_4$ described earlier, we reduce to the case that $a,b \geq 0$. Set $w = z^r g^a h^b$. If $0 \leq r \leq ab$, then $l_S(w) \leq   a+b$, by Lemma \onethr, and since $l_S(w) \geq \tilde l_S (w) \geq a+b$ in any case, we have $\tilde l_S (w) = l_S(w) = a+b$. 

If $r \notin [0,ab]$, then $w \notin S^{a+b}$ (again by Lemma \onethr), so by (a), it suffices to show that $\tilde l_S (w) \notin \brcs{a+b, a+b+1}$. It is easy to check that $\tilde l_S (w) \equiv a+b \pmod 2$; thus we reduce to showing $\tilde l_S (w) > a+b $. Otherwise, $\tilde l_S(w) = a+b$, so there would exist a positive integer $k$  \st $S^k w \subset S^{k+ a + b}$. In particular, $g^k w, h^k w \in S^{k+a+b}$. The former entails $z^r g^{a+k}h^{b} \in S^{k+a+b}$, which forces $r \geq 0$ by Lemma \onesev; the latter entails $z^{r+ak}g^a h^{b+k} \in S^{a+b + k}$, which  forces $r+ak \leq a(b+k)$, that is, $r \leq ab$.

\noindent (b) The conjugacy class of $z^r g^a h^b$ is  $\Set{z^{r+s} g^a h^b}{s \in a\Z + b\Z}$, so all the elements therein have $\tilde l_s$ value at most $|a| + |b| + 2$.  

\noindent (c) An immediate consequence of (a) or (b).

\noindent (d) The largest $\tilde l_S$ value of an element in the conjugacy class of $z^r g^a h^b$ is   $|a| + |b| +2$, and now Lemma \Bthr\  applies; this  yields local boundedness. \qed

We observe that WC fails badly---$\tilde l_S$ is bounded on every conjugacy class, and the conjugacy classes of noncentral elements are  infinite. Moreover, $\tilde l_S^{-1} (\leq 2)$ contains the centre, and a few other elements (explicitly, $\brcs{g^{-1}, h^{-1}, g^{\pm 1}, h^{\pm 1}, (gh^{-1})^{\pm 1}, (gh)^{\pm1}, (g^{-1}h)^{-1}}$; not included are $zg$, $zg^2$ or similar).

\SecT 13\ Pure traces and the Heisenberg group

Now with $G = H_3$ and $f = 1 + g + g^{-1} + h + h^{-1}$, we determine all the traces on $R_f$. First, we have the obvious ones, arising from $G \to G/[G,G] \iso \Z^2$; the faithful ones are all of this form, and their limit points yield the rest of the multiplicative traces. However, there are a lot more; explicitly, eight families of discrete traces (a trace is {\it discrete\/}) if its range is a discrete subgroup of $\R$, that is, cyclic) which overlap only trivially with the previous ones.  Proving that these consititute all of the pure traces is rather long-winded. Then we discuss the topology on the discrete traces, or at least convergent sequences of them; in this case, we rely on results of Szerkeres [Sz1, Sz2] (exposed in [Ca]) for asymptotic estimates of multiplicities in the expansion of $(g+h)^n$. 

\noindent{\it Faithful pure traces on $R_f$}. If $\tau$ is a faithful pure trace of $R_f$, then it extends to a pure trace (or anyway, part of a pure ray of traces) of $A_f$, simply by $[g,k] \mapsto \tau([e_g, \tilde l_S (g)])\tau([1,1])^{-(k-\tilde l_S(g)}$. In particular, $g \mapsto \tau([, 0])$ (where $\tau$ also denotes the extension to $A_f$) is an extremal harmonic function on the random walk. By the well-known result of Margulis, since $G$ is nilpotent, this must be a character of $G$. Hence it factors through the quotient map $G \to G/[G,G] \iso \Z^2$. Say $g \mapsto x$ and $h \mapsto y$ under this map. 

 The characters of $\Z^2$ are given by $(x,y) \mapsto (r,s)$ where $(r,s)$ is a pair of nonzero real numbers. The corresponding character is then $\Arrow \chi_{r,s};G.\R$ determined by $\chi_{r,s}(g) = r$ and $\chi_{r,s} (h) = s$ (and $z \mapsto 1$ automatically). Since traces are positive, the images of $g$ and $h$ must be positive. Hence the normalized pure faithful traces on $R_f$ are determined  by (and the same formula allows them to be extended to $A_f$)  
$$
 \tau_{r,s} \([e_w, k]\) = \frac{\chi_{r,s}(w)}{(1 + r +s +r^{-1} + s^{-1})^k }.
$$
We recognize the denominator as $(\chi_{r,s}(f))^k$ (where we have extended the character to be defined on the group ring). Knowing the pure faithful traces, we can almost finish the result in Corollary \onenin(d). 

\Lem Corollary \endsone. Let $G$ be the discrete Heisenberg group, and 
let $S = \brcs{1, g^{\pm1} , h^{\pm 1}}$.  If $w = z^r g^a h^b$ in $G$, then for all $k \geq |a| + |b| + 3$,   the endomorphism of $R_f$ (for any admissible $f$  with $\Supp f = S$) given by $\RR_{w,k}$ is  order-bounded.

\Pf Assume that $k \geq |a| + |b| + 3$. 
Then $[z^j,3] = \SS [z^j,2]\in \SS R_f$ for all integers $j$. Let $u = z^s g^c h^d$, and let $n \geq \tilde l_S (u)$. Then 
$$\eqalign{\RR_{w,k} \([u,n]\) & = [uw,k + n] \cr
& = [z^{r+s}g^c h^d g^a h^b, k+n] = [z^{r+s + da}g^{c+a}h^{d+b}, k+n ]\cr 
& = [z^{r+s + ad-bc} g^a h^b g^c h^d, k+n]\cr
& \leq [z^{r+ ad-bc} f^{|a| + |b|} z^s g^c h^d, k+n]\cr 
& = [z^{r+ ad-bc}  u, k+n - |a| - |b|] \cr
&= \RR_{z^{r+ ad-bc} , k -|a| -|b|}\( [u,n]\).
}$$ 

Let $y \in AG^+$ be central, and suppose that $[y,m] \in R_f$. From Lemma \TRthr, $\LLL_{y,m} = \RR_{y, m}$  is a bounded endomorphism with norm $\alpha = \left\| \widehat {[y,m]}\right\|$,  and for all $\epsilon > 0$ $\RR_{y,m} \leq (\alpha +\epsilon) \I$ as endomorphisms of $R_f$ (recall that $I$ is the identity operator). 

Apply this to $\RR_{z^{j},t}$ where $j = r+ad-bc$ and $t =  k - |a| - |b| \geq 3$. Now $[z^j, t] \in \SS R_f$ (since $t \geq 3$), so that for all non-faithful pure traces $\tau$, we have $\tau([z^j, t]) = 0$. Hence $\alpha  $ (which ostensibly might depend on the choice of exponent of $z$) is just the supremum of the values at faithful pure traces. As we have seen, these are given by normalized characters indexed by $(A,B)\in (\R^2)^{++}$, $[p,n]\mapsto \pi(p)(A,B)/ P(A,B)^n $ for $[p,n] \in R_f$. Applying this to $p = z^j$, we have $\pi(z^j) = 1$, so that $\alpha = \sup_{(A,B) \in (R^2)^{++}} 1/(1 + A + A^{-1} + B + B^{-1})^t = 5^{-t}$. The important thing is that this is independent of the exponent of $z$. 

Since $t \geq 3$, we have 
$$
\RR_{z^{r+ ad-bc} , k -|a| -|b|}\( [w,n]\) \leq \frac 1{120}  [w,n]   
$$
($1$ in place of $\slfrac1{120}$ would have been good enough). Thus for all $w \in G$ \st $[w,n] \in R_f$, we have $\RR_{u,k} ([w,n]) \leq [w,n]/120$, and it is immediate that $\RR_{u,k}$ is a bounded endomorphism. 
\qed

The case that $3$ is replaced by $2$ is interesting, and leads to a number of
unexpected examples. The simplest case is that of $[z^r, 2]$ where $r$ is a
nonzero integer.

 \Lem Lemma \endstwo. Assume $r$ is a nonzero integer. 
\item{(a)} $[z^r, 2] + \SS R_f$ is a positive nonzero infinitesimal in $R_f/\SS
R_f$. 
 \item{(b)} For all integers $r'$, $[z^r - z^{r'},2] \in \(\Inf R_f\) \setminus \SS
R_f$. 
 \item{(c)} $\left\| \widehat{[z^r,2]}\right\| = \slfrac 1{25}$. 

 \Rmk We know that $[w, \tilde l_S (w)]$ cannot be an infinitesimal in $R_f$ for
any group element $w$; however (a) gives an example wherein it is an
infinitesimal modulo $\SS R_f$. Condition (c) implies that for all $\epsilon >
0$, $\RR_{z^r,2} = \LLL_{z^r,2} \leq (\slfrac 1{25} + \epsilon)\I $ as
endomorphisms of $A_f$ and $R_f$, independently of $r$; in particular,
$\RR_{z^r,2}$ are (uniformly) order-bounded. 

 From this, we will complete the earlier result and obtain the following. 

 \Lem Proposition \endsthr. Let $w = z^r g^a h^b$ be an element of $H_3$. 
\item{(a)} If $k \geq |a| + |b| + 2$, then $\RR_{w,k}$ is order bounded.
 \item{(b)} Suppose that $p \in AG^+$ and let $k$ be a positive integer. Then
$\RR_{p, k}$ is order-bounded if it is locally order-bounded. 

 \Rmk In (b), order-bounded merely as an endomorphism of $A_f$ easily entails
that $[p,k] \in R_f$, and  that it must be order-bounded as an endomorphism of
$R_f$ (and with the same norm), so the apparent ambiguity is resolved. 
 
\Pf (of Lemma  \endstwo(a)). Let $k$ be any integer exceeding $2|r|$ (we will
increase $k$ without bound later on).  We wish to show that there is a positive
real-valued  function $B(k) \to \infty$ \st $B(k)z^r f^k \leq f^{k+2} + h(k)$
(coordinatewise, that is, in terms of coefficients) where $h (k) \in AG$, and
all group elements $x$ in the support of $h(k)$ satisfy $\tilde l_S (x) \leq
k+1$. This will mean that $[h(k), k+2] \in \SS R_f$, so that $B(k)[z^r,2] \leq
[1,0]$ modulo $\SS R_f$. As $B(k) \to \infty$, it follows that $[z^r,2] + \SS
R_f$ is an infinitesimal in $R_f/\SS R_f$. (We already know that it is positive
and nonzero.) It will turn out that we can take $B(k) = k/2 - |r|$.
 
 We modify our notation to make the proof less incomprehensible. For a group
element $w$ appearing in $f^n$, let $m(w,n)$ denote its coefficient ($m$ is for
{\it multiplicity\/}). For suitable (but not all) group elements $w \in \supp
f^k$, we will show  that  $m(z^r w, k+2)\geq  (k/2 - r) m(w,k)$, without
actually determining the multiplicities (which can be done, but is excruciating).
 
 Pick $w = z^s g^a h^b \in \supp f^k$. By applying the automorphism group $D_4$,
we may assume that $a,b \geq 0$; however, this comes at a slight cost, in that
$z^r$ (that we are multiplying by) might be transformed to $z^{-r}$. So we will
have two cases, $r > 0$ and $r < 0$, for which the proofs are modestly
different.  

 \noindent {\it The case that $r > 0$.} Since we (now) have $a,b \geq 0$ and $w \in \supp f^k$, it
follows that $0 \leq s \leq ab$ and $a+b \leq k$. If $\tilde l_S (z^r w) < k+2$,
then $[z^r w,k+ 2] \in \SS R_f$, so we can incorporate all those terms into
$h(k)$. So we are reduced to considering $w$ \st $z^r w = z^{r+s} g^a h^b$ has
$\tilde l_S (z^r w) = k+2$. Hence by the formula, and since $a+b \leq k$, we
must have $a+ b = k$ and $r + s > a (k-a)$. 
 
 Pick a string of elements of $\brcs {g,g^{-1}, h, h^{-1},1}$ of length  $k$
whose product in that order is $w$. The number of such strings is $m(w,k)$.
Since $a + b = k$, the string must contain   $a$ $g$'s and $k-a$ $h$s. 

 We define a substring to be a consecutive sequence inside the original string;
it is determined from the original string by specifiying the initial and
terminal positions (of the substring). We look for two types of substrings.  

 The first is a substring of the  form beginning and ending with $g$ and having
exactly $r$ $g$'s (if $r > a$, there may  not be any). The number of such
substrings, if nonzero, is at least $a-r$: at the first occurrence of $g$,
proceed along the string until the exactly $r$ $g$'s have occurred and terminate
there; proceed to form the substring beginning with the second $g$, continuing
to the $r+1$st $g$, etc. This final substring occurs when we have reached the
$a-r+1$th $g$. Hence we obtain $a-r+1$  such substrings, except if $r > a$, in
which case there are none at all. 

 For each such substring, insert $h$ immediately before the initial $g$ of the
substring, and $h^{-1}$ immediately after the terminal $g$ of the substring.
This creates a word of length $k+2$, which we can write as $w_1 h w_2 h^{-1}w_3$
where $w_i$ are the products of three substrings (the original was the product
$w_1 w_2 w_3$. We can write $w_2 $ in reduced form, say $z^p g^r h^e$; then
$hw_2 h^{-1} = z^p (hgh^{-1})^r h^e = z^{p+r} w_2$. Hence the new string yields
$z^r w$. 
 
  If $r < a$, distinct original strings  substring yield distinct strings of
length $k+2$, because of the location of the single $h^{-1}$ (and the fact that
the $r$ is fixed. It follows immediately that if $r \leq  a$, then $m(z^r w,
k+2)/m(w,k)  \geq a-r+1 $.
 
 Similarly, we can do the same exploiting $g^{-1} h^r g = z^r h^r$. In other
words, for each string that yields $w$, look for the substrings beginning and
ending with an $h$ (if $r =1$, these are just one-element substrings), and
insert $g^{-1}$ just before the initial term of the substring and $g$ just after the terminal $h$
in the substring. There are, as before, exactly $k-a -r +1$ such substrings, and
again, the location of the  unique $g^{-1}$  in the length $k+2$ position allows
to conclude that different  strings yielding $w$ yield $k-a-r +1$ different
strings of length $k+2$ yielding $z^r w$, provided $r \leq k-a$, and thus
 $m(z^r w, k+2)/m(w,k)  \geq k-a-r +1 $. 
 
  Since $r \leq k/2$, it follows that $m(z^r w, k+2)/m(w,k)  \geq \max \brcs{a-r
+ 1, k-a -r +1} \geq k/2 - r$. Thus for $r > 0$, we can take $B(k) = k/2 - r$. 
 
 {\noindent $r <0$}. Pick $w = z^s g^a h^b$; if $\tilde l_S(z^{r}w) <k +2$, then
$z^{r} w$ can be incorporated into $h(k)$. Hence $\tilde l_S(z^{r+s}g^a h^b) =
k+2$. Since $a,b \geq 0$, $0 \leq s \leq ab$, and $a+ b \leq k$, this can only
occur (since $r < 0$) if $a+b = k$ and $s+r < 0$. In particular all the strings
of elements of $\supp f$ that realize $w$ can consist only of $g$'s and $h$'s. 

  Taking the same substrings as in the previous case with $|r| = -r$ $g$s or
$|r|$ $h$'s, we exploit $gh^{-r} g^{-1} = z^{r}h^{-r}$ and $h^{-1}g^{-r} h  =
z^r g^{-r}$ by inserting $g$ immediately preceding the first $h$ in the
substring, and $g^{-1}$ immediately after ther terminal $h$ in the substring,
and similarly $h^{-1}$ and $h$ for the substrings beginning and terminating with
a $g$ and having exactly $-r$ $g$'s. 

 The same arguments as in the case $r > 0$ now yield 
$$\frac{m(z^r w, k+2)}{m(w,k) } \geq
\max \brcs{a+r + 1, k-a+r +1} \geq k/2 + r = k/2 - |r|.$$
 Thus whether   $r > 0$
or $r < 0$, we can take $B(k) = k/2 -| r|$.

 At this point, we have $(k/2 - |r|)z^r f^k  \leq f^{k+2} + h(k)$, where $[h(k),
k+2 ] \in \SS R_f$, and letting $k \to \infty$, we have that for all positive
integers $N$, $N[z^r,2] + \SS R_f \leq [1,0] + \SS R_f$, so that $[z^r,2] + \SS
R_f$ is an infinitesimal of $R_f /\SS R_f$, completing the proof of \endstwo(a).   

\noindent {\it Proof of Lemma \endstwo\paren{b,c}.} Let $\tau$ be a pure nonfaithful trace of $R_f$. By
Lemma \TRtwo, $\tau$ kills $\SS R_f$, so induces a trace on $R_f/\SS R_f$. By  \endstwo(a),
$\tau([z^r,2] = 0)$. Hence if $\tau$ is a pure trace that doesn't kill
$[z^r,2]$, then $\tau$ is a faithful trace, hence  is given by a normalized character on $G$ (by Margulis's theorem for nilpotent groups).
Thus $\tau ([z^r,2]) =  1/P(\alpha, \beta)^2$ for some $(\alpha,\beta) \in
(\R^2)^{++}$. In particular, $\tau ([z^r -z^{r'},2]) = 0$ for all pure traces,
so that $[z^r - z^{r'},2]$  is an infinitesimal of $R_f$. 
 
 If $r \neq r'$ and $[z^r - z^{r'}, 2] \in \SS R_f$, then $[z^r - z^{r'}, 1] \in
R_f$, and thus there exists $m$ and $N$ \st coordinatewise,
$$
 -N f^{m+1} \leq (z^r - z^{r'})f^m \leq f^{m+1}.
 $$
 We may suppose that $r' > r$. The coefficient of $z^{r'}g^m$  in the middle
term is $-1$ since $z^{r'-r}g^m$ cannot appear in $f^m$. On the other hand, in
order obtain a term in $z^{r'}g^m$ on the left, we must have at least $m$ $g$'s,
and the remaining term is impossible (since it can only be one of ${g^{\pm1},
h^{\pm 1},1}$), yielding (b). 
 
 Finally, 
$$\eqalign{
 \left\| \widehat{[z^r,2]}\right\| &= \sup_{\tau \in 
\partial_e S(R_f) } \tau(z^r,2]) \cr
 & = \sup_{(\alpha, \beta) \in (\R^2)^{++}} \frac 1{\(1 + \alpha + \frac 1\alpha
+ \beta + \frac 1\beta\)^2} \cr 
  & = \slfrac 1{25}.\cr 
}$$\qed 

 \noindent {\it Proof of Proposition \endsthr\paren{a}.} If $k > |a| + |b| +2$,
we are done by Corollary \endsone, so we may assume $k =
|a| + |b| + 2$. Applying our  $D_4$ group of automorphisms, we can assume that
$a, b \geq 0$ (the sign of $r$ may change, but doesn't affect anything). Now the displayed inequality in the proof of \endsone\ applies here, yielding $\RR_{w,k} \leq \RR_{z^n,2}$ for some integer $n$, thus $\RR_{w,k}$ is order bounded.

\noindent
{\it Proof of Proposition \endsthr\paren{b}.} Since $p \in ( AG)^+$, we can write $p = \sum
(p,w)w$ where the coefficients, $(p,w)$ are all positive. Hence for each $w \in
\supp p$, $\RR_{w, k}$ is bounded above by a multiple of $\RR_{p,k}$. Thus if
$\RR_{p,k}$ is locally order bounded, so is every $\RR_{w,k}$, and thus $k\geq
\sup \tilde l_S (z^r w)$ (running over the conjugates of $w$) for every $w$.
Thus (a) applies, and so $\RR_{w,k}$ is order bounded, so their linear
combination, $\RR_{p,k}$, is as well. \qed

\noindent  {\it Closure of the set of faithful pure traces} (We have already seen that this closure does not exhaust the set of traces, so there must be others.) Form the Laurent polynomial ring $A[x^{-1}, y^{-1}] = \A[\Z^2]$ (the group ring), and let $\Arrow \pi; AG . A[\Z^2]$ be the ring homomorphism induced by the quotient map $G \to G/[G,G] = \Z^2$. 
Under this map the monomial $z^r g^a h^b \mapsto x^a y^b$ and sends $f \mapsto 1 +x + y + x^{-1} +y^{-1}:= P$. Thus $\pi$ intertwines (in the obvious way), $f$ and $P$, and sends positive elements to positive elements. This induces a positive homomorphism $A_f \to A_P$ (where $A_P = \lim \Arrow P \times ; A[\Z^2]. A[\Z^2]$). This maps the positive cone onto the positive cone, and it is easy to check that the order ideal $R_f$ is sent to $R_P$, again positive and onto the positive cone (but not a quotient by an order ideal). In particular every trace on $R_P$ yields a trace on $R_f$ by composition, and moreover, pure traces correspond to pure traces. 

 Since $\Z^2$ is abelian, $R_P$ is a commutative ring [H1, H2] and its pure traces are completely known (they are multiplicative). The faithful ones are precisely the characters of $\Z^2$ (hence yield all the faithful pure traces of $R_f$), and the unfaithful ones are their limit points, which are easily described, as in [H2]. They correspond to the boundary of the  Newton polytope of $P$. The latter is the convex hull of $\brcs{\pm(1,0), \pm{0,1}}$ a lozenge; its interior points yield the faithful ones, the points on the edges other than    the vertices yield unfaithful but not discrete traces, and the four vertices yield corresponding discrete traces. The identification of the pure normalized traces of $R_P$ with the points of the Newton polytope is implemented by the moment map     [H2?]. 

 The traces corresponding to points in the relative interior of the edge in the first quadrant are given by  
$$
 \phi_t\([p, k] \) = \lim_{x \to \infty} \frac{p(x,tx)}{P^k(x,tx)}
$$
where $p \prec P^k$ and $0 < t < 1$ (that the limit exists is a consequence of l'H\^opital's rule). Traces corresponding to  other three edges  come from applying the dihedral group to the Newton polytope. Finally the trace corresponding to the vertex $(1,0)$ is a discrete trace, and given by $\Arrow \phi_{0,0}; [p,k] . (p,x^k)$,  here using inner product notation to describe the coefficients, and the other three obtained by applying the corresponding reflections. It is easy to verify that any limit point of these traces as traces on $R_f$ is a limit of the lifted traces iff the traces on $R_P$ converge to the trace on $R_P$, i.e., the topology is the same, whether as traces on $R_f$ or $R_P$, and the set of all these pure traces is compact in either case. In particular, we have described the closure of the set of pure faithful traces on $R_f$.

 Alternatively, the four discrete traces are obtained as in \thmxx, limits along a path, as we will   describe in more generality later.  

 The map $R_f \to R_P$ obtained above is {\it not\/} a quotient by an order ideal---its kernel is obtained from the augmentation ideal of the centre of $G$, and contains no positive elements. 

 \noindent {\it Remaining pure traces} We will describe eight families of uncountably many pure discrete traces (also containing the four discrete traces obtained above). Then we show that these constitute all of the remaining pure traces, and finally determine their limit points (which correspond to the traces arising from the four edges of the lozenge). The methods are brutal.

First, we observe that there is a natural action of the dihedral group $D_4$, not just the obvious one on the lozenge (the Newton polytope), but on $G$ itself (and of course inducing the action on the lozenge and on the pure traces of $R_P$).

\noindent {\it Reduction to $\overline R_f$, i.e., quotient of $R_f$ via an order ideal, via the various $\Gamma$s.} We prepare for determination of the pure traces of $R_f$ (completed in the next two sections) and in particular, determine antecedents within the Bratteli diagram of $R_f/\SS R_f$. 

\noindent  From the definitions (sections 1 and 2), if $m > \tilde l_S (w)$, then $[e_w,m] \in \SS R_f$. We have defined $\Gamma_m' = \Set{w \in G}{l_S (w) = \tilde l_S (w)}$, and it is easy to check that $\Set{[e_w,m]}{w \in S^m\setminus \brcs{S^{m-1} \cup \Gamma_m'}}$ spans $\SS R_f$. We claim there are four indecomposable  order ideals $I_{(\pm,\pm)}$---corresponding to the quadrants of $\Z^2$---\st if $\tau$ is a pure trace of $R_f$ killing $\SS R_f$, then $\ker^+ \tau$ contains one of the four order ideals. 

To explain further, an order ideal in a dimension group is {\it indecomposable\/} if it cannot be represented as an intersection of two larger order ideals (we could also call this {\it prime,} but there is potential confusion). If $\tau$ is a trace of a dimension group, $J$, then $\ker^+ \tau$ is the subgroup of $\ker \tau$ generated as an abelian group by the positive elements therein. 
It is easy to check that $\ker ^+ \tau$ is an order ideal, and is the sum of all the order ideals contained in $\ker \tau$. It is routine to prove that for dimension groups with order unit, if $\tau$ is a pure trace, then $\ker^+ \tau$ is indecomposable as an order ideal. 

Let $\Arrow \pi; H_3 . \Z^2$ be the quotient map $z^r g^a h^b \mapsto x^a y^b$ (regarding $\Z^2$ as a multiplicative group).

We describe the order ideal corresponding to the first quadrant, $I_{(+,+)}$. Define $\Gamma_m\pplus$ to be  $\Set{w \in \Gamma_m'}{\pi(w) \in \Z^2_+}$; that is, $\Gamma_m\pplus $ is $\Set{z^r g^a h^{m-a}}{a,b,m-a \geq 0; r \leq (m-a)a }$, so that $\pi(\Gamma_m\pplus)$ is in the first quadrant of $\Z^2$ (since $(a,b) \geq 0$). Define $I_{(+,+)}$ to be the span of  $\Set{[e_w,m]}{w \in \Gamma_m \setminus \Gamma_m\pplus}$. It is easy to check that $I_{(+,+)}$ is an order ideal, and the quotient, $R_f/I_{(+,+)}$ is naturally order isomorphic to the limit $\Arrow \times (g+h);\Z \Gamma_{m}\pplus . \Z \Gamma_{m}\pplus$, obtained from repeated multiplication by $\overline f:= g+h$. Denote the limit dimension group $\overline R_f$ (the overline on the $f$ is too small to be seen in $R_{\overline f}$, which is what the notation should be). This is the order ideal generated by $[1,0]$ in the dimension group $\lim \Arrow \times (g+h) ;\Z G.\Z G$, but note that $\overline f = g+ h$ is not admissible. Nonetheless, we have to deal with it. 

By applying the action of $D_4$, we obtain three more order ideals corresponding to the remaining quadrants of $\Z^2$. Now it is not difficult to show the intersection of all four is precisely $\SS R_f$. Moreover, there plenty of pure faithful traces on $\overline R_f$ (and its three automorphs), e.g., if $w = z^r g^a h^{m-a} \in \Gamma_{m}\pplus$, and $\alpha$ is a positive real number, $[e_w,m] \mapsto \alpha^a (1-\alpha)^b$ yields a pure trace. This  pulls back to a pure trace on $R_f$ which kills $I_{(+,+)}$, so is different from the faithful pure traces corresponding to eigenvectors. Hence $I_{(+,+)}$  is $\ker ^+ \tau$  for lots of pure traces $\tau$. In particular, $I_{(+,+)}$ is indecomposable, and thus so are each of $I_{(\pm,\pm)}$. It obviously contains $\SS R_f$. Now we show that if $\tau$ is a pure trace \st $\SS R_f \subset \ker \tau$, then $\ker \tau$ (and thus $\ker^+ \tau$) contains one of the four order ideals $I_{(\pm,\pm)}$. Since $\ker^+ \tau$ is indecomposable, it suffices to show the following.  

\Lem Lemma \oneffn. Any indecomposable order ideal of $R_f$ that contains $\SS R_f$ contains at least one of the four order ideals $\brcs{I_{(\pm,\pm)}}$. 

\Pf Pick $w = z^r g^a h^b \in \Gamma_m'$ and $w'= z^R g^A h^B \in \Gamma_{m'}'$ (so $a+ b = m$ and $A+B = m'$) \st $(a,b)$ do not belong to the same quadrant of $\Z^2$ (a nonzero lattice point can lie in two quadrants). By applying the appropriate automorphism from  $D_4$, we can assume that $a,b \geq 0$  and $a+ b > 0$. Consider the order ideals $J =\langle [e_w, m] \rangle$ and $J' = \langle [e_{w'},m']\rangle $. We claim that their intersection lies inside $\SS R_f$. 

First, we note that multiplication of  $w$ by  $x = 1$ or if $a> 0$ by $x=g^{-1}$,  or if $b > 0$, by  $x =  h^{-1}$,  will result in $[e_{xw},m+1] \in \SS R_f$. Similarly, multiplication of $w'$ by $x'= 1$, or if $b' < 0$, by $h$, etc, will result in $[e_{x'w'}] \in \SS R_f$. 

Suppose that the intersection does not lie inside $\SS R_f$. In a dimension group, the intersection is an order ideal, and it then follows that there exists $v \in \Gamma_n'$ and positive integer $K$ \st $[e_v,n] \leq K [e_w,m], [e_{w'},m']$ and $[e_v,n] \not\in \SS R_f$. By taking the images in subsequent rows, we may assume that additionally that $n > m,m'$, and there are paths from $w $ to $v$ (of length $n-m$) and from $w'$ to $v$ (of length $n-m'$). Moreover, since $[e_v,n] \not\in \SS R_f$, we must have that $v \in \Gamma_{n}'$. Since $(a',b')$ does not lie in any of the quadrants that contain $(a,b)$, at least one of $\brcs{\sign a,\sign a'}$ and $\brcs{\sign b, \sign b'}$ is $\brcs{+,-}$. 

This leads to a contradiction; for example, suppose $a > 0 > a'$. Then any path from $w$ to $v$ cannot contain any of $1,g^{-1}$ (else $[e_v,n]$ would belong to $\SS R_f$). Hence the exponent of $g$ in $v$ must be at least as large as $a$, in particular, must be positive. However, since $a' < 0$, no paths from $w'$ to $v$ can contain a $g$, so that the exponent of $g$ in $v$ is less than or equal to $a' <0$. 

Finally, let $L$ be an indecomposable order ideal of $R_f$ that contains $\SS R_f$. If $L$ does not contain any of the four order ideals $\brcs{I_{\pm,\pm}}$, then there exist $w \in \Gamma_m$  and $w' \in \Gamma_{m'}$ as above \st both $[e_w,m], [e_{w'},m'] \not\in L$. We may thus form the two order ideals   $L + \langle [e_w, m] \rangle$, $L + \langle [e_{w'}, m'] \rangle$ (since  these are sums of order ideals, and $R_f$ is a dimension group, they are order ideals). It is easy to check that the intersection is contained in $L$ (using Riesz decomposition), contradicting indecomposability. 
\qed

\noindent {\it Antecedents\/} 
Consider the limit, $\Z \Gamma_m\pplus \to \Z \Gamma\pplus_{m+1}$ obtained from $w \mapsto hw + gw$; this is also $R_{\overline f}$ obtained from $H_3$ from $\overline f = g+h$ (not admissible!). Because the subscripted overline is difficult to see, we will rewrite this as $\overline R_f$. We have that $\Gamma_m\pplus = \Set{z^r g^a h^{m-a}}{0 \leq a \leq m; 0 \leq r \leq a(m-a)}$. We want to determine, for $w \in \Gamma_m\pplus$ and $k < m$, what the set 
$$
\Cal A_{w,k}:= \Set{v= z^s g^c h^{k-c} \in \Gamma_k\pplus}{\exists \text{ a path $v\to w$}},
$$
 the set of antecdents, is. Explicitly, $v \in \Cal A_w$ iff there exists $u \in \Gamma_{m-k}\pplus$ \st $uv = w$. This translates to $u = wv^{-1} \in \Gamma_{m-k}\pplus$, which in turn reduces to determing the set of possible $(s,c)$ \st $z^{r-s- c\delta} g^{a-c}h^{m-k+c-a} \in \Gamma_{m-k}\pplus$, where $\delta = m-k + c-a$. 

 We thus have the following constraints: 
 \item{($-1$)} (given) $r,a$ \st $0 \leq r \leq a(m-a)$ and $0 \leq a \leq m$; 
 \item{(0)} $0 \leq c \leq k$ and $0 \leq s \leq c(k-c)$; 
 \item{(i)} $a \geq c$ and $\delta \geq 0$; \item{(ii)} $0\leq r- s - c\delta \leq (a-c)\delta$. 
 
 From (ii), we deduce $r - (m-a)(a-c) \leq s \leq r - (k-c)(a-c)$, which together with (0) yields, 
 \item{(1)}  $ 0 \vee (r - a\delta) \leq s \leq (r - c\delta)\wedge c(k-c) $
 \item{(2)}  $0 \vee (a+k-m)  \leq c \leq a \wedge k$
 
\noindent It is a trivial calculation to show that $c \leq k ,a$ entails that $r-a\delta  \leq c(k-c) $. Moreover, these conditions while obviously necessary, are also sufficient.
 
\Lem Proposition \oneten. Suppose $w  = z^r g^a h^{m-a} \in \Gamma_m\pplus$, and let $k < m$ be a positive integer. Then $\Cal A_{w,k} = \Gamma_k\pplus$ iff $k \leq a \leq m-k$ and $k(m-a) \leq r \leq a(m-a) - ak$. 
 
\Rmk The latter condition forces $k \geq m/4$, among other conditions.

 \Pf Suppose the conditions hold. Then (2)  reduces to $0 \leq c \leq k$. Moreover $r \geq k(m-a)$ entails $r - c\delta = r- c(m-a) + c(k-c)  \geq c(k-c)$, and $r \leq a(m-a-k)$ implies $r-a \delta = r - a(m-a-k) -ac \leq 0$. Thus (1) reduces to $0 \leq s \leq c(k-c)$. There are no other constraints, so all $v \in \Gamma_{k}\pplus$ can be realized. 

 The converse is routine.\qed 

 \Lem Corollary \oneele. Suppose that $M$ is an infinite subset of $\N$, and for $m \in M$, $w(m) = z^{r(m)}g^{a(m)}h^{m-a(m)} \in \Gamma_{m}\pplus$. Suppose that 
$$\min\brcs{\frac{r(m)}{m-a(m)} , m-a(m) -\frac{r(m)}{a(m)}} \to \infty \quad \text{ along $M$}.$$ If $I$ is an order ideal of $\overline R_f$ such that none of $\brcs{[e_{w(m)},m]}_{m\in M}$ belong to $I$, then $I = (0)$. 

 \Rmk The analogous result (much easier) for Pascal's triangle is that sufficient for the same conclusion is  $w(m) = (a(m),m)$ with $|a(m) - m/2| = \oh m$, where the indexing the $m$th row is by $0,1,2,\dots, m$. 
 
 \Pf Since $r \leq a(m-a)$, we have that $r(m)/(m-a(m)) \leq a(m)$; similarly, $m-a(m) \to \infty$. Thus, if we set $k (m)$ to be the greatest integer less than or equal $$
\frac{r(m)}{m-a(m)} \wedge \(m-a(m) -\frac{r}{a(m)}\) \wedge \frac m2 \wedge \frac{m-a(m)}2,
$$
then $k(m) \to \infty$. Moreover, each $k(m)$ satisfies the conditions in Proposition \oneten, so that $\Cal A_{w(m), k(m)} = \Gamma\pplus_{k(m)}$.

Suppose that $[e_u,l] \in I$; from the definition of $\overline R_f$, this entails that $u \in \Gamma_l\pplus$. Pick $m \in M$ \st $k(m)\geq l$.  Then $[e_u,l] = [f^{k(m)-l}e_u,k(m)]$, and by construction, $[f^{m- k(m)}f^{k(m)-l }u, w(m)] \neq 0$. This entails that $[e_{w(m)},m] \leq [e_u,l]$. Hence $[e_{w(m),m}] \in I$, a contradiction. \qed

\comment
odd stuff about sizes of the antecedent sets go to infinity
Now we obtain conditions to guarantee existence of $k(m)$ \st the sizes of the antecedent sets at $k(m)$ for $z^{r(m)} g^{a(m)} h^{m-a (m)} \in \Gamma\pplus_m$ go to infinity.  

this bunch of stuff may not be necessary, or it might be. Restate as theorem and prove \fillmein

 We also note that  $r - (k-c)(a-c) \geq c(k-c)$ iff $r \geq (k-c)a$, or in other words, $c \leq k - r/a$ (since all the letters represent integers, the last means $c \leq k - \flo{r/a}$). 
 
 And a necessary condition for there to be any choices for $s$ at all is that $c(k-c) \geq r - (m-a)(a-c)$. This translates to the quadratic $c^2 - c(k-m -a) + r-a(m-a) \leq 0$. As the constant term is non-positive, this occurs only if $2c \leq 2\rho:=(k-m-a) + \sqrt{(k^2+(m-a)^2 + (m-a)(4a-2k) -4r)}$. It is straightforward that at least $\rho \geq 0$, but we obtain another constraint for $c$, $c \leq \rho$. 
 
 Next we consider the possible constraints assuming $ c \leq k-r/a$. Then $r-(m-a)(a-c) \leq s \leq c(k-c)$. There are two possibilities.

\noindent $r < (m-a)(a-c)$; this yields a further restriction on $c$, namely $c < a-r/(m-a)$ (observe that $r/(m-a) \leq a$ from ($-1$)). We then have the constraints, $0 \leq c \leq (k-r/a)\wedge(a-\flo{r/(m-a)} -1)$. Set $Y =(k-r/a)\wedge(a-\flo{r/(m-a)} -1)$. Then the number of possibilities for corresponding $s$ is thus $\sum_{0 \leq c \leq Y} c(k-c)$.
 
 \noindent $r \geq (m-a)(a-c)$, which reduces to $c \geq a-r/(m-a)$. We thus have the constraints on $c$, $a-r/(m-a) \leq c \leq k - r/a$; in order for this to even occur, we must have $k -a \geq rm/a(m-a)$. Then the number of (additional) choices for 

We next  consider the case that $c > k-r/a$, so that $0 \vee(r - (m-a)(a-c)) \leq s \leq r - (k-c) (a-c)$. This entails that $r \geq (k-c)(a-c)$, that is, $c^2 - (a+k)c + ak-r \leq 0$. 

 If $ak \leq  r$, this is equivalent to $2c \leq (a+k) + \sqrt{(a-k)^2 + 4r}$, in particular, $c \leq a+k$; since we already have $c \leq k \wedge a$, this is a redundant condition. 

 On the other hand, if $a k > r$, we would have to have $(a+k) - \sqrt{(a-k)^2 + 4r} \leq 2c \leq (a+k) + \sqrt{(a-k)^2 + 4r}$; as $2c \geq 2(k - r/a)$; much easier argument in notes (somewhere) \fillmein
\endcomment

 \noindent {\it Nonnoetherianness\/} 
We show that the dimension group  $R_{\overline f}$, the order ideal of $\lim \Arrow\, (g+h)\times; \Z G.\Z G$ generated by $[g+h,1]$, is not noetherian; hence $R_f$ (original $f$) is also nonnoetherian. 

 For each positive integer $k$, define $w(k) = z^{2k} g^{2k}h^{2k+1} \in \Gamma_{4k+1}\pplus$, and $x_k = [e_{w(k)}, 4k+1] \in R_{\overline f}^+$. Now define $I_n$ to be the order ideal (of $R_{\overline f}$) generated by $\brcs{x_1, x_2, \dots, x_n}$. We will show that $x_n \not\in I_{n-1}$. 

 \Lem Lemma \onetwe. Suppose that for some integers $k,l,N > 0$ and $r,a \geq 0$, there exists $u \in \Gamma_{4l+N}\pplus$ \st $u w(k) = z^r g^{a} h^{2(k+l)+ 1} \in \Gamma_{4(k+l) + N +1}\pplus$. Then $r \geq 4kl + 2k$.
 
 \Pf %
\comment
 The element $u$ is a product of $4l + N$, $g$s and $h$s; since the exponent of $h$ in $uw(k)$ is $2(k+l) + 1$, and that in $w(k)$ is $2k+1$, $u$ must have exactly $2l$ $h$s. 
\endcomment
Since $w(k)$ has   has $2k$ $g$'s and $2k+1$ $h$'s  and $uw(k)$ has $2k+2l+1$ $h$s, $u$ must have $2l$ $h$s. But then the exponent of $z$ in the product must be at least as large as $2k + (2l)(2k)$. \qed

 \Lem Lemma \onethi. For all $n \geq 2$, $x_n \not\in I_{n-1}$. 

 \Pf Suppose $x_n \in I_{n-1}$. Then there exists a positive integer $C$ \st $x_n \leq C \sum_{i=1}^{n-1} x_i$. By Riesz decomposition, for each  $1 \leq i \leq n-1$, there exists $y_i \in \overline R_{f}^+$ \st $y_i \leq C x_i$ and $x_n = \sum_1^{n-1} y_i$. Hence there exists a positive integer $M$ \st 
$$
 \supp {\overline f}^M w(n) \subset \cup_{i=1}^{n-1} \supp {\overline f}^{M+4i}w(n-i).
$$

 Obviously $g^M w(n)  = z^{2n} g^{2n+M} h^{2n+1} \in \supp {\overline f}^M w(n) $. Hence there exists $i$ with $1 \leq i \leq n-1$ \st $g^M w(n) \in \supp {\overline f}^{M+4i}w(n-i)$. Hence there exists $u \in \supp f^{M+4i}$ \st $u w(n-i) = g^M w(n)$. It follows that $u \in \Gamma_{M+4l}\pplus$ (since $uw(n-i) \in \Gamma_{4n+1 + M}$). 
 
  Write $k = n-1$ and $l = i$. We have that $g^M w(n))  \in \Gamma_{4(k+l)+ M+1}\pplus$, and $w(n-i)) \in \Gamma_{4k+1}\pplus$. By Lemma \onetwe, the exponent of $z$ in $g^M w(n) = u w(n-i)$ has to be at least as large as $4kl + 2k$. However, $g^M w(n) = z^{2n} g^{M+2n}h^{2n+1}$, so the exponent of $z$ is $2n = 2(k+l)$. This forces $2(k+l) \geq 4kl + 2k$, that is, $2l \geq 4kl$. This is clearly impossible, as both $k$ and $l$ are at least one. \qed

 \Lem Corollary \oneftn. The dimension groups $\overline R_f$ and $R_f $ are not noetherian. 
 
\Pf By the previous two lemmas, we have  obtained an increasing sequence of order ideals $I_1 \subset I_2 \subset \dots $  in $\overline R_f$ with the property that all the inclusions are strict. This violates the definition of noetherian. Since $\overline R_f$ is a factor by an order ideal of $R_f$, the latter cannot be noetherian.\qed
 
 In this construction, the union, $\cup I_n$, is an order ideal of $\overline R_f$ with no order unit. This obviously pulls back to an order ideal of $R_f$ with the same properties.  Translated into the language of space-time cones, it corresponds to a space-time cone that requires infinitely many start-up points, inside the cone generated by $ f$ starting at $(1,0)$ (where $1$ represents the identity element of the group and $0$ is the time or level).

It should be possible to prove that for any admissible $f$ in any finitely generated nonabelian torsion-free nilpotent group, $R_f$ fails to be noetherian and  $R_f$ has countably infinitely many maximal order ideals, all of which admit order units. Both of these properties do not depend on the coefficients of $f$, merely on the support. 

 For finitely generated abelian groups, the corresponding $R_f$ is noetherian (for {\it any\/} $f \in A G^+$) as follows (eventually) from the Hilbert basis theorem (stated as, a finitely generated commutative ring is noetherian), because $R_f$ is a finitely generated commutative ring and order ideals of $R_f$ are ring ideals thereof. 

 On the other hand, if we let $G$ be the free group on two generators and $f = 1 + g + h + g^{-1} + h^{-1}$, $R_f/ \SS R_f $ is order isomorphic to $C(X,\Z)$ where $X$ is the path space of the Cayley diagram (hence is a Cantor set). It has uncountably many maximal order ideas (in contrast to countably infinite maximal order ideals for the Heisenberg group), none of which possess an order unit (in contrast to that for the Heisenberg group, where all do).

\SecT \Ends\ Ends

In this section, we define collections of pure traces arising from relatively simple paths in the Bratteli diagram for $R_f/\SS R_f$, mostly for the Heisenberg group, but we begin with general considerations. These are designed to be able to describe all the pure perfidious traces, that is, those killing $\SS R_f$. 

Suppose $B$ is a Bratteli diagram, and let $D$ be the corresponding dim group. We suppose that $D$ has an order unit, that is, for every node, there exists a path from one of the nodes at the top level. (Diagrams read down, but realizing the dimension group as a direct limit from the diagram reads across.) There may be special traces on $D$ arising from the paths. 

Suppose there is a sequence of nodes $(x_n,n)$ (where $n = 0,1,2, \dots$ represents the level, or discrete time) and $x_n $ is a node at level $n$, such that first, the sequence $((x_n,n))$ represents a path (meaning, there is at least one arrow   $(x_n,n) \to (x_{n+1},n+1)$), and second for all sufficiently large $n>0$, $(x_n,n)$ is the only node for which there is an arrow to $(x_{n+1},n+1)$ (that is, $(x_{n+1})$ has unique antecedent). If we have such a sequence $((x_n,n))$, let $c(n)$ count the number of arrows from $(x_n,n) \to (x_{n+1},n+1)$, and let $C$ denote the rank one limit group
$\lim \Arrow \times c(n); \Z.\Z$, together with the natural choice of order unit, $[(e_{x_0},0)]$. We do not require that any of the $(x_n,n)$ have a unique successor; that would be uninteresting.

Then we can define a trace $\Arrow t; D. C$ (regarding $C$ as a subgroup of $\Q$, in turn, as a subgroup of $\R$), which is onto, and whose kernel is an order ideal (this is an extremely rare property for a trace; it implies purity, but hardly any pure traces satisfy this). Simply define $t([h,k]) $ to be the $x_k$ component of $h$ divided by $\prod_{i < k} c(i)$---that is, we can write $h = \sum a(j)e_j$ where the $j$s run over the vertices at level $k$, and pick out $j = x_k$ (the vertex that hits the path), and then the value of the trace on $[h,k]$ is $a(x_k)/\prod_{i< n} c(i)$. It is easy to check (from the uniqueness of the predecessor property) that this is a pure trace.

At each level, the kernel of $t$ is generated by positive elements, so $\ker t$ is directed, and the kernel of any trace is always convex, and thus the kernel of this trace is an order ideal, and now it is clear that $D/\ker t$ is isomorphic to $C$ as ordered abelian groups.

If $G = \Z^n$ and we look at $R_f/{\SS R_f}$ (where $0 \in \supp f$ and $\supp f$ generates $G$ as a semigroup), then it is easy to check the traces that are of this form arise only from the vertices of $\Log f$, and so there are finitely many of them; in addition, every maximal order ideal of $R_f$ is the pre-image in $R_f$ of $\ker t \subset R_f/{\SS R_f}$. The latter property is lost if we go to abelian by finite groups (such as the infinite dihedral group) for every reasonable choice of $f$ (except a few degenerate ones), there are still only finitely many of them (possibly none).

When we go to nilpotent groups, even the simplest nontrivial one (the central extension of $\Z$ by $\Z^2$, with generators and relations $uv= zvu$ and $[z,u] = [z,v] = 1$), for at least one choice of $f$, there are infinitely many traces of this kind on $R_f/{\SS R_f}$, and moreover, the restrictions of the faithful pure traces fail to be dense in $\partial_e F_0$ (unlike the situation for abelian groups, and also for abelian by finite, although we haven't proved this yet).

To see what is going on, recall that if we take $\Z G$ (rather than the usual $\R G$), there is a natural Bratteli diagram associated to $R_f/{\SS R_f}$. Assume that $1 \in \supp f$ and the latter generates $G$ as a monoid. Define $\Gamma_0 = \brcs{1}$, and $\Gamma_k$ is to consist of those $g \in \supp f^k$ \st $g \not\in \supp f^j$ for all $j < k$. Then $R_f/{\SS R_f}$ is naturally isomorphic to the direct limit, $\lim \Arrow P_k; \Z \Gamma_k. \Z \Gamma_{k+1}$, where $P_k e_g = \sum_{g' \in \Gamma_{k+1}} (fe_g,g') e_{g'}$ for $g \in \Gamma_k$ (we do not have to specify the level, since $\cup \Gamma_k$ is disjoint). This can be refined, as in section 1, to restrict to $\Gamma_k'$, in case the choice of $f$ is holey.

For the Heisenberg group $G = \Z \times_{\theta}\Z^2$, as usual $f = 1 + g + g^{-1} + h + h^{-1}$. Then $\Gamma_1 = \brcs{g^{\pm 1}, h^{\pm 1}}$, $\Gamma_2 = \brcs{g^{\pm 2}, h^{\pm 2}, gh, zgh, \dots }$, and it easily follows that $\supp f$ generates $G$ as a semigroup. In $\Gamma_k$, we can obtain all words of the form $z^r g^{a}h^{k-a}$ for $0 \leq r \leq a(k-a)$, and these are the only words of degree $k$ in $\Gamma_k$. Some of these have unique predecessors, that is there is one word $w'$ in $\Gamma_{k-1}$ \st one of the words appearing in $fw'$ is $w$.

If $w \in \Gamma_k$, then $|d(w) |\leq k$. If $w' \in \Gamma_{k-1}$ and $fw'$ contains a word of degree $k$, then $w'$ has degree $k-1$, and then it follows that there are at most two possible words of degree $k$ in $fw'$, namely $uw'$ (which has to be put in reduced form) and $vw'$ (which automatically will be in reduced form). To see which   words in $\Gamma_k$ have unique predecessors (they always have at least one, just from the definition of $\Gamma_k$), consider the generic word of degree $k$ in $\Gamma_k$, $w = z^r g^a h^{k-a}$, and let $w'  = z^{r '} g^{a'} h^{k-a'-1}$, with $0 \leq r' \leq a'(k-1-a')$. Left multiplying the latter by $g$ yields $z^{r'}g^{a' +1} h^{k-a'-1}$ and left multiplying by $h$ yields $z^{r'+ a'}g^{a'}u^{k-a'}$.

Thus if $w$ has two predecessors then both $r \geq a$ and $1 \leq a \leq k-1$. Moreover, if $r > (a-1)(k-a)$, then only multiplication by $h$ can yield $w$.  

\noindent {\it Discrete traces on $\overline R_f$}
It is more convenient to work with $\overline{R}_f$, which we now do (as usual, $D_4$, the dihedral group, acts here).  Relabel the nodes at level $m$ as triples, $\brcs{(r,a, m-a}$, with $0 \leq r \leq a(m-a)$ and $m \geq a \geq 0$. Using the nodes with unique antecedent, we can construct a lot of discrete traces, in fact, a two-integer parameter family of them. In particular, $w = z^r g^a h^{m-a}$ has unique antecedent (or predecssor) in the $\overline R_f$ diagram) if and only if either $r < a$ or $r > (a-1)(m-a)$.

Let $(r,a,m-a)$  or $z^{r}g^a h^{m-a} \in\Gamma_m''$ represent a node with a unique antecedent, say with $r <a$; and its antecedent (at level $m-1$) is given by $(r, a-1,m-a-1)$; repeat this (that is, left multiplying by $g^{-1}$) until the middle coordinate equals the first, that is, $a-r$ times, to yield $(r,r, b)$. If $r = 0$, we began with $g^a h^{m-a}$, and eventually hit $h^{m-a}$---which of course does have unique antecedent, but obtained by left multiplying by $h^{-1}$. 
If $r > 0$, we arrive at $z^r g^r h^{m-a}$, which has two antecedents ($g^r h^{m-a-1}$ and $z^rg^{r-1} h^{m-a}$).

 Relabel $m-a = b$.  We obtain a path from level $r+b$ going off to infinity, such that for every level exceeding $r+b$, the corresponding node has unique antecedent: for $m \geq r+b$ the node at the  $m$th level is $((r, m-b,b)$. This allows us to define a trace,  $\tau_{r,b}$, from the general construction above. 

The multiplicities (the $c(n)$) are all $1$. The denominators in the definition are there so that it is normalized at $\1 = [\overline{f}^m,m]$. The initial point of the path, $(r,r,b$, determines the rest of it, and moreover, it follows from the unique antecedent property that $(\overline f^m, z^r g^{m-b}h^b) = (\overline f^{b+r},z^r g^r h^b)$ (or, $p(r,m-b,b) = p(r,r,b)$) for all $m \geq b+r$, that is, the successors in this particular trajectory have the same multiplicity in the corresponding power of $g+h$. 
$$
 \tau_{r,b}([w,m] )= \cases 0 & \text{if $m \geq b+r$ and $w \neq z^r g^{m-r}h^b$}\\ 
 \frac 1{(f^{b+r},z^rg^r h^b)} = \frac{1}{p(r,r,b)}& \text{if $m \geq b+r$ and $w = z^r g^{m-r} h^b$}\\ 
 \frac{(f^{b+r-m}w,z^r g^r h^b )}{(f^{b+r}, z^r g^r h^b)} & \text{if $m < b+r$.}\\ \endcases
 $$
To check the last case (that is,  in the case that $m < b+r$), we count the number of paths from $w \in \Gamma_m''$ to $z^r g^r h^b$. This is just the number of $u \in \Gamma_{b+r-m}''$ \st $uw = z^r g^r h^b$. If we write $w = z^R g^A h^{m-A}$, then $u = z^r g^r h^bw^{-1}$ simplifies to $u = z^{r-R - A(b+A-m)}g^{r-A}h^{b-m+A}$. The number of such that belong to $\Gamma_{b+r-m}''$ is precisely $p(r-R-A(b+A-m), r-A, b-m+A)$. In order for this to be nonzero, necessary and sufficient is that   $r \geq A$, $b+A \geq M$, and $0\leq r -R- A(b+A-m) \leq (r-A)(b-m+A)$ (implicitly, $A\geq 0$, $m \geq A$, and $0 \leq R \leq (m-A)A$, since $w \in \Gamma_m''$). In any event, the last line of the display can be replaced by 
$$
\tau_{r,b}([w,m]) =  \frac{p(r-R-A(b+A-m), r-A, b-m+A)}{p(r,r,b)}, 
$$
where $w = z^R g^A h^{m-A} \in \Gamma_{m}''$. 

We can similarly define the pure traces $\tau^{s,c}$ obtained from the points with unique antecedent arising from multiplication by $h$s. This is slightly more complicated. The nodes with unique antecedent arising from $h$ are of the form $z^{r}g^c h^k$ where $r > (c-1)k$. As before, we continue backwards (pre-multiplying by $h^{-1}$, until we reach either $(0,c,d)$ (that is, $g^c h^d$), which may have unique antecedent (if $c \geq 1$) via $g$, or $(0,0,0)$ which corresponds to the identity group element; otherwise, we reach $z^{(c-1)d}g^c h^d$ where $c \geq 2$ and $d\geq 1$, which has two antecedents. 

This allows us to define the traces $\tau^{c,d}$, and for variety,  we express them in terms of general elements $[q,m] \in \overline R_f$, where          $\supp q \subseteq \Gamma_m''$. 

$$
 \tau^{c,d}([q,m] )= \cases
 \frac {\(q,h^{m-c-d}\cdot z^{(c-1)d}g^c h^d\)}{p(d,c,d)} & \text{if $m \geq c+d$}\\ 
 \frac{(f^{m-c-d}q,z^{(c-1)d}g^ch^d)}{p(d,c ,d)} & \text{if $m < b+r$.}\\ \endcases
 $$
The group element $h^{m-c-d}\cdot z^{(c-1)d} g^c h^d$ is just the corresponding node at level $m$; and of course, the denominator should be $(f^{c+d}, z^{(c-1)d} g^c h^d)$, but this is $p((c-1)d,c,d) = p(d,c,d)$ by symmetry in the first variable. 

There are some restrictions on the values of $(r,b)$ and $(c,d)$ that appear as   subscripts and superscripts. For example, $\tau_{0,0}$ picks out the path $(1,g,g^2, \dots)$ (beginning at the zeroth level), and $\tau_{0,b}$ corresponds to the path $(h^b, gh^b, g^2 h^b, \dots)$ (beginning at level $b$); but $\tau_{r,0}$ is not defined if $r > 0$. The path $(zgh^b,zg^2 h^b, z g^3 h^b, \dots)$ is picked out by $\tau_{1,b}$.

On the other hand, $\tau^{0,0}$ corresponds to the path $(1,h,h^2, \dots)$, but that of  $\tau^{0,d}$ is cofinal with it for $d> 0$, so yields the same trace. 

These traces are pure  and lift to $R_f$;  the   $D_4$ action yields three more versions of each of $\tau_{r,b}$ and $\tau^{s,c}$.  Except for $\tau^{0,0}$, $\tau_{0,0}$ and their automorphs, these are different from the pure traces that factor through the abelianization $G \to G/G'$. The latter are multiplicative, even on $R_f$ (although not all are faithful). On the other hand, the kernel of a pure discrete trace is a maximal order ideal. In particular, these  (with the exceptions noted) fail to be multiplicative, and thus are not in the closure of the set of faithful pure traces ($\cup_{\lambda \geq 0} F_{\lambda}$). 

\comment

Now consider the elements $x = [vu,2] $ and $y = [uv,2] = [vuz,2]$ in $R_f^+$. Since $G$ is nilpotent, every pure trace not in $F_0$ is given by a real character, $\tau_{\chi} \in F_{1/\chi(f)}$ (here $\chi(f) = 1 +\chi(u) + \chi(u)^{-1} + \chi(v) + \chi(v)^{-1}$, and because $G/G' = \Z^2$, there are plenty of real characters of $G$. For each one of these, $\tau_{\chi}(x) = \tau_{\chi} (y)$, hence for any trace in the closure of $\cup_{\lambda} F_{\lambda}$, whether pure or not, $\tau (x-y) = 0$. However, it is easy to find (pure) traces lying in $F_0$, $\rho$ and $\sigma$, \st $\rho(x) \neq 0 = \sigma(x)$ while $\sigma(y) \neq 0 = \rho (y)$.

We just note that the words $(v^k u)$ form a path starting at level $1$, with $n=0$, of points in $\Gamma_{k+1}$, each having unique predecessor. The corresponding trace is $\rho$, since $vu $ appears among the vertices. At the other extreme, $(v^k u z^k)$ is also a path of points in successive $\Gamma_{k+1}$, each with unique predecessor, and this yields $\sigma$.

In particular, these two pure traces are not in the closure of $\cup_{\lambda >0} F_{\lambda}$, even though the faithful pure traces on $G$ all come from characters. Other (non-degenerate) choices for $f$ presumably yield similar phenomena.

In case we take the analogous $f$ for the free group $\bf F_2$, something even more interesting comes up: $R_f/{\SS R_f}$ is canonically isomorphic to $C(X,\Z)$ (as an ordered abelian group) where $X$ is the path space of the Cayley graph (tree) associated to the group.
\endcomment

\noindent {\it Traces on $\overline R_f$}
 Now we describe the pure traces on $\overline R_f$ (ultimately leading to a description of the pure traces on $R_f$), and their topology.

 The following is a very minor improvement of [GH; 2.1]. 

  \Lem Corollary \oneffn. Suppose that $H $ is an order ideal in $ G$, where both are dimension  groups such  $v$ is an order unit of $H$ and $x \in G^+$ is such that $u = v +x $ is an order unit of $G$. Then $\tau$ can be extended to a trace on $G$ iff 
$$
  \sup \Set{\tau(h)}{0\leq h \leq x; h \in H}: = \beta < \infty, 
 $$
 and when this occurs,  $\beta = \alpha-\tau(v)$. 

 \Pf Suppose $h' \in H^+$ and $h \leq u = v+ x$. Riesz interpolation yields $h' = j+ h$ where $j,h \in G^+$. and $j \leq v$, and $h \leq x$. Since $h \leq h'$, $h \in H$ (as the latter is an order ideal). Obviously $\tau(h) = \tau(h') -\tau (j)$, and thus $\beta \geq \alpha - \tau (v) $.  The reverse inequality is trivial. \qed 

 \Lem Lemma \onetto. Let $(G,u)$ be a dimension group with order unit, and let $I$  be an order ideal thereof. Suppose $\Arrow \tau ; I.\R$ is an (unnormalized) trace of $I$, and let $\brcs{I_{\alpha}}$ be a finite family of order ideals of $G$, each containing $I$, together with (unnormalized) traces $\Arrow \tau_{\alpha}; I_{\alpha}.\R$, each with the property that $\tau_{\alpha}|I = \tau$. Then there exists a trace $\psi$ on the order ideal $\sum_{\alpha} I_{\alpha}$  \st $\psi|I = \tau$. 

 \Rmk It is not generally possible to arrange that $\psi|I_{\alpha} = \tau_{\alpha}$ for all $\alpha$, because there is generally no uniqueness in the extensions of $\tau$ to $\tau_{\alpha}$, even when we assume that all are pure traces. 

 \Pf 
\comment
First, we reduce to the case that $I$ and each $I_{\alpha}$ have their own order units. We can write $I =\cup  I^j$ as the union of an ascending chain order ideals each with order unit, and similarly write each $I_{\alpha} = \cup I_{\alpha}^j$ as a union of order ideals, each having an order unit, \st $I^j \subset I_{\alpha}^j$. 
\endcomment
 Since $\alpha$ runs over a finite set, induction quickly reduces the problem to a set with two elements, $\brcs{\alpha} = \brcs{1,2}$.  Let $v_i$ be an order unit for $I_i$, and $v$  be an order unit for $I$; we may obviously arrange things so that $v \leq v_i$. By  the extension criterion, we have that $\beta(i): = \sup\Set{ \tau(a) }{ a \leq v_i \text{ and } \in I^+}$ are both finite. 

 Now we claim that  $\beta:= \sup \Set{\tau(a)} {a \in I^+ \text{ and } a \leq v_1 + v_2}$ satisfies $\beta \leq \beta_1 + \beta_2$. Suppose $a \in I^+$ and $a \leq v_1 + v_2$. Riesz decomposition yields $0 \leq a_1 \leq v_1$ and $0 \leq a_2 \leq v_2$ \st $a = a_1 + a_2$. Since $a_i \geq 0$ and are dominated by $a$, each $a_i \in I$. Obviously, $\tau(a) = \tau(a_1) + \tau(a_2) \leq \beta(1) + \beta_(2)$. 
 
Now Corollary \oneffn\ applies.\qed

 We can regard $\overline{R}_f$ as an order ideal in a  larger dimension  group, explicitly, let $\G^+$ consist of the elements  of the form $z^r g^a h^b$ where $r,a,b\geq 0$ and there is no constraint on the integer $a$. This is a subsemigroup of the Heisenberg group, and we may take the limit $J= \Z \G^+ \to \Z \G^+$ where the maps are repeated left multiplication by $f = g + h$. This contains the obvious copy $\overline R_f$, as the order ideal generated by $[1,0]$. 

We may also form the order ideal generated by $[z,0]$; it is easy to see that it is precisely $z\overline R_f$ (when viewed inside $J$), and we can do the same wit h $z^{-1}R_f$). In a dimension group, the sum of order ideals is again an order ideal, and the sum $\overline R_f + z\overline R_f$ is the order ideal generated by $\brcs{[1,0],[z,0]}$; we may also do this with the order ideal generated by $\brcs{[1,0], [z,0], [z^{-1},0]}$. We may form the union of the sums, obtaining  $\sum_{r \in \Z}z^r\overline  R_f $, an order ideal in $A_{\overline f}$. [Warning: here $1 \notin \supp \overline f$, so we cannot conclude that $\SS$, the shift, on $A_{\overline f}$, is a bounded endomorphism (it isn't), nor that $\overline R_f$ is stable under it (it also isn't).]

 \Lem Lemma \onettt. Let $\tau $ be a pure trace on $\overline R_f$. Then $\tau$ can be extended to a pure trace on $\sum_{r \in \Z}z^r\overline  R_f $. 

 \Pf Begin with the special case,  $H = \overline R_f$,  $G = \overline R_f + z \overline R_f$, with relative order units $v = [1,0]$ and $w = [z,0]$; we calculate $\beta$ (Corollary \onetty). Consider $(g+h)^m z = \sum p(r,a,m-a) {z^{r+1} g^a h^{m-r}}$. Restricting only to those $w = z^{r+1} g^a h^{m-a}$ that belong to $\Gamma_m''$ (that is, $r+1 \leq a(m-a)$), we set 
$$
\gamma_m = \sum_{w=z^rg^ah^{m-a} \in \Gamma_m'';\  r\geq 1} p(r-1, a,m-a)[ w,m],
$$ 
so that $x_m:= [\gamma_m,m] \in \overline R_f^+$ and $x_m \leq [z,0]$. 

 Now we claim that $\beta = \sup \tau(x_m)$. Select $y \in \overline R_f^+$ \st $y \leq [z,0]$. From the definition of direct limit, there exists an integer $m$ \st $y =\sum (y,w)[{w= z^rg^ah^{m-a}},m]$ with $(y,w) \in A^+$ and $\sum (y,w) z^rg^ah^{m-a} \leq (g+h)^m z$. The latter implies each $y(w) \leq p(r+1,a,m-a)$, hence $x_m \geq y$. 
 
 Now we show that $\beta \leq 2$. Since $p(r-1,a,m-a) \leq 2 p(r,a,m-a)$ (Lemma \onetwo) for $r \geq 1$, we have 
$$\eqalign{
x_m & =   \sum_{w=z^rg^ah^{m-a} \in \Gamma_m'';\ r\geq 1} p(r-1, a,m-a)[ w,m] \cr
 & \leq 2 \sum_{w=z^rg^ah^{m-a} \in \Gamma_m'';\ r\geq 1} p(r, a,m-a)[ w,m]\cr 
 & \leq 2 [\overline f^m,m] = 2[1,0]. \cr 
}$$
 Thus $\tau(x_m) \leq 2 \tau([1,0]) = 2$, so $\beta \leq 2$. 

A similar use of Lemma \onetwo\ yields the corresponding result with $z$  replaced by $z^{-1}$. Now Corollary \onetty\ and Lemma \onetto\ apply. Moreover, multiplication by $z$ is an order automorphism of $A_{\overline f}$, so we can apply the same extension result to a sum of two terms, $z \overline R_{ f} + z^2\overline  R_{f}$, with repeated application. This permits $\tau$ to extend to any finite sum, and thus  to their union. \qed

\noindent {\it Limits of pure discrete traces.}
Now we show that the set of limit points of the discrete traces is precisely the set of multiplicative traces. Let $Z$ be the set of  pure traces whose pure extension to $ \overline R_f + z \overline R_f$   kills the element $[z-1,0]$. The main theorem is that these are given by the maps $\Arrow \psi_t; \overline R_f. \overline R_f$, $[w,a+b] \mapsto  t^a (1-t)^{-b}$, where $w = z^r g^a h^b \in \Gamma_{a+b}''$, for each $t\in [0,1]$;  these are what we call the {\it multiplicative traces.}  Then $Z$ is homeomorphic to a closed interval, and the endpoints are $\tau_{0,0}$ and $\tau^{0,0}$. We will first show that the set of limit points  of $X_0:= \brcs{\tau_{r,b}} \cup \brcs{\tau^{c,d}}$ in $S(\overline R_f, [1,0])$ is precisely $Z$, and it will eventually follow that the subset $ X_0\setminus \brcs{\tau_{0,0}, \tau^{0,0}}$ is discrete in the relative topology (that is, every one of the traces, except the two indicated, are isolated points within  $\partial_e S(\overline R_f, \1)$). 

 Recall the definition of $\tau_{r,b}$. For $w \in \Gamma_m''$,
$$
 \tau_{r,b}([w,m] )= \cases 0 & \text{if $m \geq b+r$ and $w \neq z^r g^{m-r}h^b$}\\ 
 \frac 1{(f^{b+r},z^rg^r h^b)}& \text{if $m \geq b+r$ and $w = z^r g^{m-r} h^b$}\\ 
 \frac{(f^{b+r-m},z^r g^r h^b w^{-1})}{(f^{b+r}, z^r g^r h^b)}& \text{if $m < b+r$.}\\ \endcases
 $$
To check the last case (that is,  in the case that $m < b+r$), we count the number of paths from $w \in \Gamma_m''$ to $z^r g^r h^b$. This is just the number of $u \in \Gamma_{b+r-m}''$ \st $uw = z^r g^r h^b$. If we write $w = z^R g^A h^{m-A}$, then $u = z^r g^r h^bw^{-1}$ simplifies to $u = z^{r-R - A(b+A-m)}g^{r-A}h^{b-m+A}$. The number of such that belong to $\Gamma_{b+r-m}''$ is precisely $p(r-R-A(b+A-m), r-A, b-m+A)$. In order for this to be nonzero, necessary and sufficient is that   $r \geq A$, $b+A \geq M$, and $0\leq r -R- A(b+A-m) \leq (r-A)(b-m+A)$ (implicitly, $A\geq 0$, $m \geq A$, and $0 \leq R \leq (m-A)A$, since $w \in \Gamma_m''$). In any event, the last line of the display can be replaced by 
$$
\tau_{r,b}([z^R g^A h^{m-A},m]) =  \frac{p(r-R-A(b+A-m), r-A, b-m+A)}{p(r,r,b)}. 
$$
Now we notice that if $b -m + A > 0$, and either $R > 0$ or $A > 0$,  then $r-R - A(b-A-m) \leq r-A$. In that case, we can replace the numerator by $p(r-R-A(b+A-m), r-R-A(b+A-m), b-m+A)$. In view of our definitions, this is $p(r-R-A(b+A-m),  b-m+A)$ (with just two parameters). 

We also make a change of notation. We let $(\tau_{r(s),k (s)})$ (initially $s \in \N$) be  an infinite subset of $X_1 = \brcs{\tau_{r,k}}$ which has a limit point in $S(R_f,\1)$ (we do not assume that any limit point is pure; this will come out of what follows). Since $S(R_f,\1)$ is metrizable, there exists an infinite subset $S \subset \N$ together with a trace $\tau$ (normalized, but not known to be pure) \st $\(\tau_{r(s),k (s)}\)_{s\in S}$ converges to $\tau$. Now consider the sequence of real numbers $(k(s)/\sqrt {r(s)})$; there will exist an infinite subsequence of $S$ \st along this subsequence $k(s)/\sqrt{r(s)} \to \alpha \in\R^+ \cup \brcs{\infty}$; so we may assume that this convergence already takes place along $S$.  When convenient, we drop the $s$ in the notation, so that, for example, $k+r \to \infty$ means $k(s) + r(s) \to \infty$ along $S$.

With the new and improved notation, we have 
$$\eqalign{
 \tau_{r(s),k(s)}([z^R g^A h^{m-A},m] ) &= \frac {p(r(s)-R-A(k(s)-(m-A)),  k(s)-(m-A))}{p(r,k)}\cr 
& =  \frac {p(r(s)-R-A(k(s)-(m-A)),  k(s)-(m-A))}{p(r,k-(m-A))}\cdot \frac{p(r,k-(m-A))}{p(r,k)}
}\tag1$$
Sometimes we replace $ m-A$ by  $B$. Then $\tau ([w,m]) = \lim_{s \in S} \tau_{r(s),k(s)}([z^R g^A h^{m-A},m] )$, and we will obtain a limit ratio result (an easy consequence of Szekeres' asymptotic formula), which allows us to calculate the limit. 

 Define $\Arrow v; \R^{++}. \R^{++}$ implicitly by means of the equation $v(x)^2/\int_0^v t/(e^t -1)\,dt = x^2$. It is easy to verify that $v$ is a homeomorphism, and that $v$ extends uniquely to a homeomorphism $\Arrow \tilde v; \R^+ \cup \brcs{\infty}. \R^+ \cup \brcs{\infty}$, by defining $\tilde v(0) =0$ and $\tilde v(\infty) = \infty$. This arises in Szekeres' asymptotic formula for restricted partition functions, which will be recalled soon. 

The map $G \to G/G'$ yields the map on words in $\supp \overline f^m$, $w =z^r g^a h^{m-a} \mapsto x^a y^{m-a}$ and corresponding positive homomorphism $\overline R_f \to R_{x+ y}$. The latter is a ring: set $X = x/(x+y)$ (alternatively, $[x,1]$); then $R_{x+y} = A[X]$ with positive cone generated additively and multiplicatively by $\brcs{X,1-X}$. The pure traces are multiplicative, and given by $\sigma_t (F) = F(t)$ (here $F$ is a polynomial in $X$), for $0 \leq t \leq 1$. The corresponding (pure) trace on $\overline R_{f}$ is given by $\psi_t ([z^r g^a h^{m-a},m]) = t^a (1-t)^{m-a}$ (arising from the composition of maps, $[z^r g^a h^{m-a},m] \mapsto [x^a y^{m-a},m] = X^a (1-X)^{m-a} \to t^a (1-t)^{m-a}$. We refer to the $\psi_t$ as multiplicative traces on $\overline R_f$. 

Note that $\psi_0 = \tau_{0,0}$ and $\psi_1 = \tau^{0,0}$, but none of the other $\tau_{r,b}$ are multiplicative.  Define $X_1 = \brcs{\tau_{r,b}}$, $X_2 = \brcs{\tau^{c,d}}$,  $X_0 = X_1 \cup X_2$, and $Z = \brcs{\psi_t}_{t \in (0,1)}$ (the multiplicative traces). We will prove the following.

 \Lem Theorem \onetth. Let $X_1 = \brcs{\tau_{r,b}}$. 
\item{(i)} If $\tau$ is a limit point of $X_1$, then there exists an infinite subset $S \subset \N$, together with two functions $\Arrow r,s;S.\Z^+$ \st $r(s) + k(s) \to \infty$ and $\tau $ is the limit (along $S$) of the sequence $(\tau_{r(s),k(s)})_{s\in S}$. 
\item{(ii)} If $\tau = \lim_{S} \tau_{r(s),k(s)}$, then $\alpha:=\lim_{s \in S} k(s)/\sqrt{s}$ exists in $\R^+ \cup \brcs{\infty}$, and $\tau = \psi_t$ where $t = \exp(-\tilde v(\alpha))$. 
\item{(iii)}
Suppose that $S$ is an infinite subset of $\N$, and $(r(s),k(s))$ is a sequence of pairs of positive integers \st $r+k \to \infty$ along $S$, and $\lim_S k(s)/\sqrt{ r(s)}$ exists and equals $\alpha \in \R^+ \cup \brcs{\infty}$. Then $\psi:= \lim_s \tau_{r(s),k(s)}$ exists and equals $\psi_t$, where $t = \exp(-\tilde v(\alpha))$. 
\item{(iv)} $\overline X_1 = X_1 \cup Z$; $\overline X_2 = X_2 \cup Z$;  and $\overline X_0 = X_0 \cup Z$. 

 In other words, assuming $k+r \to \infty$, if $k(s) = \oh{r(s)^{1/2}}$, then the limit is $\tau_{0,0}$; if $r(s) = \oh{k(s)^{2}}$, then the limiting trace is $\tau^{0,0}$; and if $k(s)/\sqrt {r(s)}$ converges to finite nonzero $\alpha$, then the limiting trace is given by $[w= z^r g^a h^b,a+b] \mapsto t^a (1-t )^b$ where $t= \tilde v(\alpha)$ for $w \in \Gamma_{a+b}''$; this is $\psi_t$. By symmetry, there is a corresponding result for limits of sequences of the form $(\tau^{r(s),k(s)})$. 

 The following is then an easy consequence (since every limit point is a limit of a sequence in this situation). 

\Lem Corollary \onettf. Let $Z$ denote the set of multiplicative traces on $\overline R_f$, and let $X_1$ and $X_2$ denote respectively, the sets $\brcs{\tau_{r,b}}$ and $\brcs{\tau^{c,d}}$. Then $X_i, Z \subset \partial_e S(\overline R_f,\1)$ and the closure of $X_i$ in $S(\overline R_f,\1)$ is $X_i \cup Z \subset \partial_e S(\overline R_f,\1)$. If $F_i$ is the closed face of $S(\overline R_f,\1))$ generated by $X_i$, then $\partial_e F_i = X_i \cup Z$. 

\Rmk The last part is a general result for Choquet simplices: if $Y$ is a compact subset of the set of extreme points, then the closed face containing $Y$ has no new extreme point (this is part of more general separation results). It is crucial   that all limit points of $X_0$ be extreme, which is a consequence of Theorem\onetth. 
 
\noindent {\it Easy parts of the proof of Theorem \onetth.} (i) $S(\overline R_f, \1)$ is metrizable. 

\noindent (iii) implies (ii). Suppose that $\alpha $ is a limit point (possibly $\infty$) of the sequence (of real numbers with $\infty$) $(k(s)/\sqrt{r(s)})$. Compactness of the unit circle yields a subsequence $S'$ along which the sequence converges to $\alpha$. By (iii), $\lim_{s \in S'} \tau_{r(s), k(s)}$ exists in the trace space and equal $\psi_t$ where $t = \exp (-\tilde v(\alpha))$. Since $ \tilde v$ is a self-homeomorphism of $\R^+ \cup \brcs{\infty}$ and $\psi_t = \psi_{t_0}$ implies $t = t_0$, we have $\alpha$ is the only limit point of $(k(s)/\sqrt{r(s)})$, and so it converges to $\alpha$. 

\noindent (iv) This is an immediate consequence of (i--iii) (where symmetry yields the corresponding result for $X_2 = \brcs{\tau^{c,d}}$). 

\noindent {\it Proof of } (iii).
Suppose that $S$ is an infinite subset of $\Z^+$, and suppose that $\lim_{S} k(s)/\sqrt{r(s)}$ exists and equals $\alpha \in \R^+ \cup \brcs{\infty}$. We have a couple of relatively straightforward cases, and one rather interesting one. 

We have defined $p(r,a,b)$ as the coefficient of $z^r g^a h^b$   in $(g+h)^m$. Now define $p(r,b):= p(r,r,b)$ (so the same letter, $p$, but as a function of only two variables). This is the number of partitions of $r$ with a bound of $b$ on the parts (or the number of partitions of $r$ with no bound, but at most $b$ parts). 

 \noindent {(a)} {\it The sequence $(b(s))$ is  bounded.} Without loss of generality, we may assume that $b(s) = b$ is constant (by taking a further subsequence). If $(r(s))$ contains a bounded infinite subsequence, then of course there is a stationary subsequence, and so the limit point is of the form $\tau_{r,b}$ for some $r,b \in \Z_+^2$. 

 Otherwise, $r(s)  \to \infty$ (in the strong sense, that is, for all $N$, there are only finitely many $s$ in the current version of $S$ \st $r(s)< N$). Now we claim that the limit is $\gamma = \tau_{0,0}$. 

 Pick $w = z^R g^A h^{m-A}\in \Gamma_m''$. For all sufficiently large $s$, $r(s) + b > m$. For such $s$, we evaluate $\tau_{r(s),b}([w,m])$. If $m - A > b$, then there are no  paths from $w$ to $z^{r(s)}g^{r(s)}h^b$, and so in that case, $\tau_{r(s),b}([ w,m]) = 0$. If $m-A = b$, the only possibility for there to be a path (hence a nonzero value of the trace) occurs if $r(s) = R$; but this contradicts $r(s) \to \infty$. In particular, $\tau_{r(s),b}([w,m]) = 0$ whenever $m-A \geq b$. 

 Now suppose that $m - A < b$. Any path---given by $u \in \Gamma_{b-m+r(s))}''$---must contain exactly $b-(m-A)$ occurrences of $h$, and thus exactly $r(s)-A$ occurrences of $g$, and in that case, the exponent of $z$ in $uw$ is at least $R + (b-(m-A))A$. 

 Now the denominator of (1), $p(r(s),r(s),b) = p(r(s),b,r(s)) = p(r(s),b)$, and it is easy to check that the latter grows as asymptotically with a multiple of $r(s)^{b-1}$ if $b \geq 1$. As for the numerator, consider $p(x,y,z)$ where $z = b-m+A \leq b$, $y(s)= r(s) - A \leq< r(s)$, and $x(s) = r(s)-R - A(b+ A-m) \leq  r(s)$. Then $p(x,y,z)  \leq p(x(s), z)$, and this is asymptotic to a multiple of $(x(s))^{z}$ or $z = 0$.
 
  If $b = 0$, then $r(s) = 0$ (since $r \leq rb$), a contradiction. If $b \geq 1$, then $p(r(s),r(s),b) =p(r(s),b) \sim r(s)^b$, and the numerator being nonzero forces   $m=A$. But then $w = g^m$ (since $R \leq A\cdot (m-A)$). Thus the only choices for $w \in \Gamma_m''$ \st $\gamma([e_w,m]) \neq 0$ are those for which $w = g^m$, and then it is trivial to check that this forces $\gamma= \tau_{0,0}$.

 \noindent (b) {\it The sequence $(r(s))$ is bounded.} As before, we can assume that $r(s)$ is eventually stationary, with value $r$, and $b(s) \to \infty$. Then $p(r,r,b(s)) $ is just the unrestricted partition number of $r$ for all sufficiently large $s$. The only way the numerator can be nonzero is if $A = 0$ (because of the $-A(b(s) + A -m)$ term). But then $w = h^m$, and so $\gamma([ w,m]) = 0$ forces $w = h^m$, from which it follows that $\gamma = \tau^{0,0}$. 

\noindent (c) {\it Both sequences $(r(s)), (b(s))$ are unbounded.} By taking a subsequence, we may assume $r(s) \to \infty$ and by taking an infinite subsequence of that subsequence, we reduce to both $r(s), b(s) \to \infty$. We may obviously assume as well that each of $(r(s))$, $(b(s))$ is strictly increasing. 
\comment
 First suppose that $b(s)-m \geq r(s)-R$ for all sufficiently large $s$; we will show that $\gamma = \tau^{0,0}$. When  $A = 0$, then $R = 0$ and  the ratio in (1) is 
$$
 \frac{r(s), r(s),b(s)-m}{p(r(s),r(s),b(s))} = \frac{r(s), r(s)}{p(r(s),r(s))}  = 1 
 $$
 On the other hand, if $A \geq 1$, then in order for the numerator to be nonzero, we require $r(s) - b(s) -1 + m \geq 0$ (set $A = 1$), which contradicts the current hypotheses. Hence $\gamma([e_w,m])$ is nonzero iff $w = h^m$, and thus $\gamma = \tau^{0,0}$. 

 Having dealt with the possibility that $b(s) - m \geq r(s)-R$, we are thus reduced  (by considering subsequences) to  two possibilities, 
$$\eqalign{\frac{r(s)}{b(s)}  & \to \alpha \geq 0\cr
 \frac{r(s)}{b(s)} & \to \infty \cr
 }$$

 The case that $\alpha < 1$ is a subcase of $b(s)-m \geq r(s)-R$, so we obtain convergence to $\tau^{0,0}$ in this case. If $\alpha = 1$, then $|r(s)- b(s)| = \oh{r(s)}$, so that if $A = 1$ the numerator is $p(\oh{r(s)}, r(s)-A, b(s)-\dots) = p(\oh{r(s)}) = \oh{p(r(s))}$, and if $A > 1$, then the numerator is zero. Thus $\gamma([e_w,m]) = 0$ unless $A = 0$, which again forces $\gamma = \tau^{0,0}$.
 \endcomment

 Fortunately, all the hard work has been done by Szekeres [Sz1, Sz2]. This consists of asymptotic estimates for $p(r,k)$, which is  more than we need (we only require results on ratios of   restricted partition functions). We refer to the clear exposition of these results by Canfield [Ca]. We frequently drop the parameter $s$ if $S$ is understood or irrelevant; thus $k(s)/\sqrt{r(s)} \to \alpha$ might be replaced by $k/\sqrt{r} \to \alpha$. 

 \Lem Theorem \onetho. ([Sz2]; adapted from [Ca]) Let $S$ be an infinite subset of $\N$, and let $\Arrow r,k; S.\N$ be two sequences \st $r+ k \to \infty$. Then 
$$\eqalign{
 p(r,k) = \frac{f(k/\sqrt r)}r & \exp\(-\sqrt r g(k/\sqrt r) + \Oh{\frac1{r^{1/6}} \vee \frac 1k}\), \quad{\text{where}} \cr
 f(t)& = \frac{v(t)}{\sqrt 8 \pi t}\(1 - e^{-v(t)}(1+ t^2/2) \)^{-1/2},\cr 
 g(t) & = \frac{2v(t)}{t} - t \ln (1- e^{-v(t)}), \quad\text{and}\cr
v(t)  \text{ is defined implicitly via }\
v(t)^2 & = t^2 \int_0^{v(t)}\frac x{e^x-1}\,dx.
}$$

It is fortunate that all the expressions in (1) are of the form $p(r,r,b) = p(r,b)$, not the more general functions of three variables, $p(r,a,b)$. The latter's asymptotic behaviour is poorly understood (as pointed out in Vershik-Malyutin [VM]), whereas that of the former has been known for almost seven decades. 

All of the following are either easy or in [Ca], or both.

\Lem various factoids \onetht. Here $k(s)$ and $r(s)$ tend to infinity, and $u$, $R$, and $A$ are positive constants. 
 \item{(a)}   $v(t) = t^2(1-t^2/4) + \Oh{t^6}$  
\item{(b)} $\lim_{t\to \infty}\frac{v(t)}t  = \frac{\pi}{\sqrt 6}$
\item{(c)} $f\(\frac{k-u}{\sqrt r}\)/f\(\frac{k}{\sqrt r}\) \to 1$  
\item{(d)} $f\(\frac{k}{\sqrt {r-k}}\)/f\(\frac{k}{\sqrt r}\) \to 1$  (under the assumption that $r (s) - k(s)$ is unbounded on $S$)
\item{(e)} $f\(\frac{k}{\sqrt {r-R}}\)/f\(\frac{k}{\sqrt r}\) \to 1$ 
 \item{(f)} $t^2 g'(t) = tg(t) - 2 v(t)$ 
 \item{(g)} $\sqrt r\cdot (g((k-u)/\sqrt r) - g((k)/\sqrt r))  \to u \ln (1- e^{-\tilde v(\alpha)})$
 \item{(h)} $\sqrt {r-R}\cdot (g(k/\sqrt {r-R}) - g(k/\sqrt r))  \to 0$
\item{(i)}  Provided $r -Ak$ is unbounded, there exists an infinite subsequence   $S'$ of $S$, \st $\sqrt {r- Ak}\cdot (g(k/\sqrt {r-Ak}) - g(k/\sqrt r))  \to -A\tilde v(\alpha)$ along $S'$.

 \Pf (a) [Ca, p\,10]. 

\noindent (b) follows from $\int_0^{\infty} \frac x{e^x-1}\, dx = \frac{\pi^2}6$ [Ca].

 \noindent (c) Obviously, $(k-u)/\sqrt r \to \alpha \in \R^+ \cup \brcs{\infty}$ iff $k/\sqrt r \to \alpha$, so if $\alpha \in \R^{++}$, then continuity of $f$ yields the result. Otherwise, expand 
$$
\frac{f\(\frac{k-u}{\sqrt r}\)}{f\(\frac{k}{\sqrt r}\)} = 
\frac{v((k-u)/\sqrt r)}{v(k/\sqrt r)}\cdot 
\frac{k}{k-u} \cdot
\sqrt{\frac {1- e^{-v(k/\sqrt r)}(1+ k^2/2r)}{1- e^{-v((k-u)/\sqrt r)}(1+( k-u)^2/2r)}}.
  $$
 If $\alpha = 0$, then (a) can be applied, $v(t) = t^2 (1-\oh 1)$, and setting $t = k/\sqrt r$ and $(k-u)/\sqrt r$ yields the result. If $\alpha = \infty$, then $v(t)/t \to c = \pi/\sqrt 6$ as $t \to  \infty$, and the result similarly follows. 

 \noindent (d) If $k/\sqrt r \to \alpha \in \R^{++}$, then $k = \oh{r}$ (as both $k,r \to \infty$, so $k/\sqrt{r-k} \to \alpha$, and continuity of $f$ now applies). If $\alpha = 0$, then again $k = \oh{r}$, and $k/\sqrt{r-k} \to 0$. If $k/\sqrt{r} \to \infty$, then $k/\sqrt{r-k} > k/\sqrt{r} \to \infty$.
 $$
\frac{f\(\frac{k}{\sqrt {r-u}}\)}{f\(\frac{k}{\sqrt r}\)}  
= \frac{v\(k/\sqrt {r-k}\)}{v\(k/\sqrt r\)}\cdot \frac{r-k}{r} 
\cdot \sqrt{
\frac {1- e^{-v(k/\sqrt r)}(1+ k^2/2r)}{1- e^{-v(k/\sqrt {r-k)}}(1+( k^2/2(r-k))}}.
$$
 Suppose  $\alpha = 0$; from  $v(t) = t^2 + \oh {t^4}$, the first quotient goes to $1$, as does the second (as $k = \oh{r}$), and the numerator and denominator of the third quotient behave as $v(t) + \Oh{v^2(t)}$. If $\alpha = \infty$, we group the first two terms, so rewrite their product as $(v(t')/t')/(v(t)/t)$, and this tends to $1$ ($t' \to \infty$). The final quotient obviously tends to $1$. 

 \noindent (e) Follows from (d), since $r-k \leq r - R \leq r$ for almost all $s$. 
  
 \noindent (f) This is   [Ca, top p\, 4]. 

 \noindent (g)  From the derivative \wrt $k$ of $\sqrt r g(k/\sqrt r)$, we have $\Delta_u: = \sqrt r (g((k-u)/\sqrt r) - g((k/\sqrt r)) = -ug'(k^*/\sqrt{r})$ for some $k^*$ with $k-u \leq k^* \leq k$. From (f),  we have $\Delta_u = g(k^*/\sqrt r)/(k^*/\sqrt r) - 2 v(k^*/\sqrt r) ({k^*}^2/r)$, and plugging  this back into the definition of $g$, we obtain $\Delta_u =   u \ln (1- e^{- v(k*/\sqrt r)}$, and continuity of $\tilde v$ yields the result.  

 \noindent (h) The partial  derivative \wrt $r$ of $(k,r) \mapsto r^{1/2}g(k/\sqrt r)$ is 
$$ 
\frac{g(k/\sqrt r)}{2\sqrt r}- \frac{k g'(k/\sqrt r)}{2\sqrt r} = \frac{\frac{k}{\sqrt r} g(k/\sqrt r)/2 -\frac{k^2}{r} g'(k/\sqrt r))}{2k}
$$
 From (f), this simplifies to  $v(k/\sqrt r)/k$. Hence $\Delta^{(R)}:= \sqrt r\cdot (g(k/\sqrt {r-R}) - g(k/\sqrt r)) = -R v(k/\sqrt {r^*})/k$ for some $r-R \leq r^* \leq r$. If $k/\sqrt r \to \alpha \in \R^+$, then the same is true with  $r$ replaced by $r^*$, and since $v$ is continuous and $k \to \infty$, $\Delta^{(R)} \to 0$. If $k/\sqrt r \to \infty$, again, $k/\sqrt {r^*} \to \infty$, and thus $v(k/\sqrt {r^*})/(k/\sqrt {r^*}) \to c$ (a constant), so that $\Delta^{(R)} \sim c/\sqrt {r^*} \to 0$, because $r \to \infty$.   
 
\noindent  (i) Using the computation in (h) for the partial derivative, we have $\Delta^{(Ak)} : = \sqrt {r- Ak}\cdot (g(k/\sqrt {r-Ak}) - g(k/\sqrt r)) $ is $-Ak v(k/\sqrt {r^*})/k = - A v(k/\sqrt {r^*})$ for some $r^*$ with $r-Ak \leq r^* \leq r$. If $r-Ak \to \infty$ along some infinite subsequence $S'$ (that is, it is unbounded), then  $k = \oh { r}$ (along the subsequence), so that $\lim_{s\in S'} k/\sqrt{r- Ak}  = \lim_{s\in S'} k/\sqrt{r}$, and thus    $ \lim_{s\to S'} k/\sqrt{r^*} \to \alpha$. Therefore, $\Delta^{(Ak)} \to - A \tilde v(\alpha)$. \qed 

 Finally, we can prove Theorem \onetth. Recall (more for my benefit than the reader's) that we have assumed $\tau_{r(s), k(s)} \to \tau$, a normalized, but not necessarily pure trace, along an infinite subset $S \subset \N$, and by taking a suitable subsequence, we may assume that $k(s)/\sqrt r(s) \to \alpha \in \R^+\cup \infty$. We wish to show that for $w = z^r g^A h^B \in \Gamma_{A+B}''$, $\tau([w,m]) = e^{-A\tilde v(\alpha)} (1-e^{-\tilde v(\alpha)})^B$, which implies that $\tau = \psi_t$, where $t= e^{-\tilde v(\alpha)}$. We have verified this if $k(s)$ or $r(s)$ is bounded (and $r+k \to \infty$), so are reduced to the situation that both $r,k \to \infty$. Of course, the idea is to calculate the limit of the ratio in (1) using  Szekeres' asymptotic formulas.

 We need a very special case. 

 \Lem Lemma \onethe. If $k/\sqrt r \to \alpha$ and $k,r \to \infty$, and $0 \leq A \leq m$,  where $m$  and $R$ are positive constants, and $A+ R > 0$, and  $\(r(s)- R - A(k(s)-(m-A)\)$ is bounded above on an infinite subset, $S'$, of $S$, then $\alpha = \infty$ and 
$$
 \frac{p(r-R - A(k -( m-A),k)}{p(r,k)} \to 0  \quad \text{along $S'$}.
$$

 \Rmk This is a special case of the ratio tending to $e^{-A \tilde v(\alpha)}$. 

 \Pf If for  infinitely many of $s$, $r(s)- R - A(k(s)-(m-A)) < 0$, then the numerator vanishes. Otherwise, suppose $N \geq 0$, is an upper bound on $S'$. Then $\max p((r-R - A(k -( m-A),k)) \leq \pi (N)$, where $\pi$ (usually denoted $p$, but that would cause confusion) is the unrestricted partition function. Hence the numerator is bounded, but the denominator is not (since both parameters are unbounded and increasing), the limit  of the quotients is zero. \qed
 
 As $r \leq N +R +  A(k-(m-A))$, and $k \to \infty$, $k(s) > m-A = B$ for all but finitely many $s$, and thus  $ k \geq (r-N-R)/A + m-A$, so $k/\sqrt r \geq \sqrt r/A  + (m- A - (N+R)/A)r^{-1/2} \to \infty$. \qed 

 We have by (c) and (g) the right factor in (1), $p(r-B, k)/p(r,k)$ behaves as 
$$
\exp \(-B \ln \( 1- e^{- \tilde v(k/\sqrt r)} + \Oh{1/k \vee r^{-1/7}}\)\).
$$ 
Since $\Arrow \tilde v; \R^+ \cup \infty .  \R^+ \cup \infty$ is continuous and $k,r \to \infty$ (so the big Oh terms can be ignored in the limit), we see that this ratio converges to $(1- e^{-\tilde v(\alpha)})^B$. 

 Now for the other factor, $p(r-R - A(k-B), k-B)/p(r,k-B)$. By Lemma \onethe, we may assume that  $\(r(s)- R - A(k(s)-(m-A)\)$ is unbounded, so there exists a subsequence along which it goes to infinity.  With $k' = k-B$ replacing $k$, the ratio of the  values of the function $f$ tend to one if $A = 1$, by (d) and (e), and then we can iterate this for $A = 2, 3, $ etc. The upshot is that the ratio of the values of $f$ tends to $1$ in all these cases. 

 We can apply (h) with $r' = r-A(k-B)$ and $k' -B$, so the contribution from the $R$ term can be ignored. 

 By (i) applied with $k' = k-B$, we see that the ratio behaves as $\exp\( -A \tilde v(\alpha) + \Oh{1/k, 1/r^{1/7}}\)$, and since both $r,k \to \infty$, this tends to $(e^{-\tilde v(\alpha)})^A$. \qed

 \Lem Proposition \onettw. Let $\tau$ be a pure trace of $\overline R_f$, and let $\tilde \tau$ be a pure trace on $\sum z^r \overline R_f$ that extends $\tau$. 
\item{(a)} If there exists $r_1 > 0$ \st $\tilde \tau ([z^{r_1},0]) = 0$, then $\tau = \tau_{r,b}$ for some $0 \leq r \leq r_1 $ and $b$.
\item{(b)} If there exists $s_1> 0$ \st $\tilde \tau ([z^{-s_1},0]) = 0$, then $\tau = \tau^{s,c}$ for some $ 0 \leq s \leq s_1$ and $c$.

\Pf (a) Obviously $\tau$ kills $\overline R_f \cap z^{r_1} \overline  R_{f}$. Thus for every monomial $w = z^{r'} g^a h^{m-a}$ with $r' \geq r_1$ and  $r' \leq a(m-a)$, we have $\tau([w,m]) = 0$. The function induced on $A \Gamma_m''$  by $\tau$ agrees with the restriction of a positive linear combination of the restrictions of $\tau_{r,b}$ for suitable (finitely many) choices of $(r',b)$ with $r \leq r_1$, call it $\tau_m$.

Let $F_1$ denote the face of $S(\overline R_f, \1)$ generated by $\brcs{\tau_{r,b}}$ (all possible). The preceding paragraph says that the sequence $(\tau_m )$ converges (pointwise; this is the topology on the traces) to $\tau$. But then $\tau$ is an extreme point of $F_1$ (as it is an extreme point of $S(\overline R_f,\1)$ and belongs to the face $F_1$), hence is either multiplicative or of the form $\tau_{r'',b}$, by Corollary \onettf. It can't be multiplicative (since, except for $\tau_{0,0}$ and $\tau^{0,0}$, these are faithful), so it must be of the form $\tau_{r,b}$. Finally, $r \leq r_1$ just from the definitions. 

Part (b) follows from symmetry. \qed

\SecT \FTN\ Order units

Let  $X_0$ denote $\brcs{\tau_{r,b}} \cup \brcs{\tau^{s,c}} \subset \partial_e S(\overline R_f, \1)$. It is a collection of discrete pure traces, and its closure in $S(\overline R_f,\1)$ is a compact subset of $
\partial_e S(\overline R_f,\1)$ given by $X_0 \cup \brcs{\psi_t}_{0 < t < 1}$, denoted $X$.

 Define $\Gd {m}$ 
to be the set of words in $\Gamma_m''$ of the form $z^r g^a h^b$ where either $r \leq a$ or $r \geq (a-1) (m-a)$. It turns out that $\sum_{w \in \Gd m} [w,m]$ is an order unit of $R_f$. This will have a few consequences. Our immediate goal is to prove the following.

 \Lem Lemma \FTNone. Let $M $ be any integer exceeding $2m$. Then 
$$
\bigcup_{u \in \Gd{m}} \Gamma_{M-m}'' \cdot u = \Gamma_{M}''.
$$
 
In particular, if $q = \sum q(w)w$ satisfies  $\Gd m \subseteq \supp q \subset  \Gamma_m''$ and all the coefficients $q(w)$ are nonnegative, then $[q,m]$ is an order unit of $\overline R_f$. This follows because $\supp f^M q = \supp f^{M+m}$ and all coefficients are positive. 

The following is very useful. It says that any element of $\Gd{m}$ has a predecessor in $\Gd{m'}$ for every $m' < m$.

 \Lem Lemma \FTNtwo. Let $m,k$ be positive integers, and $w \in \Gd{m+k}$. Then there exists $u \in \Gamma_k''$ \st $u^{-1}w \in \Gd{m}$. 

 \Pf Alternatively,   the conclusion asserts that  we can factor $w = uv$ where $u \in \Gamma_k''$ and $v \in \Gd  m$. Write $w = z^r g^a h^{m+k-a}$, and suppose that $r \leq a$. If $k = 1$ and $a >r$, then we set $u = g$, while if $a = r$, we set $u = h$ (in this case, $v = g^a h^{m-a}$). 

 At the other end, that is, $(a-1)(m+1-a) \leq r \leq a(m+1-a)$, if the left inequality is strict, then we set $u = h$: in this case, $w$ has unique predecessor, $h^{-1}w$. The remaining case occurs when $ r = (a-1)(m+1-a)$; and set $u = g$. Then $g^{-1}w = z^r g^{a-1}h^{m-a+1}$, and $(a-2) (m+1-a) \leq (a-1)(m+1-a) = r $. 

 Hence the result is true for $k=1$. Now we can obviously continue this by induction. \qed

 \Pf (of Lemma \FTNone) Pick $w = z^R g^A h^{M-A} \in  \Gamma_{M}''$; so $0 \leq A \leq M$ and $0 \leq R \leq A(M-A)$. If $A = 0$ or $M$, then $w =g^M$ or $h^M$, and so $w$ belongs to the union on the left. So we can assume $0 < A < M$. 

\noindent {Case 1: $R \leq A(M-A)/2$.} If  $R\leq A$, then $w \in \Gd{M+m}$  and the lemma  applies. 
\comment
; then $w = g^{A-R}\cdot z^R g^R h^{M-A}$. If now $M-A +R = m$, the right factor is in $\Gd m$. If $M-A +R < m$, we have a similar factorization $w = g^{M-m}\cdot z^R g^{m+A-M}h^{M-A}$ (observe that $R < m+A-M $). If instead, $M-A +R > m$, we are done by the lemma. observe that 
$$\eqalign{w &= g^{A-R }h \cdot h^{-1}z^R g^R h^{M-A} \cr
  &= g^{A-R }h \cdot g^R h^{M-A-1} \cr 
 &= \cases g^{A-R }h g^{M-m-A+R-1} \cdot g^{m+A +1-M} h^{M-A-1} & \text{if $M \leq m+A+1$} \\  g^{A-R }h g^{R}h^{M-A} \cdot  h^{m} & \text{if $M \geq m+A+1$} \endcases
}$$
 \endcomment

 Next, suppose that $A(M-A)/2 \geq R > A$. Set $k = \flo {\frac RA}$, so that $0 \leq R- kA < A$. We note that for any $t$ with $0\leq t\leq k$ \st $M-A-t \geq 0$, $h^{-t} w =z^{R-tA} g^A h^{M-A - t} $, so belongs to $\Gamma_{M-t}''$ (it remains to verify $R-tA \leq A(M-A-t)$, but this is a consequence of $R \leq A(M-A)$), and moreover, if $t =k$, then $h^{-t}w = z^{R-kA}g^A h^{M-A -k}$, so that if $M-A-k \geq 0$ and $M-k \geq m$, we can apply the lemma.

 Now we observe that $M-A-k \geq 0$; if not, then $M-A <k \leq R/A$, yielding $R > A(M-A)$, a contradiction. Hence for all positive $t \leq k$, we have $h^{-t} w \in \Gamma_{M-t}''$. 

 Now we use $M\geq 2m$ (and $R \leq A(M-A)/2$). If $M-k < m$, then $M -m < k \leq R/A \leq (M-A)/2$, and thus $2m > M +A $, a contradiction. So $M-k \geq m$, and the lemma applies. \comment
. If $M-k = m$, we have the factorization $w = h^k \cdot z^{R-kA}g^A h^{m-A}$ (as $k \leq (M-A)/2$, we have $m \geq A$). Otherwise $M- k > m$. If $(k+1)A - R + M-A-k \leq m$, we have the factorization $w = h^k g^{(k+1 )A - R-s} \cdot z^{R-kA} g^{R-kA+ s} h^{M-A - k}$, where $s = m+ R - M- kA +k \geq 0$, and again we are done. 
 
 Otherwise $(k+1)A - R + M-A-k \leq m$ and we have the factorization $w =  h^k g^{(k+1)A-R}  \cdot   z^{R-kA}g^{R-kA} h^{M-A - k}$. This yields $w = h^k g^{(k+1)A-R} h \cdot  g^{R-kA} h^{M-A - k-1}$, and we continue in the earlier reduction to $g^c h^d$.
\endcomment 

 \noindent {Case 2: $R \geq A(M-A)/2$.} Apply the usual involution (reflection) to reduce to Case 1.  \qed 
 
 Recall the definition of $t(m)$ given in section \MI\ (just after Lemma \MIfiv).  

 \Lem Corollary \FTNthr. For   the discrete Heisenberg group and $f = g+h + g^{-1} + h^{-1} + 1$, we have $t(m) = 4(m^2 + m  -3)$  when $m \geq 5$. 

  \Rmk Thus the  growth of $t$ is roughly the square root of the growth   of the group.  

 \Pf For each integer  $m \geq 5$, set 
 $$
 x_m = \sum_{W_m} [w,m]
$$
 where $W_m$ consists of $\Gd {m}$ and its three other iterates under the action of the dihedral group  $D_4$. When $2 \leq a \leq m-a$ and $m \geq 5$, there are $(a+1) + (m-a+1)$ choices of $r$ \st $z^r g^a h^{m-a} \in \Gd m$. If $a = 0$ or $m$, there is exactly one, and if $a = 1 $ or $m-1$, there are $m+1$ choices. Hence $|\Gd m| = (m+1-4)(m+2) + 2 (m+2) = (m-1)(m+2)$. 
 
The overlap between $\Gd m$ and one of its orbits corresponding to an adjacent edge is $\brcs{g^m}$ or $\brcs{h^m}$, so the number of elements in their union is $2(m-1)(m+2)-1$, and the right half of the lozenge overlaps with those terms coming from the left in two elements. Hence $|W_m| = 2(2(m-1)(m+2)-1)-2 = 4(m^2 + m-3)$. 
 
 We verify that   $\alpha:= [x_m,m]$ is an order unit. First, let $\Arrow \pi; R_f .\overline R_f$ be the quotient map (whose kernel is an order ideal). By the previous lemma, $\pi([x_m,m])$ is an order unit (as an element of $\overline R_f$,  since $\pi(x_m)$  has only nonnegative coefficients and $\supp \overline f^m \pi(x_m)= \supp \overline f^{2m} $. The same applies to each of the three other iterates of $\pi$ under the $D_4$ action (since $x_m$ is obviously invariant). 

  If $\alpha$ is not an order unit, then (being an element of $R_f^+$), it generates a proper order ideal, and the quotient by this order ideal admits a pure trace; in particular, there thus exists a pure trace $\tau$ of $(R_f,\1)$ that kills $\alpha$. Then $\ker ^+\tau$ (the largest order ideal in the kernel) contains $\alpha$ and is indecomposable (since $\tau$ is pure).  On the other hand, every trace kills $\SS R_f$, so the intersection of the $\ker \pi$ and its three iterates is killed by $\tau$. Indecomposability of $\ker ^+ \tau$  yields that at least one of $\ker \pi$ or its iterates is contained in $\ker^+ \tau$, and thus $\tau$ induces a trace on at least one of $\overline R_f$ or its iterates that kills the image of $\alpha$---but we have already established that the image of $\alpha$ is an order unit, a contradiction. 

 Thus $t(m)\leq |W_m|$. 

 On the other hand, suppose that $y \in AG^+$ and  $\supp y \subseteq S^m = \supp f^m$, and $[y,m]$ is an order unit. Then it must be positive at all pure traces of $R_f$. If one of the coefficients of $y$ at a word in $W_m$ were zero, then there is a corresponding pure trace that kills $[y,m]$ (to see this, apply an element of $D_4$ so the word corresponds to the upper right quadrant, that is $z^r g^a h^{m-a}$ with $0 \leq a \leq m$ and $0 \leq r \leq a(m-a)$, and either $r \leq a$ or $r \geq (a-1)(m-a)$. In the former case, $\tau_{r,m-a}[y,m]= 0$ and similarly, in the latter case, $[y,m]$ will be killed by a $\tau^{s,c}$. This is contradiction, hence $W_m \subset \supp y$, and thus $t(m) \geq |W_m|$.\qed  
 
 \Lem Corollary \FTNfou. The set of maximal order ideals of $\overline R_f$  is precisely  $  \brcs{\ker \tau_{r,b}} \cup \brcs{\ker \tau^{c,d}}$. 
 
 \Rmk For this class of pure traces $\tau$, $\ker \tau = \ker^+ \tau$. 

 \Pf We already know that each of these kernels is a maximal order ideal. Suppose that $M$ is a maximal order ideal not of the indicated form. Then there exists a pure trace $\tau$ \st $M \subset \ker^+ \tau$ (the quotient, $\overline R_f/M$, admits a pure trace; it yields $\tau$). Since $M$ is maximal, $M = \ker^+ \tau$. 

 We have that $X_0 \subset 
\partial_e S(\overline R_f,\1)$, and we saw that the closure of $X_0$ is $X = X_0 \cup \brcs{\psi_t}_{0\leq t \leq 1}$, in particular, $X \subset 
\partial_e S(\overline R_f,\1)$, and $X$ is compact. Moreover, $\tau $ is a pure trace not in $X$ (the traces $\psi_t$ other than $\tau_{0,0}$ and $\tau^{0,0}$    are faithful---that is, do not kill any element of $\overline R_f^+$). 
 
 We have a standard compactness argument. Given $\phi \in X$, there exists $y_{\phi} \in M^+$ \st $\phi(y ) > 0$ (otherwise, $\ker^+ \tau = M \subset \ker^+ \phi$); if $\phi = \psi_t$, any element of $M^+ \setminus \brcs{0}$ will do. Hence there exists a neighbourhood $U_{\phi}$ of $\phi$ in $
\partial_e S(\overline R_f,\1)$ \st $\hat y_{\phi} | U_{\phi} > \epsilon_{\phi}$ for some real positive number $\epsilon_{\phi}$. 

 Obviously $\brcs{U_{\phi}}$ is a covering of $X$. So there exists a finite subcovering $\brcs{U_i}$, and corresponding $y_i \in M^+$ and $\epsilon_i$ \st $\hat y_i|U_i > \epsilon_i$. Set $ y =\sum y_i$ and $\epsilon = \min y_i$. Then $\hat y |X > \epsilon$  and $y \in M^+$. 

 We may thus write $y = [q,m]$ for some $q \in AG^+$ with $\supp q \subset \overline f^m$ for some $m$. Since all $\tau_{r,b}$ and $\tau^{c,d}$ evaluated at $y$ are positive and $y$ is positive, we have $(q,w) > 0$ (as opposed to $(q,w) = 0$) for all $w \in \Gd{m}$. Hence $y$ is an order unit (by Lemma \FTNone), contradicting $ y \in M$. \qed 

Now we show that there are no pure traces on $\overline R_f$ other than those we know about. 

\Lem Theorem \FTNnin. $\partial_e S(\overline R_f,\1) = X. $

  Let $q \in A^{\Gamma_m''}$, and let $Q \subset \Gamma_m''$. Define $q_Q = \sum_{w \in Q} q(w)w$; in other words, throw away the supporting words that are not in $Q$, and leave the rest of the coefficients untouched. For each $m$,   define $\overline f_m =( f^m)_{\Gd m}$. This is of course $\sum p(r,a,m-a) z^r g^a h^{m-a}$, the sum over the union of the two sets defined by  $ 0 \leq r \leq a \leq m$ and by $(a-1)(m-a) \leq r \leq a(m-a)$ (there is some overlap, and we do not want double counting). 

 \Lem Proposition \FTNfiv. Suppose that $\tau$ is a pure trace of $\overline R_f$. If the sequence $ \(\tau\([f_m, m]\)\)$ does not converge to zero, then $\tau \in X$.  

\Rmk One can think of this sequence of numbers as a sequence of values of measures on $\Gamma_m''$, via the measures defined by $\mu(N(m)) = \tau([f^m_{N(m)},m])$ where $N(m) \subset \Gamma_m''$. 

\Rmk Currently the proof relies on the fact that the limit points of $X_0$ are all extreme. 

\Pf For each $w \in \Gd m$, say $w = z^r g^a h^{m-a}$ with $r \leq a$ (the left half), recall that $\tau_{r,m-a}([w,m])= 1/p(r,r,m-a)  = 1/p(r,a,m-a)$, and similarly with the right half. We partition $\Gd{m}$ into its two pieces and disjointize them, $\Gd{m}_1$ consists of the $w$'s with $r \leq a$ and $\Gd{m}_2 = \Gd{m} \setminus \Gd{m}_1$, define an unnormalized  positive linear combination of the discrete traces as follows. 
$$
 \sigma_m = \sum_{w = z^r g^a h^{m-a}\in \Gd{m}_1} \tau_{r,m-a}\cdot \tau([w,m]) + \sum_{w = z^r g^a h^{m-a}\in \Gd{m}_2} \tau^{a(m-a)-r,a)}\cdot \tau([w,m]) 
$$ 
 We note that $\sigma_m$ is  a trace on $\overline R_f$, and it agrees with $\tau$ when both are restricted to the image of $A^{\Gd{m}}$ in $\overline R_f$. Moreover, on the image of $A^{\Gamma_m''}$, we have that $\sigma_m \leq \tau$. We see that $\sigma_m (\1) = \sigma_m([f^m,m]) = \tau([f^m_m,m])$. So we define $\tau_m$ to be the normalized version, $\tau_m = (1/ \tau([f^m_m,m])) \sigma_m$. 

 By hypothesis, there exists $\delta >0 $ and an infinite subset $W$ of $\N$ \st  $\tau([f^m_m,m]) \geq \delta$ for all $m \in W$. Let $\tau_0$ be a limit point in $S(\overline R_f,\1)$ of $\brcs{\tau_m}_{m \in W}$. We claim that for all $\alpha \in \overline R_f^+$, we have $\tau_0 (\alpha) \leq \tau(\alpha)/\delta$. We may write $\alpha = [q, m]$ where $q \in A^{\Gamma_m''}$ has only positive coefficients, for some $m$. Since we can always replace $[q,m]$ by $[\overline f^k q,m+k]$ for any positive integer $k$ and $W$ is infinite, we can assume that $m \in W$, and there are infinitely many such choices, say $q_m$. For each such, we have $\sigma_m (\alpha) = \sigma([q_m,m]) \leq \tau([q_m,m]) = \tau(\alpha)$. Hence $\tau_m (\alpha) \leq \tau(\alpha)/\delta$. Thus $\tau_0 (\alpha) \leq \tau(\alpha)/\delta$. Hence $\delta \tau_0 \leq \tau$. Since $\tau$ is pure, it must be a scalar multiple of $\tau_0$, and since both are pure, we must have $\tau = \tau_0$. 

 Now let $F$ be the closed face of $S(\overline R_f,\1)$ generated by $X$. Since $X$ is a compact subset of the set of pure traces, we have that $
\partial_e F = X$. The previous paragraph shows that $\tau$, being a limit of convex linear combinations of elements of $X_0$, belongs to $F$. Being pure, it belongs to $
\partial_e F = X$ (since for every closed face of a Choquet simplex, $
\partial_e F = 
\partial_e S(\overline R_f,\1) \cap F$). Hence it is in the closure of $X_0$. But Theorem \onetth\ implies it is one of $X_0 \cup \brcs{\psi_t}_{0 < t < 1}$. It cannot be in the second subset, as an easy exercise show that their values at $[f^m_m,m]$ goes to zero (for each $t$).\qed    

 \Lem Proposition \FTNsix. Suppose $\tau$ is a pure trace of $\overline R_f$ \st $\tau$ kills $z^r \overline R_f \cap \overline R_f$ for some integer $r$. If $r >0$, then $\tau = \tau_{r',b}$ for some $r' < r$, and if $r < 0$, then $\tau = \tau^{s,c}$. 

 \Pf The usual symmetry reduces to the case that $r > 0$. Pick $w = z^s g^a h^{m-a} \in \Gamma_{m}''$ with $a > 0$. We first show that for suitably large $M$ ($M > m^2/2$ is sufficient), $f^M w $ decomposes into a sum of two positive elements of $A^{\Gamma_{M+m}''}$, $f^M w =  q_1 + q_2$ where $\supp q_1 \subset \Gd{M+m}$ and every word in the support of $q_2$ has $z$ exponent at least $r$. If $s \geq r$, the result is trivial, so we may assume $r > s \geq a \geq 0$. 
 
First assume that $a > 0$. Left multiplication of $w$ by any sequence of $g$s and $h$s, with at least $t:=\ceil{(r-s)/a}$ of the latter results in a word whose $z$ exponent is at least as large as $r$. This means that  if a sequence of $M$ $g$s and $h$s results in a word with $z$ exponent less than $r$, then it must have had at most $t$ $h$s, and thus at least $M-t$ $g$s. The resulting word is of the form $z^{r'}g^{a + y} h^{M-a-y + m}$ where $y \geq M-t$. If we can guarantee that $a + y \geq r'$, then the word lies in $\Gd{M+m}$. This will occur if $a + M-t \geq r-1$, that is $M > r-1 - a + (r-s)/a$, and since $ r\leq a(m-a) \leq m^2/4$, sufficient will be $M > m^2/ 2 $. 

 This leaves the case that  $a = 0$. However, $a = 0$ entails $w = h^m$ so $s = 0$. Consider $f^m = h^m + \sum_{w\in \Gamma_m'' \setminus\brcs{h^m}} p(w) w$. Let $U(n)\subset \Gamma_n''$ denote the set of words  $z$-degree less than $r$ in $\Gamma_n'' \setminus \brcs{h^n}$. Then  
$$\eqalign{
 1 = \tau([f^m,m])  & = \tau([h^m,m]) + \tau ([f^m_{U(m)},m] \cr & = \tau([h^m,m]) +  \tau ([f^Mf^m_{U(m)},M+m]) \cr
  & =  \tau([h^m,m]) + \tau ([(f^Mf^m_{U(m)})_{U(M_m)},M+m]) \cr 
 & = \tau([h^m,m]) + \tau ([(f^Mf^m_{U(m)})_{\Gd{M+m}},M+m])  \cr 
 & \leq  \tau([h^m,m]) + \tau ([(f_{M+m},M+m]). 
}$$
  In particular, at least one of $ \tau([h^m,m]) \geq 1/2$ (which entails that $\tau([f_m,m]) \geq 1/2$) or $\tau ([f_{M+m},M+m)  > 1/2$. Since for each $m$, we can choose $M = \ceil{m^2/2} + 1$, it easily follows that there exist infinitely many $n$ \st $\tau([f^n_n,n]) \geq 1/2$, and the result follows from the previous proposition. \qed

We have shown that pure traces $\tau$ on $\overline R_f$ can be extended to pure traces on $\sum_{|i |\leq n} z^i \overline R_f$, and thus to a trace, $\tilde \tau$, on $B = \sum_{i \in \Z} z^i \overline R_f$ whose restriction to any finite sum is pure (after renormalization) or zero. Then each of $\tau \circ \RR_{z^i}$ is a trace on $B$, whose restriction to any finite sum of $z^i \overline R_f$s is pure (after renormalization) or zero. Assume that $\tau \circ \RR_{z^i}| \overline R_f$ is not zero for all $i \in \Z$, and define $\tau_i := (\tilde \tau ([z^i,0]))^{-1} \tau \circ \RR_{z^i}|\overline R_f$. These are pure normalized traces on $R_f$. We will show that $\tau = \tau_i$ for all $i$, and then it will eventually follow that $\tau$ is multiplicative, that is, there exists $t \in (0,1)$ \st $\tau = \psi_t$. 

 \Lem Lemma \FTNsev. If $M \geq (m+1)(m-2)/4$, then $g^M \Gamma_m'' \cup h^M \Gamma_m'' \subset \Gd {M+m}$. 

 \Pf Let $w = z^r g^a h^{m-a} \in \Gamma_m''$. Then $g^M w = z^r g^{M+a} h^{m-a}$, and this belongs to $\Gd{M+m}$ if $M \geq r-a$. As $0 \leq r \leq a(m-a)$, sufficient is that $M \geq a(m-a-1)$, and the latter is bounded above by $(m+1)(m-2)/4$. Similarly $h^M w = z^{r+ aM}g^a h^{M+m-a}$, so sufficient is $r+aM \geq (a-1)(M+m-a)$, that is, $M \geq (a-1)(m-a) - r$; since $r \geq 0$, sufficient for this is also $(m+1)(m-2)/4$. \qed

 \Lem Lemma \FTNeig. Suppose $\tau$ is a pure trace of $R_f$ \st $\tilde \tau ([z^i,0]) \neq 0$ for all $i$. Then $\tau_{\pm1} = \tau$. 
 
 \Pf First, we observe that $\tau([f_m,m]) \to 0$, as $\tau \not\in X_0$. Assume $\tau_1 \neq \tau$. Since $\tau$ is  thus also distinct from $\tau_{-1}$ and all three are pure traces (with the possibility that $\tau_{-1} = \tau_1$, which amounts to $\tau = \tau_2$), for all $\epsilon > 0$, there exists $\alpha \equiv \alpha(\epsilon) \in \overline R_f^+$ \st $\alpha \leq \1$, $\tau(\alpha) > 1 - \epsilon$, and $(\tau_1 + \tau_{-1} (\alpha) < \epsilon$. Let $\lambda = \tilde\tau([z,0])$ and $\lambda' = \tilde\tau([z^{1},0])$; this is positive (not zero), as $\tilde \tau\circ \RR_{z,0} \neq 0$. The second inequality translates to $\tilde \tau (z \alpha ) < \epsilon \lambda$, and if $z^{\pm1} \alpha \in \overline R_f$, then $\tau((z+ z^{-1})\alpha) < \epsilon (\lambda + \lambda')$. We work towards replacing the $\alpha(\epsilon)$ by tractible elements of $\overline R_f$. 
 
 There exists $  m(\epsilon)$ \st for all $t \geq m(\epsilon)$, $\alpha(\epsilon) = [q_t,t]$ (of course, $q_t$ depends on $\epsilon$, but the notation is becoming awkward) where $0 \leq q_t  \leq f^t$ (that is, $(q_t, w) \leq (f^t,w) = p(w)$ for all $w \in \Gamma_t''$, i.e., coordinatewise). Define for $ m \geq m(\epsilon)$,  
$$
 N(\epsilon,m) = \Set{w \in \Gamma_m'' \setminus \Gd{m}}{q(w) > p(w)/2}.
 $$ 
 We can also assume (or redefine $m(\epsilon)$ so) that $\tau([f_m,m]) < \epsilon$ for all $m \geq m(\epsilon)$. Now consider $(f^m)_{N(\epsilon,m)} = \sum_{w \in N(\epsilon)} p(w) w$. The claim is that each of $\alpha_m : = [(f^m)_{N(\epsilon,m)},m]$ (that is, replacing the coefficient of $w$ by $p(w)$ if $w \in N(\epsilon, m)$, and by $0$ otherwise)  can replace our current  $\alpha(\epsilon)$, at a slight cost in multiples of $\epsilon$. 
 
 First, at a cost of at most $\epsilon$ to its value at $\tau$, we can assume $q_m (w) = 0$ for all $w \in \Gd{m}$ (the value at $\tau_1$, small by hypothesis, becomes even smaller or remains the same). Now if $w \in \N(\epsilon)$, replacing $q_m(w) w$ by $p(w)w $ increases (or any rate doesn't decrease the value at $\tau$), but increases the value of $\tau_1$ by $(\tau_1 + \tau_{-1}) ([p(w)- q_m(w) w,m]) < (\tau_1 + \tau_{-1}) ([q_m(w) w,m])$, and thus summing over all $w \in N_{\epsilon,m}$ adds at most $(\tau_1 + \tau_{-1})(\alpha)$ to the value. Finally, for $w$ in the complement of  $N(\epsilon,m)$ and also outside $\Gd{m}$, $U$, $q_m(w) \leq p(w)/2$, we have that $\tau([\sum_U (p(w) - q(w))w,m] \geq  \tau([ \sum_U q(w)w,m])/2 $; but the sum over all of $\Gamma_m''$, $\tau([\sum_{\Gamma_m''} (p(w) - q(w))w,m] =  1 - \tau(\alpha) < \epsilon$. The upshot is that $\tau(\alpha_m) > 1 - 2\epsilon$ and $(\tau_1 + \tau_{-1}) (\alpha_m) < 2 \epsilon$. Now we further refine $N(\epsilon, m)$. 
 
 Define $N'$ (depending on $\epsilon$ and $m$) as $N(\epsilon,m) \cap z N(\epsilon,m)$. We note that 
$$\eqalign{
 \tau\([f^m_{N'},m]\) & \leq  \tau\([f^m_{zN(\epsilon,m)},m]\)  \leq 2\lambda \epsilon \cr  \tau\([f^m_{N'},m]\) &  \leq  \tau\([f^m_{N(\epsilon,m)},m]\)  \leq 2 \epsilon;\cr
 }$$ similar results hold with $z^{-1}N(\epsilon,m)$.
 Setting $N_0 \equiv N_0(\epsilon,m) = N(\epsilon,m)  \setminus \(zN(\epsilon,m) \cup z^{-1}N(\epsilon,m)\)$, we thus have 
$$\eqalign{
 \tau\([f^m_{N_0(\epsilon,m)}]\) &  \geq  1 - (2+\lambda + \lambda') \epsilon\cr 
 \tau_1\([f^m_{N_0(\epsilon,m)}]\) & \leq 1 - 2\epsilon.\cr 
 }$$ 

 Obviously, $\brcs{N_0, z N_0 , z^{-1}N_0}$ are pairwise disjoint. Pick $w \in N_0 (\epsilon, m)$, and let $M$ be a positive integer. We have 
$$\eqalign{
 \tau ([p(w)w,m]) & = p(w)\tau ([f^M w,m+M])  = \sum_{u \in \Gamma_M''}p(u)p(w) \tau([ uw, M+m]) \cr  
 \tau (p(zw)[zw,m]) & = \sum_{u \in \Gamma_M''} p(zw)\tilde \tau ([zuw,m+M]) \cr  & \geq \sum_{u \in \Gamma_M'' \cap z\Gamma_M''}p(zw) p(uz^{-1})  \tau([ uw, M+m]) \cr 
 & \geq \(\frac12\)^2 p(w) \sum p(u) \tau(uw,m+M)  - \sum_{u \in \Gamma_M'' \setminus z\Gamma_M''}p(u)p(w) \tau([ uw, M+m]); \text{ similarly, }\cr 
\tau (p(z^{-1}w)[z^{-1}w,m]) & \geq \(\frac12\)^2 p(w) \sum p(u) \tau(uw,m+M)  - \sum_{u \in \Gamma_M'' \setminus z^{-1}\Gamma_M''}p(u)p(w) \tau([ uw, M+m]) .\cr 
}$$ 
 Summing over $w \in N_0$ and combining the inequalities yields 
$$
 \tau([(z+z^{-1})f^m_{N_0}]) \geq \frac 12 \tau([f^m_{N_0},m] ) - \sum_{w \in N_0}\quad\sum_{u \in  (\Gamma_M'' \setminus z^{-1}\Gamma_M'') \cup (\Gamma_M'' \setminus z\Gamma_M''))}\tau([ uw, M+m]). \tag 3
$$
 
If the inequality $p(uz^{-1}) \geq p(u)/2$ fails, then  $u = g^a h^{M-a}$; in that case, $p(u ) = 1$. Similarly, the only way that   $p(uz) \geq p(u)/2$ fails is if $u = z^{a(M-a)} g^a h^{M-a}$. In particular, if $u \in z\Gamma_M'' \setminus z^{-1}\Gamma_M''$, then it is included in one of the sums but not the other. Explicitly, noting that $p(u) = 1$ if $u \in \(\Gamma_M'' \setminus z\Gamma_M''\) \cup \(\Gamma_M'' \setminus z^{-1}\Gamma_M'\),$
$$\eqalign{
 \sum_{w \in N_0} \sum_{u \in \(\Gamma_M'' \setminus z \Gamma_M''\)\cap z^{-1}\Gamma_M''} p(u)p(w) \tau([ uw, M+m]) 
 & = \sum_{w \in N_0} \sum_{u \in \(\Gamma_M'' \setminus z \Gamma_M''\)\cap z^{-1}\Gamma_M''} p(w) \tau([ uw, M+m])\cr 
 &\leq \sum_{w \in N_0} \sum_{u \in  z^{-1}\Gamma_M''}  p(w) \tau([ uw, M+m]) \cr 
 & = \sum_{w \in N_0} \sum_{v \in  \Gamma_M''} p(w) \tau([ zvw, M+m]) \cr 
& \leq 2 \sum_{w \in N_0} \sum_{v \in  \Gamma_M''} p(wz) \tau([ vzw, M+m]) \cr 
 & \leq 2 \lambda\tau_1 ([f^M (f^m_{N_0},m]) < 2 \lambda \epsilon. 
  }$$ 
 Similarly, 
$$
  \sum_{w \in N_0} \sum_{u \in \(\Gamma_M'' \setminus z^{-1} \Gamma_M''\)\cap z\Gamma_m''} p(u)p(w) \tau([ uw, M+m])  \leq 2 \lambda' \epsilon.
$$
 This (so far) is true for all $M$, but now we assume that $M \geq (m+1)(m-2)/4$. Then $\tau([f^{M+m}_{g^M N_0 \cup h^M N_0}]) < \epsilon$, by preceding lemma. However, $(\Gamma_M'' \setminus z\Gamma_M'') \cap (\Gamma_M'' \setminus z^{-1}\Gamma_M'')  =\brcs{g^M, h^M}$.  Hence the term $$
\sum_{w \in N_0}\sum_{u \in  (\Gamma_M'' \setminus z^{-1}\Gamma_M'') \cup (\Gamma_M'' \setminus z\Gamma_M''))}\tau([ uw, M+m])
$$ 
appearing in (2) is at most $5\epsilon$. Putting this in (3), we obtain $\epsilon \geq (1- 2(\lambda + \lambda')\epsilon)/2 - 5 \epsilon$, so if $\epsilon$ is sufficiently small, we obtain a contradiction. \qed 
 
 \Lem Corollary \FTNten. If $\tau$ is a pure trace \st $\tau \circ \RR_{z^i,0} \neq 0$ for all integers $i$, then for all $\alpha \in R_f$ \st for some integer $k$, $z^k \alpha \in R_f$, we have $\tau (z^k \alpha) = \lambda^k\tau(\alpha)$, $(\lambda')^{|k|}\tau(\alpha)$ respectively, if $k > 0$ and $k < 0$, where $\lambda = \tilde \tau ([z,0])$ and $\lambda' = \tilde \tau ([z^{-1},0]) = \lambda^{-1}$. 

 \Pf Suppose  to begin with that $k > 0$. As $\alpha, z^k \alpha \in R_f$, it easily follows that $z^i \alpha \in R_f$ for all $0 \leq i \leq k$. Then $\tau(z^k \alpha)= \tau \circ \RR_{z,0} (z^{k-1}\alpha) = \lambda \tau_1(z^{-1}\alpha)  = \lambda \tau(z^{k-1}\alpha)$, and we can continue by induction. Similarly, if $k < 0$, the same result holds but with $\lambda' $ replacing $\lambda$ and $\tau_{-1}$ replacing $\tau_1$. 

 Consider $[hg, 2] = [zgh,2]$. Then $\tau ([zgh,2]) = \lambda \tau  ([gh,2]) = \lambda \lambda' \tau([gh,2]) = \lambda \lambda' \tau_{-1} ([z gh,2]) = \lambda \lambda' \tau_{-1} ([z gh,2])$. Hence if $\tau[zgh,2]) \neq 0$, we deduce $\lambda \lambda' =1$ as desired. Otherwise, $\tau([hg,2])= \tau([zgh,2])  = 0$, and this also entails $\tau ([gh,2]) = 0$. We need two {\it Zwischenzuege\/} (Zwischenzug has been absorbed into English, so doesn't require italics, but the plural has not, so requires italics)  to complete this argument. 

\Lem Zwischenzug 1. Suppose that $v$ is a group element, and there exist  integers $M$ and $m$ \st $\supp \overline f^M v \subseteq \supp \overline f^{M+m}$. Then $v \in \supp \overline f^m = \Gamma_m''$. 

\Rmk So $\overline R_f$ satisfies the analogue of WC, whereas $R_f$ does not. 

\Pf Write $v = z^r g^a h^b$ in normal form; here $r,a,b$ can vary over all integers. If $r < 0$ or $b < 0$, then $g^M v \not\in \cup_{j \geq 1} \Gamma_j''$ for any positive $M$. Thus $r,b \geq 0$. Similarly, if $a < 0$, $h^M v$ yields a contradiction, and so $a \geq 0$. Finally, if $r > a(m-a)$, hit it (from the left) with $h^M$; the outcome is $z^{r + aM} g^a h^{b+ M} $. But now $r+aM - a(b+M) = r - ab < 0$, contradicting $h^M v \in \supp \overline{f}^{M+m}$.\qed

\Lem Zwischenzug 2. The intersection of order ideals of $\overline R_f$, $\langle [g,1] \rangle \cap \langle  [h,1]\rangle  $, is contained in $\langle [gh,2] \rangle \cap \langle  [hg,2]\rangle$.

\Rmk There are many similar such inclusions, and these can be treated much more generally. But there is no immediate need for them, except this particular one. 

\Pf For $w \in \Gamma_m''$ and $v \in \Gamma_n''$, we  have $[w,m] \in \langle [v,n]\rangle$ iff there exists a positive integer $K$ \st $[w,m] \leq K[v,n]$, and this is equivalent to $\overline f^{M+n} w \leq K \overline f^{M+m} v$ (coordinatewise) for some $M$, which in turn is equivalent to $\supp f^{M+n} w \subseteq \supp f^{M+m}v$, which we can rewrite as $\supp \overline f^{M+n}wv^{-1} \subseteq \overline f^{M+m}$. If $n=1$, by Zwischenzug 1, this entails $ wv^{-1} \in \supp \overline f^{m-1}$. Setting $v = g$ and then $h$, we have that $wg^{-1}, wh^{-1} \in \Gamma_{m-1}''$. 

Writing $w = z^r g^a h^{m-a}$ (with the usual conditions, $0 \leq a \leq m$ and $0 \leq r \leq a(m-a)$, we have $v g^{-1} = z^{r-(m-a)}g^{a-1}h^{m-a}$   and  $v h^{-1} = z^{r} g^{a} h^{m-a-1}$, and both these have to be in $\Gamma_{m-1}''$. This entails $a, m-a \geq 1$ and $m-a \leq r \leq a(m-a-1)$. Thus we can rewrite $w =z^{r-(m-a)}g^{a-1}h^{m-a-1} \cdot gh$, and we verify that $z^{r-( m-a)}g^{a-1}h^{m-a=1} \in \Gamma_{m-2}''$ from the inequalities. Thus $[w,m ] \leq [gh,2]$, so $[w,m] \in \langle [gh,2] \rangle$.\qed

\noindent {\it Rest of proof of Corollary \FTNten.}  If $\tau([gh,2]) = 0$, then $\ker^+ \tau$ contains $[gh,2]$.   As $\tau$ is pure, the order ideal $\ker^+ \tau$ is indecomposable, and thus at  least one of $\tau([g,1])$ or $\tau([h,1])$ is zero, by Zzug\, 2.  This entails  one of $\tau \circ \RR_{z^{\pm1},0}$ is zero, a contradiction. Hence $\tau([gh,2]) \neq 0$, and so $\lambda' = 1/\lambda$. \qed

 \noindent {\it Proof of theorem \FTNnin.} Let $\tau$ be a pure trace not in $X$.  First we show  that $\lambda = 1$,  and then it easily follows that both $\tau \circ \RR_{g,1} = \tau([g,1])\cdot \tau$ and  $\tau \circ \RR_{h,1} = \tau([h,1])\cdot \tau$. 

 We first note that $\tau_g:= \tau \circ \RR_{g,1}$ and $\tau_h:= \tau \circ \RR_{h,1}$ are nonzero (as in the proof above,  $\tau[g,1] = 0$ or $\tau[h,1] = 0$   implies one of $\tau\circ \RR_{z^{\pm1},0} = 0$), and thus both are (unnormalized, and not necessarily pure) traces. 

 Now we wish to show that if $\lambda \neq 0$, then $\lambda = 1$.  Assume $0 < \lambda < 1$.  Let $\pi(r)$ denote the number of (unrestricted) partitions of $r$ (this is usually denoted $p(r)$, but we are already using $p$; nor is $\pi$ likely to be confused in this context with the prime number counter). For fixed $a,m$, the sum $\sum_{0 \leq r \leq (m-a)a} \lambda^r p(r,a,m-a)$ is obviously dominated by $\sum_{r=0}^{\infty} \lambda^r \pi(r)$. As $\pi$ grows only subexponentially, the infinite sum  converges  (since $\lambda < 1$), say to $T(\lambda)$. 
\comment
, and in fact, it is easy to see that if $r_0 \geq   2\pi/3(\ln \lambda)^2$, then there exists $\lambda'$ with $\lambda < \lambda' < 1$ \st $\sum_{r \geq r_0} \lambda^r p(r) \leq (\lambda')^{r_0}/(1-\lambda') $.
\endcomment

  Thus for every $m$, 
 $$\eqalign{
 1 = \tau([f^m,m])& = \sum_{a=0}^m\sum_{0\leq r \leq a(m-a)} p(r,a, m-a) \tau([z^r g^a h^{m-a},m]) \cr
& =\sum_{a=0}^m \tau([g^a h^{m-a},m])  \sum_{0\leq r \leq a(m-a)}\lambda^r p(r,a,m-a)\cr 
&<\sum_{a=0}^m \tau([g^a h^{m-a},m])  \sum_{r=0}^{\infty}\lambda^r \pi(r)\cr 
&= \sum_{a=0}^m \tau([g^a h^{m-a},m])  T(\lambda). \cr
 }$$
Since $\brcs{g^a h^{m-a}} \subset \Gd m$, we conclude that $\tau ([f_m],m) > 1/T(\lambda)$ for all $m$. By Proposition \FTNfiv, $\tau \in X_0$,  a contradiction. Hence $\lambda \leq 1$ entails $\lambda = 1$. 

Now assume that $\lambda > 1$. We have a slightly different argument, based on the same idea, using the observation that  $p(r,a,m-a) = p(a( m-a)-r, a, m-a)$. The substitution used in line three is $s = a(m-a) -r $.
$$\eqalign{
 1 = \tau([f^m,m])& = \sum_{a=0}^m\sum_{0\leq r \leq a(m-a)} p(r,a, m-a) \tau([z^r g^a h^{m-a},m]) \cr
& =\sum_{a=0}^m \tau([z^{a(m-a)}g^a h^{m-a},m])  \sum_{0\leq r \leq a(m-a)}\lambda^{r-a(m-a)} p(a( m-a)-r, a, m-a)\cr 
& = \sum_{a=0}^m \tau([z^{a(m-a)}g^a h^{m-a},m])\sum_{0\leq s \leq a(m-a)} \frac 1{\lambda^s}p(s,a,m-a)\cr
&<\sum_{a=0}^m \tau([z^{a(m-a)}g^a h^{m-a},m])  \sum_{s=0}^{\infty}\frac 1{\lambda^s}\pi(s)\cr 
&= \sum_{a=0}^m \tau([z^{a(m-a)}g^a h^{m-a},m])  T(\lambda^{-1}). \cr
 }$$
As $\brcs{z^{a(m-a)}g^a h ^{m-a}} \subset \Gd{m}$, once again we have that the sequence $\(\tau([f_m,m])\)$ is bounded below  away from zero. This is again a contradiction, so that $\lambda \geq 1$ entails $\lambda = 1$. Thus $\lambda = 1$ in either case. 
 
Now consider $\tau_g:= \tau \circ \RR_{g,1}$ and $\tau_h:= \tau \circ \RR_{h,1}$. Each of these, if not zero,  is a trace (unnormalized, and not necessarily pure). We observed earlier  that neither $\tau([g,1])$ nor $\tau ([h,1])$ is zero; as these are the values of $\tau_g$ and $\tau_h$ at $\1$, we  conclude that both $\tau_g$ and $\tau_h$ are traces. 

For  $w= z^r g^a h^{m-a} \in \Gamma_m''$, have $wg = z^{r+ m-a}g^{a+1} h^{m-a}$, and $\tau_g ([w,m])  = \tau([g^{a+1} h^{m-a},m +1])$ (since $\lambda = 1$, and this agrees with $\tau([gw,m+1])$ (again since $\lambda = 1$). As $[gw,m+1] \leq [w,m+1]$ (straight from the definitions), we have $\tau_g \leq \tau$ as traces (it is enough to check inequality on all terms of the form $[w,m]$ for all $w \in \Gamma_m''$, for all sufficiently large $m$).

 As $\tau$ is pure, this forces $\tau_g = \gamma \tau$ for some positive scalar $\gamma$. Evaluating at $[1,0]$, we obtain $\gamma = \tau([g,1])$. Hence for all $q \in A^{\Gamma_m''}$, we have $\tau([qg,m+1]) = \tau([g,1])\tau([q,m])$. 
 
 Almost the same argument works for $\tau_h$.  
With the same $w$, $\tau_h ([w,m]) = \tau([z^r g^a h^{m+1-a},m+1]) = \tau([z^{a}g^a h^{m-a+1},m+1]) =   \tau([hg^a h^{m},m+1])  \leq   \tau([g^a h^{m},m])  = \tau([w,m])$. Thus $\tau_h \leq \tau$, so by purity of the latter, $\tau_h$ is a nonzero scalar multiple of $\tau$, and evaluating at $[1,0]$, we obtain $\tau_h = \tau([h,1])\cdot \tau$.

 This immediately yields $\tau ([z^r g^a h^{m-a},m]) = \tau([g^ah^{m-a-1},m])\tau ([h,1])$ (except if $m =a$), and inductively $\tau (z^r g^a h^{m-a},m]) = \tau([g,1])^a \tau([h,1])^b$. Set $t = \tau([g,1])$; as $[g+h,1] = \1$, we have $ 0 < t < 1$ and $\tau([h,1]) = 1-t$. So $\tau ([z^r g^a h^{m-a},m]) = t^a (1-t)^{m-a} $, that is $\tau = \psi_t$, a contradiction. (gasp)\qed  

We have a few consequences. 

\item{(1)} The trace space of $\overline R_f$, $S(\overline R_f,\1)$, is a Bauer simplex, with $\partial_e S(\overline R_f,\1) = X$. 
\item{(2)}Each point in $X_0 \setminus \brcs{\tau_{0,0},  \tau^{0,0}}$ is an isolated point of $\partial_e S(\overline R_f,\1)$.
\item{(3)} The pure trace space of $(R_f, [1,0])$ (note $R_f$, not $\overline R_f$) consists of the following.
\itemitem{(a)} the faithful pure traces; these extend to $A_f$, and are given by $\phi_s : [z^r g^a h^b, a+b] \mapsto s^a (1-s)^b$ for all $ (r,a,b ) \in \Z^3$ for each $0 < s < 1$.  These factor through the map to the abelianization, and also through $A_f \to A_P$, where $P = x+y+1 + x^{-1} + y^{-1} \in A[\Z^2]^+$. These may be identified with the interior of the Newton polytope of $P$, via the weighted moment map. 
\itemitem{(b)} non-faithful traces factoring through $R_f \to \overline R_f$ and their three other images under the action of $D_4$. These subdivide into the discrete traces on one hand, and their limits  on the other, the latter given by $\psi_t : [z^r g^a h^b,a+b] \mapsto t^a (1-t)^b$ and their $D_4$ images. The limit points correspond to the boundary of the Newton polytope, one family for each edge.

\comment
\hrule
\vskip10pt
 Other examples (comment, get integrally simplicial, get SSWC, and also $l|H$ is multiplicative on powers (i.e., $l(mx) = ml(x)$ if the group is written additively).

\SecT Appendix 3 Failures of SWC
\def\ul{{\underline l}}

Here  $G =\Z^2 \times_{\theta} C_2$, where $\theta (x) = -x$; denote the standard copy of $\Z^2$  by $H$.  We show that it fails to satisfy SWC, even though its centre is trivial. 

 We use additive notation for elements of $H$. Suppose $(l_0, W, C)$ implements SWC. By xxx, we may assume $W = \brcs{1,g}$. Let $l = l_0 |H$. Define  $\Arrow \ul ; H.\R^+$ via 
$$
 \ul (h) = \lim_{m} \frac{l(mh)}{m}
$$
 This exists, since $l$ is subadditive. Also, $\ul(h+ h') \leq \ul (h) + \ul (h')$. 

 \item{(1)} For $h \in H$, $l(h) \geq \ul (h) \geq l(h) - C$. 

 \Pf The left inequality is trivial. From $W = \brcs{1,g}$, we have $\max \brcs{l(2h), l(h-h)} \geq 2 l(h) - C$. Thus $l(2h) \geq 2l(h) - C$. By induction, $l(2^k h) \geq 2^k l(h) - (2^k-1) c$. Dividing by $2^k$, and taking limits, we deduce $\ul (h) \geq l(h) - c$. \qed
 
 \item{(2)} For all $t \in \R^+$, $\ul^{-1} (\leq t)$ is finite. 

 \Pf By (1), $\ul^{-1} (\leq t) \subset l^{-1}(\leq t + C)$, and the latter is finite. \qed

 \item{3} If $x,y \in H$ and $\ul (x+y) \geq \ul (x-y)$, then $\ul (x+y) = \ul (x) + \ul (y)$. 

 From $W = \brcs{1,g}$, for all $h,h' \in H$, 
$$
\max \brcs{l(h+h'), l(h-h')} \geq l(h) + l(h') - C \tag *
$$
  Given $\epsilon > 0$, there exists $N$ \st $n \geq N$ implies 
 $$\eqalign{
 \frac{l(2^n(x+y))}{2^n} - \frac{l(2^n(x-y))}{2^n}  &> -\epsilon, \text{ that is,}\cr
 l(2^n(x+y)) &\geq l(2^n(x-y)) - 2^n \epsilon\cr
}$$
 Applying (*) with $h = 2^n x$ and $h' = 2^n y$, we have $l(2^n (x+y)) + 2^n \epsilon \geq l(2^n x) + l(2^n y ) - C$. Dividing by $2^n$, and taking the limit, we have $\ul  (x+y) \geq \ul(x) + \ul (y) - \epsilon$. As this is true for all $\epsilon$, $\ul (x+ y) \geq \ul (x) +\ul (y)$. The reverse inequality is trivial. \qed 

 \item{(4)} For all $h,h' \in H$, $\max\brcs{\ul(h+h'),\ul (h-h')} \geq \ul (h) + \ul (h') - 4C$. 

 \Pf We have, 
 $$\eqalign{
 \max\brcs{\ul (h + h'), \ul (h-h')}
& \geq 
\max \brcs{l(h+h') - C ,l(h-h') - C }\cr 
 & \geq \max \brcs{l(h+h')  ,l(h-h')  } - C\cr 
 & \geq  l(h)  + l(h') - 2C\cr
 & \geq \ul (h) + \ul (h') - 4C 
 }$$

 \item{(5)} $G$ fails SWC. 

 Pick $x,y \in H$ \st $\brcs{x,y}$ is rationally linearly independent. Let $M,N$ be positive integers; by replacing $y$ by $-y$ if necessary, we may assume that $\ul (Mx +Ny) \geq \ul(Mx -Ny)$. 

 Thus $\ul (Mx + Ny) = \ul(Mx) + \ul (Ny)$. Moreover, we also have $\ul (kh) = k \ul (h)$ for all nonnegative integers $k$  and elements $h \in H$ (this follows from (3) applied inductively. 
 We also have $\max \brcs{\ul (2Mx), \ul (2Ny)} \geq \ul (Mx + Ny) + \ul (Mx-Ny) - 4C$. Thus 
 $$\eqalign{
 \ul (Mx-Ny) & \leq 2 \max \brcs{\ul (Mx),\ul (Ny)} - \ul (Mx +Ny) + 4C\cr
 & 2 \max \brcs{\ul (Mx),\ul (Ny)} - \ul (Mx) - \ul (Ny) + 5C\cr
 & = |\ul (Mx) - \ul (Ny)| + 5C\cr
 }$$
 
Set $x = {}1$ and $y = {}2$ in $\Z^2$, and set $t = \ul (x)$ and $u = \ul (y)$. By replacing $y $ by $-{}2$ if necessary. Now consider $t/u$. We may find infinitely many pairs of integers $(p,q)$ \st $|t/u - p/q| < 5/q^2$; thus, $|tq - pu| < 5/q$. In particular, $|l(qx) - l(py)| < 5/q$. Hence there are infinitely many choices for $M$ and $N$ \st $\ul (Mx -Ny) \leq 5C + 1$, contradicting (2). \qed

comments at end of appendix dealing with crossed product of abelian groups by finite groups.

 The argument in showing that the crossed product satisfies SWC is essentially along the following lines. There exists a subsemigroup $M$ of $H$ \st $M-M = H$,  there is a weight function on $H$ whose restriction to $M$ is additive and every element of $H$ can be moved to $M$ by the action of an element of $K$. This allows us to extend the weight function on $H$ to a weight function on $G$ with the desired properties. 

 One can ask therefore, for  an arbitrary (infinite, finitely generated, discrete) group $G$ satisfying SWC, whether there exists a subsemigroup $M_0$ of $G$ and a weight function on  $G$ that is almost subadditive on $M_0$, and $M_0 M_0^{-1} = M_0^{-1} M_0 = G$ (and we also require something analogous to the existence of a finite set of elements at least one of which can conjugate any element into $M_0$). Unfortunately, this is practically the Ore condition for groups, which would force $G$ to be almost nilpotent.  So we weaken the conditions, and ask for a family of subsemigroups $M_{n}$ whose union is $G$, the restriction of the weight function to each one is almost subadditive, and there is a finite set of elements $W$, etc.

 The definition of SWC yields a sort of locally subadditive weight function, but the presence of suitable subsemigroups does not seem to follow directly. 

 As to the converse of result (no multiplicities, etc), one part is obvious: if the trivial subrepresentation appears in $\theta \otimes 1_{\Q}$, then $\theta$ admits a nonzero fixed point, and so the centre of the crossed product is infinite, preventing SWC (or even SWC(0)). If multiplicities appear in $\theta \otimes 1_{\Q}$, the result (that SWC does not hold) is not so clear, at least in part because we do not know that much about the weight functions associated to $W$. Moreover, in at least some cases, e.g., for $p$ a prime, with  $C_p$ acting on $\Z^{p-1} \times \Z^{p-1}$ via two copies of the (essentially) irreducible action of $C_p$ on $\Z^{p-1}$, the crossed product is SWC(0). To check whether  SWC fails for $G = (\Z^{p-1} \times \Z^{p-1}) \times_{\theta} C_p$, we may assume (by xxx),  that $W = C_p $. For the natural choice of weight function ($l((x_i)) = \sum |x_i|$) and $p = 2 $ or $3$, there is a brute force argument that $(l,W)$ does not implement SWC or even any SWC($\alpha$) with  $\alpha >0$. However, this does not preclude the possibility that some other choice of weight function will yield SWC. 

\endcomment


\SecT Appendix 1 Gauging weight functions

 Let $l$ be a weight function on a group $G$. As is well known (Black's theorem), for all $g$ in $G$, the limit $\lim_{n\to \infty}l(g^n)/n$ exists; define $\Arrow \hat l ; G.\R^+$ via $g \mapsto \lim_{n\to \infty}l(g^n)/n$. Then $\hat l$ is subadditive and $\hat l (1) = 0$. However, for most groups, $\hat l^{-1}(\leq r)$ is infinite for some nonnegative real number $r$---so $\hat l$ need not be a real weight function as defined in section 3, just after Corollary \Xsev. 

 However, when $G = \Z^d$ and $l = l_S$ for some admissible $S \subset \Z^d$, there is a very nice description of $\hat l_S$, given as the gauge of the compact polytope $K = \cvx S$ restricted to $\Z^d$. In particular, there exists a positive integer $N$ \st $N\hat l_S$ is a weight function.

\Lem Lemma \Yone. Let $l$ be a weight function on a (nonabelian) group $G$. Then $\hat l$ is constant on conjugacy classes. 

\Pf If $g,h,x \in G$ and $g = xhx^{-1}$, then $|l(g) - l(h)| \leq l(x) + l(x^{-1})$. Hence if $a = x bx^{-1}$, then $a^n = x b^n x^{-1}$, so that $|l(a^n) - l(b^n)| \leq l(x) + l(x^{-1})$. Dividing by $n$ and taking limits yields $\hat l(a) = \hat l(b)$. \qed 

 \Lem Corollary \Ytwo. Let $l$ be a weight function on $G$. If $G$ is not central by finite, then there exists real $r$ \st $\hat l^{-1}(r)$ is infinite. 

 \Pf For finitely generated groups, FC is equivalent to central by finite. Hence there exists $a$ in $G$ having infinitely many conjugates, and now the preceding applies. \qed 

For the rest of this appendix, we will be discussing weight functions on $G = \Z^d$, and as a result, use additive notation (thus subadditivity is $l(g + h) \leq l(g) + l(h)$). 

 Let $L \subset \R^d$ be a compact convex set with interior (called a compact convex {\it body\/} in some references); assume in addition that $0$ belongs to the interior of $L$. Then the {\it gauge\/} associated to $L$ is the function $\Arrow \Lambda_L; \R^d.\R^+$ given by 
$$
\Lambda_L (x) = \inf \Set{\lambda > 0}{\frac x{\lambda} \in L}.  
 $$
 As is well known (and mostly obvious), for all $x,y \in \R^d$ and $r \in \R^{++}$, 
 $$\eqalign{\Lambda_L (x)  & = 0  \text{ if and only if $x = 0$}\cr 
 \Lambda_L(x+ y) &\leq \Lambda_L (x) + \Lambda_L (y)\cr 
 \Lambda_L (rx) & = r \Lambda_L (x). 
 }$$
 Moreover, $\Lambda_L$ is a norm if and only if $L = - L$. 

 As usual, $\partial L$ will denote the boundary of $L$ (that is, $L \setminus \text{int}(L)$), and $\partial_e L$ will denote the set of extreme points of $L$. The following is trivial.

 \Lem Lemma \Ythr. For $L$ a compact convex body in $\R^d$ and nonzero $x \in \R^d$, we have $\Lambda_L (x) = \lambda_0$ where $\lambda_0$ is the unique nonnegative number \st $\lambda_0 x \in \partial L$. 

 Now specialize $L$ to be a compact convex polytope in $\R^d$ with the origin in its interior. Then $\partial L$ is just the union of the facets of $L$.  For each facet $F$ there exists unique $v_F \in \R^{1\times d}$ (the normalized outward normal) \st $v_Ff = 1$ for all $f \in F$, and $v_F(x ) \leq 1$
for all $x \in L$. Then $F = v_{F}^{-1}(1) \cap L$;  also, $L = \cap_{F} v_F^{-1} (\leq 1)$ where $F$ varies over all the facets of $L$. 

 Define (for each facet $F$)  $C_F $ to be the convex hull of $\brcs{0} \cup F$. The following is straightforward. 
 
\Lem Lemma \Yfou. Let $L$ be a compact convex polytope having the origin as an interior point. \item{(i)} $L = \cup_F C_F$ where $F$ varies over the facets of $L$. 
 \item{(ii)} for all $x \in \R^d$, we have $\Lambda_L (x) = \max_F v_Fx$, and if $x \in C_{F_0}$ for some facet $F_0$, then $\Lambda_L (x) = v_{F_0}x$. 
\item{(iii)}  for a facet $F$ and $x,y \in \cup_{\lambda \in \R^+} \lambda C_F$, $\Lambda_L(x+y) = \Lambda_L (x) + \Lambda_L (y)$. 
 \item{(iv)} if $\partial_e L \subset \Q^{d}$, then $v_F \in \Q^{1\times d}$ and  $\Lambda_L (\Z^d)$ is contained in a finite union of subgroups of the form $\frac 1q \Z^+$ for various rational $q$. 

 \Pf (i,ii) Pick $x \in L\setminus \brcs{0}$; the ray from the origin through $x$ hits the boundary at $y = \lambda_0 x$ (where, as in the previous lemma, $1/\lambda_0 = \Lambda_L (x)$; obviously $\lambda_0 \geq 1$). Since the boundary is the union of the facets of $L$, there thus there exists a facet $F_0$ containing $y$. Then $x_0 = (1-1/\lambda_0)0 +  (1/\lambda_0) y $ expresses $x_0$ as a convex combination of $0$ and $y$, so $x \in C_{F_0}$, verifying (i). 

 Applying $v_{F_0}$ to $x_0=  (1-1/\lambda_0)0 +  (1/\lambda_0) y $, we obtain $v_{F_0}x = 1/\lambda_0 = \Lambda_L (x) $. For any other face $F$ of $L$, we have $v_F y \leq 1$ (with equality only if $y \in F$), and thus $v_F x \leq 1/\lambda_0$. Hence $\Lambda_L (x) = \max_F v_F x$. 

 \noindent (iii) An immediate consequence of (i) and (ii). 

 \noindent (iv) $v_F$ is the unique solution (since $F$ is a facet) to $v_F s = 1$ for all $s \in \partial_e F \subset \partial_e L \subset \Q^{1\times d}$. Any (not necessarily homogeneous) linear system  all coordinates of which (and of the nonhomogeneous part) are rational that has a real solution has a rational solution; since the solution is unique, it must have only rational coordinates.  The second part is a consequence of $v_F$ being rational and (i--iii). \qed 

 If $S$ is an admissible subset of $\Z^d$ then $K \equiv K(S): = \cvx S$ (computed within the standard copy of $\R^d$) is a compact convex polytope with the origin as an interior point, and $\partial_e K \subset S \subset \Z^d$. So all of the preceding applies to $K$. 

 \Lem Proposition \Yfiv. Let $S$ be an admissible subset of $\Z^d$, and let $K = \cvx S$. Then $\hat l_S = \Lambda_K | \Z^d$. In particular, there exists a positive integer $N$  \st $N \hat l_S$ is a weight function. 

 \Pf Pick $x \in \Z^d \setminus \brcs{0}$. There exist a positive integers $n,m$ \st $nx \in mK$, so that $y := nx/m \in K \cap \Q^d$. There exists a facet $F$ so that $y \in C_F$, hence there exists $r \geq 1$ \st $y = rf$ where $f \in F$. Then $\Lambda_K (y) = r$, and so $\Lambda_K (x) =mr/n$. By Lemma \Yfou(ii,iv), $mr/n$ is rational, so that $r$ is as well. By Carath\'eodory's theorem applied to the polytope $F$, there exists an affinely independent subset $\brcs{f_i}$ of $\partial_e F$, together with positive real numbers $\lambda_i$ \st $f = \sum \lambda_i f_i$ and $\sum \lambda_i = 1$; the $f_i$ belong to $\partial_e F \subset \partial_e K \subset \Z^d$ and $f \in \Q^d$; again, the system has a unique solution, $ (\lambda_i)$, forcing all the $\lambda_i$ to be rational. 

 We have $x = \sum_i (mr\lambda_i/n) f_i$. There exist infinitely many positive integers $M$ \st $Mmr \lambda_i/n$ are all positive integers. Thus 
$$\eqalign{
Mx& = \sum_i  \frac { Mmr \lambda_i}n f_i\cr
& \in  \sum_i  \frac { Mmr \lambda_i}n S \cr 
&= \frac{Mmr}n S. \cr
}$$
 Thus $l_S (Mx)  \leq Mmr/n$; this is true for infinitely many $M$, and so  $\hat l_S (x) \leq mr/n = \Lambda_K (x)$. 

 Now we show the reverse inequality. As above, with $\Lambda_K (x) = 1/r$, we write $x = \sum \lambda_i/r f_i$ where $\lambda_i$ are positive rationals adding to $1$ and $f_i$ are some of the extreme points of $K$. There are infinitely many positive integers $N$ \st  $N\lambda_i/r$ are integers for all $i$. We claim  for each $\epsilon > 0$ and all but finitely many positive integers $n$ that $l_S(nx) > n (1-\epsilon)/r:= t_n$. 

To the contrary, suppose $l_s (nx) \leq n(1-\epsilon)/r$ for infinitely many $n$. Then $nx \in \ceil{t_n}S \subset \ceil{t}K$, so $nx/t_n \in K$ and thus $\Lambda_K(x) \leq t_n/n$. Since this is true for infinitely many $n$, this yields $\Lambda_K (x) \leq (1-\epsilon)/r$, a contradiction. 

Hence $l_S(nx) > t_n$ for infinitely many $n$, so that $l_S(nx)/n > t_n /n$. Thus $\hat l_S (x) \geq (1-\epsilon)/r = (1-\epsilon)\Lambda_K (x)$, and this is true for all $\epsilon > 0$. 

By Lemma \Yfou(iv), there exists an integer $N$ \st $N\Lambda_K$ is integer-valued. In addition, $\Lambda_K ^{-1}(\leq n) = nK$, hence $(\Lambda|\Z^d)^{-1}(\leq n) = nK \cap \Z^d$, which is of course finite. Thus $N\hat l_S^{-1}(\leq n)$ is finite. \qed

We also have another proof that $\Z^d$ satisfies WC implicit in the following.

\Lem Lemma \Ysix. Let $S$ be an admissible subset of $\Z^d$. Then $\hat l_S = \hat {\tilde l}_S$.  

\Pf In general, $\tilde l_S \leq l_S$, so $\hat {\tilde l}_S \leq \hat l_S$. For the reverse inequality, suppose $\tilde l_S (x) = n$. Then there exists a positive integer $m$ \st $mS + x \subseteq (m+n)S$. Hence $mK + s \subseteq (m+n)K$. As earlier (where?), this forces $x \in nK$, so that $\hat l_S (x) \leq n = \tilde l_S(x)$. This yields $\hat l_S (x) \leq \hat {\tilde l}_S (x)$. \qed

 \noindent{\it A curiosity.} Let $S$ be an admissible subset of $\Z^d$. It is an interesting question to determine when $l_S = \hat l_S$; equivalently, when for all positive integers $n$ and all $x \in \Z^d$, we have $l_S (nx) = n l_S (x)$. The following result, which is practically a tautology, applies and allow us to construct lots of examples, many without any non-identity symmetries.

For $T$ a subset of $\R^d$ and $m$ a positive integer, we denote the set of sums of $m$ elements of $T$ by $mT$. If $T$ is convex, then $mT$ is both the set of sums of $m$ elements of $T$, and also the set of elements of the form $mt$ where $t \in T$, so the apparent ambiguity is illusory.

\Lem Proposition \Ysev. Let $S$ be an admissible subset of $\Z^d$, and set $K = \cvx S \subset \R^d$. Then $l_S = \hat l_S$ if and only if all of the following conditions hold.\par\item{(0)} $K \cap \Z^d = S$;
\item{(i)} for all positive integers $m$, we have $mK \cap \Z^d = m(K \cap \Z^d)$; 
\item{(ii)} for all positive integers $m$, we have $(mK \cap \Z^d )\setminus (m-1)K = \partial (mK) \cap \Z^d$.  

\Rmk   Property (i) (applied to general compact convex lattice polytopes, not necessarily those arising as the convex hull of an admissible set) is known as {\it solidity\/} [H1]; it is t here that if $K = \cvx\Log P$ is solid and $K \cap \Z^d = \Log P$, then $R_P$ (in our notation, $R_f$ for $f = P$) is integrally closed in its field of fractions. It is sufficient to verify the condition for all $m \leq d-1$. The Shapley-Folkman lemma gives a very fast argument to show that if the condition is true for all $m \leq d$, then it is true for all $m$. 

\Rmk   When $m =1$, the condition in (ii) is simply that $0$ is the only interior lattice point of $K$. To verify (ii), it is probably sufficient to verify it for $m \leq d-1$. 

\Pf Assuming the hypotheses (0--ii) hold, we show that $l_S = \Lambda_K|\Z^d $. Pick $z \in \Z^d$. If $z = 0$, then $l_S (z) = 0 = \Lambda_K (z)$. So we may assume that $z$ is nonzero. Suppose $l_S (z) = m$; necessarily $m \geq 1$. Then $z \in mS \setminus (m-1)S$. If $z \in (m-1)K$, then $z \in (m-1)K \cap \Z^d = (m-1)S$ (by (0,i)). a contradiction. Thus $z \notin (m-1)K$. Thus $z \in (mK\cap \Z^d) \setminus (m-1)K$, so by (ii), $z \in \partial (mK)$. Thus $z/m \in \partial K$, so $\Lambda_K (z) = m$. 

Now suppose that $l_S = \LL_S$. Since $K = \cvx S$, it follows that $l_S = \Lambda_K $. Suppose $x \in nK \cap \Z^d$ for some positive integer $n$. In particular, $x/n \in K$. Let $F$ run over the facets of $K$, and form $C_F = \cvx \brcs{0,F}$.  We have 
$$
\partial_e C_F = \brcs{0} \cup\partial_e F \subset \brcs{0} \cup \partial_e K \subset S .
$$
There exists a facet $F$ such that $x/n \in C_F$. By Carath\'eodory's theorem, there exists an affinely independent subset $T \subset \partial_e C_F$ \st $x/n$ is in the relative interior of $\cvx T$, i.e., there exist strictly positive $\brcs{\lambda_t}_{t\in T}$ with $\sum \lambda_t = 1$ \st $x/n = \sum \lambda_t t$. The standard uniqueness argument yields that since $x/n$ and all of the $t$ have rational coordinates, so do each of the $\lambda_t$. If one of the $t$ is the origin, then we can discard it, and in that case $\sum_{T \setminus \brcs{0}} \lambda_t < 1$; otherwise, $\sum \lambda_t =1$ and each of the $t \in F$, so that $x/n \in F$. 

In the first case, let $a = 1- \sum_{t\neq 0}\lambda_t$, so that $x/(an )\in F \subset \partial K$. Hence $l_S(x)= \Lambda_K (x) = n a$. So $na \in \N$, and $x \in naS \setminus (na-1)S$. Since $a < 1$, we have $na < n$, so that $x \in nS$ (since the origin belongs to $S$, whenever $m \leq m'$, it follows that $mS \subset m'S$). 

On the other hand, if $0 \notin t$, then $x/n \in F \subset \partial K$, and thus $\Lambda_K (x) = n$, whence $\lambda_S (x) = n$, so that $x \in nS$. 

Hence in all cases, $x \in nK \cap \Z^d$ entails $x \in nS$. When $n =1$, this yields condition (0); for all larger $n$, it yields condition (i). 

To obtain (ii), pick $x \in mK \cap \Z^d \setminus (m-1)K$. Then $m \geq \Lambda_K(x) > m-1$. As $\Lambda_K = l_S$ and the latter is integer-valued, we must have $\Lambda_K(x) = m$.  Hence $x/m \in \partial K$, so that $x \in \partial(mK)$. \qed 

\noindent{\it Examples.} We can obtain examples from orbits of $S_{n+1}$ on the fundamental weights of $SU(n)$, and more generally some orbits of reflection groups acting on $\Z^d$. Usually both conditions are easy to verify in these cases. \par
For a class of examples without any symmetry (other than the identity), we begin with right simplices in $\Z^d$. 
Let $\alpha = (a(1), a(2), \dots a(d))$ be a strictly positive $d$-tuple in $\Z^d$, and suppose that $a(1) \leq a(2)) \leq \dots \leq a(d)$. Define $L_{\alpha}$  to be the convex hull of $\brcs{0} \cup \brcs{a(i){}i}$ where ${}i$ run over the standard basis elements. 

It is easy to check that $L_{\alpha}$ has unique interior  lattice point iff $\sum 1/a(i) < 1 \leq 1/a_d + \sum 1/a(i)$, and in that case, the interior lattice point is $(1,1,\dots,1) = \sum {}i$. When this holds, define $K_{\alpha} = L_{\alpha} - (1,1,\dots, 1)$, so that the origin is the unique lattice point in the interior. Set $S \equiv S_{\alpha} = K_{\alpha} \cap \Z^d$. 

For $d = 3$, it is elementary that if $\alpha = (2,3,m)$, then $L_{\alpha}$ has unique interior point iff $7 \leq m \leq 12$, and that $K_{\alpha}$ satisfies all the conditions of Proposition \Ysev\ if $7 \leq m \leq 10$. These have no  non-identity symmetries, and although the corresponding polytopes $K$ are simple, they are not integrally simple (in contrast to the reflection group examples). 

More interestingly, when $m = 11$, conditions (0,i) hold, but (ii) fails at $m =2$, even though it succeeds for $m =1$. Presumably, even in the presence of conditions (0,i), in order that (ii) hold it is generally necessary that it apply for all $m \leq d-1$. 

One might conjecture that if $l_S = \Lambda_K|\Z^d$ and $l_S$ is $J$-stable for some action of a finite group $J$ on $\Z^d$ that acts transitively on the facets (or the extreme points), then there exists a reflection group and a corresponding action acting transitively on the facets (or the extreme points, respectively).  

\SecT Appendix 2 Questions/conjectures

A number of questions arise  concerning WC and SWC.  In all cases, $G$ is a finitely generated infinite discrete group.

 \item{(1,2)} If $H$ is a finite index subgroup of $G$, then $H$ satisfies WC iff $G$ does. [For SWC, $\Z$ is a normal subgroup of $D_{\infty}$ of index two, but only the latter satisfies SWC, not the former.]

 \item{(3)} If there exists  admissible $S$ inside $G$ \st $l_S$ satisfies WC, then $G$ satisfies WC. [If this were true, then semidirect products  of the form $G =H \times _{\theta} K$ would satisfy WC,  if $H$ satisfies WC and $K$ is  finite.] 
 
 \item{(4)} There exists an infinite group $G$ and an admissible set $S$ \st $\tilde l_S$ is identically zero. [A candidate would be a simple infinite group with every nonidentity element of order $p$.] 

 \item{(5)} Characterize when $G = \Z^d \times _{\theta} \Z$ satisfies WC, where $\Arrow \theta; \Z. \gl(d,\Z)$, and of course $A = \theta (1)$ is not of finite order. [Such crossed products appear not to satisfy SWC.]

 \item{(6)} For every $\alpha \in (0,1)$, there exists a group $G$ satisfying SWC$(\alpha)$ but not SWC$(\alpha')$ for all $\alpha' > \alpha$. [There are lots of groups satisfying SWC(0) but not SWC$(\alpha)$ for any $\alpha > 0$, e.g., $\Z^2 \times_{\theta} \Z^2$ where $\theta(k)(w) =- w$, where $k$ is the nontrivial element.]  

 \item{(7)} SWC(0) does not imply WC. \vskip 3pt

\comment
 \def\blob{{\ding something}}

 \noindent In a ring, let $r(\cdot)$ denote the right annihilator (a right ideal) of whatever is enclosed. If a ring has the ascending chain condition on right annihilators, then the sequence $r(t) \subset r(t^2) \subset r(t^3) \subset \dots$ is eventually stationary. Goldie ring satisfy this condition, so that if $\R G$ is a Goldie ring, $(r(f^n))$ is eventually stationary, and the conditions of xxx apply. If $G$ is polycyclic by finite, $\R G$ is a Goldie ring, and probably there are a lot of other conditions that guarantee it (but the group cannot contain a free subgroup of rank exceeding one). 

 But more should be true.

 \noindent {\it Questions.} If $f$ is admissible, then is $(r(f^n))$ eventually stationary (for every finitely generated group)? Is there perhaps a uniform bound on the lengths of all such sequences, as $f$ varies over admissible elements of $(\R G)$?
 
 If every admissible element of $\R G$ is not a right zero divisor, then $G$ is torsion-free. Hence we have a weakening of Kaplansky's zero divisor question.
 
 {\it Question.} If $G$ is a finitely generated torsion-free group, and $f$ is an admissible element, then $f$ is not a zero divisor in $\R G$.

Neumann and Neumann

\endcomment
\long\def\Rf[#1] #2, #3. #4\par%
{\vskip 2pt \itemitem{[#1]} #2, {\it #3,} #4\par\vskip2pt}

\SecT References

\Rf [Ca] ER Canfield, From recursion to asymptotics: on Szekeres formula for the number of partitions.
Electronic J Comb 4 (1997) number 2.

\Rf [C] IG Connell, On the group ring.
Can J Math 15 (1963) 650--685.

\Rf [GiH] T Giordano and D Handelman, Matrix-valued random walks and variations on property AT. M\"unster J Math 1 (2008) 15--72.

 \Rf [G] KR Goodearl,
Partially ordered abelian groups with interpolation.
Mathematical Surveys and Monographs 20  (1986) American Mathematical Society, Providence RI.

\Rf  [GooH] KR Goodearl \& DE Handelman, Rank functions and K${}_0$ of regular rings. 
J Pure \& App Algebra 7 (1976) 195--216.

\Rf  [GoH] KR Goodearl \& DE Handelman, Metric completions of partially ordered abelian groups. 
Indiana Univ Math J 29  (1980) 861--895.

\Rf [GrH]  R Grigorchuk \& P de la Harpe,  On problems related to growth, entropy, and spectrum in group theory. J Dynamical and Control Systems 3 (1997) 51--89.

\Rf [H1] D Handelman, Positive polynomials, convex polytopes, and a random walk problem. Lecture Notes in Mathematics 1082 (1986) Springer-Verlag 133 + x.

\Rf [H2] D Handelman, Positive polynomials and product type actions of compact groups. Memoirs of the American Mathematical Society  310 (1985)  79 + xi.

\Rf [H3] D Handelman, Matrices of positive polynomials. Electronic J of Linear Algebra  19 (2009)  2--89.

\Rf [H4] D Handelman, Eigenvectors and ratio limit theorems for Markov chains and their relatives. J d'analyse math\'ematique  78 (1999)  61--116.

\Rf [H5] D Handelman, Deciding eventual positivity of polynomials. Ergodic theory and dynamical systems  6 (1986)  57--79.

\Rf [HL] D Handelman and J Lawrence, Strongly prime rings. Trans Amer Math Soc 211 (1975) 209--223. 

\Rf [NN] B  \& H Neumann, Groups covered by permutable sets.
J London Math Soc 29 (1954) 236--248.

\Rf [P] DS Passman, The algebraic structure of group rings. Wiley-Interscience (1977) 720 + xiv.

\Rf [S] MK Smith, Central and invariant conditions in group algebras. Houston J Math 3 (1977) 125--130.

\Rf [Sz1] G Szekeres, An asymptotic formula in the theory of partitions.
Quart J Math 2 (1951) 85--108.

\Rf [Sz2] G Szekeres, Some asymptotic formul\ae\  in the theory of partitions.
Quart J Math 4 (1953) 96--111.

\Rf [VM] AM Vershik \& AV Malyutin, Asymptotic behavior of the number of geodesics in the discrete Heisenberg group. 
J Math Sci (NY), 240:5 (2019), 525--534.

{}

\vskip 10pt

Mathematics Department, University of Ottawa, Ottawa ON  K1N 6N5, Canada; dehsg\@uottawa.ca

\end

 Given $k,k' \in K$, there exists $c(k,k') \in H$ \st $E_k E_{k'} = c(k,k')E_{kk'}$ (the function $c$ is a one-cocyle, but we do not need this). Let $X = \sup_{K \times K} \brcs{l(c(k,k'))}$. Let $g = hE_k$ and $g' = h' E_{k'}$. Then 
$$
\eqalign{ \hat l(gg')& = \hat l( hE_k h' E_{k'})\cr
 & = \hat l\(h\cdot (E_k h' E_k^{-1})\cdot (E_k E_{k'})\) \cr 
 & =  \hat l\(h\cdot (E_k h' E_k^{-1})\cdot c(k.k') (E_{k k'})\) \cr 
& = l(h E_k h' E_k^{-1} c(k.k')) \leq l(h) + l(E_k h' E_k^{-1}) + l(c(k,k'))\cr
& \leq  l(h) + l(E_k h' E_k^{-1}) +X \cr 
 & = \hat l(hk) + \hat l(h'E_{k'}) + X.
}$$
 We can weaken the hypotheses considerably; for example, if $l$ is merely a semi-weight function with constant $D$, then we can replace $X$ by $X + 2D$. Moreover, we do not require full $\Cal K$-invariance of $l$, merely that $|l(E_k h E_{k}^{-1}) - l (h)|$ is uniformly bounded (as $h$ varies); that is, $l \sim l(E_k \cdot E_k^{-1})$ for each $k$.
 On the other hand, I have not been able to show that $\hat l$ is a weight function if $l$ is. 

Let $G$ be a finitely generated group. It satisfies sWC if there exists a semi-weight function $l$, a finite set $W$, and a positive constant $C$ \st for all $g,g' \in G$, $l(gWg') \geq l(g)  + l(g') - C$. It turns out that sWC implies WC, even though the former only involves semi-weights.

 \Lem Lemma. If $G$ satisfies sWC, then it satisfies WC.

 \Pf Let $S$ be an admissible set for $G$, and let $(l,W,C,D)$ implement sWC.  There exists $n$ \st $W \subset S^n$. Pick $g $ in $G$. Suppose $\tilde l_S (g) = t$. There exists a positive integer $m$ \st $S^m g \subset S^{m+t}$; we can make $m$ as large as we like, in particular, $m > n$. Select $g_1 \in S^{m-n}$ \st $l(g_1) = l(S^{m-n})$. Then $g_1 W g \subset S^{m+t}$. We have $l(g_1 W g)  \geq l(g_1) + l(g) - C = l(S^{m-n}) + l(g) - C$.

On the other hand, $l(S^{m+t}) = l(S^{m-n}S^{n+t})\leq l(S^{m-n}) + l(S^{n+t}) + D$. As $l(g_1 W g) \leq l(S^{m+t})$, we deduce $l(S^{m-n}) + l(S^{n+t}) + D \geq  l(S^{m-n}) + l(g) - C$. Thus $l (g) \leq l(S^{n+t}) + D + C \leq l(S^n) + 2D + C + l(S^t) \leq l(S^n) + 2D + C +  tl(S) + (t-1)D$.  That is, $l(g) \leq l(S^n) + D + C + t(l(S) + D)$. Since $t = \tilde l_S(g)$, we have 
$$\tilde l_S \geq  (l - K)/(l(S) +D),$$
where  $K =  l(S^n) + D + C $. It easily follows that $\tilde l_S^{-1} (\leq t)$ is finite. \qed

 \Lem Lemma. Suppose $H$ is a normal subgroup of $G$ \st $G/H:= K$ is finite. If $H$ satisfies sWC, then so does $G$. \par\par \Pf Let sWC be implemented by $(l_0,W_0,C,D)$. Let $\Cal K = \brcs{E_k}\subset G$ be a complete set of equivalence class representatives of $G/H$, and for each $k \in K$, define $l_k = l_0(E_k^{-1} \cdot E_k)$. Then each $l_k$ is a semi-weight function, and $(l_k, E_k^{-1}W E_k,C,D)$ implements sWC. We set $l = \sup_K l_k$. It is routine that $l$ itself a semi-weight function more ...

 \Lem Proposition. Let $K$ be a finite group, and $\Arrow \theta ; K.\gl(d,\Z)$ be a faithful representation over the integers. Suppose that the rational representation of $K$ given by  $\Theta:=\theta \otimes 1_{\Q}$ satisfies both of the following:
 \item{(i)} the trivial representation does not appear in $\Theta$;
\item{(ii)} $\Theta$ contains no irreducible representations with multiplicities exceeding one.{\par
\noindent} Then the crossed product $G = \Z^d \times _{\theta} K$ satisfies SWC.

 We give a concise outline of the argument. As $\Theta$ has no multiplicities,  $\Q^d$ is a cyclic $\Theta(K)$-module. It follows that by taking the transposes (the dual action), $\Theta^{\op}$ acts on $\Q^{1 \times d}$ cyclically. Hence there exists $v \in \Q^{1\times d}$ \st $v \Theta (K)$ contains a (rational) basis for $\Q^{1\times d}$. By multiplying $v$ by a sufficiently large positive integer, we may assume that $v \in \Z^{1 \times d}$.  

 Define the compact convex polytope (with rational extreme points), 
$$
 C_0 =\Set{x \in \R^d} {\theta(k)( x )\leq 1 \text{for all $k \in K$}}
$$
 Then $C_0$ is $\Theta$-invariant, and it is easy to check that it contains an open $d$-ball, and that $0$ belongs to the interior. One of the facets is given by $v^{\theta(k)}x = 1$ for some $k$, and it follows that the facets are precisely of this form, one for each element $k$. In particular, $\theta$ (or $\Theta$) acts transitively on the facets of $C_0$. 

We may multiply $C_0$  by a sufficiently large integer so that all the extreme points have integer coordinates; call the resulting lattice polytope $C$. Define (unnormalized) weight function (that is, with values a cyclic subgroup of the rationals, not necessarily the integers) on $\Z^d$ via $l_0 ((x_i)) = \max_{k \in K} v^{\theta(k)}x$. This is additive on the cone generated by any facet, and multiplication by a suitable integer yields an implementation of SWC (by $(l,K)$). \par \par Now for the details. [transcribe from notes]

There is also a corresponding result for WC, but it is extremely restrictive.

Concerning WC, there are at least three plausible stronger versions. Here $S$
denotes an admissible subset of $G$, as usual.

\noindent {\it Uniform WC.} This was discussed in xxx; $G$ satisfies {\it
uniform WC \/} if there exists a weight function $l$, together with $C \in
\Z^+$, and a finite subset $W$ of $G$, \st for all $g,h \in G$, $\max_{w \in W}
l(hwg) \geq l(h) + l(g) - C$. That this implies WC is a consequence of xxx, and
it also entails that $\tilde l_S (g) \geq l(g) - r - C$, where $r$, depending
only on the admissible set $S$ is the smallest integer \st $W \subset S^k$. It
also implies $|\tilde l_S^{-1}(k)| \leq \sum_{k-r-C \leq t \leq k} |l^{-1}(t)|$.

\noindent {\it sWC.} The pair $(G,S)$ satisfies sWC if for all $k$, there
exists $M(k) \in \Z^+$ \st for all $g \in G$,  $S^m g \subset S^{m+k}$ entails
$S^{M(k)} g \subset S^{M(k) +k}$. This is stronger than WC, but it only yields
crude  bounds, $\tilde l_S(g) \geq l_S (g) - M(\tilde l_S (g))$, and $|\tilde
l_S^{-1}(k)| \leq |S^{M(k) + k}|$.

\noindent {\it SWC.} The pair $(G,S)$ satisfies SWC if it satisfies sWC and
$\sup M(k) < \infty$. This implies $\tilde l_S (g) \geq l_S(g) - \sup M(k)$ and
$|\tilde l_S^{-1}(k)| \leq \sum_{k- \sup M(k) \leq k} |l_S^{-1}(t)|$.

For comparison, if $G = \Z^d$ and $S$ is admissible, let $K = \cvx S$ (the
convex hull in $\R^d$); then $|\tilde l_S^{-1}(k)| \leq |k K \cap \Z^d| =
\Oh{k^d}$. This is not sharp; but $|l_S^{-1}(k)| \leq |
\partial (kK) \cap \Z^d|$ (counting lattice points on the boundary of $kK$),
and this is asymptotically sharp. If $S = K \cap \Z^d$ where $K = mC$ for some
lattice polytope and $m \geq d-1$ (or more generally, if $S^m = m K \cap \Z^d$
for some  $m \geq d-1$), then $l_S = \tilde l_S$, except at the origin.

For the free group $\pmb F_k$  with generators $\brcs{g_i}$ and $S = \brcs
{g_i^{\pm 1}} \cup \brcs{1}$, it is also true that $l_S = \tilde l_S$. This
leads to an obvious question, whether for all infinite finitely generated groups
$G$, there exists an admissible set $S$ \st $l_S (g) = \tilde l_S(g)$ for all $g
\in G \setminus \brcs{1}$  (of course, $\tilde l_S(1) = 0$ while $l_S(1) = 1$. I
suspect a counter-example exists, even when $G$ is torsion-free.

So let $K$ be a finite group, and let $\Arrow \theta; K . \gl(d,)$ be a group homomorphism. Form the semidirect product, $G = \Z^d \times_{\theta} K$. We wish to investigate sufficient conditions for $G$ to satisfy  SSWC (not just SWC). Let $H$ denote the canonical copy of $\Z^d$. A necessary condition is that $ \theta$ have no fixed points  (else the semidirect product will have infinite centre, which contradicts every form of SWC). We make the following assumptions about the action, in addition to $\theta $ being one to one.

Let $S$ be an admissible subset of $H$.  Define $C$ to be the convex hull of $S$, within $\R^d$ (having as basis a basis for $H$). For each facet $F$ of $C$, define the  pyramid $C_F = \cvx \brcs{F,0}$, the set of lattice points $S_F = C_F \cap H$, and the cone generated by $\Cal S_F:= \cup_{n\geq 1} n S_F$.  (If $T$ is a subset of $H$, then $jT$ will denote the set of sums of $j$ elements of $T$, repetitions permitted.) We say $S$ is {\it $\theta$-suitable\/} (or simply {\it suitable\/}) if the following hold:
\item{(i)} $0 $ belongs to the interior of $C$; 
\item{(ii)} $S_F$ generates $ H$ as a group;
\item{(iii)} for each facet $F$ of $C$,  $F \cap S$ is a $\Z$-basis  for $H$
\item{(iv)}$\theta$ acts transitively on the set of facets of $C$.
 \comment
\item{(d)} there exists $C> 0$ \st for all facets $F$ and all  $w, w'$ in $\cup_{n\geq 0} n(S \cap F)$, $l_S(v+w) \geq l_S (v) + l_S (w) - C$. 
 \endcomment

 \Lem Lemma. If an admissible set $S$ of $\Z^d$ is $\theta$-suitable, then the semidirect product $\Z^d \times _{\theta} K$ satisfies SSWC.

\Rmk The conditions on $\theta$-suitability can be weakened for the purposes of this result. It is sufficient that (iii) and (iv) be replaced by some of their consequences, specifically, there exists a weight function $l$ and  a subsemigroup $\Cal S$ \st $l$ restricted to $\Cal S$ is additive (up to a constant $C$), and every element of $G$ can be moved to an element of $\Cal S$ by an element of $K$. But verifying the existence of an almost additive weight function (on a suitable subsemigroup) is generally difficult.

 \Pf We can write $C = \cup C_F$ (the union over all the facets); we also have that if $F$ and $F'$ are unequal facets, then $C_F \cap C_{F'} = \cvx\brcs{0,F \cap F'}$. Pick a nonzero element $w$ of $H$. Then there exists a smallest integer $n$ \st $w \in nC$. We claim that $l_S(w) = n$. 

 First,  $w \not\in (n-1)K$, implies $w \not\in (n-1)S$; hence $l_S(w) \geq n$. For the reverse inequality, if $w \in nK$, then there exists a facet $F$ \st $w \in nK_F$. Hence $w \in nK_F \cap H$; by (b),  $w \in \cup _{j=1}^n j(S \cap F) \cup \brcs{0}$. As $w$ is not zero and does not belong to $jS$ for any $j < n$ (since $jS \subset jK$). Hence $w \in n(S \cap F)$, and thus $l_S(w) \leq n$. 

 Now pick $v, w \in H$, both nonzero. By transitivity  of the action, we may choose $k \in K$ \st $\theta(k)(w)$ is in the same lattice cone as $v$. Then we note that $l_S$ is additive on points lying in the same $\Cal S_F$ (the restriction of $l_S$  to $\Cal S_F$ is given by $x \mapsto m$ where $x \in mF$).

 Now define $T = \cup_{k\in K} Sk$. This is an admissible set in $H\times_{\theta} K$, and the claim is that $l_T|H = l_S$. But this is an easy consequence of invariance of $S$. 

 Set $W = K$; we verify the SWC property for $(l_T,W)$ with $C = 1$. For $a,b \in H$, there exists $k$ \st $\theta(k)(b)$ lies in the same orthant as $a$, and thus 
$$\eqalign{
l_T(akb) &= l_T (ab^k k) \geq l_T (ab^k) -l_T(k^{-1})\cr
 & = l_S(ab^k) - 1= l_S(a) + l_S(b^k)-1\cr
 & = l_S(a) + l_S(b)  - 1 = l_T(a)+ l_T (b)-1.\cr
}$$
For general elements $ak'$ and $bk''$, pick $k \in K$ \st $b^{kk'}$ lies in the same $S_F$ as $a$; then $l_T(ab^{kk'})=l _S(ab^{kk'}) = l_S (a) + l_S (b)$ (since $S$ is invariant. We have $ak' k  b k'' = a b^{kk'}k' k k''$, and thus $l_T(ak'kbk") \geq l_T(ab^{kk'}) - 1$, and we are done. \qed 

This yields lots of examples
 
\Lem Example 1. 
 1 Let $K = \Z_2^d$ acting on $\Z^d$ by diagonal $\pm 1$ matrices. Set $S = \brcs{0, \pm {}i}$ (where ${}i$ run over the standard basis elements of $\Z^d$). 

\noindent Then $C = \cvx \brcs{\pm {}i}$, the $d$-dimensional version of the octahedron, and the facets are determined  by a  sequence of $\pm 1$ in $\Z^d_2$; transitivity of the action on the facets is clear, and  each $S_F$ is a basis of $\Z^d$. So all the hypotheses are satisfied. \qed
 
 \Lem Example 2. Let $K = P_{d+1}$ the permutation group on $d+1$ symbols, or any of its subgroups that acts transitively on the symbols. Begin with the permutation representation of $K$ on $\Z^{d+1}$, and observe that the element $v = (1,1,\dots,1)^T$ is invariant. This induces an action of $K$ on the quotient group, $\Z^{d+1}/v\Z$, which we identify with $\Z^d$. 

 \noindent Let $S = \brcs{0, {}1, \dots, {}d; -\sum {}i}$. Then it is obvious that $S$ is admissible, invariant, and $C$ is a simplex. Moreover, in this case, $K$ acts transitively on the extreme points of $C$ (unlike the situation in Example 1), and since $C$ is a simplex, it follows that $K$ acts transitively on the facets. One of the $S_F$ is just $\brcs{0,{}i}$, and thus (iii) of the definition of suitability follows from (iv). 

In particular, when $K$ is cyclic of order $d+1$, this yields the companion matrix for $x^{d} + x^{d-1} + \dots + x + 1$ as the representation. If $d+1$ is prime and the corresponding cyclotomic ring $\Z[\xi]$ is a principal ideal domain (here $\xi = \exp 2\pi i/(d+1)$), this representation is unique up to conjugacy. \qed

 With $d=1$, we obtain (only) the infinite dihedral group. With $d= 2$, faithful representations of finite groups work, with the obvious exception of $C_2$. Explicitly, we have the following. 
 
 \Lem Proposition. Let $K$ be a finite group with more than two elements, and let $\Arrow \theta; K.\gl(2,\Z)$ be a one to one homomorphism of $K$. Then $\Z^2 \times_{\theta} K$ satisfies SWC. 

 \Pf The image of $K$ is isomorphic to $K$ and a subgroup of $\gl(2,\Z)$, hence up to isomorphism must be one of $\Z_2$, $\Z_2^2$, \Z_3$, $C_3$, $C_6$, $D_3$, $D_4$, or $D_6$. First, for $K =C_6$, let $k$ be a generator. Then $\theta (k)$ has order $6$, so its set of eigenvalues must be either $\brcs{z,z^{-1}}$ or $\brcs{-z,-z^{-1}}$ where $z = (1+\sqrt{-3})/2$. Since we obtain one from the other just by multiplying by $-1 = z^3$, we can assume the spectrum is $\brcs{z^{\pm1}}$. Since $\Z[z]$ is a principal ideal domain, up to conjugacy, $\theta(k)  =\(\smallmatrix  1 & -1 \\ 1 & 0 \\ \endsmallmatrix \):= A$. Set $S $ to be the orbit of $\(\smallmatrix  1 & -1 \\ 1 & 0 \\ \endsmallmatrix \)$, together with $0$, that is, $\brcs{0} \cup \brcs{\pm {}i, \pm \(\smallmatrix  1  \\ 1  \\ \endsmallmatrix \)$. It is easy to check that the conditions (i--iv) hold; here $C$ is a hexagon. 
 
 Next, for $D_6$. Then restricted to the copy of $C_6 = \langle k\rangle$, we can assume that $\theta (k) = A$ (or $-A$, but it does not matter). Solving the equation $BAB = A^{-1}$ (that is $B A = A^{-1}B$) and $B^2 = I$ by brute force (actually quite easy), that one $\pm B$ is $\(\smallmatrix  1 & 0 \\ 1 & -1 \\ \endsmallmatrix \)$ or its transpose. In any of the four cases, the $S$ chosen above for $C_6$ is still invariant. So properties (i--iv) hold here as well. 
 
 For $K = C_3$ and $D_3 \iso S_3$, the arguments are practically the same (the ring is still $\Z[z]$). 

 This leaves the case that $K = C_2 \times C_2$. Pick $k \in K$ \st $\theta(k)$ is not scalar (this is possible since $\theta$ is one to one), and set $Q = \theta(k)$. Then $Q^2 = I$ and $Q$ is not scalar implies its spectrum is $\brcs{\pm 1}$. Using elementary simultaneous row and column operations (that is, adding a multiple of  row $i$ to row $j$ then subtracting the same multiple of column $j$ from column $i$), we reduce to the case that $Q$ is conjugate up to sign to one of $\(\smallmatrix  1 &0 \\ 0 & -1 \\ \endsmallmatrix \)$ or $B:=\(\smallmatrix  1 & 1 \\ 0 & -1 \\ \endsmallmatrix \)$ (the latter is conjugate to $\(\smallmatrix  0 & 1 \\ 1 & 0 \\ \endsmallmatrix \)$). 

Going up to M${}_2 \Q$  (the algebra over the rationals), the eigenvalues being distinct integers, the centralizer of $Q$ is just the polynomial algebra in $Q$. Hence if we let $k' \in K $\theta (k)$ must be of the form $qQ + rI$ where $q,r $ are rationals. Since the order is two, we deduce $qr = 0$, so that the only possibilities are $q = \pm 1$ and $r = 0$ or $r = \pm 1$ and $q = 0$. There are only four, so they are exhausted by the group elements. Hence $-I $ is in the image of $\theta$. Up to conjugation, the image of $\theta$ is either the set of diagonal $\pm 1$-matrices, a special case of Example 1, or the group $\brcs{B,I,-B,-I}$. In the latter case, set $S = \brcs{0, \pm {}2, \pm ({}1-{}2)}$. Then $C$ is a  is a parallelogram, and it is easy to see that (i--iv) hold. \qed

 Here is the one motivating example. Let $H$ be a normal subgroup of finite index in $G$, and set $K = G/H$. Let $\Cal K \subset G$ be a complete set of coset representatives, indexed as $\brcs{E_k}_{k \in K}$. Suppose that $l$ is a weight function on $H$ 
\st for all $k \in K$ and all $h \in H$, $l(E_k h E_k^{-1}) = l(h)$ (we are not assuming that $E_k^{-1} =E_{k^{-1}}$ or anything similar). We define the function $\Arrow \hat l; G. \Z^+$ via $\hat l(hE_k) = l(h)$. This is well-defined, and obviously $|\hat l^{-1} (t)| = |l^{-1}(t)|\cdot |K|$. Weak subadditivity  is slightly more complicated to check.

 \Lem Lemma. Suppose $G =\Z^d \times_{\theta} K$ satisfies SWC, implemented by $(l,W,C)$. Then there exist $C'\geq 0 $ and  $W' \subset K$ with  $|W'| \leq |W|$  \st $(l,W',C')$ implements SWC. 

 \Pf Set $W' = \pi (W)$ where $\Arrow \pi; G.K $ is the quotient map (that is, $k \in W'$ if there exists $v \in \Z^d$ \st $vk \in W$). We verify the property for $W'$ (with slightly bigger $C'$).

 Write $W = \brcs{v_i k_i}$ (repetitions among the $k_i$ may occur). Choose $a \in \Z^d$ and $b \in G$. There exists $i $ \st $l(a v_ik_i b) \geq l(a) + l(b) - C$. Since $av_i = v_i a$, we have 
$$\eqalign{
 l(av_i k_i b) & =  l(v_ia k_i b ) \cr 
& \leq l(ak_ib ) + l(v_i); \text{ thus,}\cr l(ak_i  b)&\geq  l(v_i a k_i b) - l(v_i); \text{ so}\cr 
l(aW' b) &\geq l(a)+ l(b) - c - \max_j \brcs{l(v_j)}.\cr
}$$
 Set $C' = C +  \max \brcs{l(w_j)}$.

 Now consider $ak \in G$ (with $a \in \Z^d$  and $k \in K$, that is, an arbitrary element of $G$). Then $ak k_j b = ak_j (k_j^{-1} k k_j  b)$, and now apply the preceding with $b$ replaced by $k_j^{-1} k k_j  b$.\qed

 In particular, when $G$ is a semidirect product of $\Z^d$ by a finite group $K $ and satisfies SWC, then we can choose $W$ so that $|W| \leq |K|$. It is easy to see that with $C_6$ acting faithfully on $\Z^2$, we can choose $W $ to be the copy of $C_3$ inside $C_6$. This of course leads to the question of minimal size of $W$. For nontrivial amalgamated products, we obtained $W$ with $|W| = 3$. I suspect that for the semidirect products considered here, there is no bound on the minimal size of $W$. Explicitly, for the diagonal action of $\Z_2^d$, it is likely that the minimal size is actually $2^d$ (achieved by the finite group itself).

\comment

Even if $K$ is cyclic and  $\theta$ is a representation on $\Z^d$, it is not clear what the right conditions should be to guarantee that $\Z^d \times_{\theta} K$ satisfies SWC. For example, if $K = C_{p^2}$ with $p$ a prime, not only are there hordes of integral representations, but it is difficult to decide which ones (if any) are suitable \wrt $\theta$. 
\Pf We extend the definition of $l$ to finite sets via $l(T) = \max_{t\in T}
l(t)$. Then it follows from (a) that for finite sets $T,U$ of $G$, $l(T\cdot U)
\leq l(T) + l(U)$ (where $T\cdot U$ is the set of all $tu$ with $t \in T$ and $u
\in U$). In particular, if $S$ is any finite set, then $l(S^{a+b}) \leq l(S^a) +
l(S^b)$.

Now let $S$ be an admissible set. There exists $n$ \st $w
\subset S^n$. Fix $k$, and suppose that $g$ is an element of $G$ \st $S^m g
\subset S^{m+k}$ for some integer $m$. This inequality also holds for any integer
$m' \geq m$ replacing $m$. Hence we can assume that $m > n$. Select $v \in
S^{m-n}$ \st $l(v) = l(S^{m-n})$. Then  $vw \in S^m$,
and thus $swg \in S^{m+n}$. Therefore
$$\eqalign{
l(S^{m+n}) & \geq  l(vwg) \geq l(v) + l(g) - C; \text{
hence}\cr
l(g) & \leq l(S^{m+k}) - l(S^{m-n}) + C \leq l(S^{n+k}) + C.\cr
}$$
The last expression does not depend on $m$. Hence $g \in \cup_{j\leq l(S^{n+k})+
C} l^{-1} (j)$, which by hypothesis (b), is finite.
\qed
\endcomment

We will show that if $f$ is admissible and $G$ is an  ICC group, then under a modest additional condition on annihilators of powers of $f$ (satisfied when $f$ is symmetric), the pair $(G,f)$  satisfies EEP. 

Recall that a group is {\it ICC\/} if the
conjugacy class of every element of $G$ other than the identity is infinite
(equivalently, the group ring $\R G$ has trivial centre). We require an amusing result about weak convergence in the W$^*$ algebra of the group, denoted W$^*(G)$. Here the inner product space is the completion of $\C G$, the inner product given by $\langle a,b\rangle$ equalling $\sum a_g b^*_{g^{-1}}$; alternatively, its trace is the linear functional arising from the augmentation map.

Let $G$ be a group, and let $\brcs{G_i}_{i=1}^n$ be a finite collection of subgroups, each of infinite index in $G$. We denote by a {\it  union of finitely many cosets of $\brcs{G_i}$}, a subset of $G$ of the form
 $$
\bigcup_{i=1}^n \(\bigcup_{j_i \in J_i} j_i G_i\),
 $$
where each $J_i$ is a finite subset of $G$.
 
 \Lem Lemma. Let $G$ be a group and let $\brcs{G_i}_{i=1}^n$ be a finite collection of subgroups, each of infinite index in $G$. Then $G$ is not a  union of finitely many cosets of $\brcs{G_i}$.
 
  \Pf If $n=1$, this follows from $[G:G_1] = \infty$. We proceed by induction on $n$. Suppose $G = \cup_{i=1^n} \(\cup_{j_i \in J_i} j_i G_i\)$. Since $[G:G_n] = \infty$, there exists $x$ in $G$ \st the coset $xG$ is disjoint from $\cup_{j\in J_n} jG_n$. Hence $xG_n \subset \cup_{i=1^{n-1}} \(\cup_{j_i \in J_i} j_i G_i\)$, whence $G_n$ is contained in a union of finitely many cosets of $\brcs{G_i}_{i \leq n-1}$. It is immediate that any union of finitely many cosets of $G_n$ is contained in a union of finitely many cosets of $\brcs{G_i}_{i=1}^{n-1}$. This forces $G$ to be a union of finitely many cosets of $\brcs{G_i}_{i=1}^{n-1}$ contradicting the induction hypothesis.\qed 
 
 Let $Z(h)$ denote the centralizer of the element $h$ of $G$; then $G$ is ICC iff for every $h \neq 1$, $Z(h)$ has infinite index in $G$.
The W*-algebra, $W^*(G)$, is a factor iff $G$ is ICC. 
 
  \Lem Lemma. Let $G$ be an ICC group. Let $Q$ be finite subset of $ G \setminus \brcs{1}$. There exists a strictly increasing nested sequence of finite subsets $\dots \subset T_n \subset T_{n+1} \subset \dots $ of $G$ \st for all $b \in Q$,
 $$
  c_n:= \frac 1{|T_n|} \sum_{g \in T_n} gbg^{-1}
 $$
  converges weakly (in W$^*(G)$) to zero.
 
  \Pf Let $Q = \Supp a \setminus \brcs{1}$. For $q \in Q$, set $G_q = Z(q)$; these are of infinite index in $G$ since $G$ is ICC. 
  
  We will inductively construct $T_n$ so that on defining $\Cal G_n = \Set{jqj^{-1}}{j \in T_n; q \in Q}$, if  $x \in \Cal G_n $  and $x = jqj^{-1}$ for some $j \in \cup_{m\geq 1} T_m$  and $q \in Q$ then $j \in T_n$. 
 
Set $T_1 = \brcs{1}$, and assume $T_n$ has been defined. We claim there exists $g \in G$ \st $gQg^{-1} \cap \Cal G_n = \emptyset$. This condition is equivalent to $g \not\in \cup_{q \in Q; j \in T_n} jZ(q)$. By the lemma, such $g $ exists. So pick one choice for $g$ and  set $T_{n+1} = T_n \cup \brcs{g}$. This yields the strictly increasing nested sequence $\brcs{T_n}$ and with $|T_n| = n$. 
 
  We verify the induction hypotheses. Suppose $jqj^{-1} = j' q'( j')^{-1}$ for some $j \in T_n$, $j' \in T_m$ with $m > n$, and $q,q' \in Q$. Let $n'$ be the smallest integer \st $j' \in T_{n'}$. If $n'\leq n$, there is nothing to do. If $n' > n$, then $j'$ was obtained as the new element \lq\lq$\brcs{g}$\lq\lq\ with the property that $j'Q(j')^{-1} \cap \Cal G_{n'-1} = \emptyset$, yielding a contradiction (since $\Cal G_n \subseteq \Cal G_{n'}$). 
 
  Now we can show that the weak limit of the $a_n$ (defined by the now-constructed sequence $\brcs{T_n}$ converges as indicated. 
 
  Let $h \in Q \cup \brcs{1}$, and define $b_{h,n} = |T_n|^{-1}\sum_{j \in T_n} jhj^{-1}$. Since $b_{1,n} =1$ for all $n$, it suffices to show that $(b_{q,n})$ converges weakly to zero for all $q$ (since the collection of finite convex linear sums of elements of $G$ is weakly dense in the unit ball of $W^*(G)$, and the unit ball is weakly metrizable). This is equivalent to showing that for all $k \in G$, $\langle b_{q,n}, k\rangle \to 0$, that is (since the inner product is given by the trace), the coefficient of $k$ in $b_{q,n}$ tends to zero as $n \to \infty$. We have 
  $$\eqalign{
 \langle b_{q,n}, k\rangle & = \frac{\#\(\Set{j \in T_n}{jqj^{-1} = k}\)}{|T_n|}. \cr
 }$$ 
  If $k \not\in \Cal G_n$ for any $n$, the numerators are all zero. If $k \in \Cal G_m$, then the numerator is simply $\#\(\Set{j \in T_m}{jqj^{-1} = k}\)$ for all $n \geq m$. Since $|T_n| \to \infty$, the limit is zero in all cases. \qed

In particular, for any finite set of elements of $\R G$, $\brcs{b_j}_{j \in J}$, there exists an increasing sequence of finite subsets of $G$, $\dots \subset T_n \subset T_{n+1} \subset \dots $ \st for all $j \in J$, the sequence 
$$
d_{n,j}:= \frac 1{|T_n|} \sum_{g \in T_n } gb_j g^{-1} \to  1\cdot (b_j, 1)
$$
weakly. This also implies that if $c \in \R G$, then $(d_{n,j}c)$ and $(cd_{n,j})$ converge weakly (in $n$) to $(a_j,1) c$ for each $j$.

\Lem Proposition. Suppose $f$ is admissible, $G$ is ICC, and there exists a nonnegative integer $n$ \st for all $a \in \R G$, $f^m a = 0$ for some $m$ implies $f^n a = 0$. Then $(G,f)$ satisfies EEP. 

\Rmk If $n = 0$, the annihilator condition is simply that $f$ is not a zero divisor in the group ring; this is generic. However, if $f_0$ is admissible and $g$ is an element of order $k > 1$, then $f : = f_0 \cdot (1 + g + g^2 + \dots + g^{k-1})$ is admissible, but is a zero divisor. So if all admissible $f$ were not zero divisors, then $G$ would be torsion-free, and we are back in the domain of the zero divisor problem. On the other hand, if the group ring were Goldie (which means $G$ must be small in some sense), then it satisfies the ascending chain condition on right annihilators, so such $n$ exists. 

\Pf Let $J$ be a finite set of nonnegative integers, and suppose there were a dependence relation of the form $\sum_{j \in J} \RR_{r_j a_j,j} = 0$, where $a_j$ are nonzero elements of $\R G$, and $r_j$ are real numbers, not all zero. Then we can obviously reduce to the case that all the $r_j$ are nonzero, and absorbing the $r_j$ into the $a_j$, we have a relation of the form $\sum_{j \in J} \RR_{a_j,j} = 0$. Let $k = \min J$, and rewrite the relation in the form $\RR_{a_k,k} +\sum_{j \in J \setminus \brcs{k}} \RR_{a_j,j} = 0$. 

For each $g$ in $G$, there exists  a positive integer $m(g)$ \st $f^{m(g)}\cdot \(\sum f^{t-j}ga_j\) = 0$, where $t = \max J$; without loss of generality, we can assume $m(g) \geq n$. The assumption on annihilators yields that for all $g \in G$, we have $f^{n}\cdot \(\sum f^{t-j}ga_j\) = 0$. Rewrite this as $f^{n+t-k}ga_k = - \sum_{J \setminus \brcs{k}} f^{n + t-j}g a_j$, for all $g$. Since $a_k$ is not zero, there exists $h \in G$ \st $(ha_k,1) \neq 0$, and of course we have (replacing $g$ by $gh$), $f^{n+t-k}g(ha_k) = - \sum_{J \setminus \brcs{k}} f^{n + t-j}g (ha_j)$ for all $g$, and in particular, $f^{n+t-k}g(ha_k)g^{-1} = - \sum_{J \setminus \brcs{k}} f^{n + t-j}g (ha_j)g^{-1}$

Set $Q = (\cup \supp a_j) \setminus \brcs {1}$.  There exists an increasing sequence of finite subsets, such that the limit of the averages converges on both sides of the equation, yielding $f^{n+t-k}\cdot (ha_k,1) = - \sum_{J \setminus \brcs{k}} f^{n + t-j}\cdot (ha_j,1)$. The left side is not zero, and its support is $\supp f^{n+t-k}$; on the other hand, the support of the left side is contained in $\cup_{j < k}\supp f^{n+t-j}$. This implies $\supp f^{M} \subseteq \supp f^m$ for some pair with $M >   m$; since $G$ is infinite and $f$ is admissible, this is impossible. 
\qed

If $f $ is symmetric (that is, for all $g$, $(f,g) = (f,g^{-1})$), or more generally if $f$ is normal (as an element of W$^*(G)$), then $n = 1$ or $0$, and so in this case, ICC implies EEP.

\comment 
\Lem (Preliminary version) Lemma. Suppose $G$ is amenable and ICC, and let $f$
be admissible and  not a zero divisor in $\R G$. If $\cap_{n\in \N} \SS^n R_f =
\brcs{0}$, then $$ \RR_{a,l} = \RR_{a',l'} \text{ implies either }\cases a= a' =
0 & \text{or}\\ a = a' \text{ and }l = l'.\endcases$$

We first examine the  condition $\cap_{n\in \N} \SS^n R_f = \brcs{0}$. I know
of no $G$ that fail to satisfy this (although if we are permitted to pick $f$
with suitably infinite support, then examples do exist with $G = \Z$). But I
cannot find any general conditions that guarantee it either except for $G$
abelian, where it is a consequence of the truth of a special case Jacobson's
conjecture for commutative noetherian rings.

\Lem Lemma. Suppose $a \in \R G$; then $[a,k] \in J:= \cap_n \SS^n R_f$ if and
only if for infinitely many (and thus all) positive integers $n$, there exist
positive integers $m(n) \to \infty$ with $f^{n+m(n)} a \prec f^{n+k}$.

\Rmk If $f^m a = 0$ for some $m$, then $[a,k] =0$; so a nonzero element in $J$
means there is always a lot of cancellation. It is easy to verify that this
cannot occur if $a \in (\R G)^+$ or more generally if $[a,k] \in R_f^+ \setminus
\brcs{0}$.

\Pf Assume the support condition holds; then $[f^{m(n)}a,k] \in R_f$, and
$$\SS^{m(n)} [f^{m(n)} a,k] = [f^{m(n)}a, k+m(n)] = [a,k].$$
Hence $[a,k] \in \SS^m R_f$ for infinitely many $m$, and since $\SS^{m-1} R_f
\subset \SS^m R_f$, we see that $[a,k] \in J$.

Conversely, if $[a,k] \in \SS^m R_f$, then we may find $b_m \in \R G$ and a
positive integer $j_m$ \st $[b_m,j_m] \in R_f$ \st $[a,k] = [b_m,j_m+m]$. This
entails the existence of a positive integer $N$ (depending on $m$) \st $f^{N+j_m
+ m}a = f^{N+k} b_m$. Since this will also hold for any larger choice of $N$, we
can also assume that $f^N b_m \prec f^{N+j_m}$. Hence $f^{N+j_m + m}a \prec
f^{N+j_m + k}$.\qed

\Pf (of preliminary version). Since $G$ is amenable, we may find a F\o olner
sequence of finite sets, $T_n$; that is, $\cup T_n = G$ and for all $\epsilon >
0$, for all $g \in G$, there exists an integer $N \equiv N(g)$ \st for all $n
\geq N(g)$,
$$ \frac{|T_n g \cap T_n|}{|T_n|} \geq 1 - \epsilon.$$

For $h \in G\setminus \brcs{1}$, define
$$
a_{h,n}  = \sum_{g \in T_n} g^{-1} h g,
$$
and now we work inside the type $II_1$ factor $W^*(G)$ (ICC guarantees that it
is a factor); its trace is obtained from the standard trace on $\R G$.

\item{(1)} {\it The sequence $(a_{h,n})_{n \in \N}$ has a strong limit point.}

\noindent As $a_{h,n}$ is a convex combination of unitaries, it is bounded in
norm by $1$. Thus the sequence has a weak limit point, and since we are in a
closed bounded convex set, every weak limit point is also a strong limit point.

\item{(2)} {\it If $z \in \C\cdot 1$ is a strong limit point of $(a_{h,n})_{n \in
\N}$, then $z = 0$.}

\noindent Let $\|\cdot \|_2$ denote the tracial 2-norm; restricted to elements
$a $ of the group ring $\|a\|_2^2$ is the absolute value of the constant
coefficient of $aa^*$. On bounded sets, convergence \wrt this norm is equivalent
to strong convergence. Hence we have $\| (a_{h,n} - z)(a_{h,n}^* - \overline
z)\|_2 \to 0$ along an infinite set of $n$s. But
$$
\| (a_{h,n} - z)(a_{h,n}^* - \overline z)\|_2 = (a_{h,n} a_{h,n}^*,1) + |z|^2 -
((z a_{h,n}^* + \overline z a_{h,n}),1).
$$
However, the constant terms in $a_{h,n}$ and its adjoint are both zero. Thus
$\| (a_{h,n} - z)(a_{h,n}^* - \overline z)\|_2 \geq |z|^2$. Since the norms go
to zero along an infinite subset, it follows that $|z|^2 = 0$, so $z = 0$.

\item{(3)} {\it Any strong limit point of $(a_{h,n})$ is central.}

\noindent Choose $\epsilon > 0$ and fix $j \in G$. For all sufficiently large
$n$, $|T_n j \cap T_n| > (1-\epsilon) |T_n|$. Hence $\| a_{h,n} -
j^{-1}a_{h,n}j \| \leq 1- 2\epsilon$ (using the usual norm, not the 2-norm); in
particular, any tracial 2-norm limit point of $(a_{j,n})$ will commute with $j$.
Since $j $ is an arbitrary of $G$, any limit point is central.

Now suppose that $ \RR_{a,l} = \RR_{a',l'} $. As in the proof of xxx, set $k =
l'-l \geq 0$, and assume it is positive; since $f $ is not a zero divisor, we
have $a' = h^{-1} f^k h a$ for all $h \in G$. Set  $H = S^k \setminus \brcs{1}$,
and write $H = \brcs{h_i}$. For $h_1$, there exists a subsequence $T_{n(j,1)}$
of the original F\o olner sequence \st $(a_{h_1,n(j,1)})$ converges strongly to
zero. But $\brcs{T_{n(i,1)}}$ is still a F\o olner sequence, and it follows
inductively that there exists a sequence $(m(i)) \to \infty$ \st
$(a_{h_j,m(i)})$ converges strongly to zero for all $j$. Thus $b_i:=
|T_m(i)|^{-1} \sum_{T_{m(i)}} g^{-1} f^k g$ converges  strongly to $(f^k,1)$ the
constant coefficient in $f_k$.

Since $1 \in \Supp f^k$ and the coefficients of $f  $ are nonnegative, we have
that $(b_i) \to r > 0$, a positive scalar.

Since $a' = b_i a$ for all $i$, we deduce that $a' = r a $ (this equality holds
in $W^*(G)$ but the map $\R G \to W^*(G)$ is an embedding. Since we also have
$a' = f^k a$ (set $h =1$), we deduce that $a = r^{-1}f^k a$. It follows from
characterization of $J$ that $[f^{sk}a,0] \in R_f$ for all positive integers
$s$, so that $[a,0] \in J$. Since our hypothesis is that $J = \brcs{0}$, there
exists $t$ \st $f^t a=0$, and since $f$ is not a zero divisor, it follows that
$a = 0$, a contradiction. Hence $k = 0$, so that $l=l'$ and $a= a'$. \qed

We can drop the amenability assumption if we could prove the following. Here
$C(h)$ will denote the centralizer of the element $h$ of $G$.

\Lem Conjecture. Suppose that $G$ is a countable ICC group, and let $H$ be a
finite subset of $G$ that does not contain $1$. Then there exists an increasing
sequence of finite subsets, $T_n \subset T_{n+1}$ \st for each $h \in H$,
$$
\frac{|gC(h) \cap T_n|}{|T_n|} \to 0 \text{ for all $g \in G$.}
$$

To see that this would be sufficient (with the other hypotheses), we note that
if $a_{h,n} = |T_n|^{-1}\sum_{j\in T_n} j h j^{-1}$, then $\| a_{h,n} a_{h,n}^*
\|_2 \to 0$ (via a direct computation), so that $a_{h,n} \to 0$ strongly in
W${}^*(G)$, and the same concluding argument in the proof above will work. There
is  no need to assume that $\cup T_n = G$. [I'm pretty much convinced that the
conjecture is correct; it is trivial for $|H| = 1$ (enumerate the infinitely
many cosets of $C(h)$, pick one representative from each, and let $T_n$ consist
of the first $n$ of them), and easy to prove for $|H| \leq 3$.]

Unfortunately, without assuming that $\R G$ is a domain, I could not find a
useful set of conditions to guarantee that $f$ is not a zero divisor
(presumably, this is generic)  or when $J = \brcs{0}$. In the latter case, it is
not clear that this depends only on $S = \Supp f$.

\endcomment

We say a group with admissible element $(G,f)$ satisfies EP, the conclusion
of xxx is satisfied, that is, for $a, a' \in \R G$,
$$
\RR_{a,l}  = \RR_{a',l'} \text{ implies either}\cases a= a' = 0 &\text{or}\\
a=a' \text{ and } l=l' .&{}
\endcases
$$
If for all admissible  $f$, the pair $(G,f)$ satisfies EP, then we simply say
that $G$ satisfies EP.

\Lem Lemma. Suppose that $G$ is finitely generated. Sufficient for $(G, f)$ to
satisfy EP are the following conditions.
\item{(a)} $f$ is not a zero divisor in $\R G$, and
\item{(b)} there exists a family of normal subgroups $N_{\alpha}$ \st $\cap
N_{\alpha} = (1)$ and each $G_{\alpha} := G/N_{\alpha}$ satisfies EP.

\Pf Suppose $G$ fails EP; then there exists an admissible $f \in \R G$,
together with $a,a' \in \R G$ \st for all $g \in G$, there exists $m(g),k \in
\N$  with $a \neq 0$ \st $f^{m(g)}(a' - gf^k g^{-1} a) = 0$. Since $f$ is not a zero
divisor, this simplifies to $a' = gf^k g^{-1}a$ for all $g \in G$ (and still
with $a \neq 0$).

Let $\Arrow \pi_{\alpha}; \R G . \R G_{\alpha}$ be the algebra homomorphism
induced by the quotient map $G \to G/N_{\alpha}$ (the kernel of $\pi_{\alpha}$
is $w(N_{\alpha})$). Then the embedding $G \to \prod G_{\alpha}$ entails that
the intersection of the kernels of all $\pi_{\alpha}$ is $\brcs{0}$. Hence there
exists $\beta$ in the index set \st $\pi_{\beta} (a') \neq 0$ (since, in
particular, $a' = f^k a$, and $f$ is not a zero divisor, so $a' \neq 0$). We
observe that $f_{\beta}:= \pi_{\beta} (f)$ is admissible as an element of $\R
G_{\beta}$, and we have $\pi_{\beta}(a') = j f_{\beta}^k j^{-1} \pi(a)$ for all
$j \in G_{\beta}$. It easily follows that the corresponding right multipliers,
$\RR_{\pi_{\beta}(a'), l+k}}$ and $\RR_{\pi_{\beta}(a),l}$ agree for any nonnegative
integer $l$. Since $k > 0$ and both $\pi_{\beta}(a) \neq 0$, we have a
contradiction to $(G_{\beta}, f_{\beta})$ satisfying EP. \qed

To obtain more examples of groups satisfying EP, we consider ICC groups. 

Let $G$ be a group, and let $\brcs{G_i}_{i=1}^n$ be a finite collection of subgroups, each of infinite index in $G$. We denote by a {\it  union of finitely many cosets of $\brcs{G_i}$}, a subset of $G$ of the form
 $$\cup_{i=1^n} \(\cup_{j_i \in J_i} j_i G_i\),
 $$
where each $J_i$ is a finite subset of $G$.
 
 \Lem Lemma. Let $G$ be a group and let $\brcs{G_i}_{i=1}^n$ be a finite collection of subgroups, each of infinite index in $G$. Then $G$ is not a  union of finitely many cosets of $\brcs{G_i}$.
 
  \Pf If $n=1$, this follows from $[G:G_1] = \infty$. We proceed by induction on $n$. Suppose $G = \cup_{i=1^n} \(\cup_{j_i \in J_i} j_i G_i\)$. Since $[G:G_n] = \infty$, there exists $x$ in $G$ \st the coset $xG$ is disjoint from $\cup_{j\in J_n} jG_n$. Hence $xG_n \subset \cup_{i=1^{n-1}} \(\cup_{j_i \in J_i} j_i G_i\)$, whence $G_n$ is contained in a union of finitely many cosets of $\brcs{G_i}_{i \leq n-1}$. It is immediate that any union of finitely many cosets of $G_n$ is contained in a union of finitely many cosets of $\brcs{G_i}_{i=1}^{n-1}$. This forces $G$ to be a union of finitely many cosets of $\brcs{G_i}_{i=1}^{n-1}$ contradicting the induction hypothesis.\qed 
 
  Recall that a group is {\it ICC\/} (infinite conjugacy classes) if every nonidentity element has infinitely many conjugates. Let $Z(h)$ denote the centralizer of the element $h$ of $G$; then $G$ is ICC iff for every $h \neq 1$, $Z(h)$ has infinite index in $G$.
 
  We now work in the W*-algebra of $G$, denoted $W^*(G)$. This is a factor iff $G$ is ICC, and its trace is the linear functional arising from the augmentation map. 
 
  \Lem Lemma. Let $G$ be an ICC group, and let $Q$ be a finite subset of $G \setminus {1}$. Let $a \in \R G$. There exists a strictly increasing nested sequence of finite subsets $T_n \subset T_{n+1}$ of $G$ \st 
 $$
  a_n:= \frac 1{|T_n|} \sum_{g \in T_n} gag^{-1}
 $$
  converges weakly (in $W^*(G)$) to the constant coefficient of $a$.
 
  \Pf Write $a = r1 + \sum_{q \in Q} r_q q$ (so $Q \cup \brcs{1}$ is the support of $a$). Set $G_q = Z(q)$; these are of infinite index in $G$ since $G$ is ICC. 
  
  We will inductively construct $T_n$ so that on defining $\Cal G_n$  by $\Set{jqj^{-1}}{j \in T_n; q \in Q}$, if  $x \in \Cal G_n $  and $x = jqj^{-1}$ for some $j \in \cup_{m\geq 1} T_m$  and $q \in Q$ then $j \in T_n$. 
 
Set $T_1 = \brcs{1}$, and assume $T_n$ has been defined. We claim there exists $g \in G$ \st $gQg^{-1} \cap \Cal G_n = \emptyset$. This condition is equivalent to $g \not\in \cup_{q \in Q; j \in T_n} jZ(q)$. By the lemma, such $g $ exists. So pick one choice for $g$ and  set $T_{n+1} = T_n \cup \brcs{g}$. This yields the strictly increasing nested sequence $\brcs{T_n}$ and with $|T_n| = n$. 
 
  We verify the induction hypotheses. Suppose $jqj^{-1} = j' q'( j')^{-1}$ for some $j \in T_n$, $j' \in T_m$ with $m > n$, and $q,q' \in Q$. Let $n'$ be the smallest integer \st $j' \in T_{n'}$. If $n'\leq n$, there is nothing to do. If $n' > n$, then $j'$ was obtained as the new element \lq\lq$\brcs{g}$\lq\lq\ with the property that $j'Q(j')^{-1} \cap \Cal G_{n'-1} = \emptyset$, yielding a contradiction (since $\Cal G_n \subseteq \Cal G_{n'}$. 
 
  Now we can show that the weak limit of the $a_n$ (defined by the now-constructed sequence $\brcs{T_n}$ converges as indicated. 
 
  Let $h \in Q \cup \brcs{1}$, and define $b_{h,n} = |T_n|^{-1}\sum_{j \in T_n} jhj^{-1}$. Since $b_{1,n} =1$ for all $n$, it suffices to show that $(b_{q,n})$ converges weakly to zero for all $q$ (since the collection of finite convex linear sums of elements of $G$ is weakly dense in the unit ball of $W*(G)$, and the unit ball is weakly metrizable). This is equivalent to showing that for all $k \in G$, $\langle b_{q,n}, k\rangle \to 0$, that is (since the inner product is given by the trace), the coefficient of $k$ in $b_{q,n}$ tends to zero as $n \to \infty$. We have 
  $$\eqalign{
 \langle b_{q,n}, k\rangle & = \frac{\#\(\Set{j \in T_n}{jqj^{-1} = k}\)}{|T_n|}. \cr
 }$$ 
  If $k \not\in \Cal G_n$ for any $n$, the numerators are all zero. If $k \in \Cal G_m$, then the numerator is simply $\#\(\Set{j \in T_m}{jqj^{-1} = k}\)$ for all $n \geq m$. Since $|T_n| \to \infty$, the limit is zero in all cases. \qed

\Lem Lemma. Let $G$ an ICC group, and let $b$ be an element of W$^*(G)$  \st the tracial moment sequence $(\tr b^n)_{n\geq 0}$ is $(1,r,r, r, \dots , r, \dots)$ for some $0 \leq r $. Then $\Spec b$ is a subset of $\brcs{0,1}$. 

 \Pf Let $B$ be the (commutative) unital Banach algebra generated by $b$, and form its Gelfand map $B \to C(\spec b)$. The moment sequence of $\hat b$ is the same one, and of course, there is a complex measure whose moment sequence is that---$(1-r) \delta_0 + r\delta_1$. Since the polynomials are dense in  $C(\Spec b)$, it follows that there is only one possible complex measure yielding the moment sequence. Hence $\Spec b \subset \brcs{0,1}$.\qed

 \Lem Proposition. Suppose the countable discrete group  $G$ is ICC, and $f$ is an admissible element of $\R G$ for which there exists an $n \geq 0$ \st the right annihilator of $f^n$  equals that of $f^{n+1}$. Then  $(G,f)$ satisfies EP. 
 
\Pf   For each $g \in G$, there exists $m(g)$ \st $f^{m(g) + k }g a = f^{m(g)}a'$, whence $f^{m(g)}  (f^k g a - g a') = 0$. Hence for all $g$, $f^{n+k} ga = f^n g a'$, and we can replace $n$ by any larger integer. This yields $g^{-1}f^{n+k}g \cdot a = g^{-1} f^n g a'$ for all $g$ and all sufficiently large $n$. 
 
  For $n$ fixed, set $J = \Supp f^{n+k}$. In particular, $J$ contains $\Supp f^n$. By xxx, there exists a strictly increasing sequence of finite sets $\subset T_j \subset T_{j+1} \subset \dots$ \st $|T_j|^{-1}\sum_{t \in T_j} t^{-1}f^{n+k}t $ converges weakly to $(f^{n+k},1)$ and $|T_j|^{-1}\sum_{t \in T_j} t^{-1}f^{n}t $ converges weakly to $(f^n,1)$. Set $r_m = (f^m,1) \neq 0$ (since $1 \in \Supp f^m$ for all positive $m$).
 
  We thus have $r_{n+k}a = r_{n}a'$, for all sufficiently large $n$. In particular, assuming $a$ (or $a'$) is nonzero, $r_{n+k}/r_n$ is constant in $n$ (for all sufficiently large $n$). Fix some sufficiently large $m$, and we  have $r_{lm} = r_{m}s^{l-1}$, where $s = r_{2m}/r_m$. Hence the tracial moment sequence for $b:= f^m/s$ is $(1,r,r,r,r \dots$, where $r = \langle f^m,1\rangle$. By   lemma xxx, $b -b^2$ is quasinilpoent. Hence, given $\epsilon > 0$, for all sufficiently large $n$, we have $\| (b-b^2)^n\|^{1/n} \leq \epsilon$, or what amounts to the same thing, $\| (b-b^2)^n\| \leq \epsilon^n$. 

 Now consider $\langle (b(1-b)^n 1,1 \rangle$; this is just $\sum_{j=0}^n \langle b^{n+j}(-1)^j {n\choose j}1,1\rangle$; but since all the tracial moments are $r$, we obtain the value $(r(1-r))^n$. It follows that $|(r(1-r))^n| \leq \| (b-b^2)^n\| \leq \epsilon^n$, and so $|r(1-r)| \leq \epsilon$. Since this is true for all $\epsilon$, either $r = 0$ or $r = 1$. 

 Finally, $r = 0$ entails that $\langle f^m,1\rangle = 0$ which is obviously false (since all the coefficients of $f$ are nonnegative, and $1$ is in the support). If $r = 1$, then $\tr f^{ml} = s^{l-1}$, but this is impossible, since $\tr f^{ml} \geq (\tr f^m)^l + \langle f^m,g \rangle \langle f^{m(l-1)},g^{-1}\rangle$ for any $g \in \Supp f^m$. For any such $g$, there exists $l$ \st $g^{-1} \in \Supp f^{m(l-1)}$ (since $f^m$ is admissible), hence the inequality is (eventually) strict, a contradiction. \qed
 
  If $f$ is not a right zero divisor and $G$ is ICC, there is a much simpler argument available; in that case the tracial moment sequence reduces to $(1,r,r^2, \dots)$ and the admissibility argument of the last paragraph precludes this.

 The condition that the right annihilators of increasing powers, $r(f^n)$, eventually become stationary, is interesting. if $f $ is symmetric (that is, $\langle f,g\rangle = \langle f,g^{-1} \rangle$ for all $g \in \Supp f$, then the right annihilator of $f$ equals that of $f^2$ (and thus of all higher powers)---more generally, this occurs whenever $f$ is a normal element of W$^*(G)$ (so this also applies if $G$ is abelian and $f$ is an arbitrary element). 
 

\SecT Quotient group $R_f/\SS R_f$

The quotient group $R_f/{\SS R_f}$ in that example is interesting. Assume that $A = \Z G$ rather than $\R G$. Define, as before $\Gamma_n$ to be the set of $g$ \st $g \in \supp f^n \setminus \cup_{i=1}^{n-1} \supp f^{n-i}$. It is easy to verify that $\Gamma_n = \brcs{x^{\pm n}, gx^{\pm(n-1)}, \sigma x^{\pm(n-1)}, g\sigma x^{\pm(n-2)}, \sigma g x^{\pm(n-2)} ,\sigma^2 g\sigma x^{\pm(n-2)} }$.

The induced maps $\Z \Gamma_n \to \Z \Gamma_{n+1}$ split in two, both given by the matrix (using a basis coming from $(1,g,\sigma g,\sigma,g\sigma, \sigma^2$),
$$
A = \(\matrix 1 & 0 & 0 & 0 & 0 & 0 \\
1 & 1 & 0 & 0 & 0 & 0 \\
0 & 1 & 1 & 0 & 0 & 0 \\
1 & 0 & 1 & 1 & 0 & 0 \\
0 & 0 & 0 & 1 & 1 & 0 \\
0 & 0 & 0 & 1 & 0& 1 \\
\endmatrix\).
$$
Thus $R_f /{\SS R_f}$ is order isomorphic to the direct sum of two copies of $H = \lim \Arrow A; \Z^6. \Z^6$ (a stationary dimension group). Each $H$ has unique trace (which corresponds to the two pure traces in $F_0$) given by the left eigenvector $(\matrix 1 & 0 & 0 & 0 & 0 & 0\endmatrix)$, the only nonnegative eigenvector of $A$, and in fact, $H$ is the lexicographic direct sum of $H_0$ with $\Z$, and $H_0$ decomposes into an order direct sum of $\Z \oplus_{\text{lex}} \Z$ and $\Z^2 \oplus_{\text{lex}} \Z$.

In contrast, if we have taken $f = (1+x+x^{-1})(1 + g + \sigma + g\sigma + \sigma g + \sigma^2)$, then the corresponding matrix is
$$
A = \(\matrix 1 & 0 & 0 & 0 & 0 & 0 \\
1 & 1 & 0 & 0 & 0 & 0 \\
1 & 0 & 1 & 0 & 0 & 0 \\
1 & 0 & 0 & 1 & 0 & 0 \\
1 & 0 & 0 & 0 & 1 & 0 \\
1 & 0 & 0 & 0 & 0& 1 \\
\endmatrix\),
$$
which yields $\Z^5 \oplus_{\text{lex}} \Z $.

Miscellaneous example

If $G$ is the infinite dihedral group, then $G$ has no real characters, but
depending on the choice of admissible $S$, has $0,1,2,3,$ or $4$ multiplicative
traces in $
\partial_e F_0$.

Other stuff

Obtain as corollary to stuff about nilpotent groups usually not satisfying the density condition

\Lem Corollary. If $G$ is nilpotent and $\cup_{\lambda > 0} \partial_e F_{\lambda}$ is dense in $S(R_f,[1,0])$ for some choice of admissible set, then $G$ is central by finite.

\SecT Other stuff


 We let $\langle \cdot \rangle$ denote the order ideal in $R_f$ generated by the contents of the angle brackets.
some sample computations

\noindent {$G = D_{\infty}$.} We write $G = \langle g \rangle \times_{\theta} \langle h \rangle$, where $h^2 = 1$ and $hgh = g^{-1}$. Then $G$ has centre $\langle g^2 \rangle$, which is of index two in $ \langle g \rangle $, and we have a normal form for elements of $G$, $g^i h^j$ where $i \in \Z$ and $j \in \brcs{0,1}$. Let $A$ be either $\Z$ or $\R$, and observe that $A G = A\langle g \rangle \oplus h A \langle g \rangle$, as right $A\langle g \langle$-modules; we identify the direct summands with the free submodules of $(A \langle g \rangle)^2 \iso A[x^{\pm1}]^2$ generated by the standard basis elements, $\( \smallmatrix 1 \\ 0 \\ \endsmallmatrix\)$ (corresponding to ${}{1_G}$) and $\( \smallmatrix 0 \\ 1 \\ \endsmallmatrix\)$ (corresponding to ${}{h}$).

Left multiplication by $g$ is thus implemented by the matrix $\( \smallmatrix x & 0 \\ 0 & x^{-1}\endsmallmatrix\)$ (coming from $g\cdot 1_G = g$ and $gh = hg^{-1}$, respectively) and left multiplication by $h$ corresponds to the matrix $\( \smallmatrix 0 & 1 \\ 1 & 0 \\ \endsmallmatrix\)$. More generally, left multiplication by the element $g^j h$ is implemented by
$$
\(\matrix 0 & x^j \\ x^{-j} & 0 \\ \endmatrix
\).
$$

Let $f = p_1(g) + p_2 (g) h$ where $p_1$ and $p_2$ are Laurent polynomials with coefficients from $A^+$. Then left multiplication by $f$ is implemented by the matrix $\Arrow M(f); A[x^{\pm 1}]^2. A[x^{\pm 1}]^2$,
$$
M(f) = \(\matrix p_1(x) & p_2(x) \\ p_2 (x^{-1}) & p_1(x^{-1}).
\\ \endmatrix
\).
$$

Then $A_f$ is identified with the direct limit (as ordered $A[x^{\pm1}]$-modules),
$$
M_f:= A[x^{\pm1}]^2 \to A[x^{\pm1}]^2 \to A[x^{\pm1}]^2 \to \dots,
$$
all maps given by $M(f)$. The element $\Bf 1$ is identified with $z:= [\( \smallmatrix 1 \\ 0 \\ \endsmallmatrix\),0]$ in this direct limit, and thus $R_f$ is identified with the order ideal in $M_f$ generated by $z$. Then $z = [c_k,k]$, where $c_k$ is the first column of $M(f)^k$.

This puts us in the situation of $[xxx]$, that is, we have a matrix-valued random walk. By xxx, the pure faithful traces are given as follows. Let $r \in (\R)^{++}$, and form the matrix $M(f) (r)$ obtained by evaluating the entries (which are Laurent polynomials in $x$) at $x\mapsto r$. Unless $f$ is really degenerate, the resulting real matrix $M(f)(r)$ will be at least primitive (hence some power will be strictly positive), hence will admit a left Perron eigenvector $v(r)$ (in fact, it is easy to check that we can find an algebraic $v$ with $r \mapsto v(r)$ real analytic), with corresponding eigenvalue $\beta (r)$. Then we can define a trace on $M_f$ (and thus on $R_f$) via
$$
\tau_r: [w,k] \mapsto \frac{v(r)w(r)}{\beta(r)^k}.
$$
Provided $M(f)$ is not upper or lower triangular or has two zero entries, these are all pure faithful traces, and every faithful pure trace is of this form [xxx]. The conditions $\cup \supp f^k$ and $1_G \in \supp f$ more than yield this.

In some cases, the resulting system is equivalent to one arising from a simple random walk on $\Z$ (equivalent, in this instance, refers to strong shift equivalence in the category of matrices with entries from $A^+$).
more examples here

The end is nigh

Here $f \in A^+$, $\supp f$ generates $G$ as a semigroup, and $1_G \in \supp f$ (or slightly weaker, $\gcd\Set{k}{1_G \in \supp f^k} = 1$. In this section, $A = \Z G$, the integral group ring. We have the shift $\Arrow \SS; A_f. A_f$ given by $\SS[a,k] = [a,k+1] = R_{(1,1)}([a,k])$. This is an order-automorphism of $A_f$, and an order endomorphism of $R_f$; in the latter case, it has the additional property (for endomorphisms) that $\SS b \in R_f^+$ implies $b \in R_f^+$.

Then ${\SS R_f}$ is just the order ideal generated by $[1_G,1]$ in $R_f$, and we also have (from somewhere) that if $\tau$ is a pure trace that is not faithful, then $\tau \circ \SS = 0$ (proof: by purity, $\tau \circ \SS = \tau([1_G,1])\tau$, and if $\tau([1_G,1]) = s \neq 0$, then $\tau$ extends to a trace on $A_f$ (necessarily on a pure ray of traces---note that $A_f$ has no order unit, so purity of traces can only be discussed \wrt rays of traces) via $\tau([a,k]) = s^{k-l} \tau([a,l])$ whenever $[a,l] \in R_f$ (this is easily checked to be well-defined, etc), and thus is automatically faithful.

Let $\tau$ be a trace on a dimension group $J$; define $\ker^+\tau$ to be the largest order ideal contained in $\ker \tau$ (in dimension groups, sums, finite intersections, and unions of chains of order ideals are order ideals, so the largest one exists and is unique); this is the order ideal generated by $\ker \tau \cap J^+$. Obviously, $\tau$ is faithful iff $\ker^+ \tau = \brcs{0}$.

It follows that if $\tau$ is a pure nonfaithful trace of $R_f$, then $S R_f \subseteq \ker ^+\tau$, and thus
$$
\SS R_f \subseteq \cap \Set {\ker^+ \tau}{\tau \text{ is a pure nonfaithful trace of $R_f$}}.
$$
In general, the inclusion is strict.

What can we say about the simple quotients of $R_f/S R_f$ (by order ideals); equivalently, the simple quotients (by order ideals) of $R_f$? We first examine a few examples, relatively well-studied. These suggest a conjecture.

\Lem Example (a) $G = \Z^d$. In this case, the Legendre transformation, here equivalent to the weighted moment map yields an identification between the pure nonfaithful traces of $R_f = R_P$ where $f = P = \sum c_w x^w$ (monomial notation; $c_w \in \Z^+$) and the boundary of the Newton polytope of $P$, the latter being $K(P)=\cvx \Set{w}{c_w \neq 0}$ in $\R^d$. The pure traces \st $\ker ^+ \tau$ is a maximal order ideal are precisely those corresponding to the vertices of $K(P)$. In particular, there are only finitely many
maximal order ideals of $R_P = R_f$, and the quotients (the factors) are precisely the rank one subgroups of $\Q$ given $\Z[1/c_w]$, one for each extreme point $w $ of $K$. Some of the $c_w$ could of course be one, and we obtain discrete pure traces in these cases.

\Lem Example. (b) $G = \pmb F_2$, and $f = 1 + g + {}{g^{-1}} + {}h + {}{h^{-1}}$.

In this case, it is easy to verify that the pure traces of $R_f/{\SS R_f}$ correspond precisely to the set of ends of $G$, and in fact $R_f/S R_f$ is just $C(X, \Z)$ where $X$ is the path space of the obvious tree. This is one of the few cases where $S R_f \subseteq \cap \Set {\ker^+ \tau}{\tau \text{ is a pure nonfaithful trace of $R_f$}}$ can be shown.

However, if we change the support of $f$, relatively drastic things can happen, cf Example xxx (section on maximal order ideals).

\comment
\Lem Lemma. Suppose $g $ is an element of $G$
\st $[g,n] \in R_f$. Then $\Arrow \phi_g:= \RR_{(g,n)}; R_f .R_f$ is a positive, one to one endomorphism of $R_f$ \st $\phi_g^{-1}(R_f^+) = R_f^+$; moreover, it is an order-isomorphism, $R_f \to \langle [g,n] \rangle$.

\Pf We have that $\Arrow \phi_g; R_f .R_f$ is an order-preserving endomorphism.
\noindent {(0)} {$\it \phi_g$ is one to one and $\phi^{-1}(R_f^+) \subset R_f^+$.} There exists $m$ \st $g^{-1} \in \Gamma_m$, and it follows that $ \phi_{g^{-1}} \phi_g = \phi_g\phi_{g^{-1}} = \RR_{1,m+n} = S^{m+n}$. As $S$ is one to one, so $S^{m+n}$, and thus so is $\phi$. If $\phi_g(x) \in R_f^+$, then $S^{m+n}(x) \in R_f^+$, and if we write $x = [a,k]$, then $S^{m+n}(x) = [a,{m+n + k}]$; as this is in $R_f^+$, there exists $t$ \st $f^t a \in A^+$, and therefore $[a,k] \in A_f^+\cap R_f = R_f^+$.

\noindent {(1)} {\it $\phi_g (R_f) = \langle [g,n]\rangle$.} For any $k$ and $a \in A$ \st $\supp a \subseteq \supp f^k$, $\phi_g ([a,k]) = [ag, k+n]$; as $\supp ag =( \supp a)\cdot g$, and therefore $\supp ag \subset (\supp f^k )g$. Thus there exists $K$ \st both $\pm[ag,k+n] \leq K[f^k g, k+n] = K[g,n]$. Hence $\phi_g ([a,k]) \in \langle [g,n]\rangle$, so $\phi_g(R_f) \subseteq \langle [g,n]\rangle$.

Select $[a,k] \in \langle [g,n]\rangle$; we can assume that $\supp a \subset f^k$. This entails there exists a positive integer $K$ \st both $\pm [a,k] \leq K [g,n]$. This in turn implies there exists $m$ \st $K f^{m+k}g \pm f^{m+n} a \in A^+$, and this yields $\supp (f^{m+n}a) \subset \supp (f^{m+k}g) = (\supp f^{m+k})g$. Set $b = a{}{g^{-1}}$, so $\supp f^{m+n}b \subset \supp f^{m+k}$. Thus $[f^{m+n}b, m+k] = [f^n b,k] \in R_f$, and moreover, $\phi_g([f^n b,k]) = [f^n b g, k+n] = [f^n a, k+n] = [a,n]$. Hence $\phi(R_f)$ contains $\langle [g,n]\rangle$ in its image, and thus (1) holds.

Now (0) and (1) together yield the conclusion.
\qed
\endcomment

Only rarely is $\RR_{g,n}$ an automorphism of $R_f$.

Let $N$ be an ordered module over $R$, itself a commutative ordered ring with $1$ as an order unit. First, we want conditions on limits to be rank 1 over their endomorphism ring.

Let $M$ be a primitive matrix over $A = \Z[x_i^{\pm1}]$, with a projectively faithful assumption ([Hxxx], see stuff about faithfulness). Form $G_M$; under what conditions is it a rank one module over $E(G_M)$, essentially $C_Z (M)[\hat M^{-1}]$? (The latter is not commutative, but is generically.) Next, look at orders, and the fact that under modest conditions, $I\cdots E(G_M) = G_M$ and the bounded subring is an order in $E(G_M)$.

If $E(G_M)$ is a commutative domain (as occurs precisely when the characteristic polynomial of $M$ (as a polynomial with coefficients from $A$) is irreducible over $A$, modulo powers of $X$, then it should be very easy to show that $G_M$ is rank one over $E$.

Let $\tau$ be a trace on $(I ,w)$, an order ideal in $G_M$, having its own order unit, and let $E_b$ be the bounded subring of $E(G_M)$; then, $E_b$ is also the bounded subring of $E(I)$, and so $I$ is an $E_b$-module. Hence $\phi(\tau) = L_{\tau}$, defined by $L_{\tau} (e) = \tau (ea)$ is a normalized trace on $E_b$, and if $\tau$ is pure, then $L_{\tau}$ is multiplicative.

We would like something like $\tau$ is order unit good (as a trace on $I$) iff $L_{\tau}$ is order unit good, particularly if $\tau$ is pure (so $L_{\tau}$ will then be pure too). In one direction, we note that $\ker L \cdot I \subset \ker \tau$. Thus $I/\ker \tau $ is a quotient of $I_L:= I /(\ker L \cdot I)$, and the latter is an ordered $E_b/\ker L$-module (so is $I /\ker \tau$, but we have less information about it). If $L$ is an order unit good trace of $E_b$, then $E_b / \ker L := R_L$, with the quotient ordering is an ordered subring of the reals (\wrt the relative ordering). If now $I$ is known to be rank one over $E_b$, then $I /(\ker L \cdot I)$ should be rank one over $R_L$ (this may require using the $\hat M^{-1}$ stuff).

Anyway, if the rank one thing holds, then up to order isomorphism, $R_L \subseteq I_L \subseteq K \subset \R$, where $K$ is the field of fractions of $R_L$ (inside the reals). This has no proper quotients which have a partially ordered structure; it would then follow that $\ker L \cdot I = \ker \tau$. It should then follow that the closure of $\ker \tau$ in the affine representation of $(I,w)$ should be (a) a real vector space, and (b) dense in $\tau^{\vdash}$.

This would give $L_{\tau}$ order unit good implies $\tau$ is order unit good (in the case that of suitable $G_M$, and presumably more generally when the module is rank one. This apparently requires that the centralizer be a domain. Perhaps if the characteristic polynomial is not irreducible, then the large eigenvalue function corresponds to just one of the factors, and we work over the thing module whatever stuff there is ...?

danger, danger, danger. If $L_{\tau} = L_{\tau'}$, then we are sunk (at least with this method), since $\tau$ and $\tau'$ will both kill $\ker L \cdot M$. However, for faithful pure traces $\tau$, it's probably OK (uniquely determined by $L$, need a bit more ...)

Really need I iso to R-ideal, as there is a problem with general result about rank one for quotients ...

For this, $P$ is a prime ideal in $R$.
Or, write $M = \cup M_i$ where each $M_i$ is cyclic over $R$ (and free?). Then $MP = \cup M_i P$, so ?? $M/MP \iso \lim M_i/M_iP$, a limit of cyclic $R/P$-modules, therefore rank one. Replace rank one by limit of free of rank one modules.

Did we get, if at each vertex, there exists exactly one nonnilpotent block, then the map on traces should be an affine homeomorphism?

Easy to construct thingies with $E_b$ not having faithful traces dense, e.g.,
$$
M = \( \matrix 2x + 1 & 1 \\ 1 & 1+x \\
\endmatrix\)
$$
Then the left eigenvector is given by $(\beta_0, 1)$ where $\beta_0 (x) = x/2 + \sqrt{x^2 +4}/2$. Then $\beta_0(x)/(1+x) \to 1$ (as $x \to \infty$), so the faithful pure traces (corresponding to $x \mapsto r \in \R^{++}$ as $r \to \infty$) converge to the trace given by $(1, 0)$ on $\lim_{x\to \infty }M/P = \diag(2,1)$, which yields one of the pure traces on the dimension group, but not the other one, corresponding to $(0,1)$. The latter also yields a pure trace, as is easy to see, and is not a limit of pure faithful ones. This is a pretty generic example. But it is likely that as in the previous question, if there is just one block per face (or merely per vertex? probably not), the faithful pure traces are dense in the pure trace space.

^^
This is an improvement on the earlier result for ideals, and makes the result transparent.

All modules are unital (meaning $m\cdot 1 = m$ for $m$ in the module, and $1$ the identity of the coefficient ring).

\Lem Lemma. Let $R$ be a partially ordered commutative ring with $1$ as order unit. Let $M$ be a partially ordered $R$-module. If $R$ is approximately divisible, then $M$ is nearly divisible.

\Pf Suppose $m \in M^+\setminus \brcs{0}$. We may find order units $a$ and $b$ \st $1 = 2a + 3b$. In particular, there exists $K$ \st $1 \leq Ka, Kb$ in $R$. Then we have $m = 2(ma) + 3(mb)$ and $m \leq K(ma), K(mb)$, so $M$ is nearly divisible.
\qed

--------------------------
more other stuff

Let $N$ be an ordered module over $R$, itself a commutative ordered ring with $1$ as an order unit. First, we want conditions on limits to be rank 1 over their endomorphism ring.

Let $M$ be a primitive matrix over $A = \Z[x_i^{\pm1}]$, with a projectively faithful assumption ([Hxxx], see stuff about faithfulness). Form $G_M$; under what conditions is it a rank one module over $E(G_M)$, essentially $C_Z (M)[\hat M^{-1}]$? (The latter is not commutative, but is generically.) Next, look at orders, and the fact that under modest conditions, $I\cdots E(G_M) = G_M$ and the bounded subring is an order in $E(G_M)$.

If $E(G_M)$ is a commutative domain (as occurs precisely when the characteristic polynomial of $M$ (as a polynomial with coefficients from $A$) is irreducible over $A$, modulo powers of $X$, then it should be very easy to show that $G_M$ is rank one over $E$.

Let $\tau$ be a trace on $(I ,w)$, an order ideal in $G_M$, having its own order unit, and let $E_b$ be the bounded subring of $E(G_M)$; then, $E_b$ is also the bounded subring of $E(I)$, and so $I$ is an $E_b$-module. Hence $\phi(\tau) = L_{\tau}$, defined by $L_{\tau} (e) = \tau (ea)$ is a normalized trace on $E_b$, and if $\tau$ is pure, then $L_{\tau}$ is multiplicative.

We would like something like $\tau$ is order unit good (as a trace on $I$) iff $L_{\tau}$ is order unit good, particularly if $\tau$ is pure (so $L_{\tau}$ will then be pure too). In one direction, we note that $\ker L \cdot I \subset \ker \tau$. Thus $I/\ker \tau $ is a quotient of $I_L:= I /(\ker L \cdot I)$, and the latter is an ordered $E_b/\ker L$-module (so is $I /\ker \tau$, but we have less information about it). If $L$ is an order unit good trace of $E_b$, then $E_b / \ker L := R_L$, with the quotient ordering is an ordered subring of the reals (\wrt the relative ordering). If now $I$ is known to be rank one over $E_b$, then $I /(\ker L \cdot I)$ should be rank one over $R_L$ (this may require using the $\hat M^{-1}$ stuff).

Anyway, if the rank one thing holds, then up to order isomorphism, $R_L \subseteq I_L \subseteq K \subset \R$, where $K$ is the field of fractions of $R_L$ (inside the reals). This has no proper quotients which have a partially ordered structure; it would then follow that $\ker L \cdot I = \ker \tau$. It should then follow that the closure of $\ker \tau$ in the affine representation of $(I,w)$ should be (a) a real vector space, and (b) dense in $\tau^{\vdash}$.

This would give $L_{\tau}$ order unit good implies $\tau$ is order unit good (in the case that of suitable $G_M$, and presumably more generally when the module is rank one. This apparently requires that the centralizer be a domain. Perhaps if the characteristic polynomial is not irreducible, then the large eigenvalue function corresponds to just one of the factors, and we work over the thing module whatever stuff there is ...?

danger, danger, danger. If $L_{\tau} = L_{\tau'}$, then we are sunk (at least with this method), since $\tau$ and $\tau'$ will both kill $\ker L \cdot M$. However, for faithful pure traces $\tau$, it's probably OK (uniquely determined by $L$, need a bit more ...)

Really need I iso to R-ideal, as there is a problem with general result about rank one for quotients ...

For this, $P$ is a prime ideal in $R$.
Or, write $M = \cup M_i$ where each $M_i$ is cyclic over $R$ (and free?). Then $MP = \cup M_i P$, so ?? $M/MP \iso \lim M_i/M_iP$, a limit of cyclic $R/P$-modules, therefore rank one. Replace rank one by limit of free of rank one modules.

Did we get, if at each vertex, there exists exactly one nonnilpotent block, then the map on traces should be an affine homeomorphism?

Easy to construct thingies with $E_b$ not having faithful traces dense, e.g.,
$$
M = \( \matrix 2x + 1 & 1 \\ 1 & 1+x \\
\endmatrix\)
$$
Then the left eigenvector is given by $(\beta_0, 1)$ where $\beta_0 (x) = x/2 + \sqrt{x^2 +4}/2$. Then $\beta_0(x)/(1+x) \to 1$ (as $x \to \infty$), so the faithful pure traces (corresponding to $x \mapsto r \in \R^{++}$, $r \to \infty$) converge to the trace given by $(1, 0)$ on $\lim_{x\to \infty }M/P = \diag(2,1)$, which yields one of the pure traces on the dimension group, but not the other one, corresponding to $(0,1)$. The latter also yields a pure trace, as is easy to see, and is not a limit of pure faithful ones. This is a pretty generic example. But it is likely that as in the previous question, if there is just one block per face (or merely per vertex? probably not), the faithful pure traces are dense in the pure trace space.

^^
This is an improvement on the earlier result for ideals, and makes the result transparent.

All modules are unital (meaning $m\cdot 1 = m$ for $m$ in the module, and $1$ the identity of the coefficient ring).

\Lem Lemma. Let $R$ be a partially ordered commutative ring with $1$ as order unit. Let $M$ be a partially ordered $R$-module. If $R$ is approximately divisible, then $M$ is nearly divisible.

\Pf Suppose $m \in M^+\setminus \brcs{0}$. We may find order units $a$ and $b$ \st $1 = 2a + 3b$. In particular, there exists $K$ \st $1 \leq Ka, Kb$ in $R$. Then we have $m = 2(ma) + 3(mb)$ and $m \leq K(ma), K(mb)$, so $M$ is nearly divisible.
\qed

---------------------------

\SecT Infinitesimals

When is $\Inf R_f$ nonzero? Here $1 \in \supp f$, $\supp f$ is finite and generates $G$ as a semigroup.

\Lem Examples.

If $G$ is a finite group, then $\Inf R_f =0$ implies $f$ is a scalar multiple of $\sum_{g} g$.

If $G = \Z^d$, then $\Inf R_f = \brcs{0}$ (much more is known of course, but in this case, it's just a consequence of the fact that a polynomial in $d$ variables that vanishes on an open set in $\R^d$ is constant.

If $G$ is the infinite dihedral group or $\Z \times \Z_2$, there exist some choices of $f$ for which $\Inf R_f$ is nonzero. This holds more generally for any $\Z^d$ by finite group, but currently at least, the $f$ has to be chosen so that the characteristic polynomial of the corresponding matrix-valued random walk (a matrix with entries in the Laurent polynomial ring in $d$ variables) is reducible (after ignoring zero eigenvalues) over the Laurent ring; this is far from generic. If is possible that if the characteristic polynomial is irreducible over the Laurent polynomial ring (or irreducible after throwing out any zero eigenvalues), then $\Inf R_f$ is trivial. Of course, irreducibility is generic in more than one variable.

Much more surprisingly, if $G$ is the class two nilpotent group with generators and relations $uv=zvu$ and $z$ is central, and $f = 1 + u + v + u^{-1} + v^{-1}$, then $x = [{}{uv},3] - [{}{vu},3] $ is a nonzero element of $\Inf R_f$. It is easily seen to be nonzero; since $G$ is nilpotent, all pure faithful traces are given by real characters, so that $\tau (x) = 0$ for all such $\tau$. The remaining pure traces kill ${\SS R_f}$, and since both $uv$ and $vu$ appear in $\supp f^2$, we can write $x = Sy$ for $y = [{}{uv}-{}{vu},2] \in R_f$, whence $x$ is killed by the remaining pure traces. This generalizes to practically all finitely generated nilpotent groups. The only property used of nilpotent groups is that every eigenfunction arises from a real character.

\Lem Lemma. Let $G$ be a nilpotent group, and let suppose that $1 \in \supp f$ and $\supp f$ generates $G$ as a semigroup. Let $g,h$ be elements of $G$ \st $gh \neq hg$ and {\it either\/} of the following two conditions hold.
\item{(a)} The commutator $ghg^{-1}h^{-1}$ has infinite order, or
\item{(b)} $f$ is not a zero divisor in the group ring $\R G$.{\par}\noindent
If $gh,hg \in \supp f^j$, then for all $k > j$, the element $[{}{gh}- {}{hg},k] \in \Inf R_f$ and is nonzero.

\Rmk For nilpotent groups, the existence of such $g,h$ is extremely likely.

\Pf Since $\supp f$ generates $G$ as a semigroup and $1 \in \supp f$, there exists $j$ \st $gh,hg \in \supp f^j$. Thus $[{}{gh} - {}{hg},j] \in R_f$, and so for all $k > j$, it follows that $y_k:= [{}{gh}-{}{hg},k] \in {\SS R_f}$. Hence for all $\sigma \in \partial_e F_0$, $\sigma (x_k) = 0$. On the other hand, for every pure $\tau \in F_{\lambda}$ with $\lambda \neq 0$, a theorem of Margulis asserts that $\tau = \tau_{\chi}$ for some positive real character $\chi$ of $G$. Obviously $\tau_{\chi}(y_k) =0$. Hence $x_k$ is killed by every pure trace of $R_f$, so that $y_k \in \Inf R_f$.

To check that all $y_k$ are nonzero, $y_k = 0$ would entail that there exists $m$ \st $f^m* (gh-hg) = 0$ in $\R G$. Since $gh \neq hg$, (b) yields a contradiction, while if (a) holds, we can rewrite the equation as $-f^m * (1-ghg^{-1}h^{-1})* hg = 0$, whence $f^m*(1-z) = 0$ for $z$ a nonzero element of $G$ having infinite order. But it is well known that in a group ring, such $1-z$ is not a zero divisor. This would force $f^m = 0$ which is clearly impossible (as all of its nonzero coefficients are positive).
\qed

Curiously, if $G = \pmb F_2$, the free group on $\brcs{g,h}$, and $f = 1 + g + g^{-1} + h + h^{-1}$, then $R_f$ has no infinitesimals. To see this, suppose that $p \in \R G$ and $p \prec f^k$ but $p \not\prec f^{k-1}$. There thus exists a reduced word $w$ of length $k$ with nonzero coefficient. There exists a path in the Cayley graph of $G$ hitting $w$ at level $k$, so there exists a corresponding pure trace in $F_0$ sending $[{}w,k] \to 1$, but for all other words of length $k$, $v$, $[{}v,k] \to 0$; moreover, if $v \in \supp p$ and has length less than $k$, $[{}v,k] = S[{}v,k-1]$, so $[{}v,k] \to 0$ (as $\tau \in F_0$). Hence the value of $[p,k]$ under this trace is exactly the coefficient of $w$ in $p$, which by hypothesis is nonzero. So the general element $[p,k]$ is not an infinitesimal.

This remains true if we replace $f$ by $f'$ where $\supp f' = \supp f^k$ for some $k$, but for other replacements, it is not clear what can happen.

  It is reasonable to conjecture that for all finitely generated groups $G$, there is a bound on the length of strictly increasing sequences of annihilators, $(r(f^n))$ when $f$ is admissible, or even a uniform bound over all groups $1$?). 

 Notes to me: work out $\R S_4$ if there exists admissible thingy with length 2. 

Some properties discussed here are unlikely to be residual, but are almost so.

There is a corresponding limited residual result for the property that $\cap
\SS^n R_f = (0)$. If this holds for the pair $(G,f)$ (where $f$ is admissible,
then we say $(G,f)$ satisfies IP (intersection property), and if it holds for
all admissible $f$ in $\R G$, then we say that $G$ satisfies IP. We have already
seen that unique product groups satisfy IP. The following extends this,
somewhat.

\Lem Lemma. Suppose $G$ is a finitely generated group. Sufficient for $G$ to
satisfy IP is the following condition.

\item{} There exists a family of normal subgroups $\brcs{N_{\alpha}}$ \st $\cap
N_{\alpha} = (1)$, $G_{\alpha}:= G/N_{\alpha}$ satisfies IP for all $\alpha$, and there
exists a positive integer $K$ \st for each $G_{\alpha}$  and
admissible $f_{\alpha} \in \R G_{\alpha}$, whenever $f_{\alpha}^k c = 0$ for
some $c \in \R G$, then $f_{\alpha}^K c = 0$.

\Rmk For example, if $\R G_{\alpha}$ all have no zero divisors, then we can set
$K = 0$. But in general, the existence of such a $K$ is difficult to guarantee, hence it is only a
limited residual theorem. A sufficient condition would be that lengths of chains of right annihilators are uniformly bounded over all $\R G_{\alpha}$.

\Pf Suppose there exists admissible $f \in \R G$ \st $\cap \SS^n R_f$ is nonzero. By xxx,
there exists $a \in \R G$ \st for infinitely many $n$, there exists $m(n) \in
\N$ \st $f^{m(n) + n}a \prec f^{m(n)}$ and $f^n a \neq 0$. Let $\pi_{\alpha}$ be
the algebra map $\R G \to \R G_{\alpha}$, and set $f_{\alpha} = \pi_{\alpha}
(f)$. Then $\supp \( f_{\alpha}^{m(n) + n}\pi_{\alpha} (a)\) \subset \supp
\pi_{\alpha}(f)^{m(n)}$ for infinitely many $n$. By xxx and $G_{\alpha}$
satisfying IP, $f_{\alpha}^{m(n)} \pi_{\alpha}(a) = 0$. Hence $f_{\alpha}^K
\pi_{\alpha} (a)= 0$, thus $\pi_{\alpha} \(f^K a\) = 0$ for all $\alpha$; so
$f^K a = 0$, a contradiction (choose $n > K$). \qed

 Let us call the property of the group $G$ that all admissible elements are not zero divisors in $\R G$, $(\blob)$. Then $(\blob)$ is a residual property via an easy proof (unlike the case of the group ring having no zero divisors).

 \Lem Lemma. Suppose that $G$ is a finitely generated group, and there exists a family of normal subgroups $\brcs{N_{\alpha}}$ \st each $G_{\alpha}:= G/N_{\alpha}$ satisfies $(\blob)$, and $\cap N_{\alpha} = \brcs{1}$. Then $G$ satisfies $(\blob)$. 

 \Pf Let $\Arrow \pi_{\alpha}; \R G . \R G_{\alpha}$ denote the onto algebra homomorphism induced by the quotient map $G \mapsto G/N_{\alpha}$. Let $f \in \R G$ be admissible; then $f_{\alpha}:= \pi_{\alpha}(f)$ is admissible (\wrt $G_{\alpha}$). So if for some $a \in \R G$, the product $fa = 0$, then $f_{\alpha}\pi_{\alpha}(a) = 0$. By hypothesis, this forces $\pi_{\alpha}(a) = 0$. Since $G \to \prod G_{\alpha}$ is one to one, so is the map on group rings, $\R G \to \R\prod G_{\alpha}$, and it follows immediately that $a = 0$.\qed 

In contrast, the following is  easy. It yields immediately that all indicable groups (thus including all infinite abelian groups---recall, we are under the blanket assumption that all groups discussed here are finitely generated) satisfy IP for all choices of admissible $f$.

   false: only get modulo augmentation ideal

 \Lem Lemma. Let $N$ be a normal subgroup of $G$. If $(G,f)$ fails to satisfy IP, then  $(G/N, \overline f)$ also fails to satisfy IP. 

 \Pf Let $\Arrow \pi; G.G/N$ be the quotient map; this induces the obvious algebra/ring homomorphism  $AG \to A[G/N]$, denoted $a \mapsto \overline a$. As $f$ is admissible, then $\overline f$ is admissible (\wrt $G/N$). If $(G/N, \overline f)$ satisfied IP, then given $a \in AG$ and nonnegative integer $k$, there would exist $n$ \st the equation (in $R_{\overline f}$, in the variables $b \in AG$ and integer $m \geq n$) $[\overline a,k] = [\overline b,m]$ and $\supp \overline b \subseteq \supp \overline f^{m-n}$ cannot be solved (for $(b,m)$); but this clearly means that neither can $[a,k] = [b,m]$ with $m \geq n$ and $\supp  b \subseteq \supp f^{m-n}$ (we can always reduce to $\supp \overline b \subseteq \supp \overline f^{m-n}$, since the requirement is that $[ b,m-n] \in R_{f}$, and we can just replace $b$ by $f^N b$ (and $m$ by $m+n$) for suitably large $N$). \qed


Powers of the $f$ yield isomorphic things?? Different $\Gamma$s

Yes, really looking at two invariants: intersections of the positive kernels of perfidious traces, and intersections of maximal order ideals (different, even for $G = \Z^2$.

Note that $[g, n] \in R_f$ iff there exists $M$ \st $(\supp f^M)g \subseteq \supp f^{M+n}$; in general, $g \in \Gamma_n$ is sufficient, but not necessary.

It is relatively easy to characterize $\brcs{[g,n]}$ in $R_f$ \st all perfidious traces kill it; this is equivalent to its image in $R_f/\SS R_f$ being a positive infinitesimal, hence for every $N$, there exists $K \equiv K(N)$ \st
$N[g,n] \leq \Bf 1 + K [{}{1_G},1]$. This is equivalent to the existence of $M \equiv M(N)$ \st
$N f^{M+1} g \leq f^{M+N} + K f^{M}$; in other words, $\supp_{-} (f^{M+N} - N f^{M}) \subset \supp f^m$, where $\supp_{-} h$ is the part of the support of $h$ having negative coefficients.

But really looking at intersection of maximal order ideals. I suppose it has something to do with
$N[g,n]-\Bf 1$ belonging to a proper order ideal (or perhaps the negative part of it) (for then, $[g,n]$ is an order unit modulo some order ideal, therefore modulo a maximal order ideal). This amounts to looking at $S:= \supp f^{M+n} \setminus (\supp f^M)g $ (assuming $g \in \supp f^n$), and seeing whether this generates the improper ideal, in the sense that there exists $L$ \st $(\supp f^L)\cdot S = \supp f^{M+n +L}$.

\fillmein Unlikely, but plausible: if $G$ is a noetherian group (all increasing chains of subgroups are eventually stationary; equivalently, all subgroups are finitely generated ) or simply that $\Z G$ is a noetherian ring, then $R_f$ has only finitely many maximal order ideals.
(True for $\Z^d$, all choices of $f$---maximal order ideals correspond to extreme points of Newton polyhedron of $f$; but for the free group, get infinitely many, etc). Maybe noetherian should be replaced by finitely generated (as all these groups are) and amenable.

\fillmein Where is the result about separating maximal order ideals (for all $M \neq M'$, there exists $k'$ \st for all $k\geq k'$, the sets $\Gamma_{k,M} $ and $\Gamma_{k,M'}$ are disjoint. This yields $s(k) = N$ for all sufficiently large $k$ if $R_f$ has exactly $N$ maximal order ideals. 

Pick a positive integer $k$, and let $w = z^k gh^k = h^kg \in \Gamma_{k+1}''$. Then $\tau([w,k+1]) = \lambda^k \tau{[gh^k,k+1]} \leq \lambda ^k \tau([h^k])$. On the other hand, $\tau([h^kg,k+1]) = \tau_g ([h^k,k]) = \tau([g,1])\tau ([h^k,k])$. Hence $\lambda^k \geq \tau([g,1])$ for all $k \geq 1$, provided $\tau([h^k,0]) \neq 0$. Since $\tau([g,1]) > 0$, this forces $\lambda \geq 1$, yielding (since we assumed $\lambda \leq1$) $\lambda = 1$.

 \Lem Zwischenzug. If $\lambda < 1$, then $\tau([h^i,i]) > 0$ for all $i \geq 1$. 

 \Pf Assume not.  Let $k$ be the smallest positive integer \st $\tau([h^i,i]) = 0$. First, we show that $\tau([h,1]) > 0$.  We have $i < k$ entails $\tau([h^i,i]) > 0$ (if $k =1$, then $h^0 =1$, so $\tau([h^0,0]) = \tau([\overline f,1]) =1$ anyway), and if $i > k$, $)$[h^i,i] \leq [h^k,k]$, so $\tau ([h^i,i])  0$. Set $t = \tau([g,1]) > 0$, so that $\tau([h,1]) = 1-t < 1$, and this is also $\max \tau([h^i,i])$. 

 For any positive $m$ and  element $w = z^r g^a h^{m-a} \in \Gamma_m''$, we have 
$$\eqalign{
 \tau([w,m])  & = \lambda^{r-a(m-a)} \tau([z^{a(m-a)}g^a h^{m-a}, m])\cr 
 & = \lambda^{r-a(m-a)} \tau([h^{m-a}g^a , m]) =\lambda^{r-a(m-a)} \tau_g([h^{m-a}g^{a-1} , m-1]) \cr
& = \dots = t^a \lambda^{r-a(m-a)}\tau([h^{m-a},m-a]).\cr 
 }$$
In particular, if $m - a \geq k$, then $\tau([w,m]) = 0$. Equivalently, $\tau(h^{m-a},m-a]) > 0$ entails $a \geq m-k+1$. If $\tau([h,1])= 0 $, then for any $m$, $1 =\sum_{w}\tau([w,m) p(w)= \tau(g^m,m) $ (the only term without an $h$ in $f^m$), so that $\tau = \tau_{0,0}$, a contradiction. Hence $k \geq 2$.  

 Next, we show that $t/\lambda^{k-1} < 1$. For any $\epsilon > 0$, for all sufficiently large $m$, we have that $\tau([f_m,m]) < \epsilon$. In particular, $\sum_{0\leq a \leq m} \tau([g^a h^{m-a}m)},m]) < \epsilon$. We have 
$$\eqalign{
 \epsilon >  \sum_{0\leq a \leq m} \tau([g^a h^{m-a},m]) & = \sum_{0\leq a \leq m} \lambda^{-a(m-a)} \tau([h^{m-a}g^a ,m]) \cr 
 & = \sum_{0\leq a \leq m} \lambda^{-a(m-a)} t^a \tau([h^{m-a} ,m]) \cr  
 & = \sum_{m-k+1\leq a \leq m} \frac{t}{\lambda^{(m-a)} }^a \tau([h^{m-a} ,m-a]) \cr 
 &   = \sum_{j=0}^{k-1} \( \frac{t}{\lambda^{j}}  \)^{m-k-j-1} \tau([h^{j},j]). \cr
 }$$
Since each of the $\tau([h^{k-1},k-1]) > 0$, allowing $m \to \infty$, we must have $t/\lambda^{k-1} < 1$.

 Hence for all sufficiently large  $m$, 
$$\eqalign{
 1 = \tau([\overline f^m,m]) & = \sum_{0\leq a \leq m} \sum_{0\leq r \leq a(m-a)} p(r,a,m-a)\tau([z^r g^a h^{m-a}])\cr 
 & = \sum_{m-k+1\leq a \leq m}  \lambda^{-a(m-a)}t^a\tau([h^{m-a},m-a])\sum_{0\leq r \leq a(m-a)} \lambda^r p(r,a,m-a)\cr
& \leq \sum_{m-k+1\leq a \leq m} \(\frac{t}{ \lambda^{m-a}} \)^{a}(\tau([h^{m-a},m-a]) {m \choose a} \cr 
 & \leq (1-t)\sum_{m-k+1\leq a \leq m} \(\frac{t}{ \lambda^{m-a}} \)^a {m \choose a} \cr
}$$
  The sequences (one for each $m$) $a\mapsto (t/\lambda^{m-a})^a {m \choose a}$ are  strongly unimodal (also known as log concave), at least for $m-k +1 \leq a \leq m$; we compute the ratio of the $s$th term to the $s+1$st term for $m-k+1 \leq s \leq m-1$:
$$
 \rho_{m,s}  = \frac{1}{t\lambda^{2s-1 - m}}\frac{a+1}{m-s}.
$$ 
 With $m$ fixed, increasing $s$ increases the ratio (since $\lambda < 1$), and so each sequence is strongly unimodal. Moreover, if we define $\rho_m = \rho_{m,m-k+1}$ (the left most ratio), we see that $\rho_m =  ((a+1)/(k+1) (t/\lambda^{m-2k+1})$, which for all sufficiently large $m$ is (extremely) large. Hence  
 $$\eqalign{
(1-t)\sum_{m-k+1\leq a \leq m} \(\frac{t}{ \lambda^{m-a}}\)^a  {m \choose a} & \leq (1-t) \(\frac{t}{ \lambda^{k-1}}\)^{m-k+1} \frac {m^k}{k!} \(1 + \frac{1}{\rho_m} + \frac1{\rho_m^2}+ \dots + \) \cr 
   & \leq (1-t)\(\frac{t}{ \lambda^{k-1}}\)^{m-k+1} \frac {m^k}{k!}  \frac{\rho_m}{\rho_m-1}.
}$$ 
 As $\rho_m \to \infty$, and the logarithm of the rest is  $(m-k+1)\ln \gamma  + k \ln m + \oh{1}$ where $\gamma =  {t}/{ \lambda^{k-1} <1$, we see that the sum goes to $0$, a contradiction. \qed 

 So with the Zwischenzug, we see that $\lambda \leq 1$ entails $\lambda = 1$.

 Associated to $A_f$ is the (partially ordered) ring $\Cal E_f$ generated by the positive endomorphisms that commute with $\SS$. In the  case that $G = \Z^d$, this is not very exciting---it is simply the group ring with $f$ inverted, $A[\Z^d][1/f]$, and $R_f$ is a subring and an order in $A_f$. However, with a modicum of noncommutativity, there is a drastic change. For example, if $f = f^*$ (that is, $f$ is symmetric \wrt $g \mapsto g^{-1}$) and $G$ is ICC (all non-identity conjugacy classes are infinite), then $\Cal E_f$ is enormous, which in this context means that the natural  map $AG \times \Z^+ \to \Cal E_f$ is almost one to one; a lesser form of enormousness can be proved when $G$ is any nonabelian group whose group ring has no zero divisors.